\input amstex

\newif\ifproofmode                      
\proofmodefalse				

\newif\ifforwardreference		
\forwardreferencefalse			

\newif\ifchapternumbers			
\chapternumbersfalse			

\newif\ifcontinuousnumbering		
\continuousnumberingfalse		

\newif\iffigurechapternumbers		
\figurechapternumbersfalse		

\newif\ifcontinuousfigurenumbering	
\continuousfigurenumberingfalse		

\font\eqsixrm=cmr6			
\def\marginstyle{\eqsixrm}		

\newtoks\chapletter			
\newcount\chapno			
\newcount\eqlabelno			
\newcount\figureno			

\chapno=0
\eqlabelno=0
\figureno=0


\def\chapfolio{\ifnum \chapno>0 \the\chapno \else \the\chapletter \fi}


\def\bumpchapno{\ifnum \chapno>-1 \global \advance \chapno by 1
	\else \global \advance \chapno by -1 \setletter\chapno \fi
	\ifcontinuousnumbering \else \global\eqlabelno=0 \fi
	\ifcontinuousfigurenumbering \else \global\figureno=0 \fi}

%


%

\def\tempsetletter#1{\ifcase-#1 {}\or{} \or\chapletter={A}\or\chapletter={B}
  \or\chapletter={C} \or\chapletter={D} \or\chapletter={E}
  \or\chapletter={F} \or\chapletter={G} \or\chapletter={H}
  \or\chapletter={I} \or\chapletter={J} \or\chapletter={K}
  \or\chapletter={L} \or\chapletter={M} \or\chapletter={N}
  \or\chapletter={O} \or\chapletter={P} \or\chapletter={Q}
  \or\chapletter={R} \or\chapletter={S} \or\chapletter={T}
  \or\chapletter={U} \or\chapletter={V} \or\chapletter={W}
  \or\chapletter={X} \or\chapletter={Y} \or\chapletter={Z}\fi}

%

\def\chapshow#1{\ifnum #1>0 \relax #1%
   \else {\tempsetletter{\number#1}\chapno=#1 \chapfolio} \fi}

%
\def\today{\ifcase\month\or
January\or February\or March\or April\or May\or June\or
July\or August\or September\or October\or November\or December\fi
\space\number\day, \number\year}

%

\def\initialeqmacro{\ifproofmode
 \headline{\tenrm \today\hfill \jobname\ --- draft\hfill\folio}
     \hoffset=-1cm \immediate\openout2=allcrossreferfile \fi
 \ifforwardreference \input labelfile
     \ifproofmode \immediate\openout1=labelfile \fi \fi}


%

\def\chaplabel#1{\bumpchapno \ifproofmode \ifforwardreference
   \immediate\write1{\noexpand\expandafter\noexpand\def
   \noexpand\csname CHAPLABEL#1\endcsname{\the\chapno}}\fi\fi
   \global\expandafter\edef\csname CHAPLABEL#1\endcsname
   {\the\chapno}\ifproofmode\llap{\hbox{\marginstyle #1\ }}\fi\chapfolio}

%
\def\eqnum{\global\advance\eqlabelno by 1
   \eqno(\ifchapternumbers\chapfolio.\fi\the\eqlabelno)}

\def\eqlabel#1{\global\advance\eqlabelno by 1 \ifproofmode\ifforwardreference
 \immediate\write1{\noexpand\expandafter\noexpand\def
 \noexpand\csname EQLABEL#1\endcsname{\the\chapno.\the\eqlabelno?}}\fi\fi
 \global\expandafter\edef\csname EQLABEL#1\endcsname
 {\the\chapno.\the\eqlabelno?} \eqno(\ifchapternumbers\chapfolio.\fi
 \the\eqlabelno)\ifproofmode\rlap{\hbox{\marginstyle #1}}\fi}

\def\leqlabel#1{\global\advance\eqlabelno by 1 \ifproofmode\ifforwardreference
 \immediate\write1{\noexpand\expandafter\noexpand\def
 \noexpand\csname EQLABEL#1\endcsname{\the\chapno.\the\eqlabelno?}}\fi\fi
 \global\expandafter\edef\csname EQLABEL#1\endcsname
 {\the\chapno.\the\eqlabelno?} \leqno(\ifchapternumbers\chapfolio.\fi
 \the\eqlabelno)\ifproofmode\rlap{\hbox{\marginstyle #1}}\fi}

\def\eqalignnum{\global\advance\eqlabelno by 1
   &(\ifchapternumbers\chapfolio.\fi\the\eqlabelno)}

\def\eqalignlabel#1{\global\advance\eqlabelno by1 \ifproofmode
 \ifforwardreference\immediate\write1{\noexpand\expandafter\noexpand\def
 \noexpand\csname EQLABEL#1\endcsname
     {\the\chapno.\the\eqlabelno?}}\fi\fi
 \global\expandafter\edef\csname EQLABEL#1\endcsname
 {\the\chapno.\the\eqlabelno?}\ifchapternumbers\chapfolio.\fi
 \the\eqlabelno\ifproofmode\rlap{\hbox{\marginstyle
 #1}}\fi}

\def\eqref#1{(\ifundefined{EQLABEL#1}***\ifproofmode\ifforwardreference)%
   \else \write16{ ***Undefined Equation Reference #1*** }\fi
   \else \write16{ ***Undefined Equation Reference #1*** }\fi
   \else \edef\LABxx{\getlabel{EQLABEL#1}}%
   \def\LAByy{\expandafter\stripchap\LABxx}\ifchapternumbers
   \chapshow{\LAByy}.\expandafter\stripeq\LABxx
   \else\ifnum \number\LAByy=\chapno \relax\expandafter\stripeq\LABxx
   \else\chapshow{\LAByy}.\expandafter\stripeq\LABxx\fi\fi)\fi
   \ifproofmode\write2{Equation #1}\fi}

%

\def\fignum{\global\advance\figureno by 1 \relax
   \iffigurechapternumbers\chapfolio.\fi\the\figureno}\

\def\figlabel#1{\global\advance\figureno by 1\relax
 \ifproofmode\ifforwardreference
 \immediate\write1{\noexpand\expandafter\noexpand\def
 \noexpand\csname FIGLABEL#1\endcsname{\the\chapno.\the\figureno?}}\fi\fi
 \global\expandafter\edef\csname FIGLABEL#1\endcsname
 {\the\chapno.\the\figureno?}\iffigurechapternumbers\chapfolio.\fi
 \ifproofmode$^{\hbox{\marginstyle #1}}$\relax\fi\the\figureno}

\def\figref#1{\ifundefined{FIGLABEL#1}!!!!\ifproofmode\ifforwardreference)%
   \else \write16{ ***Undefined Equation Reference #1*** }\fi
   \else \write16{ ***Undefined Equation Reference #1*** }\fi
   \else \edef\LABxx{\getlabel{FIGLABEL#1}}%
   \def\LAByy{\expandafter\stripchap\LABxx}%
   \iffigurechapternumbers\chapshow{\LAByy}.\expandafter\stripeq\LABxx
   \else\ifnum\number\LAByy=\chapno \relax\expandafter\stripeq\LABxx
   \else\chapshow{\LAByy}.\expandafter\stripeq\LABxx\fi\fi
   \ifproofmode\write2{Figure #1}\fi\fi}

%

%

\def\getlabel#1{\csname#1\endcsname}
\def\ifundefined#1{\expandafter\ifx\csname#1\endcsname\relax}
\def\stripchap#1.#2?{#1}
\def\stripeq#1.#2?{#2}

\figurechapternumberstrue  

\chapternumberstrue        

\def\thmlbl#1{\figlabel{#1}}
\def\thmref{\figref}
\def\eqnlbl#1{\leqlabel{#1}}
\def\eqnalignlbl#1{\eqalignlabel{#1}}
\def\eqnref#1{\eqref{#1}}
\def\sectionnumber{\chapno}
\def\theoremnumber{\figureno}
\def\equationnumber{\eqlabelno}

\documentstyle{amsppt}
\overfullrule=0pt
\magnification =\magstep1
\baselineskip=18pt
\vcorrection{-.33truein}
\pageheight{9.0truein}
\input amssym.def
\def\BU{\bar U}
\def\tG{\tilde G}
\def\ratio{\bar \alpha}
\def\Tpsi{{\tilde{\psi}}}
\def\tvar{{\tilde{\phi}}}
\def\bvarphi{{\bar{\varphi}}}

\def\balpha{{\bar{\alpha}}}
\def\ratio{\bar \alpha}
\def\etaT{{\widetilde{\eta}}}
\def\Tmu{{\widetilde{\mu}}}

\def\tildeM{{\widetilde{M}}}
\def\tildeW{{\widetilde{w}}}
\def\tildeC{{\widetilde{C}}}

\def\BbbR{\Bbb{R}}
\def\BbbC{\Bbb{C}}

\def\BbbA{\Bbb{A}}

\def\BbbJ{\Bbb{J}}
\def\BbbK{\Bbb{K}}

\def\CalD{\Cal{D}}
\def\CalH{\Cal{H}}

\def\CalT{\Cal{T}}

\def\CalW{\Cal{W}}
\def\bU{{\bar{u}}}
\def\bv{{\bar{v}}}
\def\R{\roman{Re\ }}
\def\I{\roman{Im\ }}
\def\sgn{\text{\rm sgn\ }}

\def\det{\text{\rm det\ }}

\def\Range{\text{\rm Range\ }}

\def\Residue{\text{\rm Res}}
\def\Span{\text{\rm Span\ }}
\def\stab{{($\Cal{D}$)}}

\def\factorial{!}
\def\dim{\text{\rm dim \ }}
\def\adj{\text{\rm adj\ }}
\def\diag{\,\text{\rm diag}\,}
\def\Re{\text{\rm Re\ }}
\def\errfn{\text{\rm errfn\ }}
\def\Im{\text{\rm Im }}
\def\deltdot{{\cdot \atop {\raise8pt\hbox{$\delta$}}}}
\def\myqed{\vrule height3pt depth2pt width3pt \bigskip}
\def\newsection{\centerline}
\def\cal{\Cal}
\def\exp{{\text{\rm exp\ }}}
\def\Loc{\text{\rm Loc\ }}

\def\bfO{{\Cal{O}}}

\def\zigzag{\buildrel > \over <}


\def\Fyx{{\Cal F}^{y\to x}\,}
\def\Fzx{{\Cal F}^{z\to x}\,}
\def\Fzy{{\Cal F}^{z\to y}\,}
\def\e{\varepsilon}

\pagewidth{6truein}
\pageheight{9truein}
\topmatter
\title
{Pointwise Green's function bounds and
stability of relaxation shocks}
\endtitle
\leftheadtext{Stability of Relaxation Shocks}
\rightheadtext{Corrado Mascia and Kevin Zumbrun}
%
\thanks
The first author thanks Indiana University and Scholarship 
CNR n.203.01.68 (07/05/97) for their hospitality and for
making possible the postdoctoral visit
(January--August 1999) during which 
this work was initiated and partly carried out.
The second author thanks 
Instituto per le Applicazioni del Calcolo 
``M. Picone'' (CNR) and European TMR Project ``Hyperbolic Systems of 
Conservation Laws'' for their hospitality and for making possible
the visit (April 10--May 10, 2000) 
during which the work was completed.
Research of the first author was supported in part by 
European TMR Project ``Hyperbolic Systems of Conservation Laws''.
Research of the second author was supported
in part by the National Science Foundation under Grants No. DMS-9107990
and DMS-0070765.
\endthanks

\abstract
We establish sharp pointwise Green's function bounds 
and consequent linearized stability
for smooth traveling front solutions,
or relaxation shocks, of general hyperbolic
relaxation systems of dissipative type,  under the necessary assumptions 
([G,Z.1,Z.4])
of spectral stability, i.e., stable point spectrum of the linearized 
operator about the wave,
and hyperbolic stability of the corresponding ideal shock 
of the associated equilibrium system,
with no additional assumptions on the structure or strength of the shock.
Restricting to Lax type shocks,
we establish the further result of nonlinear stability with respect to small 
$L^1\cap H^2$ perturbations, 
with sharp rates of decay in $L^p$, $2\le p\le \infty$,
for weak shocks of general simultaneously symmetrizable systems;
for discrete kinetic models, 
and initial perturbation small in $W^{3,1}\cap W^{3,\infty}$,
we obtain sharp rates of decay in $L^p$, $1\le p\le \infty$, for
(Lax type) shocks of arbitrary strength.
This yields, in particular, nonlinear stability
of weak relaxation shocks of the 
discrete kinetic Jin--Xin and Broadwell models, for which
spectral stability has been established in [HL,JH] and [KM], respectively.

Our analysis follows the basic pointwise semigroup approach introduced 
by Zumbrun and Howard [ZH] for the study 
of traveling waves of parabolic 
systems; however, significant extensions are required to deal with
the nonsectorial generator and more singular short-time behavior of 
the associated (hyperbolic) linearized equations.
Our main technical innovation is a systematic method for 
refining large-frequency (short-time) estimates on the resolvent kernel, 
suitable in the absence of parabolic smoothing.
This seems particularly interesting from the viewpoint of general
linear theory, replacing the zero-order estimates of existing theory with a 
series expansion to arbitrary order.
The techniques of this paper should have further application in the 
closely related case of traveling waves of systems with partial 
viscosity, for example in compressible gas dynamics or MHD.
\endabstract
\author
{Corrado Mascia and Kevin Zumbrun}
\endauthor
\date{July 11, 2001; Revised: May 2, 2002}
\enddate
\address
Dipartimento di Matematica ``G. Castelnuovo'',
Universit\`a di Roma ``La Sapienza'',
P.le Aldo Moro, 2 - 00185 Roma (ITALY)
\endaddress
\email
mascia\@ mat.uniroma1.it
\endemail 
\address
Department of Mathematics,
Indiana University,
Bloomington, IN  47405-4301
\endaddress
\email
kzumbrun\@indiana.edu
\endemail 
\endtopmatter
\document
\newsection {\bf Section 1. Introduction}
\sectionnumber=1 
\theoremnumber=0
\equationnumber=0
\smallskip
\TagsOnLeft

A variety of nonequilibrium processes in continuum mechanics can be
modeled by {\it hyperbolic  relaxation systems} of general form
$$
 \pmatrix u\\ v \endpmatrix_t
 + \pmatrix f(u,v)\\ g(u,v) \endpmatrix_x
 = \pmatrix 0\\ \tau^{-1} q(u,v)
 \endpmatrix,
 \eqnlbl{general}
$$
$u$, $f\in \Bbb{R}^n$,
$v$, $g$, $q \in \Bbb{R}^r$,
where 
$$
\text{\rm Re }\sigma\big(q_v(u,v^*(u))\big)<0
\eqnlbl{qv}
$$
along a smooth {equilibrium manifold} defined by 
$$
q(u,v^*(u))\equiv 0,
\eqnlbl{eq}
$$
and $\tau$ is a (usually small) parameter determining relaxation time:
for example, non-thermal equilibrium gas dynamics 
[W,MR,Lev], traffic dynamics [LW,AR,Li],
and multiphase flow [BGDV,BV,NT,MaP]; 
for a general discussion of relaxation models, we refer the reader to
[L.2,CLL] or to
the surveys [N,Yo.4].

A particularly natural class of relaxation systems are those
obtained by discretizing continuous
kinetic models such as Boltzmann or Vlasov--Poisson equations,
for which the lefthand side of \eqnref{general} reduces to a simple
linear transport operator: 
for example, Broadwell and other lattice gas models [PI].
Here, loosely following Natalini [N], we define the class of
{\it discrete kinetic models} as
systems of form \eqnref{general} for which $f$ and $g$ are linear, 
constant coefficient.  These include in particular the subclass
of semilinear relaxation schemes, such as the Jin--Xin [JX] 
and Natalini [N]
models, that have enjoyed considerable recent popularity as a
method for numerical approximation of hyperbolic systems.

The relaxation system \eqnref{general} can be viewed ``to zeroth order'' as a 
regularization of the associated {\it equilibrium}, or ``relaxed'' system
of conservation laws
 $$
  u_t + f^*(u)_x=0,
  \eqnlbl{hyp}
 $$
$f^*(u):= f(u,v^*(u))$, 
which system may or may not be hyperbolic (see [JX,BGDV,BV]
for examples of nonhyperbolicity).
To ``first order,'' it may be approximated, at least formally,
by an associated parabolic system of conservation laws
 $$
  u_t + f^*(u)_x=\tau (B^*(u)u_x)_x,
  \eqnlbl{parab}
 $$
where 
 $$
  \eqalign{
   B^*(u)
    &:= -f_v q_v^{-1}(g^*_u - v^*_u f^*_u) \cr
    &= -f_v q_v^{-1}(g_u - g_vq_v^{-1}q_u
     + q_v^{-1}q_u(f_u - f_v q_v^{-1} q_u), \cr
           }
  \eqnlbl{B}
 $$ 
$g^*:=g(u,v^*(u))$, is determined by Chapman--Enskog expansion as described
in, e.g.,  [Wh,L.2,Z.1]. 
Here, the order referred to is with respect to relaxation time
$\tau$; however, as pointed out by Liu in his seminal work [L.2], 
this formal approximation
has a more readily justifiable interpretation in terms of {\it large time},
rather than small $\tau$ behavior.
It is this point of view that is relevant for the present work,
in which we investigate large time behavior for a {\it fixed} 
relaxation time.  
From now on, we take without loss of generality
$\tau=1$,  and suppress the parameter $\tau$.

Note that, in the frequently occurring case that $r<n$, the $n\times n$ matrix
$B^*$ is necessarily singular, since it factors through the $r\times r$ matrix
$q_v^{-1}$.  Thus, the proper analogy is to ``real viscosity,'' or
partially parabolic systems of conservation laws.
Indeed, this analogy holds even in the case that $B^*$ is nonsingular,
since the high frequency behavior of the original system \eqnref{general}
is by necessity hyperbolic, irregardless of formal expansion;  see
[Ro,Sl.1--2,Ma] for interesting discussions on the
validity/meaning of the Chapman--Enskog expansion in the high-frequency regime.
The analogy is in fact much deeper than this. 
As pointed out by Zeng [Ze.2],
the two problems (relaxation/real viscosity) are essentially dual in 
the linear constant--coefficient case.  A similar duality,
at the structural level,
may be seen in the variable-coefficient case [Z.4]

A particularly interesting phenomenon, suggested by \eqnref{hyp},
is the existence of smooth {\it traveling front solutions}
$$
\eqalign{
(u,v)(x,t)&= (\bar u, \bar v)(x-st) \cr
\lim_{z\to \pm \infty}
(\bar u, \bar v)&= (u_\pm , v_\pm),
}
\eqnlbl{traveling wave}
$$
of \eqnref{general}, where by necessity 
$v_\pm=v^*(u_\pm)$ and $u_\pm$  corresponds to a shock solution of \eqnref{hyp}.
See [L.2] for a treatment of existence in the general case $n=r=1$, 
and [YoZ,FZe] for generalizations in case $n$ or $r>1$;
further results/discussion are given in Section 1.1, below.
Such traveling waves are known as {\it relaxation shocks} or
{\it relaxation profiles}.
The question of their stability was investigated in [L.2]
in the general case $n=r=1$, for which the equilibrium system is scalar,
under the assumption of ``weak,'' or small-amplitude shocks.
This analysis has been generalized to strong shocks in the special
case of the $2\times 2$ Jin--Xin model [MaN].
For the $3\times 3$ Broadwell model, a standard discrete kinetic model
for which $n=2$, $r=1$,  Szepessy and Xin [S,SX.2] 
have announced the result of
linear and nonlinear stability of weak shocks, extending partial
results of [KM,CaL]:
to our knowledge, the only such result for
the system case $n>1$ occurring in most physical applications.
However, to date there is no complete analysis of stability for any other case,
in particular for the case $n>1$, $r>>1$ arising through
approximation of Boltzmann
or Vlasov--Poisson equations by moment closure or discretization.
For strong shocks, there is no complete analysis for {\it any} system
such that $n$ or $r>1$.

Useful necessary conditions for stability have been obtained for 
shock profiles of the Jin--Xin relaxation model in [Z.1,G] 
and for shock profiles of general relaxation 
and real viscosity models in [Z.4], generalizing corresponding
results of [GZ,BSZ,ZS] in the strictly parabolic case.
Strengthened versions of the one-dimensional hyperbolic stability 
criterion of Erpenbeck--Majda [Er,M.1--3],
these conditions yield interesting results of {\it spectral instability}
(hence linearized and nonlinear instability), similarly as in the
strictly parabolic case (see e.g. [GZ,FreZ,ZS,Z.4]). 
%
Notably, the above-mentioned analyses apply to shocks of arbitrary strength
and type (i.e., Lax, under-, or overcompressive).

However, to our knowledge, for $n>1$, the only positive stability
result for relaxation shocks other than the ones of [KM,CaL,S]
is a partial result recently established by H. Liu [Liu]
of linear and nonlinear {\it zero-mass stability}
of weak shocks of the Jin--Xin system,
i.e., stability with respect to perturbations whose integral in
the $u$ variable is zero; see also a subsequent and closely related proof 
by Humpherys [Hu] of {\it spectral stability} of small-amplitude 
Jin--Xin shocks, defined as stability of the
point spectrum of the linearized operator about the wave.
(Note: spectral stability is implied by zero-mass
stability, see discussion, [ZH] or [Hu]).
Being based on energy methods, 
none of these previous analyses for systems ([KM,S,Liu,Hu]) 
generalize to the case of strong (i.e., large-amplitude) shocks.
Moreover, they each depend to some degree on special structure of
the respective systems treated;
hence, the extent to which they generalize to other systems is not clear.

In the present paper, following the general philosophy set out in [ZH]
for the study of viscous shock waves,
we investigate the linearized and nonlinear stability
of relaxation shocks of general systems \eqnref{general},
with arbitrary $n$, $r$, 
but {\it assuming spectral stability.}
As in the parabolic case, what makes this enterprise nontrivial
is that the essential spectrum of the linear generator accumulates
at the imaginary axis, so that linearized decay cannot immediately
be concluded from spectral stability, and moreover can be at most 
at time-algebraic and not exponential rate.  Indeed, much of the
technical difficulty of the problem lies in this usually trivial
step.  Spectral stability can then be investigated separately as
an ODE problem, or checked numerically as in [Br.1--2].
An important advantage of this approach is that it 
is completely general, applying in principle to shocks {\it of all types}
(i.e., Lax, under-, or overcompressive) {\it and strengths}.
\footnote{
However, see Remark 2 just below, regarding subshocks.}

More precisely, we study linearized stability for 
systems \eqnref{general} satisfying appropriate dissipation
conditions as in [Yo.4] or [Ze.2],
and nonlinear stability for weak Lax shocks of (general) symmetrizable systems
(defined Section 7.2),
or strong Lax shocks of discrete kinetic models.
Our main example will be the Jin--Xin model,  for which calculations
are particularly transparent. 
In this case, combining our results with the existing
spectral stability results of [Liu,Hu],
we obtain the result that Lax-type Jin--Xin shocks of sufficiently weak
strength are nonlinearly stable under small $L^1\cap L^\infty$ perturbations,
with no additional hypotheses:
the first such result for a relaxation shock in the case $n$, $r>1$.
Likewise, using the zero-mass stability
result of [KM] (recall: zero-mass implies spectral stability), 
we recover a more detailed version of the result
of [S,SX.2] that weak Broadwell shocks are nonlinearly stable,
with the additional information of sharp rates of decay in all $L^p$.
(Note: No rates of decay are given in [S,SX.2].)
In the course of our analysis, we develop extremely detailed 
bounds on the Green's function of the linearized evolution equations 
about the traveling wave, of interest in their own right.
These generalize bounds obtained in 
[ZH,Z.2--3] for shock profiles of strictly parabolic systems,
and in [Ze.1,LZe], [Ze.2] for constant solutions of partially parabolic
and relaxation systems, respectively.

Our results are of interest from the physical point of view
as some of the first complete shock stability results for 
systems with realistic dissipative mechanism such as relaxation
or partial viscosity: 
in particular the first positive stability results
for strong relaxation shocks (namely, the reduction to the numerically
verifiable spectral stability condition ($\Cal{D}$)),
the first for quasilinear systems,
and the first results for relaxation shocks of any strength 
to give rates of decay.
From the technical side, they are of interest
as a blueprint for the application of the pointwise semigroup
methods of [ZH] to the case that the associated semigroup is $C^0$ 
and not analytic: in particular, when the Green's function $G(x,t;y)$
is a distribution, lying outside the class of smooth functions for $t>0$.

With regard to the latter, technical issue, 
we mention related analyses carried out by
Dodd [Do] and Howard--Zumbrun [HZ.2] for shock profiles of
dissipative--dispersive conservation laws, for which the
semigroups associated with the linearized equations are likewise
merely $C^0$.  The analysis of [HZ.2] in particular proceeds
from a somewhat similar point of view to that followed here,
and could be considered as an intermediate point between the
analysis of [ZH] and that carried out here.
The main technical innovation distinguishing the present treatment
is a systematic method for refining high-frequency bounds on
the resolvent kernel to any desired accuracy;  for comparison,
the zeroth order step in this algorithm gives bounds at the level 
of [ZH], and the first order approximately at the level of [HZ.2]
\footnote{
In fact, we make also key improvements generalizing the approach of [HZ.2]
to the system case, and to operators with nonconstant principal part, 
see Remarks 8.16 and 4.10 below.
}, 
whereas here we require bounds to the second order.
In effect, we generate an arbitrary order asymptotic expansion,
limited only by smoothness of the coefficients of the linearized operator, 
of the resolvent about an essential singularity at infinity, 
with strong error bounds
valid almost up to the boundary of analyticity/essential spectrum.
This, and other techniques of the paper
should apply also to systems with incomplete 
parabolicity, such as the equations of compressible gas dynamics or MHD.

\medskip
We now describe our results in somewhat greater detail.
Consider a smooth traveling wave solution \eqnref{traveling wave} of
relaxation system \eqnref{general}.
Besides \eqnref{qv}--\eqnref{eq}, we make the following, rather
standard assumptions (see, e.g., [Yo.1--3,YoZ,Z.4]):

\medskip
\noindent(H0) \quad $f$, $g$, $q \in C^3$ in the general case,
$f$, $g$, $q \in C^2$ for discrete kinetic models.

\smallskip
\noindent(H1) \quad $\displaystyle{
\sigma\left( (df,dg)^t (u,v)\right)}$ real, semi-simple with
constant multiplicity, and different from $s$.

\smallskip

\noindent(H2) \quad $\displaystyle{
 \sigma\left( df^*(u_\pm)\right)}$ real, distinct, and different from $s$. 
\smallskip

\noindent(H3) \quad $\displaystyle{
\R \sigma\left(  i\xi\binom{df}{dg} (u_\pm, v_\pm) +
\binom 0{dq} (u_\pm, v_\pm)\right) \leq \frac{-\theta |\xi|^2}{1+|\xi|^2}}$,
\quad $\theta>0$,
for all $\xi\in \BbbR$.
\smallskip

\noindent(H4) \quad
The set of solutions of \eqnref{traveling wave} forms a smooth
manifold $(\bU^\delta,\bv^\delta),\ \delta\in \Cal{U} \in \BbbR^\ell$.

\medskip
With the possible exception of the ``noncharacteristic condition,''
$s\not\in \sigma (df,dg)^t$, these are precisely analogous 
to the hypotheses given in [ZH,Z.4] for the viscous case.
Conditions (H1)--(H3), with $\theta$ set to zero in (H3),
and relaxing the requirement of distinct eigenvalues in (H1),
are essentially the standard set of hypothesees proposed 
by W.-A. Yong for general relaxation systems [Yo.1--3];
as discussed in [Y.4], these are satisfied for most known relaxation
sytems, and indeed may be considered as criteria for a ``reasonable'' model.
A slight difference here is that we have assumed strict 
hyperbolicity/nonsonicity of the
equilibrium system, (H2), and stability, (H3), only at the endpoints 
$(u_\pm,v_\pm)$ of the shock, 
and not everywhere along the profile; this allows applications to 
interesting nonhyperbolic situations such as are discussed in [JX,BGDV,BV]:
in particular, shocks of nonclassical, under- and overcompressive type 
are admitted in our treatment, along with the classical, Lax type.
Existence of weak (necessarily Lax-type) relaxation profiles has 
been shown under the 
hypotheses of [Y.1--3] 
together with some additional mild nondegeneracy conditions
(satisfied for most physical systems)
by Yong and Zumbrun [YoZ]; see [FreZ,BR] 
for related studies.
In Appendix A1, we give a simplified treatment of existence,
showing that existence of weak relaxation profiles is
in fact guaranteed by (H0)--(H3) alone.

The slightly strengthened condition (H3) ($\theta>0$) may be recognized as 
the {\it dissipativity} condition of Zeng--Kawashima [Ze.2,Kaw],
and is likewise satisfied for a variety of physical models,
though not all.
The degenerate case $\theta=0$ is quite delicate, and requires 
additional hypotheses on the nonlinear structure of \eqnref{general}; for 
a discussion of this interesting situation and related phenomena, see [Ze.2].

{\bf Remarks.} 1. The noncharacteristic condition 
$s\not\in \sigma (df,dg)^t$ is necessary
in order that the traveling wave ODE be of nondegenerate
type, a natural assumption in the context of stability 
of relaxation profiles.
(See Section 1.1, below, for a statement of the traveling wave ODE).

2. In (H0)--(H4), we have made no assumptions on the structure or amplitude 
of the profile $(\bar u,\bar v)$.
However, note that the assumption of a smooth profile is in effect
a limitation on the strength of the shock, since, as is well known,
profiles of sufficiently large amplitude will develop ``subshocks,'' 
or jump discontinuities [Whi,L.2].
More precisely, this is a restriction on the amplitude 
with respect to the number of resolving modes in the system [BR,MR]:
For example, the allowable amplitude for a moment closure system
is a strictly increasing function of the number of moments,
apparently approaching infinity in the limiting (infinite-dimensional)
case of the Boltzmann equations.
The stability of profiles containing subshocks is a very interesting
problem for future investigation,  
and likewise the stability (or instability) of strong Boltzmann profiles.
Of course, stability of even weak Boltzmann profiles is 
a fundamental open problem.
\medskip
{\bf Example \thmlbl{JX}.} 
The Jin--Xin model,  
$$
\aligned
u_t + v_x&=0,\\
v_t + a^2 u_x&= h(u)-v,\\
\endaligned
\eqnlbl{JX}
$$
$u$, $v\in \BbbR^n$, was introduced in [JX] as a numerical
scheme approximating solutions of the equilibrium system $u_t+h(u)_x=0$.
System \eqnref{JX} is in general nonstrictly hyperbolic,
with $\sigma (df,dg)^t$ consisting of the two eigenvalues $\pm a$,
each with constant multiplicity $n$.
It is readily checked that (H0)--(H3) hold if and only if:
(i) $h\in C^3$;
(ii) $a\ne s$;
(iii) the equilibrium system $u_t+h(u)_x=0$ is (nonstrictly) hyperbolic
at the endstates $u_\pm$,
with $s\not \in \sigma\big(dh(u_\pm)\big)$; and,
(iv) there holds {\it at the endpoints of the shock} 
the subcharacteristic condition $a^2> \sigma \big(dh(u_\pm)^2\big)$.

Likewise, one finds that the traveling wave ODE for a relaxation shock
of speed $s$ reduces to that for an ``effective'' viscous system,
$$
u_t+ h(u)_x= (a^2-s^2)u_{xx},
\eqnlbl{JXode}
$$
with constant, scalar viscosity $a^2-s^2$, so that (H4) reduces to
the corresponding requirement in the viscous case.
Note that (iv) above implies that the effective viscosity $a^2-s^2$ is 
positive whenever a profile exists, since in that case $s$ must lie 
between the extreme values of the spectra of $dh(u_\pm)$.
This furnishes a rich set of examples of relaxation shocks
falling under our framework, through the existence theory
for viscous traveling waves.
Note that \eqnref{JX} is a discrete kinetic model, since
its lefthand side is linear.
\myqed
\medskip

{\bf 1.1. Structure/Classification of Profiles.}
We begin by presenting some basic, but apparently new results
relating the structure of relaxation profiles to the characteristic
structure of the corresponding equilibrium shock in the relaxed system
\eqnref{hyp}.  
These generalize observations of Majda and Pego 
in the viscous, strictly parabolic case [MP];
analogous results hold for viscous profiles with
real (i.e., partially parabolic) viscosity, as can be seen by
essentially the same arguments used here (see, for example,
the unified framework for relaxation/real viscosity
models given in [Z.4], Appendices A1/A2).

Without loss of generality
taking $s=0$, by a linear change in the flux functions if necessary,
we may write the traveling wave ODE in the form:
$$
\aligned
f(u,v)'&=0,\\
g(u,v)'&= q(u,v),\\
\endaligned
\eqnlbl{zero ode}
$$  
where $'$ denotes $d/dx$.
From \eqnref{zero ode}, we immediately see that $q(u_\pm,v_\pm)=0$
and $f(u)\equiv {\roman constant}$, so that endstates 
$(u_\pm,v_\pm)$ lie on the equilibrium manifold, and satisfy
the Rankine--Hugoniot jump condition
$$
f^*(u_+)-f^*(u_-)=
f(u_+,v_+)-f(u_-,v_-)=0
\tag{RH}
$$
for the equilibrium system \eqnref{hyp}: that is, $(u_-,u_+)$
corresponds to a shock solution of \eqnref{hyp}.
Let us now consider the behavior of solutions of \eqnref{zero ode} near the
rest points $(u_\pm,v_\pm)$,
for definiteness $(u_+,v_+)$.

By assumption (H1) together with the Inverse Function Theorem,
$(u,v)\to (f,g)$ may be inverted in a neighborhood of $(u_+,v_+)$,
reducing \eqnref{zero ode} to the $r$-dimensional ODE
$$
g'=q\big(u(f_+,g),v(f_+,g)\big),
\eqnlbl{gode}
$$
where $f$ is held constant at value $f_+$.  Linearizing about the rest point
$g_+$, we obtain after a brief calculation the linearized equation
$$
g'=(q_u,q_v)
\pmatrix 
f_u & f_v\\
g_u & g_v\\
\endpmatrix ^{-1}
\pmatrix 
0\\
I_r\\
\endpmatrix g,
\eqnlbl{relaxtype}
$$
where the righthand side is understood to be evaluated at $(u_+,v_+)$.

\proclaim{Lemma \thmlbl{structure}}
Under assumptions (H0)--(H3), 
the matrix
$$
(q_u,q_v)
\pmatrix 
f_u & f_v\\
g_u & g_v\\
\endpmatrix ^{-1}
\pmatrix 
0\\
I_r\\
\endpmatrix _\pm
\eqnlbl{redcoeff}
$$
has no center subspace, i.e., $(u_\pm,v_\pm)$
are hyperbolic rest points of the reduced traveling
wave ODE \eqnref{gode}.
In particular, traveling wave solutions
\eqnref{traveling wave} satisfy
$$
|(d/dx)^k \big((\bar u(x),\bar v(x))-(u_\pm,v_\pm)\big)|
\le Ce^{-\theta |x|}, \quad  k=0,\dots,4,
\eqnlbl{expdecay}
$$
as $x\to \pm \infty$.
\endproclaim

Lemma \thmref{structure} verifies the hypotheses of
the Gap Lemma, [GZ,ZH]\footnote{
See also the version established in [KS],
independently of and simultaneously to that of [GZ].
},
on which we shall rely to obtain the basic ODE estimates
underlying our analysis in the low-frequency/large-time regime.
A proof of Lemma \thmref{structure} is given in Appendix A1,
a description and proof of the Gap Lemma in Appendix A3.

We shall classify relaxation shocks mainly according to their
{\it hyperbolic type}, i.e., the type of the corresponding hyperbolic
shock for the equilibrium system \eqnref{hyp}. 
Specifically, let $i+$ denote the
dimension of the stable subspace of $df^*(u_+)$,
$i_-$ denote the dimension of the unstable subspace of $df^*(u_-)$,
and $i:=i_+ + i_-$.  Indices $i_\pm$ count the number of incoming
characteristics from the right/left of the shock, while $i$ counts
the total number of incoming characteristics toward the shock.
Then, the relaxed (hyperbolic) shock $(u_-,u_+)$ is called:
$$
\cases
\hbox{Lax type} & \hbox{if} \quad i=n+1,\\
\hbox{Undercompressive (u.c.)} &\hbox{if} \quad i \le n,\\
\hbox{Overcompressive (o.c.)} &\hbox{if} \quad i \ge n+2.
\endcases
$$
In case all characteristics are incoming on one side, i.e.,
$i_+=n$ or $i_-=n$, a shock is called  {\it extreme}.

Similarly as in the parabolic case, [MP,ZH,ZS,Z.4], there is
a close connection between the hyperbolic type of a shock 
profile, and the nature of the corresponding connection
in the phase portrait of the $r$-dimensional 
reduced traveling wave ODE \eqnref{gode} on the $n$-dimensional submanifold 
$$
f(u,v)-su \equiv f(u_-,v_-)-su_-
$$
of $\BbbR^{n+r}$, obtained by integration of the 
$u$-equation $\big((f(u,v)-su\big)\equiv 0$ in the full traveling wave ODE.
Define $d_+$ to be the dimension within this submanifold
of the stable manifold at $(u+,v_+)$, and $d_-$ the dimension
of the unstable manifold at $(u_-,v_-)$, and $d:=d_-+d_+$.
Then, we have:

\proclaim{Lemma \thmlbl{connection}}
Under assumptions (H0)--(H3), 
there holds
$$
d-r=i-n.
\eqnlbl{indexrelation}
$$
\endproclaim

Lemma \thmref{connection} generalizes an analogous relation pointed
out in the viscous case by Majda and Pego [MP]
(see, e.g., Lemma 1.1 [Z.4]).  However, the proof is
different from the essentially algebraic one given in [MP], 
following instead an alternative argument suggested in Remark 2.3, [ZH],
that is more closely linked to the dynamics of the traveling wave;
for an exposition of this proof, see Appendix A1.

A complete description of the connection, of course,
requires the further index $\ell$ defined in (H4) as the 
dimension $\ell$ of the connecting manifold
between $(u_\pm,v_\pm)$ in the traveling wave ODE.
Generically, one expects that $\ell$ should be equal to the surplus $d-r=i-n$.
In case the connection is ``maximally transverse'' in this sense, i.e.:
$$
\ell = \cases
1 &\hbox{undercompressive or Lax case},\\
i-n &\hbox{overcompressive case},
\endcases
\eqnlbl{pure}
$$
we call the shock ``pure'' type, and classify it according to 
its hyperbolic type;
otherwise, we call it ``mixed'' under/overcompressive type.
Throughout the paper, we shall assume that \eqnref{pure} holds,
so that all relaxation profiles are of {\it pure, hyperbolic type}.
This holds in particular for the weak profiles whose
existence was established in [YoZ], which are always of pure Lax type.

{\bf Remark.}
In the case of the Jin--Xin model \eqnref{JX}, relation 
\eqnref{indexrelation} follows directly from
the corresponding parabolic result of [MP], see discussion,
Example \thmref{JX}.

\medskip

{\bf 1.2. Spectral vs. Linearized Stability.}
We next discuss spectral and linearized 
stability of profiles, and the relation between them:
specifically, we show that spectral stability, appropriately
defined, is necessary and sufficient for linearized stability.
Linearizing \eqnref{general} about the stationary solution 
$(\bar u,\bar v)$, we obtain the linearized equations
$$
U_t= LU:=
-(AU)_x  + QU,
\eqnlbl{linearized}
$$
where
$$
A := \binom{df}{dg} (\bU,\bv),\quad Q:=
\binom{0}{dq(\bU,\bv)}, 
\eqnlbl{coefficients}
$$
and
$$
U  = (u,v)^t,\ u\in\BbbR^n,\ v\in \BbbR^r.
\eqnlbl{U}
$$

The appropriate notion of stability of shock profiles 
is {\it orbital stability}, or convergence to the 
manifold of stationary solutions $\{(\bar u^\delta,\bar v^\delta)\}$ 
defined in (H4).
Likewise, the relevant notion of linearized stability is
that of {\it linearized orbital stability}, defined as
convergence to the {\it tangent stationary manifold}, 
consisting of stationary solutions
$\{(\partial/\partial \delta_j)(\bar u^\delta,\bar v^\delta)\}$
of the linearized equations \eqnref{linearized} [ZH].
An evident necessary condition for linearized orbital stability,
(and thus for nonlinear orbital stability, appropriately defined) is 
\medskip
(D1) \qquad $\sigma (L) \cap \{\R \lambda \ge 0\}=\{0\},$
\medskip
\noindent where $\sigma$ refers to the spectrum of $L$ with respect to some
$L^p$, say $L^2$ for simplicity.  We shall refer to condition (D1) as
{\it strong spectral stability}.

A standard result of Henry [He] gives that the rightmost boundary
of the essential spectrum $\sigma_{\roman ess}(L)$ is given by
the right envelope of the union of the (purely essential) spectra 
of the limiting, constant-coefficient
operators $L_\pm$ as $x\to \pm \infty$, which by (H3) intersect
the non-strictly-stable half-plane
$\R \lambda \ge 0$ precisely at the origin $\lambda=0$.
Thus, (D1) is equivalent (as often the case for stability of traveling waves) 
to nonnegativity of the point spectrum: 
$$
\R \sigma_p (L)
\cap \{\R \lambda \ge 0\}=\{0\}.
\eqnlbl{nonpos}
$$
One might conjecture that sufficient conditions for linearized
orbital stability would be \eqnref{nonpos} augmented with the requirement
that $\lambda=0$ be an $\ell$-fold, semi-simple eigenvalue of $L$, 
where $\ell$ is the dimension of the tangent stationary manifold.
However, the lack of spectral separation
associated with accumulation of the essential spectrum of $L$ 
at the stationary eigenvalue $\lambda=0$ prevents us from
drawing immediate conclusions in this direction.
Indeed, it is not a priori clear how one should properly define eigenspace
at a point of accumulation of essential spectrum.

As described in [ZH], such technical difficulties may be conveniently
resolved using an extended, or ``effective'' spectrum based on the
{\it Evans function} of [E.1--4,AGJ,GZ], defined precisely in 
Section 3, below.  An analytic function playing a role for
differential operators analogous to that played for finite-dimensional
operators by the characteristic polynomial, the Evans function
$D(\lambda)$ has the property that, on the resolvent set of $L$,
its zeroes agree (in both location and multiplicity) with the
eigenvalues of $L$ [AGJ,GJ.1--2].
Taking this correspondence as the definition of effective point spectrum, 
we can thus extend the notion of eigenvalue past the essential spectrum
boundary of $L$, via analytic extension of the Evans function [GZ,KS].
The associated effective eigenprojection may then be defined
in terms of the induced meromorphic extension of the resolvent kernel,
via calculus of residues [ZH]. 
For details, see Section 5.4, below.

The first main result of this paper, generalizing the corresponding
result established for viscous, strictly parabolic, shocks in [ZH], is:

\proclaim{Theorem \thmlbl{D}}
Under assumptions (H0)--(H3), 
shocks of ``pure'' hyperbolic type are $L^1\cap L^p\to L^p$
linearly orbitally stable for $p>1$ if and only if
$L$ has precisely $\ell$ effective eigenvalues in $\{\R \lambda \ge 0\}$
(necessarily at $\lambda=0$), or, equivalently,
there holds the Evans function condition:
\medskip
($\CalD$)\qquad
$D(\cdot)$ has precisely $\ell$ zeroes in $\{\R \lambda \ge 0\}$
(necessarily at $\lambda=0$).
\endproclaim

As pointed out in [ZH], the apparently simple criterion
($\CalD$) in fact hides considerable subtlety in behavior, 
in particular in the critical case that (D1) holds, but ($\CalD$) does not:
see [ZH], section 12 for examples of possible behaviors in this case.

Evidently, ($\Cal D$) is equivalent to the pair of conditions:
\medskip
($\CalD$1) \qquad $D(\cdot)$ has no zeroes in 
$\{\R \lambda \ge 0\}\setminus \{0\}$.

($\CalD$2)\quad \qquad \, \, \, (ii) $(d/d\lambda)^\ell D(0)\ne 0$.

\medskip
\noindent The first condition may be recognized as (D1), above.
The second has been shown in [Z.4], Appendix A1 to be equivalent
to
\medskip
(D2)\qquad $\bar \Delta :=\gamma \Delta \ne 0$,
\medskip
\noindent where $\gamma$ is a coefficient measuring transversality of the
connecting manifold between $(u_\pm,v_\pm)$ (i.e., $\gamma\ne0$
equivalent to transversality of the intersection of unstable/stable
manifolds at $(u_-,v_-)$/$(u_+,v_+)$), and $\Delta$ is the one-dimensional
``hyperbolic stability coefficient'' defined in [ZS,Z.4], 
appropriate to the type (Lax, undercompressive, or overcompressive) 
of the wave.  
For example, in the Lax case, $\Delta$ is just the Liu-Majda determinant
$$
\Delta:=\det (r_1^{*-},\dots,r_{n-i_-}^{*-},r_{i_++1}^{*+},\dots, r_n^{*+},[u]),
\eqnlbl{liumajda}
$$
for which $\Delta\ne 0$ is
equivalent to hyperbolic stability of the corresponding ideal
shock of the equilibrium system \eqnref{hyp} [M.1--3];
here $r_j^{*\pm}$ are the eigenvectors of $df^*(u_\pm)$, 
in order of increasing eigenvalue, and $i_\pm$ are as defined above 
(Section 1.1).
For the generalizations of \eqnref{liumajda} in the undercompressive and
overcompressive cases, and their interpretations with respect to
hyperbolic theory, see [ZS] or [Z.4]. (Note: the overcompressive version is 
cited also in (6.54)--(6.55) below).  
In the physical literature,
transversality, $\gamma\ne 0$, is referred to as {\it structural}
stability (i.e., stability of inner shock structure), and evolutionary
(hyperbolic) stability $\Delta \ne 0$ as {\it dynamical} stability [BE].

\proclaim{Corollary \thmlbl{D12}}
Under assumptions (H0)--(H3), 
linear $L^1\cap L^p\to L^p$ orbital stability of pure shocks, for $p>1$,
is equivalent to (D1)--(D2), or equivalently:
strong spectral stability plus classical structural and dynamical stability.
\endproclaim

{\bf Remark \thmlbl{nec}.}  The necessary conditions obtained in
[Go,Z.1,Z.4] were {\it signed conditions}
$$
\gamma \Delta \ge 0,
\eqnlbl{sign}
$$
with appropriate normalization of $\gamma$, $\Delta$.
For $\gamma \Delta\ne 0$, the sign of $\gamma \Delta$ 
is associated with a {\it stability
index} determining the parity of the number of unstable eigenvalues
$\R \lambda >0$ of $L$, hence yields information intermediate 
to that given in (D1) and (D2).
The advantage of \eqnref{sign} is that in the extreme, Lax shock case,
the sign of the quantity on the lefthand side is explicitly 
computable as the product of \eqnref{liumajda} and a Wronskian determinant
evaluated entirely at one infinity: {two linear-algebraic quantities}.
For history/applications/further discussion of the stability index, 
see e.g. [GZ,BSZ,ZS,Z.4].

Note, with the result of Theorem \thmref{D}, that \eqnref{sign}
is strengthened to
$$
\gamma \Delta > 0.
\eqnlbl{strictsign}
$$
\medskip

{\bf Remark \thmlbl{scales}.} 
The twin conditions (D1)--(D2) reflect
multiple scales present in the problem, (D1) corresponding to high-frequency
and (D2) to low-frequency scales.
As discussed in [Z.4], introduction, this dichotomy is natural
from the point of view of formal, matched asymptotic expansion.
For a generalization of (D1)--(D2) to multi-dimensions, see [ZS,Z.4].
\medskip

{\bf Remark \thmlbl{strongspec}.}
In the case discussed in [YoZ], of weak shock profiles arising in regions 
where the equilibrium system is strictly hyperbolic, profiles are necessarily
of Lax type, and transversal, $\gamma\ne0$;  moreover,
the Liu-Majda condition $\Delta\ne 0$ is automatically satisfied, from
standard hyperbolic considerations [M.1--3].
Thus, {\it linearized stability in this case is equivalent to strong
spectral stability, (D1).}
As remarked earlier, strong spectral stability is equivalent to
zero mass stability, so this is indeed a substantial reduction of
the problem.  For example, zero-mass stability has been shown in the case
of weak Jin--Xin profiles [Liu,Hu], whereas linearized stability has not.
The verification of (D1) for weak relaxation profiles 
in this general setting we regard as an extremely interesting
open problem in the one-dimensional theory.
\medskip

{\bf 1.3. Pointwise Green's Function Bounds.}  
Theorem \thmref{D} is obtained as a consequence of
detailed, pointwise bounds on the Green's function 
$G(x,t;y)$ of the linearized evolution equations 
\eqnref{linearized} (more properly speaking, a distribution),
which we now describe. 
For simplicity, we first restrict to the case that $(df, dg)^t$
is strictly hyperbolic, indicating afterward
the extension to the general case.

Let $a_j(x)$, $j=1,\dots (n+r)$ denote the eigenvalues of
$A(x)=(df,dg)^t(\bar u(x))$, and $l_j(x)$ and $r_j(x)$ 
be smooth families of associated left and right eigenvectors, respectively,
normalized so that $l_j^t r_k=\delta^j_k$. These are well-defined, 
by strict hyperbolicity.
Eigenvalues $a_j(x)$, and eigenvectors $l_j$, $r_j$
correspond to hyperbolic characteristic speeds, and modes of propagation 
of the relaxation model \eqnref{general}.
Likewise, let $a^{*\pm}_j$, $j=1,\dots,n$ 
denote the eigenvalues of $df^*(u_\pm)$,
and $l^{+\pm}_j$, $r_j^{*\pm}$ associated left and right eigenvector,
corresponding to hyperbolic characteristic speeds and modes of propagation
of equilibrium system \eqnref{hyp} at the endpoints $u_\pm$ of the
relaxation profile.

Define local, {\it scalar dissipation rates}
$$
\eta_j(x):= -l_j^t Q r_j (x), \quad j=1,\dots,n+r, 
\eqnlbl{eta}
$$
and asymptotic, {\it scalar diffusion rates}
$$
\beta_j^{*\pm}:= \left(l_j^{*t} B^* r_j^{*}\right)_\pm, \quad j=1,\dots, n,
\eqnlbl{beta}
$$
where $B^*_\pm:=B^*(u_\pm)$ is the effective diffusion
defined in \eqnref{B}, predicted by formal, Chapman--Enskog expansion.
As described in Appendix A2, these quantities arise in a natural 
way, through Taylor expansion of the (frozen-coefficient) Fourier symbol 
$$
-i\xi A(x)+Q(x)
\eqnlbl{frozenx}
$$
of the linearized operator $L$ about $\xi=\infty$ and $\xi=0$, respectively.
As a consequence of dissipativity, (H3), we have that
$$
\eta_j^\pm:=\eta_j(\pm \infty)>0, \quad
j=1,\dots,n+r,
\eqnlbl{goodeta}
$$
and
$$
\beta_j^{*\pm}>0, \quad 
j=1,\dots,n.
\eqnlbl{goodbeta}
$$
Let us define also $a_j\pm:=a_j(\pm\infty)$.

The important quantities $\eta_j$, $\beta_j^*$ were identified by
Zeng [Z2] in her study of decay to {\it constant} 
(necessarily equilibrium) {\it solutions} 
$(\bar u, \bar v) \equiv (u_\pm,v_\pm)$ 
of relaxation systems, corresponding at the linearized level
to the study of the limiting equations 
$$
U_t=L_\pm U:= -A_\pm U_x+ Q_\pm U
\eqnlbl{limiting}
$$
as $x\to\pm \infty$ of the linearized evolution equations \eqnref{linearized}.
Using Fourier Transform techniques, Zeng established the following 
sharp, pointwise bounds on the Green's function of \eqnref{limiting},
stated for reference in our own somewhat different notation:

\proclaim{Proposition \thmlbl{zeng}(Ze.2)}
Assuming (H0)--(H3), 
plus strict hyperbolicity of $(df, dg)^t$,
the Green's function $G(x,t;y)$ associated with
the linearized, constant-coefficient evolution equations 
\eqnref{limiting} may be decomposed as
$$
G(x,t;y)= H + S + R,
\eqnlbl{yannidecomp}
$$
where 
$$
H(x,t;y):=\sum_{j=1}^{n+r} e^{-\eta_j^\pm t}\delta(x-y-a_j^\pm t)
r_j^{\pm } l_j^{\pm t}
$$
and
$$
S(x,t;y):= \sum_{j=1}^{n} \big(4\pi \beta_j^{*\pm}(1+t)\big)^{-1} 
e^{|x-y-a_j^{*\pm}t|^2/4\beta_j^{*\pm}(1+t)}
R_j^{*\pm} L_j^{*\pm t},
$$
denote hyperbolic and scattering terms, respectively, 
and $R$ denotes a faster decaying residual.
Here, $a_j^\pm$, $a_j^{*\pm}$, $l_j^\pm$, $l_j^{*\pm }$,
$r_j^\pm$, $r_j^{*\pm}$, $\eta_j^\pm$ and $\beta_j^{*\pm}$ are
as defined just above, and $R_j^{*\pm}$, $L_j^{*\pm}$ by
$$
R_j^{*\pm}:=\pmatrix r_j^{*\pm} \\ -q_v^{-1}q_u r_j^{*\pm})
\endpmatrix, \qquad
L_j^{*\pm}:= \pmatrix l_j^{*\pm t} \\0 \endpmatrix.
\eqnlbl{RL}
$$
\endproclaim

Note, as suggested by the derivation of constants $\eta_j^\pm$,
$\beta_j^\pm$ in Appendix A2, that short-time ($\sim$ large-frequency) behavior
is dominated by hyperbolic term $H$, while large-time
($\sim$ small-frequency) behavior is dominated by scattering term $S$.
Moreover, the latter term is precisely that
identified in [LZe] as describing large-time behavior
for the Chapman--Enskog approximation \eqnref{parab},
in agreement with the general philosophy of [L.2].

Our second main result, and the principal result of this paper,
is the following, nonconstant-coefficient analog of Proposition \thmref{zeng}:

\proclaim{Theorem \thmlbl{greenbounds}}
Assuming (H0)--(H4),
plus strict hyperbolicity of $(df, dg)^t$, 
and the spectral stability criterion ($\CalD$)
(equivalently, (D1)--(D2)),
the Green's function $G(x,t;y)$ associated with
the linearized evolution equations \eqnref{linearized} 
may in the Lax or overcompressive case be decomposed as
$$
G(x,t;y)= H + E+  S + R,
\eqnlbl{ourdecomp}
$$
where, for $y\le 0$:
$$
H(x,t;y):=
\sum_{j=1}^{n+r} e^{-\bar \eta_j t}\delta_{x-\bar a_j t}(-y)
r_j(x)l_j^t(y),
\eqnlbl{H}
$$
$$
\aligned
E&(x,t;y):=\\
&\sum_{a_k^- > 0}
[c^{j,0}_{k,-}]
\frac{\partial}{\partial \delta_j}
\pmatrix
\bar u^\delta(x) 
\\
\bar v^\delta(x)
\endpmatrix
L_k^{*-t}
\left(\errfn\left(\frac{y+a_k^{*-}t}{\sqrt{4\beta_k^{*-}t}}\right)
-\errfn \left(\frac{y-a_k^{*-}t}{\sqrt{4\beta_k^{*-}t}}\right)\right),
\endaligned
\eqnlbl{E}
$$
and
$$
\aligned
S(x,t;y)&:=
\chi_{\{t\ge 1\}} 
\sum_{a_k^{*-}<0}R_k^{*-}  {L_k^{*-}}^t
(4\pi \beta_k^-t)^{-1/2} e^{-(x-y-a_k^{*-}t)^2 / 4\beta_k^{*-}t} 
\\
&+ 
\chi_{\{t\ge 1\}} 
\sum_{a_k^{*-} > 0} R_k^{*-}  {L_k^{*-}}^t
(4\pi \beta_k^{*-}t)^{-1/2} e^{-(x-y-a_k^{*-}t)^2 / 4\beta_k^{*-}t}
\left({e^x \over e^x+e^{-x}}\right)\\
&+ 
\chi_{\{t\ge 1\}}
\sum_{a_k^{*-} > 0, \,  a_j^{*-} < 0} 
[c^{j,-}_{k,-}]R_j^{*-}  {L_k^{*-}}^t
(4\pi \bar\beta_{jk}^{*-} t)^{-1/2} e^{-(x-z_{jk}^{*-})^2 / 
4\bar\beta_{jk}^{*-} t} 
\left({e^{ -x} \over e^x+e^{-x}}\right),\\
&+ 
\chi_{\{t\ge 1\}}
\sum_{a_k^{*-} > 0, \,  a_j^{*+} > 0} 
[c^{j,+}_{k,-}]R_j^{*+}  {L_k^{*-}}^t
(4\pi \bar\beta_{jk}^{*+} t)^{-1/2} e^{-(x-z_{jk}^{*+})^2 / 
4\bar\beta_{jk}^{*+} t} 
\left({e^{ x} \over e^x+e^{-x}}\right)\\
\endaligned
\eqnlbl{S}
$$
denote hyperbolic, excited, and scattering terms, respectively, 
and $R$ denotes a faster decaying residual (rates given in
Proposition 6.1, below).
Symmetric relations hold for $y\ge 0$.

Here, the averaged dissipation and convection rates
$\bar \eta_j(y,t)$ and $ \bar a_j(y,t)$ in \eqnref{H} denote the 
time-averages over $[0,t]$ 
of $\eta_j(x)$ and $a_j(x)$, respectively,
along characteristic paths $z_j=z_j(y,t)$ defined by
$$
dz_j/dt= a_j(z_j), \quad z_j(0)=y,
\eqnlbl{char}
$$
while, in \eqnref{S},
$$
z_{jk}^{*\pm(y,t)}:=a_j^{*\pm}\left(t-\frac{|y|}{|a_k^{*-}|}\right)
\eqnlbl{zjk}
$$
and
$$
\bar \beta^{*\pm}_{jk}(x,t;y):= 
\frac{|x^\pm|}{|a_j^{*\pm} t|} \beta_j^{*\pm}
+
\frac{|y|}{|a_k^{*-} t|} 
\left( \frac{a_j^\pm}{a_k^{*-}}\right)^2 \beta_k^{*-},
\eqnlbl{barbeta}
$$
represent, respectively, approximate scattered characteristic
paths and the time-averaged diffusion rates along those paths.
In all equations, $a_j$, $a_j^{*\pm}$, $l_j$, $L_j^{*\pm }$,
$r_j$, $R_j^{*\pm}$, $\eta_j$ and $\beta_j^{*\pm}$ are
as defined just above, 
and scattering coefficients $[c_{k,-}^{j,i}]$, $i=-,0,+$, are constants, 
uniquely determined by
$$
\sum_{a_j^{*-} < 0} [c_{k, \, -}^{j, \, -}]r_j^{*-} +
\sum_{a_j^{*+} > 0} [c_{k, \, -}^{j, \, +}]r_j^{*+} +
\sum_{j=1}^\ell [c_{k,-}^{j,0}]
\int_{-\infty}^{+\infty}(\partial/\partial \delta_j)
\bar u^\delta(s) ds
= r_k^{*-}
\eqnlbl{scattering}
$$
for each $k=1,\dots n$, and satisfying
$$
\sum_{a_k^->0} [c_{k,-}^{j,0}] L_k^{*-}
=
 \sum_{a_k^+<0} [c_{k,+}^{j,0}] L_k^{*+} 
=\pi_j
\eqnlbl{pi}
$$
for each $j=1,\dots,\ell$, where the constant vector $\pi_j$
is the left zero effective eigenfunction associated
with the right eigenfunction $(\partial/\partial \delta_j)(\bar u^\delta,
\bar v^\delta)$. 
Similar, but more complicated formulae hold in the undercompressive
case (see Remark 6.9, below).

In addition, $G$ vanishes identically outside the 
hyperbolic domain of influence
$$
z_1(y,t)\le x \le z_{n+r}(y,t),
\eqnlbl{finiteprop}
$$
$z_j$ defined as above.
\endproclaim

Theorem \thmref{greenbounds} will be established by Laplace Transform 
(i.e., semigroup) techniques generalizing the Fourier Transform 
approach of [Ze.1--2,LZe].
The extension of Fourier Transform techniques to the variable-coefficient 
case is of course a central issue in mathematical physics and analysis, 
associated with such topics as semigroup, Sturm--Liouville, 
and scattering theory [Ti,RS,LP]; however, many unresolved issues
remain, particularly in connection with higher order
or nonnormal differential operators.
From this, general theoretical point of view,
the results of Theorem \thmref{greenbounds}, and the methods
used to obtain them, seem of wider interest,
independent of their application to stability of shock profiles.

{\bf Remarks.}
Just as in the constant--coefficient case, short-time/large-frequency 
behavior is dominated by hyperbolic term $H$, while large-time/small-frequency
behavior is dominated by terms $E+S$ agreeing with those predicted
for the partially parabolic Chapmann--Enskog approximation \eqnref{parab},
or, equivalently, the ``effective,'' strictly parabolic, system 
defined by Kawashima [Kaw].

Terms $E$ and $S$ are essentially those derived for strictly parabolic
systems in [ZH,Z.2--3], combined with projection onto the equilibrium
manifold.
As discussed in those references, they may be understood
by the caricature of Gaussian signals propagating for $x\gtrless 0$
according to the limiting, constant-coefficient equations at $\pm\infty$
(of the associated effective, strictly parabolic system),
interacting with an ``ideal'' shock layer at $x=0$ according to the
rules encoded by scattering coefficients $[c_{k,\pm}^{j,i}]$,
$i=-,0,+$: That is, a signal initiating as a delta-function at $y<0$
first decomposes into $n$ Gaussian signals traveling in each of the equilibrium
characteristic modes $r_k^{*-}$ with approximately equilibrium 
characteristic speed $a_k^{*-}$; each of these signals then travels until
it reaches the shock layer $x=0$, whereupon it is reflected/transmitted
in each of the outgoing characteristic modes $r_j^{*\pm}$, thereafter
moving with approximate equilibrium characteristic speed 
$a_j^{*\pm}\gtrless 0$ along path $z_{jk}^{*\pm}$; 
at the same time, the incoming signal excites stationary modes,
or ``bound states,''
$(\partial/\partial \delta_j)(\bar u^\delta ,\bar v^\delta )$.
Incoming and outgoing Gaussian signals comprise the scattering term $S$,
excited stationary modes the excited term $E$.

The ``parabolic'' terms $E+S$ are considerably
more complicated in the variable-coefficient than in the 
constant--coefficient case;
by contrast, the hyperbolic term $H$,
since localized, has an appealingly simple generalization 
to the variable-coefficient case, consisting in the
``frozen-coefficient'' evolution described by \eqnref{char}.

\medskip

{\bf The nonstrictly hyperbolic case.}  The result of Theorem
\thmref{greenbounds} holds also in the case that $(df,dg)^t)$
has multiple characteristics, with the following
modification.  If $a_j$ has multiplicity $m_j$, $j=1,\dots, J$,
then $l_j(x)$, $r_j(x)$
are defined instead as $(n+r)\times m_j$ matrices of left and
right eigenvectors of $(df,dg)^t(\bar u, \bar v)(x)$, normalized
so that $l_j^tr_j\equiv I_{m_j}$.
In addition to this ``static'' normalization, 
we make the ``dynamic'' normalization 
$$
l_j^t(r_j)'(x)\equiv 0.
\eqnlbl{dnorm}
$$
This is always possible, as shown in Lemma 4.9 below; for
discrete kinetic models, it may be accomplished simply by taking
$l_j$, $r_j$ constant.
Then, defining the $m_j\times m_j$ matrix
$$
\eta_j(x):= -l_j^t Q r_j(x),
\eqnlbl{multeta}
$$
and augmenting characteristic equation \eqnref{char} 
with the {\it dissipative flow} 
$$
d\zeta_j/dt= -\eta_j(z_j)\zeta_j, \quad \zeta_j(y)=I_{m_j},
\eqnlbl{diss}
$$
governing the evolution of dissipation matrix
$\zeta_j=\zeta_j(y,t)\in \BbbR^{m_j\times m_j}$, we have:

\proclaim{Theorem \thmlbl{auxgreenbounds}}
Assuming (H0)--(H4) and the spectral stability criterion ($\CalD$)
(equivalently, (D1)--(D2)),
decomposition \eqnref{ourdecomp} holds also in the
general (nonstrictly hyperbolic) case, with \eqnref{H} replaced by 
$$
\aligned
H(x,t;y)&:=
\sum_{j=1}^{J} 
r_j(x) \zeta_j(y,t)
\delta_{x-\bar a_j t}(-y)
l_j^t(y)\\
&=
\sum_{j=1}^{J} r_j(x) \Cal{O}(e^{-\eta_0 t})
\delta_{x-\bar a_j t}(-y)
l_j^t(y),
\endaligned
\eqnlbl{multH}
$$
for some uniform $\eta_0>0$, with $l_j$, $r_j$, and $\zeta_j$ as
defined just above.
\endproclaim

That is, the trivial, {\it scalar flow} 
$ \zeta(t)=e^{\int_0^t -\eta_j(z(s)) ds} =e^{\bar -\eta_j t} $
must be replaced in the case of corresponding characteristics
by genuine, system dynamics \eqnref{diss}.
In the strictly hyperbolic case, \eqnref{multH} reduces to
\eqnref{H}, as it must.
Likewise, condition \eqnref{goodeta} must be replaced by
$$
\R \sigma (\eta_j^\pm)>0, \quad
j=1,\dots,J.
\eqnlbl{nsgoodeta}
$$

{\bf Remark \thmlbl{multeq}.} Similarly, we may relax the
requirement of strict hyperbolicity of $df^*(u_\pm)$, under
a mild additional assumption.
If $a_j^{*\pm}$ has multiplicity $m^{*\pm}_j$, then, likewise,
we define $l_j^{*\pm}$ and $r_j^{*\pm}$ as $n\times m_j^{*\pm}$
blocks, and 
$$
\beta_j^{*\pm}:= l_j^{*\pm t} B^{*\pm}r_j^{*\pm}.
\eqnlbl{multbetaj}
$$
Make the additional assumption that each $\beta_j$ be
{\it diagonalizable}, without loss of generality diagonal.
Then, \eqnref{ourdecomp} holds with \eqnref{E}--\eqnref{S} 
replaced by matrix versions as in \eqnref{multH} above.
Diagonalizability of $\beta_j^{*\pm}$ follows, for example,
from simultaneous symmetrizability of $A$, $Q$, which holds
for most physical systems of interest; this is satisfied, in 
particular, for our main example, the Jin--Xin system [Yo.4].
\medskip

{\bf 1.4. Nonlinear Stability.}  
Finally, we discuss our results on nonlinear stability.
In this part of the analysis, we restrict our attention to weak Lax
shocks of {\it simultaneously symmetrizable models} 
(defined precisely in Section 7.2, below),
or strong Lax shocks of {\it discrete kinetic models} 
(defined as above as the subclass of models for which $f$, $g$ 
are linear, constant-coefficient).
Simultaneously symmetrizable models of the dissipative type considered here
admit favorable energy estimates [Kaw,Ze.2] allowing us to bound
nonlinear derivative terms and close our nonlinear iteration;
however, these energy estimates appear to require the weak shock assumption.
Discrete kinetic models, on the other hand, do not have any such 
derivative terms to begin with, since nonlinearities enter only through 
the collision term $q$.

With these restrictions, we obtain our third main result:

\proclaim{Theorem \thmlbl{nonlin}}
Let $(\bar u(x),\bar v(x))$ be a Lax type relaxation profile
of a general relaxation model \eqnref{general}, satisfying
(H0)--(H3) and the spectral stability criterion ($\CalD$) 
(equivalently, (D1)--(D2)).
If, also, $(\bar u,\bar v)$ is sufficiently weak, in the sense
that $|(\bar u',\bar v')|_{L^\infty}$ is sufficiently small
relative to the parameters of \eqnref{general},
and \eqnref{general} is simultaneously symmetrizable (with
coefficients $q\in C^3$ now also in the case of discrete kinetic models),
then, provided that initial perturbation $(u_0,v_0)-(\bar u,\bar v)$
is bounded by $\zeta_0$ in $L_1$ and $H^2$, for $\zeta_0$
sufficiently small,
the solution $(u,v)(x,t)$ of \eqnref{general} with
initial data $(u_0,v_0)$ satisfies
$$
|(u,v)(x,t)-(\bar u,\bar v)(x-\delta(t))|_{L^p}\le C \zeta_0
(1+t)^{-\frac{1}{2}(1-1/p)},
\eqnlbl{3.31}
$$
for all $2\le p\le \infty$, for some $\delta(t)$ satisfying
$$
|\dot \delta (t)|\le C \zeta_0 (1+t)^{-\frac{1}{2}}
\eqnlbl{3.32}
$$
and
$$
|\delta(t)|\le C \zeta_0.
\eqnlbl{3.33}
$$
In particular, $(\bar u,\bar v)$ is nonlinearly orbitally stable from
$L^1\cap H^2$ to $L^p$, for all $p\ge 2$.
%
%
For weak shocks of simultaneously symmetrizable discrete kinetic models
($q\in C^2$) and data merely small in $L^1\cap H^1$, 
we obtain \eqnref{3.31} for $p=2$, but only boundedness 
$|(u,v)(x,t)-(\bar u,\bar v)(x-\delta(t))|_{L^\infty}\le C \zeta_0$,
for $p=\infty$, and interpolated rates for $p$ between $2$ and $\infty$.

For discrete kinetic models,
not necessarily simultaneously symmetrizable, with $q\in C^3$, 
and data small in $W^{3,1}\cap W^{3,\infty}$, 
we obtain the rates \eqnref{3.31} for all $1\le p\le \infty$,
for shocks of arbitrary strength.
\endproclaim

Theorem \thmref{nonlin} is obtained using a modified version
of the argument used in [Z.2] to treat the strictly parabolic case,
new complications arising from the need to control 
contributions of the
singular component $H$ of the Green's function.
Similar results should be possible for nonclassical,
over- and under-compressive shocks by more detailed, pointwise
analyses in the spirit of [L.1,ZH] and [Z.5], respectively.

{\bf Remark \thmlbl{stronglim}.}
As pointed out by W.-A. Yong [Yo.4], the condition of 
simultaneously symmetrizability 
is satisfied for most systems of interest in applications: for example,
the discrete kinetic models of Platkowski--Illner [PI]; 
the numerical approximation schemes of Jin--Xin [JX] and
their generalizations by Natilini [N]; the BGK models of
Bouchut [B]; and, most notably, perhaps, for {\it all} of the 
extended thermodynamic models in the moment closure heirarchies 
of Levermore or Dreyer [Lev,Dre]; further examples may be found
among the general relaxation models described in [CLL,Yo.1--3].
\medskip

{\bf Remark \thmlbl{stronglim}.}
It is an interesting question whether the limitation to
weak shocks for general relaxation models in Theorem \thmref{nonlin}
is only a technical artifact of our analysis, or a genuine difference 
in nonlinear behavior between quasilinear and semilinear models.
We point out that this limitation arises very late in our argument,
and at the (nonlinear) level of energy estimates only.
On the other hand, discrete kinetic models are in some
sense closer to original kinetic models such as Boltzmann
equations than are quasilinear models based on moment closure 
assumptions,  hence it is conceivable that they 
truly give a more faithful approximation of behavior in the 
large-amplitude regime.

\medskip

{\bf Remark \thmlbl{weakJX}.}
Theorem \thmref{nonlin} in particular shows that {\it strong spectral
stability implies nonlinear stability}, for weak shocks of general
simultaneously symmetrizable models, or for strong shocks of discrete kinetic models.
As noted above, strong spectral stability follows for
weak relaxation profiles of the Jin--Xin and Broadwell models
from the results of [KM] and [Liu,Hu], respectively, thus yielding
complete nonlinear stability results for these two cases. 
A very interesting open problem is spectral stability of weak 
shocks of general dissipative relaxation systems, 
as identified by Zeng/Kawashima [Kaw,Ze.2], or, likewise, spectral
stability of strong shocks of special discrete kinetic models.

\medskip
{\bf Remark \thmlbl{deltadecay}.}
By a more refined, pointwise analysis, it should be possible to show
that $\delta(t)$ approaches a limit as $t\to \infty$, which may then be
explicitly determined using conservation of mass; however, it is clear
that there can be no uniform rate of convergence for general $L^1\cap L^\infty$
perturbations, since interaction with the shock layer (and thus formation
of the final asymptotic state) may be delayed for
arbitrarily long time by moving the perturbation to spatial infinity.
Likewise, it should be possible to show that the perturbed shock
converges in $L^1$ to a translated profile superimposed with
diffusion waves in the far field: see, for example, related results
obtained in [L.1,ZH] for localized data in the strictly parabolic case.

As pointed out in [Z.2],
the latter results, expressed in terms of the
perturbation $U(x,t)-\bar U (x-\delta(+\infty))$ from
the time-asymptotic translate $\bar U(x-\delta(+\infty))$
rather than the instantaneous translate $\bar U(x-\delta(t))$,
contain an additional error term 
$$
\bar U(x-\delta(+\infty))
-\bar U(x-\delta(t))
\sim
(\delta(t)-\delta(+\infty))\bar U'(x)
\eqnlbl{coupled}
$$
not found here, the so-called ``coupled viscous wave'' of [SX.1,S].
As observed from a different point of view in [S], the shape $\bar U'$ 
of this coupled wave depends on a ``dynamical'' effective viscosity 
induced by the traveling wave ODE, 
which is in general different from either of the effective viscosities 
$B^{*}(U_\pm)$ predicted by Chapman--Enskog expansion at the endstates 
$U_\pm$ of the profile: 
$a^2-s^2$, for example, vs. $a^2-dh^2(u_\pm)$, in the case of the Jin--Xin
model \eqnref{JX}.
As discussed in [S], this is related to the modified Chapman--Enskog
expansion of [L.2], in which derivatives are assumed small in directions 
of propagation of the solution, i.e., shock speed $s$ in the shock layer
vs. equilibrium characteristic speeds $a_j^{*\pm}$ in the far fields:
in the case of the Jin--Xin model, $s$ vs. $\sigma(dh(u_\pm))$. 
For the class of initial perturbations $U_0-\bar U\sim \zeta_0 (1+|x|)^{-3/2}$
considered in [L.1,ZH], it can be shown that 
$
(\delta(t)-\delta(+\infty))\sim \zeta_0 (1+t)^{-1/2},
$
recovering the results of those papers.

\medskip
{\bf Plan of the paper.}
In Section 2, we establish the basic spectral representation
of $G$ via inverse Laplace Transform of the resolvent kernel.
In Section 3, we make a careful construction of the resolvent
kernel; using this representation, we make careful resolvent
kernel estimates in Sections 4 and 5, in high- and low-frequency
regimes, respectively.
Combining these results, we derive in Section 6 the Green's function
bounds that comprise the main results of the paper,
establishing Theorems \thmref{greenbounds}--\thmref{auxgreenbounds}
and a partial generalization yielding necessity of ($\Cal{D}$)
for linearized stability.
Finally, in Section 7, we establish
sufficiency of ($\Cal{D}$) for linearized stability,
completing the proof of Theorem \thmref{D},
and carry out the nonlinear iteration establishing
nonlinear stability, Theorem \thmref{nonlin}.
For completeness, we carry out in Appendices A1--A3 various auxiliary
calculations/lemmas used throughout the paper.

\bigskip
{\bf Note}:  
Subsequent to the completion of our analysis, 
F. Rousset [R] has announced the partial result of zero-mass stability 
of Jin--Xin relaxation shocks 
satisfying the spectral stability condition of [G,Z.1,Z.4]
(the same one assumed here): that is, the extension
to strong shocks of the zero-mass result of H. Liu [Liu].
However, his analysis appears to be incomplete in the 
high-frequency regime;
specifically, the Gap Lemma of [GZ,ZH]
is incorrectly cited at a key step in the proof of his
high-frequency resolvent estimates,
see Remark 8.5, Appendix A3.1, below.
This gap is remedied by our own, more precise bounds.

%
\bigskip

\newsection {\bf Section 2. The Spectral Resolution Formula.}
\sectionnumber=2 
\theoremnumber=0
\equationnumber=0
\smallskip
\TagsOnLeft

The starting point for our analysis is the spectral resolution,
or inverse Laplace Transform formula for the Green's function
$G(x,t;y)$ associated with the linearized evolution equations
\eqnref{linearized}, defined by:

(i) $(\partial_t -L_x)G=0$ in the distributional sense, for all $t>0$;
and, 

(ii) $G(x, t;y)\rightharpoondown \delta(x-y)$ as $t\to 0$.

\noindent Here, $G$ is to be interpreted in (i) as a distribution
in the joint variables $(x,y,t)$, and in (ii) as a distribution in $(x,y)$, 
continuously parametrized by $t$.  We shall see that $G$, uniquely
so defined, by uniqueness of weak solutions of \eqnref{linearized}
within the class of test function initial data, in fact lies
in the space of measures and not functions, so should 
properly be called a Green's distribution rather than a Green's function.
Following the convention of [LZe,Ze.2], we shall use the two
terms interchangeably. 

The Green's distribution may be conveniently constructed using
standard semigroup theory, considering the linearized operator
$L$ as a closed, densely defined operator on $L^2$, with domain 
$\CalH:=\{U:\, U\in L^2, \, AU\in H^1\}$.

\proclaim{Lemma \thmlbl{garding}}
If $A(x)$ is smooth and symmetrizable, in the sense that there exists
a smooth, invertible $S(x)$ such that $SAS^{-1}$ is symmetric,
and $Q(x)$ is bounded, 
then operator $L:\, LU:= -(AU)_x+QU$ satisfies the bound
$$
|U|_{L^2}\le |\lambda-\lambda_*|^{-1}|(L-\lambda) U|_{L^2}
\eqnlbl{garding}
$$
for all $U\in \CalH$, and real $\lambda$ greater than some value $\lambda_*$.
\endproclaim

{\bf Proof.}  Without loss of generality, we may take $A$ symmetric,
since we may achieve this by the coordinate change $U\to SU$, without
affecting boundedness of $Q$.
Setting
$$
(L-\lambda) \phi=f,
\eqnlbl{semi2}
$$
for test functions $\phi \in C^\infty_0$,
and performing an elementary energy estimate,
we obtain
$$
\langle \phi',A\phi\rangle + \langle \phi,Q\phi\rangle = \langle \phi,f \rangle
+ \lambda |\phi|^2,
\eqnlbl{energy}
$$
where $\langle \cdot,\cdot \rangle$ and $|\cdot|$ refer to 
complex $L^2$ inner product and norm, respectively.
Taking real parts in \eqnref{energy}, we have the G\"arding type estimate
$$
-\frac{1}{2}\langle \phi,A'\phi\rangle +
Re \langle \phi,Q\phi\rangle = Re \langle \phi,f \rangle
+ Re \lambda |\phi|^2,
\eqnlbl{semi3}
$$
or, rearranging:
$$
(Re \lambda-\lambda_*) |\phi|^2\le |\phi||f|,
\eqnlbl{semi4}
$$
for some real $\lambda_*>0$ sufficiently large.
But, this gives the claimed bound
$$
|\phi|\le |\lambda-\lambda_*|^{-1}|f|
=|\lambda-\lambda_*|^{-1}|(L-\lambda)\phi|
\eqnlbl{semi5}
$$
for all $\lambda>\lambda_*$ on the real axis, 
for $\phi\in C^\infty_0$ within the class of test functions.
Observing that $C^\infty_0$ is dense within $\CalH$
(as may be seen by the straightforward computation
$|(AU)^\varepsilon - AU^\varepsilon|_{H^1}\to 0$ as $\varepsilon\to 0$,
where $f^\varepsilon := f*\eta^\varepsilon$, and $\eta^\varepsilon\in C^\infty_0$ 
denotes the standard mollifier, 
together with the standard facts that $|(AU)^\varepsilon - AU|_{H^1}\to 0$
and $|U^\varepsilon - U|_{L^2}\to 0$),
we may pass to the limit to obtain the same bound for $\phi\in \CalH$.
\myqed

\proclaim{Corollary \thmlbl{semigroup}}
If the operator $L:\, LU:= -(AU)_x+QU$ of Lemma \thmref{garding}
has in addition the properties that $A$ is nonsingular, 
and $A$ and $Q$ are asymptotically
constant as $x\to \pm \infty$, then $L$
generates a $C^0$ semigroup $e^{Lt}$
on $L^2$, satisfying $|e^{Lt}|_{L^2}\le Ce^{\omega t}$ for some real $\omega$.
\endproclaim

{\bf Proof.}
The bound \eqnref{semi4} applies also to the limiting, constant-coefficient
operators $L_\pm:= -A_\pm \partial_x + Q_\pm$ as $x\to \pm \infty$, whence 
the spectra of these operators is confined to $\R \lambda \le \lambda_*$.
Under the assumption that $A$ is nonsingular, we can apply a
standard result of Henry ([He], Lemma 2, pp. 138--139) to find that
the essential spectrum of $L$ is also confined to this set.
Since \eqnref{garding} precludes point spectrum for real $\lambda>\lambda_*$,
we thus find that all such $\lambda$ belong to the resolvent set $\rho(L)$,
with the resolvent bound
$$
|(L-\lambda)^{-1}|_{L^2}\le |\lambda-\lambda_*|^{-1}.
\eqnlbl{resbound}
$$
But, this is a standard sufficient condition that a closed, 
densely defined operator $L$ generate a $C^0$ semigroup, 
with $|e^{Lt}|_{L^2}\le Ce^{\omega t}$ 
for all real $\omega >\lambda_*$
(see e.g. [Pa], Theorem 5.3, or [Fr,Y]).
\myqed

An immediate consequence (see, e.g. [Pa], Corollary 7.5) is:

\proclaim{ Corollary \thmlbl{inverseLT}}
For $L$, $e^{Lt}$ as in Corollary \thmref{semigroup},
there holds for all real $\eta$ greater than some $\eta_0$
the inverse Laplace Transform (spectral resolution) formula
$$
e^{Lt}\phi=
P.V. \int_{\eta -i\infty}^{\eta+i\infty} e^{ikt}(L-\lambda)^{-1}\phi \, d\lambda,
\eqnlbl{iLT}
$$
for all $\phi \in D(L^2)\subset H^2$, where domain
$D(L^2)$ 
is defined as in [Pa], p. 1. 
\footnote{
Note: From the definition in [Pa], we find that the domain
$D(L)$ of $L$ satisfies $H^1\subset \CalH \subset D(L)$.
Likewise, the domain $D(L^2)$, consisting of
those functions $V$ such that $LV\in D(L)$, satisfies
$H^2\subset D(L^2)$.
}
\endproclaim

{\bf Remark.}  Symmetrizability, in the one-dimensional
context, is essentially
equivalent to hyperbolicity: $A$ has real,
semisimple eigenvalues; in particular, it is implied by
either strict hyperbolicity, or else nonstrict hyperbolicity with
constant multiplicity of eigenvalues, so that all cases considered
here fall into this class. 
Thus, this is no real requirement
other than the usual one of hyperbolicity of the principal
part of the symbol (Note: this may be satisfied even when the 
subcharacteristic condition is violated, or when the 
reduced system is not hyperbolic,
as in applications of the Jin--Xin model [JX] to phase transition problems).
\medskip

With the result of  Corollary \thmref{semigroup},
we can express the Green's distribution simply as
$$
G(x,t;y):=e^{Lt}\delta(x-y),
\eqnlbl{g}
$$
where the evolution $e^{Lt}\phi(x,y)$ of a distribution $\phi(x,y)$
with the property that, for $k \in C^\infty_0$, the distributional pairing
$\langle \phi(\cdot,y),k(y)\rangle_y$ is in $L^2$ for each fixed $x$, 
is defined through its action on $C^\infty_0$ test functions $h(x)k(y)$:
$$
\langle e^{Lt}\phi(x,y),h(x)k(y) \rangle_{x,y}
:=
\langle h(x), e^{Lt} \langle \phi(x,y),k(y) \rangle_{y} \rangle_x,
$$
or
$$
\langle G(x,t;y),h(x)k(y) \rangle_{x,y}
:=
\langle h(x), e^{Lt}k(x) \rangle_x.
\eqnlbl{distg}
$$
The defining properties (i)--(ii) of the Green's distribution
are then consequences of the corresponding properties defining
a $C^0$ semigroup (see, e.g., [Pa], p. 1--4).

Defining the {\it resolvent kernel}, similarly, through
$$
G_\lambda(x,y):= (L-\lambda)^{-1}\delta(x-y),
\eqnlbl{glambda}
$$
where the righthand side is defined again through its action on test functions,
$$
\langle (L-\lambda)^{-1}\phi(x,y),h(x)k(y) \rangle_{x,y}
:=
\langle h(x), (L-\lambda)^{-1} \langle \phi(x,y),k(y) \rangle_{y} \rangle_x,
$$
or
$$
\langle G_\lambda(x,y),\psi(x,y) \rangle_{x,y}
:=
\langle h(x), (L-\lambda)^{-1} k(x)\rangle_x,
$$
we find through \eqnref{iLT}, as a tautology, the fundamental:

\proclaim{Proposition \thmlbl{gres}}
For $L$, $e^{Lt}$ as in Corollary \thmref{semigroup},
and $G$, $G_\lambda$ as in \eqnref{g}--\eqnref{glambda},
there holds the inverse Laplace Transform (spectral resolution) formula
$$
G(x,t;y)=
P.V. \int_{\eta -i\infty}^{\eta+i\infty} e^{ikt} G_\lambda(x,y)\, d\lambda,
\eqnlbl{giLT}
$$
valid for all $\eta$ greater than some $\eta_0$.
\endproclaim

It is this relation that we shall use to compute the time-evolutionary
Green's distribution $G$.  Note that, on the resolvent set $\rho(L)$, 
the resolvent kernel $G_\lambda$ likewise
has an alternative, intrinsic characterization as the unique distribution
satisfying 
$$
(L_x-\lambda)G_\lambda(x,y)=\delta(x,y)
\eqnlbl{dreskernel}
$$
and taking $f\in L^2$ to $\langle G_\lambda(x,y),f(y)\rangle_y \in L^2$.
This is the characterization that we shall use in computating $G_\lambda$.

{\bf Remark \thmlbl{solnformula}.}  Note that a semigroup on even a still more restricted
function class than $L^2$ would have sufficed in this construction, 
since we are constructing objects in the very weak class of distributions.
Later, by explicit computation of $G$, we will verify that $L$
in fact generates a $C^0$ semigroup in any $L^p$.
Note also that \eqnref{distg} implies the expected, standard solution formula
$e^{Lt}f=\langle G(x,t;y),f(y)\rangle_y$, or, formally:
$$
e^{Lt}f=\int G(x,t;y)f(y) dy,
\eqnlbl{stdsoln}
$$
for $f$ in any underlying Banach space on which $e^{Lt}$ is defined,
in this case on $L^p$, $1\le p\le \infty$.
\medskip

{\bf Remark \thmlbl{altpf}.}
For general interest, we point out another, more concrete route
to \eqnref{giLT}, \eqnref{stdsoln}, generalizing the approach
followed in [ZH] for the parabolic case.
Namely, we may observe that, in each finite integral in the approximating
sequence defined by the principal value integral \eqnref{giLT},
we may exchange orders of integration and distributional differentiation,
using Fubini's Theorem together with the fact,
established in the course of our analysis, that $G_\lambda$ 
is uniformly bounded on the contours under consideration
(in fact, we establish the much stronger result 
that $G_\lambda$ decays exponentially in $|x-y|$, with uniform rate),
to obtain
$$
\aligned
(\partial_t - L_x)
\int_{\eta-iK}^{\eta+iK} e^{\lambda t}G_\lambda(x,y)\, d\lambda
&=
\delta(x-y)\int_{\eta-iK}^{\eta+iK} e^{\lambda t}\, d\lambda\\
& \rightharpoondown \delta(x-y)\delta(t)\\
\endaligned
\eqnlbl{formalcalc}
$$
as $K\to \infty$, for all $t\ge 0$. 
So, all that we must verify is: (i') the Principal Value integral \eqnref{giLT} 
in fact converges to some distribution $G(x,t;y)$ as $K\to \infty$; and, 
(ii') $G(x,t;y)$ has a limit as $t\to 0$.
For, distributional limits and distributions {commute} 
(essentially by definition), so that (i')--(ii') together with
\eqnref{formalcalc} imply (i)--(ii) above.
Facts (i')--(ii') will be established by direct calculation
in the course of our analysis, thus verifying formula \eqnref{giLT} 
at the same time that we use it to obtain estimates on $G$.
Note that we have made no reference in this argument 
to the semigroup machinery cited above.
\bigskip

\newsection {\bf Section 3. Construction of the Resolvent kernel.}
\sectionnumber=3 
\theoremnumber=0
\equationnumber=0
\smallskip
\TagsOnLeft

We next derive explicit representation formulae for the resolvent kernel
$G_\lambda$, using the classical construction (see, e.g. [CH]) 
for the Green's function of an ordinary differential operator, 
in terms of solutions of the homogeneous eigenvalue equation
$$
0=(L-\lambda)W=
-(AW)'+ (Q -\lambda)W,
\eqnlbl{evalue}
$$
matched across the singularity $x=y$ by appropriate jump conditions,
in the process obtaining standard decay estimates on the resolvent kernel
(see Proposition 3.4, below)
suitable for analysis of intermediate frequencies $\lambda$.
Here, $A$ and $Q$ are as defined as in \eqnref{linearized}, 
$W\in \BbbR^N:=\BbbR^{n+r}$, and $'$ as usual denotes $d/dx$.
In the details of our treatment, we follow [ZH,Z.4]; see also
related analyses of [AGJ,K.1--2] regarding solution of
the resolvent equation by variation of constants.

{\bf 3.1. Basic construction.}
First, consider the limiting, constant-coefficient eigenvalue equations
$$
0=(L_\pm-\lambda)W:=
-A_\pm W'+ (Q_\pm -\lambda)W.
\eqnlbl{limitevalue}
$$
Define 
$$
\Lambda:=\cup \Lambda_j^\pm,\quad j=1,\dots,n+r,
\eqnlbl{Lambda}
$$
where $\Lambda_j^\pm$ denote the open sets bounded on the left by 
the algebraic curves $\lambda_j^\pm(\xi)$
determined by the eigenvalues of $ i\xi\binom{df}{dg} (u_\pm, v_\pm) +
\binom 0{dq} (u_\pm, v_\pm)$ as $\xi$ is varied along the real axis.

\proclaim{Lemma \thmlbl{frozen}}
The limiting equations \eqnref{limitevalue}
have no center manifold on the set $\Lambda$.
Moreover, assuming (H0)--(H1), (H3),  we have:

(i)
$
\Lambda\subset
\{\lambda: \R \lambda >
- \eta |\I \lambda|/(1 + |\I \lambda|),
$ 
$\eta>0$; and 

(ii) the stable manifold at $+\infty$ and 
the unstable manifold at $-\infty$ have dimensions summing to the
full dimension $N=n+r$.
\endproclaim

{\bf Proof.}
The fundamental modes of \eqnref{limitevalue} are of form $e^{\mu x}V$,
where $\mu$, $V$ satisfy the {\it characteristic equation}
$$
(-\mu A_\pm + Q_\pm -\lambda)V=0.
\eqnlbl{char2}
$$
The existence of a center manifold thus corresponds with existence of
solutions $\mu=i\xi$, $V$ of \eqnref{char2}, $\xi$ real, i.e., solutions
of the {\it dispersion relation}
$$
(-i\xi A_\pm + Q_\pm -\lambda)V=0.
\eqnlbl{disp}
$$
But, $\lambda\in \sigma (-i\xi A_\pm + Q_\pm)$ implies, by definition
\eqnref{Lambda}, that $\lambda$ lies outside of $\Lambda$,
establishing the first claim.
Bound (i), likewise, follows easily from the bound (H3) on the
dispersion curves $\lambda_j^\pm$.
Finally, nonexistence of a center manifold, together with 
connectivity of $\Lambda$, implies that the dimensions of 
stable/unstable manifolds at $+\infty$/$-\infty$ 
are constant on $\Lambda$.
Taking $\lambda\to +\infty$ along the real axis, we find from
\eqnref{char2} that these dimensions are those of the unstable/stable
subspaces of $A_+$/$A_-$; but, these must sum to $N$, since, by (H1),
the dimensions of unstable/stable subspaces of $A(x)$ are independent
of $x$.  This establishes (ii), completing the proof.
\myqed

For asymptotically constant-coefficient ordinary differential
operators, the rightmost component ($\Lambda$, in this case)
of the set of $\lambda$ such 
that the limiting, constant-coefficient eigenvalue equations  
$(L_\pm -\lambda)W=0$ are of hyperbolic type, with dimensions 
of unstable/stable manifolds summing to that of the full phase 
space is  called in the Evans function literature the {\it domain of 
consistent splitting} [AGJ].
Its significance (see, e.g. [He], Lemma 2, pp. 138--139) 
is that it lies in the complement of the 
essential spectrum of the variable-coefficient operator $L$,
with its boundary lying in the essential spectrum: that is,
it is a maximal domain in the essential spectrum complement,
consisting of the component containing real, plus infinity.
Moreover, provided that the coeffients of $L$
approach their limits at integrable rate, it can be shown
(see, e.g. [AGJ,GZ,ZH]) that each connected component consists either
entirely of eigenvalues, or else entirely of {\it normal points},
defined as resolvent points or isolated eigenvalues of finite
multiplicity.  The latter, and for our purposes more significant
fact will be seen directly in the course of our construction. 

Recall, \eqnref{expdecay}, that coefficients of $L$ converge exponentially
to their limiting values in $L_\pm$, in particular at {\it integrable rate}.
In standard fashion
(see, for example, [Co,AGJ,ZH,Z.1,ZS,Z.4]), we may thus conclude:

\proclaim{Proposition \thmlbl{bases}}
If (H0)--(H1), (H3) hold, then, about any $\lambda$
in the domain of consistent splitting $\Lambda$,
there exist locally analytic (in $\lambda$) bases
$\phi_1^+, \dots, \phi_k^+(x;\lambda)$ and $\phi_{k+1}^-,\dots, \phi_{n+r}^-(x;\lambda)$
respectively spanning the stable manifold at $+\infty$
and the unstable manifold at $-\infty$ of solutions of
the variable-coefficient equation \eqnref{evalue}.
Moreover, these manifolds are tangent as $x\to -\infty$, $x\to +\infty$,
respectively, to their constant-coefficient counterparts, decaying
exponentially for $x\ge 0$, $x\le 0$, respectively, with uniform rate
$Ce^{-\eta |x|}$, $\eta>0$.

Likewise, there exist bases
$\psi_1^+, \dots, \psi_k^+(x;\lambda)$ and $\psi_{k+1}^-,\dots, \psi_{n+r}^-(x;\lambda)$
respectively spanning the unstable manifold at $+\infty$
and the stable manifold at $-\infty$ of solutions of
\eqnref{evalue}, that are locally analytic in $\lambda$,
tangent as $x\to -\infty$, $x\to +\infty$,
respectively, to their constant-coefficient counterparts, and
growing exponentially for $x\ge 0$, $x\le 0$, respectively, with uniform rate
$Ce^{\eta |x|}$, $\eta>0$.
\endproclaim

{\bf Proof.}  By spectral separation of stable/unstable subspaces,
we may conclude the existence of globally analytic bases of solutions
of the limiting equations \eqnref{limitevalue},
using a standard result of Kato ([Kat], pp. 81--83).
The result for the variable-coefficient equations \eqnref{evalue}
then follows by standard, asymptotic theory of ODE [Co,Z.1],
using the integrable rate of convergence of coefficients to
their limiting values.  Alternatively, using the fact that 
convergence is in fact exponential, we may apply the results
of Appendix A2.1 (specifically, Remark 8.4) to obtain a self-contained
treatment.
\myqed

{\bf Observation \thmlbl{strat}.} 
To show that $\lambda$ is in the resolvent set
of $L$, 
with respect to any $L^p$, $1\le p\le \infty$,
it is enough to:
(i) construct a resolvent kernel $G_\lambda(x,y)$
satisfying \eqnref{glambda} and obeying a uniform decay estimate
$$
|G_\lambda(x,y)|\le Ce^{-\eta|x-y|};
\eqnlbl{udecay}
$$
\noindent and,
(ii) show that there are no $L^p$ solutions of $(L-\lambda)W=0$, i.e., 
$\lambda$ is not an eigenvalue of $L$.
(For $p<\infty$, these are necessary as well as sufficient, 
though we do not show it here [Z.1]).
For, then, the Hausdorff--Young inequality yields that
the distributional solution formula
$(L-\lambda)^{-1}f:=\int G_\lambda(x,y) f(y)dy$ yields a bounded right
inverse of $(L-\lambda)$ taking $L^p$ to $L^p$, while nonexistence of
eigenvalues implies that this is also a left inverse.
On the domain of consistent splitting, $\Lambda$, we shall see that (i)--(ii)
are equivalent.
\medskip

Following the strategy of Observation \thmref{strat}, we are forced to
to choose $G(x,y)$, for fixed $y$ and $x\gtrless y$ from among
the decaying solutions at $+\infty$/$-\infty$, respectively,
of the homogeneous eigenvalue equations \eqnref{evalue}, in such
a way the $G_\lambda$ satisfies at $x=y$ the {\it jump condition}
 $$
  \bigl[\,G_\lambda\,\bigr]_{x=y}:= G_\lambda(y+0,y)-
   G_\lambda(y-0,y)=-A^{-1}(y),
 \eqnlbl{jump}
 $$
where $[\,h\,]_{x=y}$ denotes the jump of the function $h$ at $x=y$:
that is, according to the classical construction of [CH].

Implementing this strategy on the domain $\Lambda$ of consistent
splitting, we seek a solution of form
 $$
  G_\lambda(x,y)=\cases
   \Phi^+(x;\lambda)N^+(y;\lambda)\,\qquad &x>y,\\
   \Phi^-(x;\lambda)N^-(y;\lambda)\,\qquad &x<y,\\
                 \endcases 
\eqnlbl{gref}
 $$
where
 $$
  \Phi^+(x;\lambda)=(\,\phi^+_1(x;\lambda)\, \dots\, \phi^+_{k}(x;\lambda)\,)
   \in {\Bbb R}^{N\times k}
  \eqnlbl{Phi+}
 $$
and
 $$ 
  \Phi^-(x;\lambda)=(\,\phi^-_{k+1}(x;\lambda)\, \dots\, \phi^-_{N}(x;\lambda)\,)
  \in {\Bbb R}^{N\times (N-k)},
  \eqnlbl{Phi-}
 $$
$\phi_j^\pm$ as defined in Proposition \thmref{bases}.
Here, $N^+=N^+(y;\lambda)$ and $N^-=N^-(y;\lambda)$ are $k\times N$ 
and $(N-K)\times N$ matrices, respectively, $N:=n+r$,
to be determined by imposing the jump condition \eqnref{jump},
or equivalently
 $$
  (\, \Phi^+(x;\lambda) \quad \Phi^-(x;\lambda)\,)
  \pmatrix N^+(y;\lambda)\\ -N^-(y;\lambda) \endpmatrix=-A^{-1}(y).
  \eqnlbl{jump2}
 $$

Inverting \eqnref{jump2}, we obtain the expression 
 $$
  \pmatrix N^+(y;\lambda)\\ -N^-(y;\lambda) \endpmatrix =
  - 
\pmatrix \Phi^+ & \Phi^- \endpmatrix^{-1} 
(y;\lambda) \, A^{-1}(y).
  \eqnlbl{Nsoln}
 $$
for the matrices $N^\pm$, yielding
a final expression for the resolvent kernel $G_\lambda$,
in terms of the matrices $\Phi^\pm$, of
 $$
  G_\lambda(x,y)=\cases
   -(\,\Phi^+(x;\lambda)\quad 0\,)\,
\pmatrix \Phi^+ & \Phi^- \endpmatrix^{-1} 
    (y;\lambda)\, A^{-1}(y)\qquad &x>y,\\
   (\, 0\quad \Phi^-(x;\lambda)\,)\,
\pmatrix \Phi^+ & \Phi^- \endpmatrix^{-1} 
(y;\lambda)\, A^{-1}(y)\qquad &x<y,\\ \\
                 \endcases 
  \eqnlbl{2.12}
 $$
valid when the matrix
$$
\Phi:=\pmatrix \Phi^+ & \Phi^-\endpmatrix
\eqnlbl{Phi}
$$
is invertible.  If $\Phi$ is singular, on the other hand, then
$\lambda$ is clearly an eigenvalue, since there must be some
function lying in the span of both $\Phi^+$ and $\Phi^-$, hence
decaying exponentially at both $\pm \infty$; thus, we see that
(i) implies (ii) on $\Lambda$, partially verifying our claim that
they are equivalent.

It will at times be useful to represent the resolvent kernel
in a more intrinsic, geometric fashion.  Observing that
 $$
  (\, 0\quad \Phi^-(x;\lambda)\,)= 
  \Phi (x;\lambda)\,\pmatrix 0& 0\\ 0 &I_n \endpmatrix
  = \Phi (x;\lambda) \, \Phi^{-1} (y;\lambda) \, 
  (\,0\quad \Phi^-(y;\lambda)\,),
  \eqnlbl{geometric}
 $$
and, similarly,
$$
(\,\Phi^+(x;\lambda) \quad 0\,)=\Phi (x;\lambda) \, \Phi^{-1} (y;\lambda)\,
(\, \Phi^+(y;\lambda) \quad 0\,),
\eqnlbl{}
$$
we get the expression
 $$
  G_\lambda(x,y)=\cases
   -\Fyx\Pi^+_y\, A^{-1}(y)
                                  \qquad &x>y,\\
   \Fyx\Pi^-_y\, A^{-1}(y)
                                  \qquad &x<y,\\ \\
                        \endcases 
  \eqnlbl{basicflow}
 $$
where we have defined the solution operator from $y$ to $x$ of
\eqnref{evalue}, denoted by $\Fyx$, as
 $$
  \Fyx= \Phi(x;\lambda)\, \Phi^{-1}(y;\lambda)
  \eqnlbl{Fyx}
 $$ 
and the projections $\Pi^\pm_y$ on the stable manifolds at $\pm\infty$ as
 $$
  \Pi^+_y=(\,\Phi^+(y;\lambda)\quad 0\,)\,\Phi^{-1}(y;\lambda) 
   \quad\text{and}\quad
  \Pi^-_y=(\,0\quad \Phi^-(y;\lambda)\,)\,\Phi^{-1}(y;\lambda). 
  \eqnlbl{Pi}
 $$
Note that in \eqnref{geometric}, there is no longer any reference
to the specific choice of bases functions $\phi_j^\pm$, but only
to the (uniquely defined) stable and unstable manifolds at $+\infty$,
$-\infty$, respectively.

With the above representations in hand, we readily obtain the 
following basic, intermediate-frequency result.
More careful analyses will be required in the large- and small-freqency limits.

\proclaim{Proposition \thmlbl{intbounds}}
With respect to any $L^p$, $1\le p\le \infty$, the domain of consistent 
splitting, $\Lambda$, consists entirely of normal points of $L$, i.e., 
resolvent points, or isolated eigenvalues of constant multiplicity.
Moreover, the resolvent kernel is meromorphic on $\Lambda$, 
satisfying either of formulae
\eqnref{2.12},\eqnref{geometric}.
Finally, on any compact subset $K$ of $\rho(L)\cap \Lambda$, 
there holds the uniform decay estimate
$$
|G_\lambda(x,y)|\le Ce^{-\eta |x-y|},
\eqnlbl{uKdecay}
$$
where $C>0$ and $\eta>0$ depend only on $K$, $L$.
\endproclaim

{\bf Proof.}  As noted above, eigenvalues coincide with
zeroes of the (locally) analytic function
$$
D(\lambda):=\det \Phi,
\eqnlbl{localevans}
$$
which are either isolated of finite multiplicity
or else fill all of $\Lambda$.
For the remaining claims, it is sufficient to verify
\eqnref{uKdecay} holds for the function $G_\lambda$
defined by formulae \eqnref{2.12},\eqnref{geometric}.
For, it then follows that this function, off of the
set of isolated eigenvalues, indeed represents the
resolvent kernel, and that all noneigenvalue points
are in the resolvent set.  By inspection, \eqnref{2.12}
is meromorphic, with poles of finite order less
than the multiplicity of the corresponding zeroes of $D$;
it follows using the spectral resolution formula 
(see, e.g., [Kat,Y,ZH]) that the eigenprojection 
of $L$ at each eigenvalue $\lambda_0$, defined as the residue at 
$\lambda_0$ of the resolvent $(L-\lambda)^{-1}$, is
well-defined, with range contained in the $L^p$ kernel
of $(L-\lambda)^k$, where $k$ is the order of the pole.
But, the solution set of $(L-\lambda)^kW=0$ is at most
$nK$, hence finite.

Exponential decay is best seen using the geometric formula
\eqnref{geometric}.  First, observe that, for $D(\lambda)$ bounded
away from zero, the manifold of decaying solutions at $-\infty$
is uniformly transverse at $x=0$ to the manifold of decaying solutions
at $+\infty$, hence is uniformly exponentially growing as $x\to +\infty$:
that is, away from eigenvalues of $L$, this manifold is {\it globally}
uniformly exponentially decaying as $x$ decreases, on the whole
real line $-\infty\le x\le +\infty$.  Likewise, the manifold of
decaying solutions at $+\infty$ is globally uniformly exponentially
decaying as $x$ increases.  A side-consequence is that the two
manifolds are uniformly transverse as $x\to \pm\infty$,
approaching the unstable/stable manifolds of the limiting coefficient
matrix of \eqnref{evalue}; by continuity, they are uniformly transverse,
then, on the whole real line.  Thus, projections $\Pi_y^\pm$ are bounded
in \eqnref{geometric}, and flows $\Cal{F}^{y\to x}$ are uniformly
exponentially decaying, giving the result \eqnref{uKdecay}.
\myqed

The transmission function $D(\cdot)$ defined in \eqnref{localevans} is
known as an {\it Evans function}. 
It is possible, with more care, to define a globally analytic 
Evans function on $\Lambda$, and thereby to determine information
on location of eigenvalues via topological considerations,
see e.g. [E.1--4,J,AGJ,PW,GZ].  
However, for our purposes local analyticity is sufficient.
Though there exist as many versions of the Evans function
as there are analytic bases $\phi_j^\pm$, it is evident that
these differ, locally, only by a nonvanishing analytic factor.
Thus, location and multiplicity of zeroes agree for all versions
of the Evans function, as the following proposition of Gardner
and Jones shows that they must.

\proclaim{Proposition \thmlbl{gj}[GJ.1--2]}
On $\Lambda$, the zeroes of (any version of) the 
Evans function $D(\cdot)$ agree in both location and multiplicity
with the eigenvalues of $L$.
\endproclaim

{\bf Proof.} It is evident that they agree in location,
whence the result follows if one can show that zeroes
and eigenvalues both split under regular perturbation of $L$;
the existence of such a perturbation (a nontrivial assertion)
was shown by Gardner and Jones in [GJ.1].  Alternatively, they showed
in [GJ.2] that the correspondence may be seen by explicit computation 
in a convenient Jordan basis.
\myqed

Note that Proposition \thmref{gj} gives an alternative
proof that eigenvalues of $L$ are of finite multiplicity on $\Lambda$.

{\bf Remark \thmlbl{cc}.}
In the constant-coefficient case, the matrices $\Phi^\pm$ above can be 
calculated (almost) explicitly.
To accomplish this aim, we have to find the decaying modes at $\pm\infty$,
that is we have to find independent solutions to
 $$
  - A W' + B W - \lambda W = 0\qquad\iff\qquad
   W' = A^{-1} (B - \lambda I) W,
  \eqnlbl{}
 $$
such that $W(\pm\infty)=0$.
So, let $\mu_i^\pm=\mu_i^\pm(\lambda)$ with $i=1,\dots,N$ be the solutions of
 $$
  p(\mu)=\det\left(A^{-1}(B-\lambda I)-\mu I\right)=0,
  \eqnlbl{}
 $$
with $\Re\mu_i^+<0<\Re\mu_j^-$ for any $i=1,\dots,k$ and $j=k+1,\dots,N$.
Then, assuming $\mu_i^\pm\neq \mu_j^\pm$ for any $i\neq j$, 
$\phi^\pm_i=r^\pm_i e^{\mu^\pm_i x}$ where $r^\pm_i$ is a right
eigenvector of the matrix  $A^{-1}(B-\lambda I)$ relative to the
eigenvalue $\mu_i^\pm$, so that
 $$
  \Phi(x;\lambda)=(\,
   r^+_1 e^{\mu^+_1 x} \quad \dots\quad r^+_k e^{\mu^+_k x}\quad 
   r^-_{k+1} e^{\mu^-_{k+1} x} \quad \dots\quad r^-_N e^{\mu^-_N x}\,)
  \in{\Bbb R}^{N\times N}. 
  \eqnlbl{}
 $$

\noindent{\bf Example \thmlbl{JXN}.} 
For the Jin--Xin model, the equations linearized at the stationary 
solution $(\bar u,\bar v)$ are
 $$
  \cases
   u_t - s   u_x +   v_x = 0,\\ 
   v_t + a^2 u_x - s v_x = \alpha(x) u - v,
         \endcases
  \eqnlbl{JX2}
 $$
where $\alpha(x):=dh(\bar u(x))$.
Hence
 $$
  W=\pmatrix     u\\ v            \endpmatrix,\quad
  A=\pmatrix -s &1\\ a^2 &-s      \endpmatrix,\quad
  B=\pmatrix  0 &0\\ \alpha(x) &-1\endpmatrix,
  \eqnlbl{}
 $$
(here any element in the matrices $A$ and $B$ corresponds to an $n\times n$
block).

If $\alpha(x)\equiv$ constant, then
$p(\mu)=\Pi_{j=1}^n
\bigl[(a^2-s^2)\mu^2+[(1+2\lambda)s-\alpha_j]\mu-\lambda(1+\lambda)
\bigr]$, 
where $\alpha_j$ denote the eigenvalues of $\alpha$,
so that there are $n$ pairs of eigenvalues 
 $$
  \mu_j^\pm(\lambda)=\frac{\alpha_j-(1+2\lambda)s\mp\sigma_j(\lambda)}{2(a^2-s^2)},
  \eqnlbl{}
 $$
where $\sigma_j(\lambda)=\sqrt{[\alpha_j-(1+2\lambda)s]^2+
4\lambda(1+\lambda)(a^2-s^2)}$. 
If the subcharacteristic condition holds, then $a^2-s^2>0$
and for any $\lambda>0$ there holds $\mu_j^+(\lambda)<0<\mu_j^-(\lambda)$.

The corresponding eigenvectors (for $\lambda>0$) are
 $$
  W_j^\pm = C \left(\,\left[(1+\lambda)s+(a^2-s^2)\mu_j^\pm\right]\,r_j\,,
\,(\alpha_j s-a^2\lambda)\,r_j\,\right)^t\qquad\qquad(C\neq 0),
 $$ 
where $r_j$ denote right eigenvectors of $\alpha$ associated with $\alpha_j$. 
Hence,
 $$
  \Phi(x;\lambda)=
  \Bigl( W_1^+\, e^{\mu_1^+ x}\quad\cdots\quad W_n^+\, e^{\mu_n^+ x}\quad
         W_1^-\, e^{\mu_1^- x}\quad\cdots\quad W_n^-\, e^{\mu_n^- x}\Bigl),
   \eqnlbl{JX_decaying}
 $$
and the resolvent kernel is given by
 $$
  G_\lambda(x,y)=\cases
    -\sum_{i=j}^n W^+_j V^+_j A^{-1} e^{\mu_j^+(x-y)}\quad &x<y,\\
     \sum_{j=1}^n W^-_i V^-_j A^{-1} e^{\mu_j^-(x-y)}\quad &x>y,
                 \endcases 
  \eqnlbl{}
 $$
where $V^\pm_j$ are row vector defined by
 $$
  {\bigl(\, W_j^+ \quad W_i^-\,\bigr)}^{-1}=
   \pmatrix V^+_j\\ V^-_j\\ \endpmatrix.
  \eqnlbl{inverse_coeff}
 $$
\myqed


{\bf 3.2. Dual formulation.}
For use in our later treatment of low-frequency behavior,
we develop a modified representation of the 
resolvent kernel involving solutions of the adjoint eigenvalue equation
$$
(L^*-\lambda^*){\tilde W}=0,
\eqnlbl{dualevalue}
$$
where
$$
L^*{\tilde W}:=A^*{\tilde W}'+Q^*{\tilde W}
\eqnlbl{Ladjoint}
$$
denotes the adjoint operator for $L$, and $^*$ for a matrix
or vector denotes matrix adjoint, or conjugate transpose.
Following [ZH,Z.4], we first point out:

\proclaim{Lemma \thmlbl{dual}} For any ${\tilde W}$, $W$ such that
$(L^*-\lambda^*){\tilde W}=0$ and $(L-\lambda)W=0$, there holds
$$
\langle {\tilde W},AW\rangle\equiv\text{\rm constant}.
\eqnlbl{duality}
$$
\endproclaim

\noindent{\bf Proof.}
 By direct calculation, 
 $$
\aligned
\langle {\tilde W},AW\rangle'&=
\langle A^*{\tilde W}',W\rangle +
\langle {\tilde W},(AW)'\rangle \\
&=
\langle -(Q^*-\lambda^*){\tilde W},W\rangle +
\langle {\tilde W},(Q-\lambda)W\rangle =0.
\endaligned
  \eqnlbl{calc}
 $$
\myqed

{}From \eqnref{calc}, it follows that if there are $k$ independent solutions
$\phi^+_1, \dots, \phi^+_{k}$ of $(L-\lambda I)W=0$ decaying at $+\infty$, 
and $N-k$ independent solutions $\phi^-_{k+1}, \dots, \phi^-_N$ of the same 
equations decaying at $-\infty$, then there exist $N-k$ independent 
solutions $\tilde\psi^+_{k+1}, \dots, \tilde\psi^+_{N}$ of $(L^*-\lambda^* I){\tilde W}=0$ 
decaying at $+\infty$, and $k$ independent solutions 
${\tilde\psi}^-_{1}, \dots, \tilde\psi^-_{k}$ decaying at $-\infty$.
More precisely, setting
 $$
  \Psi^+(x;\lambda)=
  \pmatrix \psi^+_{k+1}(x;\lambda) & \cdots & 
           \psi^+_N(x;\lambda)\endpmatrix
   \in {\Bbb R}^{N\times (N-k)},
  \eqnlbl{Psi+}
 $$
$$
  \Psi^-(x;\lambda)=
  \pmatrix\psi^-_1(x;\lambda) & \cdots & 
          \psi^-_{k}(x;\lambda)\endpmatrix
  \in {\Bbb R}^{N\times k},
  \eqnlbl{Psi-}
 $$
and
 $$
  \Psi (x;\lambda) = 
  \pmatrix \Psi^-(x;\lambda) &
           \Psi^+(x;\lambda)\endpmatrix
  \in {\Bbb R}^{N\times N},
  \eqnlbl{Psi}
 $$
where $\psi_j^\pm$ are the exponentially growing solutions described
in Proposition \thmref{bases}, we may define dual exponentially decaying
and growing solutions $\tilde \psi_j^\pm$ and $\tilde \phi_j^\pm$ via
$$
\pmatrix \tilde \Psi & \tilde \Phi \endpmatrix_\pm^* A
\pmatrix \Psi & \Phi \endpmatrix_\pm\equiv I.
\eqnlbl{dualbase}
$$

Recalling (see, e.g. [ZH], Lemma 4.2) the classical duality principle 
that the transposition ${G}_{\lambda}^*(y,x)$ of the Green's function
$G_\lambda(x,y)$ associated with operator $(L-\lambda)$ should
be the Green's function for the adjoint operator $(L^*-\lambda^*)$, i.e.,
$$
(L^*_y-\lambda^*)G_\lambda^*(x,y)=\delta(x-y),
$$
we may seek, alternatively, to represent the resolvent kernel in the form
 $$
  G_\lambda(x,y)=\cases
   \Phi^+(x;\lambda)M^+(\lambda)\tilde\Psi^{-*}(y;\lambda)\,\qquad &x>y,\\
   \Phi^-(x;\lambda)M^-(\lambda)\tilde\Psi^{+*}(y;\lambda)\,\qquad &x<y,\\
                 \endcases 
  \eqnlbl{dualglambda1}
 $$
where $M^+=M^+(\lambda)$ and $M^-=M^-(\lambda)$ are $k\times k$ 
and $(N-k)\times (N-k)$ matrices, respectively, to be determined by 
imposing the jump condition \eqnref{jump}.

Solving, as before, we obtain
 $$
  M(\lambda):=\pmatrix -M^+(\lambda) &0\\ \\
           0 & M^-(\lambda)\endpmatrix=
  \Phi^{-1}(z;\lambda)A^{-1}(z)\tilde\Psi^{-1*}(z;\lambda),
  \eqnlbl{M}
 $$
where 
$$
\tilde \Psi:=\pmatrix \tilde \Psi^- & \tilde \Psi^+\endpmatrix.
\eqnlbl{tildePsi}
$$
Note that the independence of the righthand side with respect to $z$ 
is a consequence of the previous lemma.
Thus,
 $$
  G_\lambda(x,y)=\cases
   -(\,\Phi^+(x;\lambda)\quad 0\,)\, M(\lambda) \,
    (\,\tilde\Psi^-(y;\lambda)\quad 0\,)^*\qquad &x>y,\\
   (\,0\quad \Phi^-(x;\lambda)\,)\, M(\lambda) \,
    (\,0\quad \tilde\Psi^+(y;\lambda)\,)^*\qquad &x<y,\\
                 \endcases 
  \eqnlbl{dualglambda}
 $$
with $M$ as in \eqnref{M}.
It is readily seen [ZH,Z.4] by duality, 
\eqnref{duality},  that $\tilde\Psi$ is nonsingular if 
and only if $\Phi$ is nonsingular, so that formula \eqnref{dualglambda} 
is indeed valid on the same region as were the previous ones; see also
Proposition 3.11, below.

Finally, let us note as before that it is possible to 
represent the matrix $G_\lambda$ by means of intrinsic objects 
as solution operators and projections on stable manifolds
 $$
  G_\lambda(x,y)=\cases
   -\Fzx\Pi^+_z\, A^{-1}(z)\,\tilde\Pi^-_z\,\tilde\Fzy
                                  \qquad &x>y,\\
   \Fzx\Pi^-_z\, A^{-1}(z)\,\tilde\Pi^+_z\,\tilde\Fzy
                                  \qquad &x<y,\\ \\
                        \endcases 
  \eqnlbl{dualgeometric}
 $$
where 
 $$
  \tilde\Fzy= \tilde\Psi^{-1}(z;\lambda)\, \tilde\Psi(y;\lambda)
  \eqnlbl{}
 $$ 
and 
 $$
  \tilde\Pi^+_z= 
    \tilde\Psi^{-1}(z;\lambda)  
    \pmatrix 0 \\ \tilde\Psi^+(z;\lambda) \endpmatrix
   \quad\text{and}\quad
  \tilde\Pi^-_z=
    \tilde\Psi^{-1}(z;\lambda)  
    \pmatrix \tilde\Psi^-(z;\lambda) \\ 0 \endpmatrix.
  \eqnlbl{}
 $$
\medskip

{\bf Remark \thmlbl{decaypf}.} The exponential decay asserted in Proposition 
\thmref{intbounds} is somewhat more straightforward to see in the
dual formulation, by judiciously choosing $z$ in \eqnref{dualgeometric}
and using the exponential decay of forward and adjoint flows 
on $x\gtrless 0$ alone.
\medskip

{\bf 3.3. Spectral decomposition.}
The matrix $M(\lambda)$ in \eqnref{M} depends on the choice of 
functions $\phi_j^\pm$ and $\psi_j^\pm$. 
For example, formula \eqnref{dualglambda} may be considerably simplified
in the constant-coefficient case by the choice of bases
$\Psi^\pm=\Phi^\mp$, which yields
 $$
  \tilde\Psi(x;\lambda)A(x)\Phi(x;\lambda)=I
   \eqnlbl{}
$$
and $M(\lambda)=I$, or equivalently 
 $$
  G_\lambda(x,y)=\cases
   -\sum_{j=k+1}^{N} \phi_j^+(x;\lambda)\tilde \phi_j^{+*}(y;\lambda)
\qquad &x>y,\\
   \sum_{j=1}^k \phi_j^-(x;\lambda)\tilde \phi_j^{-*}(y;\lambda)
    \qquad &x<y,\\
                 \endcases 
  \eqnlbl{cc2}
 $$
where $\phi_j^\pm$, $\tilde {\phi}_j^\pm$ may be taken as pure exponentials
$$
\phi_j^\pm(x)\tilde \phi_j^{\pm*}(y)=e^{\mu_j^\pm(\lambda)(x-y)}
V_j^\pm(\lambda) \tilde{V}_j^{\pm*}(\lambda).
\eqnlbl{spectral}
$$
This appealing formula may be viewed as a {\it generalized spectral
decomposition}.
Note that, moving individual modes $\phi_j^\pm\tilde \phi_j^\pm$
in the spectral resolution formula \eqnref{giLT}, using
Cauchy's Theorem, to contours $\mu_j(\lambda)\equiv i\xi$ lying along
corresponding dispersion curves $\lambda=\lambda_j(\xi)$, we obtain
the standard decomposition of $e^{Lt}$ into eigenmodes of continuous
spectrum: in this (constant-coefficient) case, just the usual representation
obtained by Fourier transform solution.

{\bf Example \thmlbl{JXcc}.} 
For the Jin--Xin model \eqnref{JX2}, the dual bases of 
$\phi^{\pm}_i$ can be obtained by inverting the relation
$\tilde\Psi(y;\lambda)A\Phi(y;\lambda)=I$.
Hence, recalling formula \eqnref{JX_decaying} and definition of $V^\pm_j$ in
\eqnref{inverse_coeff}, we obtain
 $$
  \tilde\Psi(y;\lambda)= 
    \pmatrix 
   V^+_1\,e^{-\mu_1^+ y}\,\\ \vdots \\ V^+_n\,e^{-\mu_n^+ y}\,\\ 
   V^-_1\,e^{-\mu_1^- y}\,\\ \vdots \\ V^-_n\,e^{-\mu_n^- y}
    \endpmatrix 
    A^{-1}
 $$
 $$
   \frac1{(a^2\lambda-\alpha s)\sigma(\lambda)}
   \pmatrix
   (a^2\mu_+-\alpha)e^{-\mu_+ y} & (\mu_+s-\lambda)e^{-\mu_+ y} \\ \\
   (\alpha-a^2\mu_-)e^{-\mu_- y} & (\lambda-\mu_-s)e^{-\mu_- y}
  \endpmatrix.
  \eqnlbl{}
 $$
(we used $\sigma(\lambda)=(a^2-s^2)(\mu_- - \mu_+)$).
As noted before, for such a choice of a basis for the space of 
solutions to the adjoint equation, $M(\lambda)\equiv I$ and 
the resolvent kernel $G_\lambda$ can be represented by 
 $$
  G_\lambda(x,y)=\cases
   (\,0\quad \Phi^-(x;\lambda)\,)\,
    (\,0\quad \tilde\Psi^+(y;\lambda)\,)^*\qquad &x<y,\\
   -(\,\Phi^+(x;\lambda)\quad 0\,)\,
    (\,\tilde\Psi^-(y;\lambda)\quad 0\,)^*\qquad &x>y,
                 \endcases 
  \eqnlbl{cc1}
 $$
or, expanding, \eqnref{cc2}.
\myqed

In the asymptotically constant-coefficient case, the notion
of continuous spectrum is in general not particularly useful
outside of the self-adjoint case, see e.g.  discussion of [OZ].
However, as pointed out in [ZH], the more primitive relation
\eqnref{cc2} persists in the form of a {\it scattering decomposition}
encoding interactions of scalar modes at $\pm \infty$.
This decomposition gives a natural generalization 
of the Fourier transform to the asymptotically constant-coefficient case,
one of the main contributions of [ZH].

The following result of [Z.4], a slight extension of Proposition
7.1 of [ZH], generalizes \eqnref{cc2} to the asymptotically 
constant-coefficient case.

\proclaim{Proposition \thmlbl{kochelformulae} [Z.4]}  
On $\Lambda \cap \rho(L)$, there hold
$$
G_\lambda(x,y) =
\sum_{j,k} M^+_{jk}(\lambda ) \phi^+_j (x;\lambda) \Tpsi^-_k(y;\lambda)^* 
\eqnlbl{Mplusrep}
$$
for $y\le 0\le x$, 
$$
G_\lambda(x,y) =
\sum_{j,k} {d}^+_{jk}(\lambda ) \phi^-_j (x;\lambda) 
\Tpsi^-_k(y;\lambda)^* 
-\sum_k \psi^-_k (x;\lambda)
\Tpsi^-_k (y;\lambda)^*
\eqnlbl{dplusrep}
$$
for $y\le x\le 0$, and 
$$
G_\lambda(x,y) =
\sum_{j,k} {d}^-_{jk}(\lambda ) \phi^-_j (x;\lambda) 
\Tpsi^-_k(y;\lambda)^* +
\sum_k \phi^-_k (x;\lambda)
\tilde \phi^-_k (y;\lambda)^*
\eqnlbl{dminusrep}
$$
for $x\le y\le 0$, with
$$
M^+=
(-I,0) 
\pmatrix
\Phi^+ & \Phi^-\\
\endpmatrix^{-1}
\Psi^-
\eqnlbl{Mplus}
$$
and
$$
d^\pm= 
\big( 0 , I\big)
\pmatrix
\Phi^+ & \Phi^-\\
\endpmatrix^{-1 }
\Psi^-.
\eqnlbl{dplusminus}
$$
Symmetric representations hold for $y\ge 0$.
\endproclaim

{\bf Proof.}  
The matrix $M^+$ in \eqnref{M} may be expanded using duality relation 
\eqnref{duality} as 
$$
\aligned
M^+=
(-I,0) 
\pmatrix
\Phi^+ & \Phi^-\\
\endpmatrix^{-1}
&A^{-1} 
\pmatrix
\tilde \Psi^- & \tilde \Phi^-\\
\endpmatrix^{*\, -1} 
\pmatrix
I\\
0
\endpmatrix_{|_{z}}\\
&= (-I,0)
\pmatrix
\Phi^+ & \Phi^-\\
\endpmatrix^{-1}
\pmatrix
\Psi^- & \Phi^-\\
\endpmatrix  
\pmatrix
I\\
0
\endpmatrix_{|_{z}}\\
&=
(-I,0) 
\pmatrix
\Phi^+ & \Phi^-\\
\endpmatrix^{-1}
\Psi^-_{|_z},\\
\endaligned
\eqnlbl{4.29}
$$
yielding \eqnref{Mplus} for $x\ge y$, in particular for $y\le 0\le x$.

Next, expressing $\phi_j^\pm(x;\lambda)$ as a linear combination of
basis elements at $-\infty$, we obtain the preliminary representation
$$
G_\lambda(x,y) =
\sum_{j,k} {d}^+_{jk}(\lambda ) \phi^-_j (x;\lambda) 
\Tpsi^-_k(y;\lambda)^* +
\sum_{j,k} e^+_{jk}\psi^-_j (x;\lambda)
\Tpsi^-_k (y;\lambda)^*,
\eqnlbl{predplus}
$$
valid for $y\le x\le 0$.
Duality, \eqnref{duality}, with \eqnref{basicflow},
and the fact that $\Pi_+=I-\Pi_-$, gives
$$
\aligned
-\pmatrix d^+ \\  e^+ \\ \endpmatrix&=
\pmatrix \tilde \Phi^{-} & \tilde  \Psi^{-}\endpmatrix^*
A \Pi_+ 
{\Psi^-}_{|_{x}}\\
&=
-\pmatrix
\Phi^- & \Psi^-\\
\endpmatrix^{-1}
[I-
\pmatrix 0 & \Phi^-\\ \endpmatrix
\pmatrix
\Phi^+ & \Phi^-\\
\endpmatrix^{-1}]
\Psi^{-}
\\
&=
\pmatrix 0 \\ -I_{k} \\ \endpmatrix 
-
\pmatrix 0 & I_{N-k}\\
0&0\\ \endpmatrix
\pmatrix
\Phi^+ & \Phi^-\\
\endpmatrix^{-1}
\Psi^-,
\endaligned
\eqnlbl{dpluscalc}
$$
yielding \eqnref{dplusrep} and \eqnref{dplusminus} for $y\le x\le 0$. 
Relations \eqnref{dminusrep} and \eqnref{dplusminus} follow for
$x\le y\le 0$ in similar,
but more straightforward fashion from \eqnref{duality} and \eqnref{basicflow}.
\myqed

{\bf Remarks}.  
1. Note, in the various scattering decompositions 
given in Proposition \thmref{kochelformulae}, that
all functions $\phi_j^\pm$, $\tilde \phi_j^\pm$,
$\tilde \phi_j ^\pm$, $\tilde \psi_j^\pm$ are evaluated
on the intervals $[0,\pm \infty]$ on which their behavior
is known.  Indeed, by Proposition \thmref{bases}, they
are well-approximated by corresponding solutions
of the limiting, constant-coefficient equations at $\pm \infty$,
i.e., {\it exponential modes} generalizing the continuous
eigenmodes obtained by Fourier transform in the constant-coefficient
case.

\medskip
2. Letting $x$, $y \to -\infty$ in \eqnref{dminusrep},
we find that
$$
G_{\lambda }
\sim \sum_k \phi^-_k (x) \tvar^-_k (y)^*
\sim \sum_k \bar \phi^-_k (x) \bar{\tvar}^-_k(y)^* =
\bar{G}_{\lambda }^-,
$$
where $\bvarphi_k^-$, $\bar{\tvar}_k^-$ denote solutions of the asymptotic
constant-coefficient ODE of the eigenvalue and adjoint eigenvalue
equations at $-\infty$ and $\bar{G}_{\lambda }^-$ denotes the
associated (constant-coefficient) resolvent kernel.  Thus,
Proposition \thmref{kochelformulae} has the important and 
intuitively appealing implication that far-field behavior 
reduces to that of the constant coefficient case.

\medskip
3.  The fact that $d^+=d^-$ may be deduced without calculation 
from the fact that the constant-coefficient formula \eqnref{cc2}
satisfies the resolvent kernel equation \eqnref{glambda}, but not the 
correct boundary conditions of decay at $\pm \infty$.  
The continuous $d^\pm$ correction satisfies the {homogeneous}
eigenvalue equation, serving only to enforce decay at $+\infty$.
As noted already, this term decays at $-\infty$, as do the remaining,
constant-coefficient terms $\sum_k \phi_k^-\tilde \phi_k^-$.
\bigskip

\newsection {\bf Section 4. High frequency expansion.}
\sectionnumber=4 
\theoremnumber=0
\equationnumber=0
\smallskip
\TagsOnLeft

We now turn to the crucial estimation of the resolvent kernel in
the high frequency regime, as $|\lambda|\to 0$ with $\R \lambda \ge -\eta$.
According to the usual duality between spatial and frequency variables,
large frequency $\lambda$ corresponds to small time $t$, 
or {\it local behavior} in space and time.
Therefore, we begin by examining the frozen-coefficient equations
at each value $x_0$ of $x$:
$$
U_t=L_{x_0}U:=
-A(x_0)U_x+ Q(x_0)U.
\eqnlbl{Lx0}
$$
As usual, we shall for clarity state first the result in the case that
$(df,dg)^t$ is strictly hyperbolic, 
indicating afterward the extension to the general case.

\proclaim{Lemma \thmlbl{infexpansion}}
Let $(df,dg)^t(\bar u,\bar v)(x_0)$ have distinct, real eigenvalues.
Then, for $|\lambda|$ sufficiently large, the eigenvalue equation
$(L_{x_0}-\lambda)W=0$ associated with the frozen-coefficient
operator $L_{x_0}$ at $x_0$ has a basis of $N=n+r$ solutions
$$
\{\bar{\phi}_{1}^+, \dots, \bar{\phi}_{k}^+,
\bar{\phi}_{k+1}^-,\dots, \bar{\phi}_{N}^-\}(x;\lambda, x_0) 
$$
that are analytic in $1/\lambda$, satisfying 
$\bar{\phi}_j^\pm= e^{\mu_j(\lambda,x_0)x}V_j(\lambda,x_0)$, with
$$
\aligned
\mu_j(\lambda,x_0)&= -\lambda/a_j(x_0) -\eta(x_0)/a_j(x_0) 
+ \Cal{O}(1/\lambda),\\
V_j(\lambda,x_0)&= 
r_j(x_0)/a_j(x_0) + \Cal{O}(1/\lambda),
\endaligned
\eqnlbl{frozenexp}
$$
$a_j(\cdot)$, $\eta_j(\cdot)$, and $r_j(\cdot)$ as defined in 
Section 1.3.
Likewise, the adjoint eigenvalue equation 
$(L_{x_0}-\lambda)^*\tilde W=0$ has a basis of solutions
$
\{\bar{ \tilde \phi}_{1}^-, \dots, \bar{ \tilde \phi}_{k}^-,
\bar{ \tilde \phi}_{k+1}^+,\dots, \bar{ \tilde \phi}_{N}^+\}(x;\lambda,x_0)
$
satisfying
$\bar{\tilde \phi}_j^\pm= 
e^{-\mu_j(\lambda,x_0)x}\tilde{V}_j(\lambda,x_0)$, with
$$
\tilde{V}_j(\lambda,x_0)= l_j(x_0) + \Cal{O}(1/\lambda),
\eqnlbl{frozendualexp}
$$
$l_j(\cdot)$ as defined in Section 1.3.

The expansions \eqnref{frozenexp}--\eqnref{frozendualexp} hold
also in the nonstrictly hyperbolic case,
with $\bar \phi_j$, $\bar{\tilde \phi}_j$ and $r_j$, $l_j$ now
denoting $n\times m_j$ blocks, and $\mu_j$, $\eta_j$ denoting
$m_j\times m_j$ matrices, $j=1,\dots, J $,
where $m_j$ is the multiplicity of eigenvalue $a_j$ of $A(x_0)$.
\endproclaim

{\bf Proof.}
Similarly as in the proof of Lemma \thmref{frozen},
this follows by inversion of the corresponding 
expansions about $\xi=\infty$ of the eigenvalues $\lambda_j(\xi)$
and eigenvectors $V_j(\xi)$ of the frozen-coefficient Fourier symbol
\eqnref{frozenx}, using the basic relation $\mu_j=i \xi$ relating
characteristic and dispersion equations \eqnref{char2} and \eqnref{disp}.
The expansions for the Fourier symbol are carried out in Appendix A2.
\myqed

{\bf Example \thmlbl{JXexp}.} 
(Jin--Xin model) For model \eqnref{JX}, everything is
almost explicit and expansion of the eigenvalues and eigenvectors 
can be computed directly.
Recalling that $\mu^\pm=\frac\alpha{2a^2}\mp\frac{\sigma(\lambda)}{2a^2}$
where $\sigma(\lambda)=\sqrt{\alpha^2+4a^2\lambda+4a^2\lambda^2}$,
$\alpha=dh(\bar u(x_0))\in \BbbR^{n\times n}$,
we find that
 $$
  \mu^\pm(x;\lambda)=
  \mp\frac \lambda a + \frac {\alpha(x)\mp a}{2a^2}
   + O\left(\frac 1{\lambda}\right)\qquad \lambda\to+\infty,
 $$
where $\mu_\pm$, $\eta_\pm$ are $n\times n$ blocks.
Similar explicit expansions for $V_j$, $\tilde V_j$ can be computed 
using the formulae
$V_\pm=(\mu_\pm, -\lambda)^*$ and 
$\tilde V_\pm=(\pm 1/\lambda\sigma(\lambda))
(-a^2\mu_\pm,  \lambda)^*$
of Examples \thmref{JXN} and \thmref{JXcc}.
\myqed

\medskip

By \eqnref{cc2}, Lemma \thmref{infexpansion} 
gives an expression for the constant-coefficient resolvent kernel 
$G_{\lambda; x_0}$ at $x_0$ of, to lowest order:
$$
G_{\lambda; x_0}\sim
\cases
   \sum_{j=1}^k 
a_j(x_0)^{-1}
e^{(-\lambda/a_j -\eta_j/a_j)(x_0)(x-y)}r_j(x_0)l_j^{t}(x_0)
    \qquad &x<y,\\
   -\sum_{j=k+1}^{N} 
a_j(x_0)^{-1}
e^{(-\lambda/a_j -\eta_j/a_j)(x_0)(x-y)}r_j(x_0)l_j^{t}(x_0)
\qquad &x>y,
                 \endcases 
\eqnlbl{glambdax0}
 $$
\medskip

Our main result of this section, and perhaps the main technical
challenge of this paper,
will be to establish on an appropriate unbounded 
subset of the resolvent set $\rho(L)$
that the variable coefficient solutions $\phi_j^\pm$,
$\tilde \phi^\pm_j$ ``track'' their 
frozen coefficient analogs to order $1/|\lambda|$ in both
direction and (exponential) growth rate, 
yielding order $1/|\lambda|$ estimates for the resolvent kernel.
More precisely, define 
$$
\Omega:=
\{\lambda: -\eta_1 \le  \R \lambda \},
\eqnlbl{Omega}
$$
with $\eta_1>0$ sufficiently small that 
$\Omega \setminus B(0,r)$ is compactly contained 
in the set of consistent splitting $\Lambda$
(defined in \eqnref{Lambda}),
for some small $r>0$ to be chosen later;
this is possible, by Lemma \thmref{frozen}(i).

\proclaim{Proposition \thmlbl{highbounds}}
Assume that there hold (H0)--(H4), plus strict hyperbolicity 
of $(df,dg)^t$.
Then, for any $r>0$, and $\eta_1=\eta_1(r)>0$ chosen sufficiently
small, there holds $\Omega \setminus B(0,r)\subset \Lambda \cap \rho(L)$.
Moreover, for $R>0$ sufficiently large, there holds on
$\Omega \setminus B(0,R)$ the decomposition
$$
G_{\lambda}(x,y)= H_\lambda(x,y) +\Theta_\lambda(x,y),
\eqnlbl{highfreq}
$$
$$
H_{\lambda}(x,y):=
\cases
   -\sum_{j=k+1}^{N} 
a_j(y)^{-1}
e^{\int_y^x(-\lambda/a_j -\eta_j/a_j)(z)\, dz}r_j(x)l_j^{t}(y)
\qquad &x>y,\\
   \sum_{j=1}^k 
a_j(y)^{-1}
e^{\int_y^x(-\lambda/a_j -\eta_j/a_j)(z)\, dz}r_j(x)l_j^{t}(y)
    \qquad &x<y,\\
                 \endcases 
\eqnlbl{Hlambda}
 $$
$$
\Theta_\lambda(x,y)= 
\lambda^{-1}B(x,y;\lambda)+ 
\lambda^{-1}(x-y)C(x,y;\lambda)+ 
\lambda^{-2}D(x,y;\lambda),
\eqnlbl{errorterms}
$$
where
$$
B(x,y;\lambda)=
\cases
   \sum_{j=k+1}^{N} 
e^{-\int_y^x\lambda/a_j(z) \,dz }b^+_j(x,y)
\qquad &x>y,\\
   \sum_{j=1}^k 
e^{-\int_y^x\lambda/a_j(z) \,dz }b^-_j(x,y)
    \qquad &x<y,\\
                 \endcases 
\eqnlbl{errortermsB}
$$
$$
\aligned
C(x,y&;\lambda)=\\
&\cases
   -\sum_{i,j=k+1}^{N} 
\text{\rm mean}_{z\in [x,y]}
e^{-\int_y^z\lambda/a_i(s) ds
-\int_z^x\lambda/a_j(s) ds }c_{i,j}^+(x,y;z),
\qquad &x>y,\\
   \sum_{i,j=1}^k 
\text{\rm mean}_{z\in [x,y]}
e^{-\int_y^z\lambda/a_i(s) ds
-\int_z^x\lambda/a_j(s) ds }c_{i,j}^-(x,y;z),
    \qquad &x<y,\\
                 \endcases \\
\endaligned
\eqnlbl{errortermsC}
$$
with
$$
|b^\pm_j|, \, |c_{i,j}^\pm|\le  C e^{-\theta |x-y|}
\eqnlbl{unidecay}
$$
and
$$
D(x,y;\lambda)=
\cases
\Cal{O}(
e^{-\int_y^x\R\lambda/a_{k+1}(z) \,dz }
e^{-\theta |x-y|})
\qquad &x>y,\\
\Cal{O}(
e^{-\int_y^x\R\lambda/a_k(z) \,dz }
e^{-\theta |x-y|})
    \qquad &x<y,\\
                 \endcases 
\eqnlbl{errortermsD}
$$
for some uniform $\theta>0$ independent of $x$, $y$, $z$.

Likewise, there hold derivative bounds
$$
\aligned
(\partial/\partial x)\Theta_\lambda(x,y)&= 
\big(B^0_x(x,y;\lambda)+ 
(x-y)C^0_x(x,y;\lambda)\big) + 
\lambda^{-1} \big( B^1_x(x,y;\lambda)\\
&\quad+ 
(x-y)C^1_x(x,y;\lambda)
+(x-y)^2D^1_x(x,y;\lambda)
\big)
+ \lambda^{-2}E_x(x,y;\lambda)
\endaligned
\eqnlbl{xder}
$$
and
$$
\aligned
(\partial/\partial y)\Theta_\lambda(x,y)&= 
\big(B^0_y(x,y;\lambda)+ 
(x-y)C^0_y(x,y;\lambda)\big) 
+\lambda^{-1} \big( B^1_y(x,y;\lambda)\\
&\quad+ 
(x-y)C^1_y(x,y;\lambda)
+ (x-y)^2D^1_y(x,y;\lambda)
\big)
+ \lambda^{-2}E_y(x,y;\lambda),
\endaligned
\eqnlbl{yder}
$$
where $B^\alpha_\beta$ and $C^\alpha_\beta$ satisfy bounds 
of form \eqnref{errortermsB} and \eqnref{errortermsC},
$D^1_\beta$ now denotes the iterated integral
$$
\aligned
D^1_\beta(x,y&;\lambda)=\\
&\cases
   -\sum_{h,i,j=k+1}^{N} 
\text{\rm mean}_{y\le w\le z\le x}
e^{-\int_y^w\lambda/a_h(s) ds
-\int_w^z\lambda/a_i(s) ds 
-\int_z^x\lambda/a_j(s) ds }\\
\qquad \qquad \times d_{h,i,j}^{\beta,+}(x,y;z),
\qquad &x>y,\\
   \sum_{h,i,j=1}^k 
\text{\rm mean}_{x\le z\le w\le y}
e^{-\int_y^w\lambda/a_h(s) ds
-\int_w^z\lambda/a_i(s) ds 
-\int_z^x\lambda/a_j(s) ds }\\
\qquad \qquad \times
d_{h,i,j}^{\beta,-}(x,y;z),
    \qquad &x<y,\\
                 \endcases \\
\endaligned
\eqnlbl{errortermsD1}
$$
with
$ |d_{h,i,j}^{\beta,-}|\le C e^{-\theta |x-y|}$,
and $E_\beta$ satisfies a bound of form \eqnref{errortermsD}.

The bounds \eqnref{errorterms}--\eqnref{yder}
hold also in the nonstrictly hyperbolic case,
with \eqnref{Hlambda} replaced by
$$
H_{\lambda}(x,y)=
\cases
   -\sum_{j=K+1}^{J} 
a_j(y)^{-1}
e^{-\int_y^x \lambda/a_{j} (z)\, dz}r_j(x) \tilde \zeta_j(x,y) l_j^{t}(y)
\qquad &x>y,\\
   \sum_{j=1}^K
a_j(y)^{-1}
e^{-\int_y^x \lambda/a_j (z)\, dz}r_j(x) \tilde \zeta_j(x,y) l_j^{t}(y)
    \qquad &x<y,\\
                 \endcases 
\eqnlbl{nshighfreq}
 $$
where $\tilde\zeta_j(x,y)\in \BbbR^{m_j\times m_j}$
denotes a dissipative flow similar to $\zeta_j$ in \eqnref{diss}, 
but with respect to variable $x$: i.e.,
$$
d\tilde\zeta_j/dx=\eta_j(x) \tilde\zeta_j(x,y) /a_j(x), 
\quad \tilde \zeta_j(y,y)=I,
\eqnlbl{zetax}
$$
or, equivalently, $\tilde \zeta_j(x,y)=\zeta_j(y,\tau)$ for $z_j(y,\tau)=x$,
where $z_j$ as in \eqnref{char} denotes the characteristic 
path associated with $a_j$;
$r_j$, $l_j$ now
denote $n\times m_j$ blocks; and $\mu_j$, $\eta_j$ denote
$m_j\times m_j$ matrices, $j=1,\dots, K $,
$j=K+1,\dots, J$,
where $m_j$ is the multiplicity of eigenvalue $a_j$ of $A(x_0)$,
with
$$
a_1\le \cdots\le a_K<0<a_{K+1}\le \cdots\le a_J.
$$
(Note that we have not here assumed ($\Cal{D}$)).
\endproclaim

{\bf Remark.}
The three error terms in \eqnref{errorterms} may be recognized
as roughly analogous to the terms that would arise in the
constant-coefficient case through expansions 
\eqnref{frozenexp}--\eqnref{frozendualexp}; see also the
related Fourier transform analysis of [LZe,Ze.1--2].
The complicated ``mixed'' form of the second, averaged term 
(corresponding to $\Cal{E}_j$ in Proposition 8.15, Appendix A3.2.3) 
comes from our search for minimal regularity hypotheses on $f$, $g$, $q$.
At the expense of a further order of differentiability, allowing
us to approximately diagonalize the eigenvalue equation to one higher order, 
we could obtain instead a term of form 
$\lambda^{-1}(x-y)\tilde C(x,y;\lambda)$, where, for example, 
$$
\aligned
\tilde C(x,y;\lambda)&=
   \sum_{j=1}^k 
e^{\int_y^x(-\lambda/a_j -\eta_j/a_j)(z)\, dz}\\
&\times
\big(e^{\tilde c_j^-(x,y)(x-y)/\lambda}-1\big)/\big((x-y)/\lambda \big)
r_j(x)l_j^t(y)\\
\endaligned
\eqnlbl{altform}
$$
for $x<y$ and symmetrically for $x>y$,
precisely mimicking that of the constant-coefficient case. However,
this distinction is unimportant for our analysis, in which the precise
form of these terms plays no role; indeed, the (variational) form of 
\eqnref{errortermsC} is somewhat more convenient for our purposes.
\medskip

The proof of Proposition \thmref{highbounds}, and its nonstrictly
hyperbolic analog (stated at the end of the section)
will be the work of the remainder of the section.

\medskip

{\bf 4.1. General Framework.}  
We have to solve, approximately, the eigenvalue equation
$$
(AW)'= (Q-\lambda)W
$$
as $|\lambda|\to \infty$ within the region $\Omega$ defined above:
or, setting $Z=AW$, the equation
$$
Z'=(Q-\lambda)A^{-1}Z.
\eqnlbl{Zeqn}
$$
(As noted in [ZH], such an invertible transformation is a trick to
reveal the reduced regularity requirements on coefficients of 
divergence-form operators).

Following standard procedure (see, e.g., [AGJ,GZ,ZH,Z.4]),
we perform the rescaling $\tilde x = |\lambda|x$ (suggested by 
the leading order behavior of $\mu(x_0)$, \eqnref{frozenexp})
to obtain the perturbation equation
$$
Z'=(A_0(\tilde x) + |\lambda|^{-1} A_1(\tilde x))Z,
\eqnlbl{pertproblem}
$$
where
 $$
\aligned
  A_0 &:= -\tilde \lambda A^{-1}(|\lambda|^{-1}\tilde x),\\
  A_1 (\tilde x) &:= QA^{-1} (|\lambda|^{-1}\tilde x),\\
\tilde \lambda &:=\lambda/|\lambda|.\\
\endaligned
  \eqnlbl{coeffs}
 $$
This expansion is to be considered as a
continuous family of one-parameter perturbation equations,
indexed by $\tilde \lambda \in S^1$.
We shall have to take some care, therefore, to ensure that all derived
bounds depend continuously on the index $\tilde \lambda$; 
this accounts for the rather meticulous description of error terms 
in the relevant Proposition 8.15 of Appendix A3.2.3, below.

We seek to develop corresponding expansions of the 
resolvent kernel in orders of $|\lambda|^{-1}$, 
or, equivalently,
by representation \eqnref{basicflow}, expansions of 
stable/unstable manifolds of \eqnref{pertproblem}, and the
reduced flows therein.
This we will accomplish by a two-stage process. 
First, we will reduce by a series of coordinate changes 
to an approximately block-diagonal system segregating spectrally 
separated (in particular, stable/unstable) modes, with formal
error of a suitably small order.
Second, we will convert the formal error into rigorous error bounds
using a refined version of the Tracking Lemma of [GZ,ZH,Z.1,Z.4].

Both of these steps will be carried out in a general
setting suitable for applications to it arbitrary
ordinary differential operators of nondegenerate type (i.e., for
which the principal part of the symbol has fixed order and
type), even including operators that do not have asymptotically
constant, but only smooth coefficients. 
Together, they give an algorithm for the systematic estimation 
of high-frequency behavior, to an order of accuracy depending
only on the smoothness of the coefficients of the operator in question.

\medskip
{\bf 4.2. Formal diagonalization procedure.}
More generally 
(fixing any indexing parameters such as $\tilde \lambda$ above), 
consider an ODE
 $$
  W'=A(\e x,\e) W\qquad (x\in\BbbR),
  \eqnlbl{generalode}
 $$
with small parameter $\varepsilon$, where $A$ has formal Taylor expansion
 $$
  A(y,\e)=\sum_{k=0}^p \e^k A_k(y) + O(\e^{p+1}).
  \eqnlbl{seriesA}
 $$
(Note that the case considered in \eqnref{pertproblem}
corresponds to linear dependence with respect to $\e$.)
Assume:

(h1) \quad $A$ is $C^{p+1}$ in $y,\e$, or, if
$A_0,\dots,A_r$ are constant, then $(\partial/\partial \e)^{r+1} A$ is 
$C^{p-r}$ in $y,\e$.  In either case, all derivatives considered, in
both $y$ and $\e$, are uniformly bounded for all $y\in\BbbR$, $\e\le \e_0$, 
for some $\e_0>0$.

(h2) \quad $A_0(y)$ is block-diagonalizable to form
$$
D_0(y)=\text{\rm diag }\{d_{0,1},\dots,d_{0,s}\},
\eqnlbl{blockform}
$$
$d_j\in \Bbb{C}^{n_j\times n_j}$, 
with spectral separation between blocks 
(i.e., complex distance between their eigenvalues) 
uniformly bounded below by some
constant $\gamma >0$, for all $\e\le \e_0$, $x\in \BbbR$.

We begin by showing that \eqnref{generalode} may be approximately
block-diagonalized to order $\varepsilon^{p}$.

\proclaim{Proposition \thmlbl{diag}}
Given (h1)--(h2), there exists a uniformly well-conditioned
change of coordinates $W=T \tilde W$ such that
 $$
  \tilde W' = D(\e x,\e)\tilde W + O(\e^{p+1})\tilde W
  \eqnlbl{diagform}
 $$
uniformly for $x\in \BbbR$, $\e\le \e_0$, with 
 $$
  D(y,\e)=\sum_{k=0}^p \e^k D_k(y) + O(\e^{p+1})
  \eqnlbl{seriesD}
 $$
and
 $$
  T(y,\e)=\sum_{k=0}^p \e^k T_k(y) + O(\e^{p+1}),
  \eqnlbl{seriesT}
 $$
where $D_j$, $T_j\in C^{p+1-j}(y)$ for all $0\le j \le p$.
and each $D_j$ has the same block-diagonal form
$$
D_j(y)=\text{\rm diag }\{d_{j,1},\dots,d_{j,s}\},
\eqnlbl{jblockform}
$$
$d_{j,k}\in \Bbb{C}^{n_k\times n_k}$, as does $D_0$. 
\endproclaim

Our proof of this result relies on the following linear algebraic fact.

\proclaim{Lemma \thmlbl{commutator}}
Let $d_1 \in \BbbC^{m_1\times m_1}$ and
$d_2 \in \BbbC^{m_2\times m_2}$ have norm
bounded by $C_1$ and respective spectra
separated by $1/C_2>0$.  Then, the matrix
commutator equation
$$
d_1 X - X d_2= F,
\eqnlbl{meqn}
$$
$X\in \BbbC^{m_1\times m_2}$, 
is soluble for all $F\in \BbbC^{m_1\times m_2}$,
with $|X|\le C(C_1,C_2) |F|$. 
\endproclaim

{\bf Proof.}
Consider \eqnref{meqn} as a matrix equation $\Cal{D}\Cal{X}=\Cal{F}$
where $\Cal{X}$ and $\Cal{F}$ are vectorial representations of
the $m_1\times m_2$ dimensional quantities $X$ and $F$, and $\Cal{D}$
is the matrix representation of the linear operator (commutator)
corresponding to the lefthand side.
The spectrum of $\Cal{D}$ is readily seen to be the difference
$$
\{\alpha_1-\alpha_2: \, \alpha_j \in \sigma(d_j)\},
$$
between the spectra of $d_1$ and $d_2$,
with associated eigenvectors given by $r_1 l_2^*$, where
$r_1$ is the right eigenvector associated with $\alpha_1$ and
$l_2$ is the left eigenvector associated with $\alpha_2$.
(Here as elsewhere, $*$ denotes adjoint, or conjugate transpose of
a matrix or vector).
By assumption, therefore, the spectrum of $\Cal{D}$ has modulus uniformly
bounded below by $1/C_2$, whence the result follows.
\myqed

{\bf Proof of Proposition \thmref{diag}.}
Substituting $W=T\tilde W$ into \eqnref{diagform} and rearranging, we obtain
$$
\aligned
\tilde W'&= (T^{-1}AT - T^{-1}T')\tilde W\\
&=
(T^{-1}AT - \varepsilon T^{-1}T_y)\tilde W,\\
\endaligned
$$
yielding the defining relation
$$
(T^{-1}AT - \varepsilon T^{-1}T_y)= D
\eqnlbl{defrel}
$$
for $T$.

By (h2), there exists a uniformly well-conditioned family of matrices
$T_0(x)$ such that $T_0^{-1}A_0 T_0=D_0$, $D_0$ as in \eqnref{blockform};
moreover, these may be chosen with the same regularity in $y$ as $A_0$
(i.e., the full regularity of $A$).
Expanding 
$$
T^{-1}(y,\e)=\left(I - \big(\sum_{j=2}^p \e^j  T_0^{-1} T_j\big)
+ \big(\sum_{j=2}^p \e^j T_0^{-1} T_j\big)^2 - \dots \quad
\right)T_0^{-1}
\eqnlbl{Neumann}
$$
by Neumann series, and matching terms of like order $\varepsilon^j$,
we obtain a hierarchy of systems of linear equations of form:
$$
D_0 T_0^{-1} T_j- T_0^{-1}T_j D_0 - F_j= D_j,
\eqnlbl{hierarchy}
$$
where $F_j$ depends only on $A_k$ for $0\le k\le j$ and $T_k$, 
$(d/dy)T_k$ for $0\le j-1$.

On off-diagonal blocks $(k,l)$, \eqnref{hierarchy} reduces to
$$
d_l [T_0^{-1} T_j]^{(k,l)}- [T_0^{-1}T_j]^{(k,l)} d_k = F_j^{(k,l)},
\eqnlbl{offdiag}
$$
uniquely determining $[T_0^{-1} T_j]^{(k,l)}$, by Lemma 
\thmref{commutator}.  On diagonal blocks $(k,k)$,  we are free to
set $[T_0^{-1} T_j]^{(k,k)}$, whereupon \eqnref{hierarchy}
reduces to
$$
[D_j]^{(k,k)}= - [F_j]^{(k,k)},
\eqnlbl{diag}
$$
determining the remaining unknown $[D_j]^{(k,k)}$.

Thus, we may solve for $D_j$, $T_j$ at each successive stage in a
well-conditioned way.  Moreover, it is clear in the general
case that $A_0\not \equiv \text{\rm constant}$ that the regularity
of $D_j$, $T_j$ is as claimed, since we lose one order of regularity
at each stage, through the dependence of $F_j$ on derivatives of
lower order $T_k$, $k<j$.
In the case that $A_0,\dots,A_s$ are constant in $y$, we may
choose also $T_0,\dots,T_s$ to be constant in $y$, again recovering
the claimed regularity.
\myqed

{\bf Remark \thmlbl{oldway}.}
Proposition \thmref{diag} may be verified in more indirect fashion
without using Lemma \thmref{commutator}, by
expressing $T$ instead as a product $T_0 T_1 \cdots T_p$
and solving the resulting succession of nearby diagonalization problems.
\medskip

{\bf Remark \thmlbl{regreq}.}
For relaxation systems, we will need to expand
to formal error order $\varepsilon^2$, i.e., $p=1$.
Thus, we require $A$, $Q\in C^2$ in the general case, and $C^1$
for discrete kinetic models, in which case $A_0$ above is constant;
this translates to our requirements of $C^3$/$C^2$ on $f$, $q$ in (H0).
\medskip

{\bf Refined diagonalization procedure.}
Now, in addition to (h1)--(h2), assume:

(h3) \quad  $A_0(\cdot)$ decays exponentially in $y$.

\noindent 
(Note: this is the case for operators whose coefficients approach constants
at exponential rate as $y=\e x\to \pm \infty$, in particular for
relaxation profiles).
In this case, we point out that, by judicious choice of the
initial diagonalizing transformation $T_0$, we may arrange that
$D$ is given to first order by the ``frozen-coefficient'' approximation
$T_0^{-1}AT_0$.

\proclaim{Proposition \thmlbl{refined}}
Let (h1)--(h3) hold.  Then, 
in Proposition \thmref{diag}, it is possible to
choose $T_0(\cdot)$ in such a way that, also, $D_0 + \e D_1$
is determined simply by the block-diagonal part of
$$
T_0^{-1}(A_0 + \varepsilon A_1)T_0.
\eqnlbl{firstorder}
$$
\endproclaim

Reviewing \eqnref{defrel} and \eqnref{diag}, we find that \eqnref{firstorder}
is equivalent to the requirement that 
$ T_0^{-1}(d/dy)T_0 $ vanish on diagonal blocks.
The matrix $T_0$ is given by
 $$
  T_0= \pmatrix R_1, \dots, R_K \endpmatrix,
  \eqnlbl{R}
 $$
where the blocks $R_j$ are made up of columns spanning eigenspaces
of $A_0$ corresponding to the various blocks.  Let us write suggestively
 $$
  T_0^{-1}=\pmatrix
          L_1 \\
          \vdots \\
          L_K
         \endpmatrix,
  \eqnlbl{L}
 $$
where $L_j R_k$ is a zero block if $j\neq k$ and an identity block
if $j=k$.  
Then, Proposition \thmref{refined} follows from:

\proclaim{Lemma \thmlbl{goodman}}
Let (h1)--(h3) hold.
Then, given any initial prescription of $R_j(0)$, $L_j(0)$ at $y=0$,
there exists a unique choice of $R_j$, $L_j$ in \eqnref{R}--\eqnref{L} 
such that
 $$
  L_j R_j'= L_j'R_j \equiv 0,
  \eqnlbl{zero}
 $$
which, moreover, is uniformly smooth and bounded to the same order
of differentiability as is $A_0$.
\endproclaim

{\bf Proof.}
Let $\tilde R_j$, 
$\tilde L_j$ denote any smooth choice, and define
 $$
  R_j:=\tilde R_j \alpha_j \qquad
  L_j:=\alpha_j^{-1} \tilde L_j,
  \eqnlbl{}
 $$
where $\alpha_j$ are nonsingular $k_j\times k_j$
matrices to be determined,  $k_j$ denoting the dimension of the
$j$th block (i.e., $R_j$ dimension $N\times k_j$).
Defining $\alpha_j$ to solve the linear ODE
 $$
   \alpha_j'  = - \tilde L_j \tilde R_j'\, \alpha_j,\qquad\qquad
   \alpha_j(0) = I_{k_j},
  \eqnlbl{}
 $$
where $'$ in this single instance denotes $d/dy$,
we obtain \eqnref{zero} by direct calculation:
 $$
  L_j R'_j = L_j (\tilde R_j ' \alpha_j + \tilde R_j \alpha'_j )= 
  L_j (\tilde R_j ' \alpha_j - \tilde R_j \tilde L_j \tilde R'_j \alpha_j )=
  L_j ( I - \tilde R_j \tilde L_j) \tilde R'_j \alpha_j = 0.
  \eqnlbl{}
 $$
The other relation comes from $L_j R_j = I$, from which 
$L'_j R_j + L_j R'_j= 0$.
Uniform global bounds and smoothness then follow by linear
theory, and the fact that $\tilde R_j$, $\tilde L_j$ could
be chosen in the first place so that 
$\tilde L_j \tilde R_j' $ decays exponentially in $y$
at spatial infinity.
\myqed

{\bf Remark \thmlbl{goodcite}.}
The key result of Lemma \thmref{goodman} 
was first pointed out in the diagonalizable case by Goodman, in [Go.1].
There, the careful choice of diagonalizing coordinates
was essential in his analysis of stability by energy methods.
It is essential for us here, as well, in this different context.
Note that we do not need this result in the case that $A_0$
is identically constant, corresponding to the case that the
principal part of the symbol of the associated operator $L$ has
constant coefficients, for which we may simply choose $L$, $R$
constant as well.  This was done, for example, in the treatment
of scalar dispersive--diffusive systems in [HZ.2], for which the
dispersion coefficient was assumed constant, and can be done also
in the case of discrete kinetic (relaxation) models.
{\it Thus, the role of Lemma \thmref{goodman} in our analysis
is to give sharp estimates in the case of 
a variable-coefficient principal part.}
\medskip

{\bf Remark \thmlbl{reg2}.}
An important observation for us is that,
in both Propositions \thmref{diag} and \thmref{refined}, 
any additional regularity of $A$
with respect to $\e$ translates directly to regularity of
the error term $\Cal{O}(\e^{p+1})$.
More precisely, if $A\in C^{q}(\e \to C^{p+1}(y))$, $q\ge p+1$, then
the error term has regularity $C^{q-p-1}(\e \to C^0(y))$.
In case $A\in C^{\omega}(\e \to C^{p+1}(y))$, $q\ge p+1$, as often happens
(for example, for relaxation profiles),
the error term has regularity $C^{\omega}(\e \to C^0(y))$.
This observation is used in conjunction with Proposition 8.15 below;
specifically, it gives regularity with respect to $\delta$
(in this case $\e^k$) in the error term $\Theta_{11}$.

\medskip
{\bf Remark \thmlbl{analytic}.}
A second important observation is that even though we have
used in \eqnref{pertproblem} the (standard) homogeneous
rescaling $\tilde x =|\lambda| x$, which introduces terms such
as $\tilde \lambda= \lambda/|\lambda|$ that are singular at $\lambda=\infty$,
any associated singularities at infinity in the resulting expansion
are only apparent.  That is, they are {\it coordinate singularities}
that disappear when the diagonalized system is converted back
to original coordinates.  
More precisely, since $|\lambda| A$ is analytic at $\infty$ with 
respect to $\lambda^{-1}=\varepsilon/\tilde \lambda$, we may choose
the initial diagonalizing transformation $T_0$ to be analytic
in $\lambda^{-1}$ as well (indeed, it may be chosen independent
of $\lambda$; see \eqnref{coeffs}, with the effect that $|\lambda|D_0$
is also analytic in $\lambda^{-1}$.
If we make this choice, then it is readily verified that all 
further $T_j$, $D_j$ also enjoy these properties, since they are 
preserved by the (uniquely determined, from this point) 
diagonalization procedure. 
Clearly, this observation is not unique to the case of relaxation systems, 
but holds for general equations and rescalings.
An immediate consequence is that the rigorous bounds obtained 
from the formal diagonalization series
by the fixed point iteration described in Appendix A3.2, when
expressed in original coordinates, must also be analytic in $\lambda^{-1}$;
for, the fixed point iteration, being scale invariant, may just as well
be carried out in the original (i.e., unrescaled) coordinates to begin with,
whereupon we obtain analyticity as usual from 
the uniform convergence of iterates.

We point out that we could alternatively have performed the 
formal diagonalization to begin with in variable $\e:=\lambda^{-1}$, using
the modified perturbation formulation
$$
W'=\tilde \lambda (-A^{-1} + \lambda^{-1}QA^{-1}),
$$
and taking $A_j$ to depend on $(\tilde \lambda^{-1} \e x, \e)$.

\medskip

{\bf 4.3. Error bounds/application to relaxation systems.}  
We now return to our original problem, that of approximating
the stable/unstable manifolds at $+\infty$/$-\infty$ of ODE
\eqnref{generalode}, and the reduced flows therein.
In Appendix A3.2, we present a systematic method for
generating such approximations, given a system in the approximately diagonal
form \eqnref{diag}.
The main requirement in the theory is that the approximate flows 
$\Cal{F}_j^{y\to x}$ generated
by the blocks $d_j:=D^{(j,j)}$ be ``exponentially separated''
to some order $\e^k$, $k\le p$, in the sense that, as $x$ increases,
they are either uniformly exponentially decaying 
or else uniformly exponentially 
growing, i.e.: 
$$
\cases
|\Cal{F}_j^{y\to x}|\le Ce^{\e^k\theta_1 (x-y)},\\
|\Cal{F}_j^{x\to y}|\le Ce^{\e^k\theta_2 (y-x)},\\
\endcases
\eqnlbl{separation}
$$
for $x\ge y$, with
$$
\theta_1< 0 < \theta_2.
\eqnlbl{expsep}
$$

{\bf Observation \thmlbl{pertevalues}.} 
In the simplest case, that $A_0$ is diagonalizable, the verification
of \eqnref{expsep} is straightforward.  For, in this case, the block
flows are scalar, and exponential separation is equivalent to the
integral spectral condition
$$
\cases
\int_y^x \R d_j(z) \, dz
\le \e^k \theta_1 (x-y) + \Cal{O}(1),\\
\int_y^x \R d_j(z) \, dz
\ge \e^k \theta_2 (x-y) + \Cal{O}(1),\\
\endcases
\eqnlbl{integral}
$$
or, roughly, the satisfaction on average of a spectral
gap $\R d_j <\e^k\theta_1$/$\R d_j > \e^k\theta_2$.  Moreover, this
is clearly equivalent to the satisfaction of the corresponding condition
$$
\cases
\int_y^x \R(\sum_{l=0}^k \e^l d_{l,j}(z) \, dz
\le \e^k \theta_1 (x-y) + \Cal{O}(1),\\
\int_y^x \R(\sum_{l=0}^k \e^l d_{l,j}(z) \, dz
\ge \e^k \theta_2 (x-y) + \Cal{O}(1),\\
\endcases
\eqnlbl{kstage}
$$
at the $k$th stage of the diagonalization process; thus, we need only
check \eqnref{kstage} for each $1\le k\le p$.
In particular, it is sufficient that there hold the pointwise spectral gap
$$
\cases
\R(\sum_{l=0}^k \e^l d_{l,j}
\le \e^k\theta_1,\\
\R(\sum_{l=0}^k \e^l d_{l,j}
\ge \e^k\theta_2,\\
\endcases
\eqnlbl{pgap}
$$
on all except a set of measure $1/\alpha$, where it is violated by at
most $\alpha$; for, the corresponding error can be absorbed in the multiplying
constant $C$ of \eqnref{integral}--\eqnref{kstage}.

In the case of strictly hyperbolic relaxation systems satisfying
(H0)--(H4), \eqnref{pgap} holds on $\Omega$ for $k=p=1$
for all $|y|\ge M$, or equivalently $|x|\ge M/\e$,
by Lemma \thmref{refined}
together with Lemma \thmref{infexpansion}, \eqnref{goodeta}, and
the exponential approach in $y$ of coefficients to their constant
states, Lemma \thmref{structure},
while on $|y|\le M$, or $|x|\le \e$, it is violated by at most order $\e$.
Thus, the entire theory of Appendix A3.2 applies, yielding the 
error bounds reported below.
\myqed

\medskip

{\bf Observation \thmlbl{nspertevalues}.} 
When $A_0$ is only block-diagonalizable, the verification
of \eqnref{expsep} is more involved.  
Since $d_j$ is now a matrix, \eqnref{integral} is neither
necessary nor sufficient.
Moreover, the expansion \eqnref{seriesD}
must be considered as a matrix perturbation problem,
and so the dependence of spectrum on coefficients $D_k$
is not necessarily monotone in $k$. 
Thus, \eqnref{separation}--\eqnref{expsep} are
in general not equivalent (as in the scalar case) 
to the corresponding conditions on the $k$th order approximate flow,
unless there hold additional properties such as diagonalizability
of $\sum_{l\le k-1} \e^l d_{l,j}$ (as implied, in particular, by {\it symmetry})
or special structure of $d_{l,j}$ relative to $d_{k-1,j}$, 
for all $k\le l\le p$.

On the other hand, suppose as in the present case that $d_{l,j}$ 
are {\it scalar} until step $k$, where splitting occurs. 
Then, the contribution of these lower order terms to the approximate
flow induced by $d_{k,j}$ may be split off as an exponential factor
$$
e^{\int_y^x \sum_{l=0}^{k-1} \e^l d_{l,j}(z)\, dz}
\eqnlbl{factor}
$$
times the approximate flow induced by $\e^k d_{k,j}$ alone, simplifying
the computation somewhat.
In particular, for nonstrictly hyperbolic relaxation systems satisfying
(H0)--(H4), with $k=p=1$, we obtain just
$e^{-\int_y^x \tilde \lambda/a_j(z)\, dz}$ times the flow induced
by $\e \eta(z)$, $\eta$ as defined in Section 1.3.
The latter flow may be seen to have an order $\e$ spectral separation,
by \eqnref{nsgoodeta} plus exponential convergence of $\eta_j$ to
$\eta_j^\pm$ as $y\to \pm \infty$ 
\footnote{
For example, we may apply the
Gap Lemma of Appendix A3.1 separately to cases $x$, $y\in (-\infty,0]$
and $x$, $y\in [0,+\infty)$, treating the case $y<0<x$ by the semigroup
property as the concatenation of cases $y<0=x$ and $y=0<x$. 
Alternatively, we may proceed as in Appendix A3.2 by 
applying a change of coordinates $w_j= P \tilde w_j$, 
$P=I + \Cal{O}(e^{-\theta \e x})$,
converting $\eta_j$ to positive/negative definite form for all $x\ge M/\e$).}, 
while the former has exponential
separation (in agreeing sense) of order $\R \tilde \lambda/\min_{j,y} |a_j(y)|$.
Combining these observations, we find that there is spectral separation
on all of $\Omega$.
Thus, the theory of Appendix A3.2 applies again, yielding the 
desired error bounds.
\myqed

{\bf Remark \thmlbl{generalexp}.}
More generally, we point out, in the case relevant to traveling waves 
that all $A_j$ converge exponentially in $y$ to constants $A_j^\pm$
as $y\to \pm \infty$,
that, provided that $d_{1,j}$ are scalar, and there is
neutral separation to zeroth order
$$
\cases
\R d_{1,j}\ge 0,\\
\R d_{1,j}\le 0,\\
\endcases
$$
for all $y$ (as is then a {\it necessary} condition for spectral separation
along rays pointing into the domain of consistent splitting),
then a sufficient condition for spectral separation of order 
$\e^k$ is {\it spectral gap} of order $\e^k$ in the limiting blocks
$d_j(\pm \infty)$ (also a necessary condition, by the Gap Lemma of Appendix
A3.1).  For, the diagonalizing construction described
in Proposition \thmref{diag} clearly preserves the property of
exponential decay to constant states.
Thus, making a constant change of coordinates such that $d_j^\pm$
becomes positive/negative definite of order $\e^k$, we have that
$d_j$ is positive/negative definite to order $\e^k +\Cal{O}(\e e^{\e x})$,
readily yielding \eqnref{separation}.
More precise accounting of the separation, as in the two observations just
above, requires more detailed consideration of the flows of different orders.
We remark that $d_{1,j}$ scalar is roughly equivalent to 
diagonalizability/constant multiplicity of $A_0$, 
since it may then be arranged that each block $d_j$ correspond to
a single eigenvalue of $A_0$.
{\it Thus, the requirement of spectral separation is in practice 
no restriction for applications to stability of traveling waves.}
\medskip

{\bf Proof of Proposition \thmref{highbounds}.}
We will carry out in detail the basic estimate \eqnref{errorterms},
indicating the extension to derivative bounds \eqnref{xder}--\eqnref{yder}
by a few brief remarks.
First, observe that in the modified coordinates $Z=AW$, 
\eqnref{basicflow} just becomes
 $$
  G_\lambda(x,y)=\cases
   A^{-1}(x)\Fyx\Pi^-_y
                                  \qquad &x<y,\\ \\
   -A^{-1}(x)\Fyx\Pi^+_y
                                  \qquad &x>y,
                        \endcases 
  \eqnlbl{zbasicflow}
 $$
where $\Pi_\pm$, $\Fyx$ now denote the corresponding objects
in $Z$-coordinates.

Applying the formal diagonalization procedure of
Proposition \thmref{diag} to \eqnref{pertproblem}, 
and choosing the special initialization of Lemma \thmref{goodman},
at the same time taking care as described in Remark \thmref{analytic}
to preserve analyticity with respect to $\lambda^{-1}$ in original
coordinates, we obtain the approximately diagonalized system
$$
\aligned
z_j'&=(\tilde \lambda a_j + \e \eta_j)z_j 
+ \e^2 \sum_{j,k}\tilde \Theta_{jk} z_k \\
&=|\lambda|^{-1}\left(
(\lambda a_j + \eta_j)z_j 
+ \lambda^{-1}\sum_{j,k} \Theta_{j,k}z_k \right),
\endaligned
$$
where error terms $\Theta_{j,k}$ are bounded and analytic 
in $\lambda^{-1}$ as $\lambda$ goes to infinity,
or, converting back to the original spatial coordinates:
$$
z_j' =\left( (\lambda a_j + \eta_j)z_j 
+ \lambda^{-1}\sum_{j,k} \Theta_{j,k}z_k \right).
\eqnlbl{originaldiag}
$$
Moreover, the approximatedly diagonalizing transformation $T$ is given by
$$
\aligned
T&=T_0 + \lambda^{-1} T_1\\
&=(r_1,\dots,r_J) + \lambda^{-1}T_1,
\endaligned
$$
where $r_j$ are as described in Proposition \thmref{highbounds}.

As described in Observations \thmref{pertevalues}--\thmref{nspertevalues}
just above, the approximate, block-diagonal flows 
$z_j' =( (\lambda a_j + \eta_j)z_j $ are exponentially separated,
with uniform rates
$$
\cases
C^{\int_y^x \R \lambda/a_K(z)\, dz} e^{-\theta |x-y|},\\
Ce^{-\int_y^x \R \lambda/a_{K+1}(z)\, dz} e^{-\theta |x-y|},\\
\endcases
\eqnlbl{srates}
$$
where $a_K<0<a_{K+1}$.
Applying the fixed point construction of Appendix A3.2.2(ii), we
find that the stable/unstable manifolds of \eqnref{originaldiag}
are {\it uniformly transverse} for $\lambda \in \Omega$, $|\lambda|$
sufficiently large, with associated projections $\Pi_\pm$ given by
$$
\Pi_+= \sum_{j=1}^{K} r_j l_j^* + \lambda^{-1} E_+ + \Cal{O}(\lambda^{-2})
\eqnlbl{Piplus}
$$
and
$$
\Pi_-= \sum_{j=K+1}^{J} r_j l_j^*
+ \lambda^{-1} E_- + \Cal{O}(\lambda^{-2}),
\eqnlbl{Piminus}
$$
for some continuous matrices $E_\pm(x)$.

From transversality, we obtain immediately that 
$\Omega \setminus B(0,R) \subset \rho(L)\cap \Lambda$, for $R$ sufficiently
large, as claimed.
Moreover, approximating the reduced flows on stable/unstable manifolds
using Proposition 8.15 of Appendix A3.2.3, and recalling representation
\eqnref{zbasicflow}, we obtain after rearrangement the claimed
estimates on the resolvent kernel $G_\lambda$, where error terms
$B$ correspond to error terms $E\pm$ in \eqnref{Piplus}--\eqnref{Piminus},
$C$ to term $\Cal{E}_1^\pm$ in (8.61)--(8.62) of Proposition 8.15,
and $D$ contains all remaining, $\Cal{O}(\lambda^{-2})$ errors:
for example, the third, $\Cal{O}(\delta^2/\eta^2)= 
\Cal{O}(\lambda^{-2})$ term in the righthand side of (8.61), the
$\Cal{O}(\lambda^{-2})$ terms in \eqnref{Piplus}--\eqnref{Piminus},
and cross-terms coming from the product of different order $\lambda^{-1}$
terms.

{\it Derivative estimates.}
To obtain derivative estimates \eqnref{xder}--\eqnref{yder}, we may
follow precisely the same procedure, but carrying out
to one further order the fixed point iteration at the stage of the 
reduced flow estimate described in Proposition 8.15:
that is, to order $\lambda^{-2}$, with error term $\Cal{O}(\lambda^{-3}$. 
(Note: this leads to the second-order iterated integral
\eqnref{errortermsD1}). This yields a Taylor expansion to the same
higher order for the $x$ and $y$ derivatives of the reduced flow,
within the {\it rescaled} spatial coordinates, and thus, in original
coordinates, to the desired $\lambda^{-1}$ order, with
$\Cal{O}(\lambda^{-2})$ error. 
\myqed

\bigskip


\newsection {\bf Section 5. Low frequency bounds.}
\sectionnumber=5 
\theoremnumber=0
\equationnumber=0
\smallskip
\TagsOnLeft

Our next object is the estimation of the
resolvent kernel in the remaining critical regime
$|\lambda|\to 0$,
corresponding to large time behavior of the Green's function $G$, 
or {\it global behavior} in space and time.
Here, the details of the system are not important,
and we may follow essentially the same treatment as in
the strictly parabolic case, combining aspects of [ZH,Z.2,Z.4].
From a global perspective, the structure of the linearized
equations is that of two nearly constant-coefficient regions
separated by a thin shock layer near $x=0$; accordingly, behavior is
essentially governed by the two, limiting far-field equations:
$$
U_t=L_\pm U,
\eqnlbl{Lpm}
$$
coupled by an appropriate transmission relation at $x=0$.

Similarly as in the previous section,
we begin by examining these limiting systems.
For simplicity of exposition, we make the provisional assumption:
\medskip

(A1) \quad 
The $r\times r$ matrices $dq A^{-1}_\pm(0,I_r)^t$ have distinct eigenvalues.

\medskip

\noindent
(Recall, \eqnref{relaxtype}, that these matrices govern behavior near
$(u_\pm,v_\pm)$ in traveling wave ODE \eqnref{zero ode}.)
Since none of our estimates will depend on the distance between
eigenvalues, this assumption is easily removed by a limiting argument.

\proclaim{Lemma \thmlbl{zeroexpansion}}
Assume that there hold (H1)--(H3) and (A1).  
Then, for $|\lambda|$ sufficiently small, the eigenvalue equation
$(L_\pm-\lambda)W=0$ associated with the limiting, constant-coefficient
operator $L_\pm$ has a basis of $n+r$ solutions
$$
\bar{W}_j^\pm= e^{\mu_j^\pm(\lambda)x}V_j(\lambda), 
\eqnlbl{barphi}
$$
$\mu_j^\pm$, $V_j^\pm$, analytic in $\lambda$, consisting of
$r$ ``fast'' modes
$$
\aligned
\mu_j^\pm&= \gamma^\pm_j + \Cal{O} (\lambda), \\
V_j^\pm&= A^{-1}_\pm S^\pm_j + \Cal{O} (\lambda), \\
\endaligned
\eqnlbl{fast}
$$
$S_j^\pm=(0,s_j^{\pm t})^t$
where $\gamma^\pm_j$, $s_j^\pm$ are eigenvalues and associated right 
eigenvectors of $dq A^{-1}_\pm(0,I_r)^t$ (equivalently, $\gamma_j^\pm$,
$S_j^\pm$ are nonzero eigenvalues and associated right eigenvectors
of $QA^{-1}$),
and $n$ ``slow'' modes
$$
\aligned
\mu_{r+j}^\pm (\lambda):=
-\lambda/a^{*\pm}_j + \lambda^2 \beta^{*\pm}_j/{a^{*\pm}_j}^3 + 
\Cal{O}(\lambda ^3),\\
 V_{r+j}^\pm(\lambda) := 
R^{*\pm}_{j}
+  \Cal{O} (\lambda),\\
\endaligned
\eqnlbl{slow}
$$
where $a_j^{*\pm}$, $R_j^{*\pm}$, and $\beta_j^*\pm$ are as 
in Proposition \thmref{zeng}.
Likewise, the adjoint eigenvalue equation 
$(L_\pm-\lambda)^*Z=0$ has a basis of solutions
$
\bar{ \tilde W}_{1}^\pm =
e^{-\mu_j^\pm(\lambda)x}\tilde{V}_j(\lambda)$, with
$$
\tilde{V}_{j}(\lambda)= \tilde T_j^{\pm} + \Cal{O}(\lambda)
\eqnlbl{dualfast}
$$
and
$$
\tilde{V}_{r+j}(\lambda)= L_j^* + \Cal{O}(\lambda),
\eqnlbl{dualslow}
$$
$\tilde V$ analytic in $\lambda$,
where $\tilde T_j$ are the left eigenvectors of $QA^{-1}_\pm$
associated with the nonzero eigenvalues $-\mu_j^\pm$, and
$L_j^{*\pm}$ are as in Proposition \thmref{zeng}.
\endproclaim

{\bf Proof.}  Equivalently, we must show that $\mu_j^\pm$, $V_j^\pm$
describe the eigenvalues and associated eigenvectors of
coefficient matrix $A^{-1}(Q-\lambda)_\pm$, and 
(see Lemma \thmref{dual})
that $\tilde V_j^\pm$ describe the left eigenvectors of $(Q-\lambda)A^{-1}_\pm$.
At $\lambda=0$, there are $r$ distinct nonzero eigenvalues 
$\mu_j^\pm=\gamma_j^\pm$, by assumption (A1), 
and an $n$-fold eigenvalue $\mu=0$ corresponding
to the $n$-dimensional kernel of $Q$.  
By spectral separation, the former vary analytically as $\lambda$ is varied, 
with Taylor expansions \eqnref{fast} and \eqnref{dualfast}.

The latter may be seen to vary analytically, with expansions
\eqnref{slow} and \eqnref{dualslow}, by inversion of the expansions
$$
\lambda_j(i\xi)=-i a_j^{*\pm}\xi -\beta_j^{*\pm}\xi^2 + \Cal{O}(\xi^3)
$$
carried out in Appendix A2 for the the dispersion curves near $\xi=0$, 
together with the fundamental relation $\mu=i\xi$ between roots
of characteristic and dispersion equations 
(recall \eqnref{char2}--\eqnref{disp} of Section 3.1).
Alternatively, one can carry out directly the associated
matrix bifurcation problem as in [ZH] to obtain the same result.
\myqed

Our main result of this section is then:

\proclaim{Proposition \thmlbl{lowbounds}}
Assume that there holds (H0)--(H4) and ($\Cal{D}$) (equivalently, (D1)--(D2)).
Then, for $r>0$ sufficiently small, the resolvent kernel $G_\lambda$ has
a meromorphic extension onto $B(0,r)$, 
which may in the Lax or overcompressive case be decomposed as
$$
G_\lambda=E_\lambda + S_\lambda + R_\lambda
\eqnlbl{ldecomp}
$$
where, for $y\le0$:
$$
E_\lambda(x,t;y)
:= \lambda^{-1}\sum_{a_k^{*-} > 0}
[c^{j,0}_{k,-}]
\frac{\partial}{\partial \delta_j}
\pmatrix
\bar u^\delta(x) \\
\bar v^\delta(x)
\endpmatrix
L_k^{*-t}
e^{(\lambda/a^{*\pm}_k - \lambda^2 \beta^{*\pm}_k/{a^{*\pm}_k}^3 )y},
\eqnlbl{Elambda}
$$
$$
S_\lambda(x,t;y):=
\sum_{a_k^{*-} > 0, \,  a_j^{*+} > 0} 
[c^{j,+}_{k,-}]R_j^{*+}  {L_k^{*-}}^t
e^{(-\lambda/a^{*+}_j + \lambda^2 \beta^{*+}_k/{a^{*+}_k}^3 )x
+(\lambda/a^{*-}_k - \lambda^2 \beta^{*-}_k/{a^{*-}_k}^3 )y}
\eqnlbl{Slambda1}
$$
for $y\le 0\le x$,
$$
\aligned
S_\lambda&(x,t;y):=\\
&\sum_{a_k^{*-}>0}R_k^{*-}  {L_k^{*-}}^t
e^{ (-\lambda/a^{*-}_k + \lambda^2 \beta^{*-}_k/{a^{*-}_k}^3 )(x-y)}
\\
&+ 
\sum_{a_k^{*-} > 0, \,  a_j^{*-} < 0} 
[c^{j,-}_{k,-}]R_j^{*-}  {L_k^{*-}}^t
e^{(-\lambda/a^{*-}_j + \lambda^2 \beta^{*-}_j/{a^{*-}_j}^3 )x
+(\lambda/a^{*-}_k - \lambda^2 \beta^{*-}_k/{a^{*-}_k}^3 )y}
\\
\endaligned
\eqnlbl{Slambda2}
$$
for $y\le x\le 0$,
and
$$
\aligned
S_\lambda&(x,t;y):=\\
&\sum_{a_k^{*-}<0}R_k^{*-}  {L_k^{*-}}^t
e^{ (-\lambda/a^{*-}_k + \lambda^2 \beta^{*-}_k/{a^{*-}_k}^3 )(x-y)}
\\
&+ 
\sum_{a_k^{*-} > 0, \,  a_j^{*-} < 0} 
[c^{j,-}_{k,-}]R_j^{*-}  {L_k^{*-}}^t
e^{(-\lambda/a^{*-}_j + \lambda^2 \beta^{*-}_j/{a^{*-}_j}^3 )x
+(\lambda/a^{*-}_k - \lambda^2 \beta^{*-}_k/{a^{*-}_k}^3 )y}
\\
\endaligned
\eqnlbl{Slambda3}
$$
for $ x\le y \le 0$,
and $R_\lambda=R_\lambda^E+ R_\lambda^S$
is a faster-decaying residual satisfying 
$$
\aligned
(\partial/\partial x)^\gamma &
(\partial/\partial y)^\alpha R^E_\lambda= 
\Cal{O}(e^{-\theta|x-y|})+
\sum_{a_k^{*-} > 0}
e^{-\theta|x|}
e^{(\lambda/a^{*\pm}_k - \lambda^2 \beta^{*\pm}_k/{a^{*\pm}_k}^3 )y}\\
&\qquad\times
\left( \lambda^{\alpha-1}\Cal{O}(e^{\Cal{O}(\lambda^3)y}-1)+ \lambda^{\alpha-1}
\Cal{O}(e^{\Cal{O}(\lambda^3)x}-1)+ 
\Cal{O}(\lambda^{\alpha}) 
\right),\\
\endaligned
\eqnlbl{RElambda}
$$
$$
\aligned
(\partial/&\partial x)^\gamma 
(\partial/\partial y)^\alpha R^S_\lambda(x,t;y)=\\
&\sum_{a_k^{*-} > 0, \,  a_j^{*+} > 0} 
e^{(-\lambda/a^{*+}_j + \lambda^2 \beta^{*+}_k/{a^{*+}_k}^3 )x
+(\lambda/a^{*-}_k - \lambda^2 \beta^{*-}_k/{a^{*-}_k}^3 )y}\\
&\times
\left( \lambda^{\gamma+\alpha}\Cal{O}(e^{\Cal{O}(\lambda^3)y}-1)+
\lambda^{\gamma+\alpha}\Cal{O}(e^{\Cal{O}(\lambda^3)x}-1)
+\lambda^{\alpha}\Cal{O}(e^{-\theta|x|})
+ \Cal{O}(\lambda^{1+\gamma}) 
\right)\\
\endaligned
\eqnlbl{RSlambda1}
$$
for $y\le 0\le x$,
$$
\aligned
(\partial&/\partial x)^\gamma 
(\partial/\partial y)^\alpha R^S_\lambda(x,t;y)=
\sum_{a_k^{*-}>0}
e^{ (-\lambda/a^{*-}_k + \lambda^2 \beta^{*-}_k/{a^{*-}_k}^3 )(x-y)}\\
&\times
\left( \lambda^{\gamma+\alpha}\Cal{O}(e^{\Cal{O}(\lambda^3)y}-1)+
\lambda^{\gamma+\alpha}\Cal{O}(e^{\Cal{O}(\lambda^3)x}-1)
+\lambda^{\alpha}\Cal{O}(e^{-\theta|x|})
+ \Cal{O}(\lambda^{1+\gamma}) 
\right)\\
&+ 
\sum_{a_k^{*-} > 0, \,  a_j^{*-} < 0} 
e^{(-\lambda/a^{*-}_j + \lambda^2 \beta^{*-}_j/{a^{*-}_j}^3 )x
+(\lambda/a^{*-}_k - \lambda^2 \beta^{*-}_k/{a^{*-}_k}^3 )y}\\
&\times
\left( \lambda^{\gamma+\alpha}\Cal{O}(e^{\Cal{O}(\lambda^3)y}-1)+
\lambda^{\gamma+\alpha}\Cal{O}(e^{\Cal{O}(\lambda^3)x}-1)
+\lambda^{\alpha}\Cal{O}(e^{-\theta|x|})
+ \Cal{O}(\lambda^{1+\gamma}) 
\right)\\
\endaligned
\eqnlbl{RSlambda2}
$$
for $y\le x\le 0$,
and
$$
\aligned
(\partial&/\partial x)^\gamma 
(\partial/\partial y)^\alpha R^S_\lambda(x,t;y)=
\sum_{a_k^{*-}<0}
e^{ (-\lambda/a^{*-}_k + \lambda^2 \beta^{*-}_k/{a^{*-}_k}^3 )(x-y)}
\\
&\times
\left( \lambda^{\gamma+\alpha}\Cal{O}(e^{\Cal{O}(\lambda^3)y}-1)+
\lambda^{\gamma+\alpha}\Cal{O}(e^{\Cal{O}(\lambda^3)x}-1)
+\lambda^{\alpha}\Cal{O}(e^{-\theta|x|})
+ \Cal{O}(\lambda^{1+\gamma}) 
\right)\\
&+ 
\sum_{a_k^{*-} > 0, \,  a_j^{*-} < 0} 
e^{(-\lambda/a^{*-}_j + \lambda^2 \beta^{*-}_j/{a^{*-}_j}^3 )x
+(\lambda/a^{*-}_k - \lambda^2 \beta^{*-}_k/{a^{*-}_k}^3 )y}
\\
&\times
\left( \lambda^{\gamma+\alpha}\Cal{O}(e^{\Cal{O}(\lambda^3)y}-1)+
\lambda^{\gamma+\alpha}\Cal{O}(e^{\Cal{O}(\lambda^3)x}-1)
+\lambda^{\alpha}\Cal{O}(e^{-\theta|x|})
+ \Cal{O}(\lambda^{1+\gamma}) 
\right)\\
\endaligned
\eqnlbl{RSlambda3}
$$
for $ x\le y \le 0$,
$\alpha+ \beta \le 1$,
with each $\Cal{O}(\cdot)$ term separately analytic in $\lambda$,
and where $[c^{j,0}_{k,\alpha}]$, $\alpha=-,0,+$, 
are constants to be determined later.
Symmetric bounds hold for $y\ge 0$.
Similar, but more complicated formulae hold in the undercompressive case
(see Remark 6.9, below). 
\endproclaim

\medskip
{\bf 5.1.  Normal modes (behavior in $x$).}
We begin by relating the normal modes of the variable coefficient
eigenvalue equation \eqnref{evalue} to those of the limiting, 
constant-coefficient equations \eqnref{Lpm}.

\proclaim{Lemma \thmlbl{gaprel}} 
Let (H0)--(H3) hold.  Then,
for $\lambda\in B(0,r)$, $r$ sufficiently small,
there exist solutions $W^{\pm}_j (x; \lambda)$ of
\eqnref{evalue}, $C^{2}$ in $x$ and analytic in $\lambda$,
satisfying  
$$
\eqalign{
W^{\pm}_j (x; \lambda)&= V^{\pm}_j(x;\lambda) e^{\mu^\pm_j x}\cr
(\partial/\partial \lambda)^k V^\pm_j(x;\lambda)&= 
(\partial/\partial \lambda)^k 
V^\pm_j(\lambda) +{\Cal O}(
e^{-\tilde \theta |x|} 
| V^\pm_j(\lambda)|), \quad  x \gtrless  0, 
}
\eqnlbl{gapbounds}
$$
for any $k\ge 0$ and $0<\tilde \theta < \theta$,
where $\theta$ is the rate of decay given in \eqnref{expdecay},
$\mu_j^\pm(\lambda)$, $V_j^\pm(\lambda)$ are as in 
Lemma \thmref{zeroexpansion} above, and $\Cal{O}(\cdot)$
depends only on $k$, $\tilde \theta$.
\endproclaim

{\bf Proof.}
With the exception of the regularity $C^2$ asserted in $x$, this
is a direct consequence of Lemma \thmref{zeroexpansion} together with
Proposition 8.1 of Appendix A3.1 (the Gap Lemma).
The additional regularity of $C^2$ vs. the $C^1$ asserted
in the Proposition is a consequence of the divergence form of
the eigenvalue equation, and may be deduced as in the previous
section by the change of coordinates \eqnref{Zeqn}.
\myqed

The bases $\phi_j^\pm$, $\psi_j^\pm$ 
of Section 3.1 may evidently be chosen from among $W_j^\pm$,
yielding an analytic choice of bases in $\lambda$, with
the detailed description \eqnref{gapbounds}.
It follows that the dual bases $\tilde \phi_j^\pm$, $\tilde \psi^\pm_j$ 
defined in Section 3.2 are also analytic in $\lambda$ and satisfy 
corresponding bounds with respect to the dual solutions 
\eqnref{dualfast}--\eqnref{dualslow}.
With this observation, we have immediately,
from Definition \eqnref{localevans} and Proposition \thmref{kochelformulae}:

\proclaim{Corollary \thmlbl{extensions}}  Let (H0)--(H3) hold.  Then, 
the Evans function $D(\lambda)$ admits an analytic extension onto $B(0,r)$,
for $r$ sufficiently small, and the resolvent kernel $G_\lambda(x,y)$
admits a meromorphic extension, with an isolated pole of finite multiplicity
at $\lambda=0$.
\endproclaim

Moreover, we have close control on the $(x,y)$ behavior of $G_\lambda$
through the spectral decomposition formulae of Proposition 
\thmref{kochelformulae}.  
For {\it slow, dual modes},
the bounds \eqnref{gapbounds} 
(in particular, the consequent bounds on first spatial derivatives) 
can be considerably sharpened,
provided that we appropriately initialize our bases at $\lambda=0$.
This observation will be quite significant in the Lax or overcompressive case.
Likewise, {\it fast-decaying forward modes} 
can be well-approximated near $\lambda=0$
by their representatives at $\lambda=0$, using only the basic bounds
\eqnref{gapbounds}.
These two categories comprise the modes determining behavior of
the Green's function to lowest order.

Specifically, due to the special, partially conservative structure 
of the underlying evolution equations, the adjoint eigenvalue equation 
\eqnref{dualevalue} at $\lambda=0$ admits an $n$-dimensional 
subspace of constant solutions 
$$
\tilde W=(\tilde w_1^t,\tilde w_2^t)^t 
\equiv(\text{\rm constant}, 0)^t;
\eqnlbl{constsoln}
$$
this is equivalent to the observation [L.2] that integral 
quantities in variable $u$ are conserved under time evolution 
for relaxation systems of form \eqnref{general}.
Thus, at $\lambda=0$, we may choose, by appropriate change of coordinates
if necessary, that slow-decaying dual modes ${\tilde \phi}^\pm_j$ and
slow-growing dual modes ${\tilde \psi}^\pm_j$ be identically constant, 
Note that this does not interfere with our previous choice in 
Lemma \thmref{gaprel},
since that concerned only the choice of limiting solutions 
${\overline W}_j^\pm$ of the asymptotic, constant coefficient 
equations at $x\to \pm \infty$,
and not the particular representatives $W_j^\pm$ that approach them
(which, in the case of slow modes, are specified only up to the addition
of an arbitrary fast-decaying mode).

{\bf Remark \thmlbl{fastgrowth}.}  The prescription of constant dual
bases just described requires that we reverse our previous approach, 
choosing dual bases first using the Gap Lemma, then defining forward 
bases using duality.
Alternatively, rewriting the forward eigenvalue equation at $\lambda=0$ as
$$
\pmatrix z_1\\z_2\\
\endpmatrix'=QA^{-1}Z=\pmatrix
0 \\ dq A^{-1}Z\\
\endpmatrix,
\eqnlbl{zode}
$$
where $Z=\pmatrix z_1\\z_2\endpmatrix:=AW$, 
we may observe that fast-decaying modes satisfy 
$$
\aligned
z_1&\equiv 0\\
z_2'&= dq A^{-1}\pmatrix 0\\z_2\endpmatrix,\\
\endaligned
\eqnlbl{linode}
$$
i.e., the linearization of the reduced traveling wave ode \eqnref{gode}.
By duality relation \eqnref{duality}, requiring slow dual
modes to be constant is equivalent to choosing fast-growing (forward) modes 
as well as fast-decaying modes from among the solutions of \eqnref{linode}.
(Recall: though fast-decaying modes are uniquely determined as a subspace,
fast-growing modes are only determined up to the addition of faster decaying
modes.)
\medskip

\proclaim{Lemma \thmlbl{refinedderiv}}  Assuming (H0)--(H3),
with the above choice of bases at $\lambda=0$, 
and for $\lambda \in B(0,r)$
and $r$ sufficiently small, 
slow modes ${\tilde W}^\pm_j$ satisfy
$$
{\tilde W}^\pm_j(y;\lambda) = e^{-\mu(\lambda) y}\bar{\tilde V}^\pm_j(0) 
+  \lambda \tilde \Theta_j^\pm(y; \lambda),
\eqnlbl{slowtaylor}
$$
where
$$
\aligned
|\tilde \Theta_j^\pm|&\le
C|e^{-\mu(\lambda) y}|\\
|(\partial/\partial y)\tilde \Theta_j^\pm|
&\le C|e^{-\mu(\lambda) y}|(|\lambda|+ e^{-\theta|y|}),
\endaligned
\eqnlbl{refderiv}
$$
$\theta>0$,
as $y\to \pm \infty$, and $\bar {\tilde V}_j^\pm\equiv \text{\rm constant}$.
\footnote{
In particular, this includes 
$|(\partial/\partial y)\tilde{W}_j^\pm(y;\lambda)|\le 
C|\lambda e^{\mu_j^\pm(\lambda)y}|$
in place of the general spatial-derivative bound 
$
|(\partial/\partial y)\tilde{W}_j^\pm(y;\lambda)|\le
C(|\mu \tilde{V}_j^\pm|+|(\partial/\partial y)\tilde{V}_j^\pm|) e^{\mu_j^\pm(\lambda)y}
\sim
C|e^{\mu_j^\pm(\lambda)y}|
$.
}
Similarly, fast-decaying (forward) modes $W_j^\pm$ satisfy
$$
W^\pm_j(x;\lambda) = W^\pm_j(x;0) + \lambda \Theta_j^\pm(x;\lambda),
\eqnlbl{fasttaylor}
$$
where
$$
|\Theta_j^\pm|, \,
|(\partial/\partial x)\Theta_j^\pm|
\le C e^{-\theta|x|}
\eqnlbl{refderiv2}
$$
as $x\to \pm \infty$,
for some $\theta>0$.
\endproclaim

{\bf Proof.}  
Applying the Gap Lemma to the augmented variables
$$
\aligned
\tilde{\Bbb{W}}_j^\pm(y;\lambda):=\pmatrix \tilde{W}_j^\pm 
\\ \tilde{W}_j^{\pm'} \endpmatrix(y;\lambda)
&=: e^{-\mu(\lambda) y}\tilde{\Bbb{V}}_j^\pm(y;\lambda)
\pmatrix \tilde{W}_j^\pm 
\\ \tilde{W}_j^{\pm'} \endpmatrix(y;\lambda)\\
&= e^{-\mu(\lambda) y}
\pmatrix \tilde{V}_j^\pm 
\\ -\mu_j^\pm\tilde{V}_j^{\pm} 
 +\tilde{V}_j^{\pm'} 
\endpmatrix(y;\lambda)\\
\endaligned
$$
and
$$
\aligned
\Bbb{W}_j^\pm(x;\lambda):=\pmatrix W_j^\pm \\ {W_j^\pm}' \endpmatrix(x;\lambda)
&=: e^{\mu(\lambda) x}\Bbb{V}_j^\pm(x;\lambda)\\
&= e^{-\mu(\lambda) y}
\pmatrix \tilde{V}_j^\pm 
\\ \mu_j^\pm{V}_j^{\pm} 
 +{V}_j^{\pm'} 
\endpmatrix(y;\lambda)\\
\endaligned
$$
we obtain bounds
$$
\eqalign{
\tilde{\Bbb{W}}^{\pm}_j (x; \lambda)&= 
\tilde{\Bbb{V}}^{\pm}_j(x;\lambda) e^{-\mu^\pm_j( \lambda)x}\cr
(\partial/\partial \lambda)^k \tilde{\Bbb{V}}^\pm_j(x;\lambda)&= 
(\partial/\partial \lambda)^k 
\tilde{\Bbb{V}}^\pm_j(\lambda) +{\Cal O}(
e^{-\tilde \theta |x|} 
| \tilde{\Bbb{V}}^\pm_j(\lambda)|), \quad  x \gtrless  0, 
}
\eqnlbl{auggapbounds}
$$
and
$$
\eqalign{
{\Bbb{W}}^{\pm}_j (x; \lambda)&= 
{\Bbb{V}}^{\pm}_j(x;\lambda) e^{\mu^\pm_j(\lambda) x}\cr
(\partial/\partial \lambda)^k {\Bbb{V}}^\pm_j(x;\lambda)&= 
(\partial/\partial \lambda)^k 
{\Bbb{V}}^\pm_j(\lambda) +{\Cal O}(
e^{-\tilde \theta |x|} 
| \tilde{\Bbb{V}}^\pm_j(\lambda)|), \quad  x \gtrless  0, 
}
\eqnlbl{auggapbounds}
$$
$\tilde \theta>0$, analogous to \eqnref{gapbounds}, valid
for $\lambda \in B(0,r)$,  where
$$
{\tilde{\Bbb{V}}}_j^\pm(\lambda)=
\pmatrix {\tilde{V}}(\lambda)\\
-\mu_j^\pm {\tilde{V}}(\lambda)\\
\endpmatrix
$$
and
$$
{{\Bbb{V}}}_j^\pm(\lambda)=
\pmatrix {{V}}(\lambda)\\
\mu_j^\pm {{V}}(\lambda).\\
\endpmatrix
$$

By Taylor's Theorem with differential
remainder, applied to $\tilde{\Bbb{V}}$, we have:
$$
{\tilde \Bbb{W}}_j^\pm(y,\lambda)=
e^{-\mu_j^\pm (\lambda)y}
\left({\tilde \Bbb{V}}_j^\pm(y;0) + 
\lambda (\partial/\partial \lambda)\tilde {\Bbb{V}}_j^\pm(y;0) 
+\frac{1}{2}
\lambda^2 (\partial/\partial \lambda)^2\tilde {\Bbb{V}}_j^\pm(y;\lambda_*) 
\right),
\eqnlbl{expB}
$$
for some $\lambda_*$ on the ray from $0$ to $\lambda$, where, recall,
$(\partial/\partial \lambda)\tilde {\Bbb{V}}_j^\pm(y;\cdot)$ and
$(\partial/\partial \lambda)^2\tilde {\Bbb{V}}_j^\pm(y;\cdot)$ are 
uniformly bounded in $L^\infty[0,\pm\infty]$ for $\lambda \in B(0,r)$. 
Together with the choice ${\tilde V}_j^\pm(y,0)\equiv \text{\rm constant}$,
this immediately gives the first (undifferentiated) bound in \eqnref{refderiv}.

Applying now the bound \eqnref{auggapbounds} with $k=1$, we may expand the
second coordinate of \eqnref{expB} as
$$
\aligned
(\partial/\partial y) {\tilde {W}}_j^\pm(y,\lambda)&=
e^{-\mu_j^\pm (\lambda)y}
\Big(
-\mu_j^\pm\tilde {V}_j^\pm(y;0) 
+\tilde {V}_j^{\pm'}(y;0) \\
&\qquad -\lambda\big(
 (\partial/\partial \lambda)(\mu_j^\pm
\tilde {{V}}_j^\pm)(0) 
+ \Cal{O}(e^{-\theta|y|})
\big)
+ \Cal{O}( \lambda^2)
\Big)\\
&=
e^{-\mu_j^\pm (\lambda)y}
\left(
-\lambda\big(
 (\partial/\partial \lambda)\mu_j^\pm(0)
\tilde {{V}}_j^\pm(0) 
+ \Cal{O}(e^{-\theta|y|})
\big)
+ \Cal{O}( \lambda^2)
\right),\\
\endaligned
\eqnlbl{expB}
$$
and subtracting off the corresponding Taylor expansion 
$$
\aligned
(\partial/\partial y)&\big( e^{-\mu_j^\pm(\lambda)y}{\tilde {V}}_j^\pm(y,0)
\big)\\
&=
\mu_j^\pm(\lambda) e^{-\mu_j^\pm(\lambda)y}{\tilde {V}}_j^\pm(y,0)\\
&=
e^{-\mu_j^\pm (\lambda)y}
\left(-\mu_j^\pm(0)\tilde {V}_j^\pm(y;0) 
-\lambda (\partial/\partial \lambda)\mu_j^\pm(0)
\tilde {{V}}_j^\pm(y;0) 
+ \Cal{O}( \lambda^2)
\right)\\
&=
e^{-\mu_j^\pm (\lambda)y}
\left( -\lambda (\partial/\partial \lambda)\mu_j^\pm(0)
\tilde {{V}}_j^\pm (y;0) 
+ \Cal{O}( \lambda^2)
\right),\\
\endaligned
\eqnlbl{expapprox}
$$
we obtain
$$
\aligned
(\partial/\partial y) {\tilde {\Theta}}_j^\pm(y,\lambda)&=
e^{-\mu_j^\pm (\lambda)y}
\left(
\lambda \Cal{O}(e^{-\theta|y|})
+ \Cal{O}( \lambda^2)
\right),\\
\endaligned
\eqnlbl{expB}
$$
as claimed.

Similarly, we may obtain \eqnref{fasttaylor}--\eqnref{refderiv2}
by Taylor's Theorem with differential remainder applied to 
$W=e^{\mu(\lambda) x}$, and the Leibnitz calculation
$$
(\partial/\partial \lambda)W=(d\mu/d\lambda)x e^{\mu x}V+ e^{\mu x}
(\partial/\partial \lambda)V,
$$
together with the observation that $|xe^{-\theta|x|}\le Ce^{-\theta|x|/2}$.
\myqed

{\bf Remark \thmlbl{hoz.4}.}
For a more explicit derivation of \eqnref{slowtaylor}--\eqnref{refderiv}
in the scalar case, see Lemma 2.2, (2.6)--(2.7) in [HoZ.2].
(Note: in [HoZ.2], the one-dimensional case corresponds to $\xi\equiv 0$).
\medskip

This leaves only the problem of determining
behavior in $\lambda$ through the study of coefficients $M$, $d^\pm$.
To this end, we make the following further observations in the
Lax or overcompressive case, generalizing the corresponding observation
of Lemma 4.30, [Z.4], in the strictly parabolic case:

\proclaim{Lemma \thmlbl{4.7}}  For transverse ($\gamma \neq 0$) Lax and
overcompressive shocks, with the above-specified choice of basis at $\lambda=0$,
fast-growing modes $\psi^+_j$, $\psi^-_j$
are fast-{\it decaying} at $-\infty$, $+\infty$, respectively.  
Equivalently, fast-decaying modes $\tilde \psi^+_j$, 
$\tilde \psi_j^-$ are fast-growing at $-\infty$, $+\infty$:
i.e., the only bounded solutions of the adjoint eigenvalue equation
are constant solutions.
\endproclaim

{\bf Proof.} Noting that the manifold of solutions of the
$r$-dimensional ODE \eqnref{linode} 
that decay at either $x \to +\infty$ or $x \to -\infty$ 
is by transversality, together with \eqnref{indexrelation},
exactly $d_+ + d_- -\ell = r$, we see that 
{\it all} solutions of this ODE in fact decay at least at one infinity.
This implies the first assertion, by the alternative characterization
of our bases described in Remark \thmref{fastgrowth} above; the second
follows by duality, \eqnref{duality}.
Since the manifold of fast-decaying dual modes is uniquely determined,
independent of the choice of basis, this implies that 
nonconstant dual modes that are bounded at one infinity must blow up 
at the other, hence bounded solutions must be constant as claimed.
\myqed

{\bf Remark \thmlbl{invariants}.}
As noted in [LZ.2,ZH,Z.2,Z.4], the property of Lax and overcompressive
shocks that the adjoint eigenvalue
equation has only constant solutions has the interpretation
that the only $L^1$ time-invariants of the evolution of the one-dimensional
linearized equations about $(\bU,\bar v)(\cdot)$ are given by {\it
conservation of mass}, see discussion [LZ.2].  This distinguishes
then from undercompressive shocks, which {\it do} have additional
$L^1$ time-invariants.  The presence of additional
time-invariants for undercompressive shocks 
has significant implications for their behavior, 
(see discussions [LZ.2] Section 3 and [ZH], Section 10):
in particular, the time-asymptotic location of a perturbed
undercompressive shock (generically) evolves nonlinearly, 
and is not determinable by any linear functional 
of the initial perturbation.
By contrast, the time-asymptotic location 
of a perturbed (stable) Lax or overcompressive shock may be 
determined by the mass of the initial perturbation alone.
\medskip

{\bf 5.2. Scattering coefficients (behavior in $\lambda$).}
We next turn to the estimation of scattering coefficients
$M_{jk}$, $d^\pm_{jk}$.
Consider coefficient $M_{jk}$.
Expanding \eqnref{4.29} using Kramer's rule, and setting
$z=0$, we obtain 
$$
M^\pm_{jk}= D^{-1} C^\pm_{jk},
\eqnlbl{4.33}
$$
where
$$
C^+:=
(I,0) 
\pmatrix
\Phi^+ & \Phi^-\\
{\Phi^+}{'}& {\Phi^-}{'}
\endpmatrix^{\adj}
\pmatrix
{\Psi^-}\\
{\Psi^-}'\\
\endpmatrix_{|_{z=0}} 
\eqnlbl{4.34}
$$
and a symmetric formula holds for $C^-$.
Here, $P^{\adj}$ denotes the {\it adjugate matrix} of a matrix $P$, 
i.e. the transposed matrix of minors. 
As the adjugate is polynomial in the entries of the
original matrix, it is evident that
$|C^\pm|$ is uniformly bounded and therefore
$$
|M_{jk}^\pm|\le C_1|D^{-1}| \le
C_2 \lambda^{-\ell}
\eqnlbl{4.35}
$$
by ($\Cal{D}$), where $C_1$, $C_2>0$ are uniform constants.

However, the crude bound \eqnref{4.35} hides considerable
cancellation, a fact that will be crucial in our analysis.  
Relabel the $\{\varphi^\pm\}$ so that, at $\lambda =0$,
$$
\varphi^+_{N-j+1} \equiv \varphi^-_j = (\partial /\partial \delta_j)\pmatrix
\bar u^\delta\\ \bar v^\delta
\endpmatrix, \quad j=1,\cdots, \ell.
\eqnlbl{4.38}
$$
With convention \eqnref{4.38}, we have the sharpened bounds:

\proclaim{Lemma \thmlbl{4.5}} 
Let (H0)--(H4) and ($\Cal{D})$ (equivalently, (D1)--(D2))  hold,
and let $\phi_j^\pm$ be labeled as in \eqnref{4.38}.
Then, for $|\lambda|$ sufficiently small, there hold
$$
|M_{jk}^\pm|, |d_{jk}^\pm|, |\bar d_{jk}^\pm|
\le C
\cases
\lambda^{-1} &\hbox{for}\ j=1,\cdots,\ell,\\
1 &\hbox{otherwise,}
\endcases
\eqnlbl{4.40a}
$$
where $M^\pm$, $d^\pm$ are as defined in Proposition \thmref{kochelformulae}.
Moreover, 
$$
\text{Residue }_{\lambda= 0} \, M^+_{N-j+1,k}=
\text{Residue }_{\lambda= 0} \, d^\pm_{j,k},
\eqnlbl{res1}
$$
for $1\le j\le \ell$, all $k$.
\endproclaim

That is, blowup in $M_{jk}$ occurs to order $\lambda^{\ell-1}|D^{-1}|$ 
rather than $|D^{-1}|$, and, more importantly, only in 
(fast-decaying) {\it stationary modes} 
$(\partial/\partial \delta_j) (\bar u^\delta,\bar v^\delta)$.
Moreover, each stationary mode $(\partial/\partial \delta_j)(\bar u^\delta,
\bar v^\delta)$ has consistent scattering coefficients with dual mode
$\tilde \psi_k^-$, to lowest order, across all of the various representations
of $G_\lambda$ given in Proposition \thmref{kochelformulae}.

{\bf Proof.}  
Formula \eqnref{4.34} may be rewritten as
$$
C_{jk}^+=
\det 
\pmatrix
\varphi^+_1,\cdots, \varphi^+_{j-1},  \psi^-_k, \varphi^+_j,
\cdots, \varphi^+_n,
&\Phi^-\\ 
\\
{\varphi^+_1}{'},\cdots, {\varphi^+_{j-1}}{'},
 {\psi^-_k}{'}, {\varphi^+_j}{'}, 
\cdots, \varphi^+_n,
&{\Phi^-}{'}
\endpmatrix_{|_{z=0}}, 
\eqnlbl{4.42}
$$
from which we easily obtain the desired cancellation in $M^+=C^+D^{-1}$.  
For example, for $j>\ell$, we have
$$
\aligned
C_{jk}^+
&= \det 
\pmatrix
\varphi^+_1 + \lambda  \varphi^+_{1_\lambda } +\cdots, &\cdots,
&\varphi^+_n +\lambda
\varphi^+_{n_\lambda } + \cdots\\
\\
\cdots, &,\cdots, &{\varphi^+_n}' +\lambda {\varphi^+_{n_\lambda }}{'}
+\cdots
\endpmatrix \\
&= \Cal{O}(\lambda ^\ell) \le C|D|,
\endaligned
\tag\eqnalignlbl{4.43}
$$
yielding $|M_{jk}|=|C_{jk}||D|^{-1}\le C$ as claimed, by elimination
of $\ell$ zero-order terms, using linear dependency among fast modes
at $\lambda=0$.  For $1\le j\le \ell$, there is only an $\ell -1$ order
dependency, and we obtain instead the bound $|C_{jk}|\le C|\lambda|^{\ell-1}$,
or $|M_{jk}|=|C_{jk}||D|^{-1}\le C|\lambda|^{-1}$. 
The bounds on $|d_{jk}|$ and $|\bar d_{jk}|$ follow similarly;
likewise, \eqnref{res1} follows easily from the observation
that, since columns $\psi_j^+$ and $\psi_{N-j+1}^-$ agree,
and $\tilde \psi_j^-$ is held fixed, the expansions of the 
various representations by Kramer's rule yield determinants
which at $\lambda^{-1}$ order differ only by a transposition
of columns.
\myqed

In the Lax or overcompressive case, we can say a bit more:

\proclaim{Lemma \thmlbl{4.8}}  
Let (H0)--(H4) and ($\Cal{D}$) (equivalently, (D1)--(D2)) hold,
with $|\lambda|$ sufficiently small.
Then, for Lax and overcompressive shocks, 
with appropriate basis at $\lambda=0$ (i.e. slow dual modes
taken identically constant), 
there hold 
$$
 |M_{jk}|, |d_{jk}|, |\bar d_{jk}| \leq C
\eqnlbl{4.53} 
$$
if $\Tpsi_k$ is a fast mode, and
$$
 |M_{jk}|, |d_{jk}|, |\bar d_{jk}| \leq 
C |\lambda|
\eqnlbl{slow-fast} 
$$
if, additionally, $\varphi_j$ is a slow mode.
\endproclaim

As we shall see in the following section, this result has the
consequence that only slow, constant dual modes play a
role in long-time behavior of $G$.

{\bf Proof.} Transversality, $\gamma \neq 0,$ follows from
(D2), so that the results of Lemma \thmref{4.7} hold.  
We first establish the bound \eqnref{4.53}.
By Lemma \thmref{4.5}, we need only consider $j=1,\dots,\ell$,
for which $\phi_j^-=\phi_j^+$.
By Lemma \thmref{4.7}, all fast-growing modes $\psi^+_k$, $\psi_k^-$
lie in the fast-decaying manifolds $\Span
(\varphi_1^-, \cdots, \varphi_{i_-}^-)$, $\Span (\varphi_1^+, \cdots,
\varphi_{i_+}^+)$ at $-\infty$, $+\infty$ respectively.  It follows that in
the righthand side of \eqnref{4.42}, there is at $\lambda=0$ 
a linear dependency between columns 
$$
\varphi^+_1, \cdots, \varphi_{j-1}^+, \psi^-_k, \varphi^+_{j+1},
\varphi_{i_+}^+ \ \hbox{and}\  \varphi^-_j = \varphi_j^+
\eqnlbl{4.54}
$$
i.e. an $\ell$-fold dependency among columns
$$
\varphi^+_1, \cdots, \psi^-_k, \cdots, \varphi_{i +}^+ \ \hbox{and}\
\varphi^-_1, \cdots, \varphi^-_\ell.
\eqnlbl{4.55}
$$
It follows as in the proof of Lemma \thmref{4.5} that
$|C_{jk}|\leq C \lambda^\ell $ for $\lambda$ 
near zero, giving bound \eqnref{4.53} for $M_{jk}$.
If $\varphi_j^+$ is a slow mode, on the other hand, then the same
argument shows that there is a linear dependency in columns
$$
\varphi^+_1, \cdots, \varphi_{j-1}^+, \psi^-_k, \varphi^+_{j+1},
\varphi_{i_+}^+ 
\eqnlbl{4.54z}
$$
and an $(\ell+1)$-fold dependency in \eqnref{4.55},
since the omitted slow mode $\varphi_j^+$ plays no
role in either linear dependence; thus,
we obtain the bound \eqnref{slow-fast}, instead.
Analogous calculations yield the result for $d^\pm_{jk}$,
$\bar d_{jk}$ as well.
\myqed 

\medskip

{\bf 5.3. Proof of Proposition \thmref{lowbounds}.} 
The proof of Proposition \thmref{lowbounds} is now just a matter
of collecting the bounds of Lemmas \thmref{gaprel}--\thmref{4.8},
and substituting in the representations of Proposition \thmref{kochelformulae}.
More precisely, approximating fast-decaying dual modes 
$\phi_1^+,\dots,\phi_\ell^+$ and $\phi^-_{N-\ell+1},\dots,\phi_N^-$
by stationary modes $(\partial/\partial \delta_j)(\bar u^\delta,\bar v^\delta)$,
and slow dual modes by $e^{-\mu(\lambda)y}\bar {\tilde{V}}(\lambda)$
as described in Lemma \thmref{refinedderiv}, truncating $\mu(\lambda)$
at second order, and keeping only
the lowest order terms in the Laurent expansions of scattering coefficients
$M^\pm$, $d^\pm$ we obtain $E_\lambda$ and $S_\lambda$, respectively,
as order $\lambda^{-1}$ and order one terms, except for negligible
$\Cal{O}(e^{-\theta |x-y|})$ terms which we have accounted for in
error term $R^E_\lambda$.

Besides the latter, we have accounted in
term $R^E_\lambda$ for: (i) the truncation errors involved in the aforementioned
approximations (difference between stationary modes and actual fast modes,
and between constant-coefficient approximant and actual slow dual modes,
plus truncation errors in exponential rates); and, (ii)
neglected fast-decaying forward/slow-decaying dual mode combinations
(of order $\lambda$, by Lemma \thmref{4.5}).
A delicate point is the fact that term
$\Cal{O}(e^{-\theta|x-y|})\Cal{O}(e^{-\theta|y|})$
arising in the dual truncation of $E$ through 
the $\lambda\Cal{O}(e^{-\theta|y|})$ portion
of estimate \eqnref{refderiv} for 
$(\partial/\partial y)\tilde{\Theta}_j^\pm$, absorbs in the
(leading) $\Cal{O}(e^{-\theta|x-y|})$ term of \eqnref{RElambda}.
It was to obtain this key reduction that we carried out expansion
\eqnref{refderiv} to such high order.

In $R^S_\lambda$, 
we have accounted for truncation errors in the slow/slow pairings
approximated by term $S_\lambda$ (difference between
constant-coefficient approximants and actual slow forward modes,
corresponding to an additional $e^{-\theta|x|}$ term in the truncation
error factor,  and between constant-coefficient approximants and actual 
slow dual modes, plus truncation errors in exponential rates).
We have also accounted for the remaining, slow-decaying
forward/fast-decaying dual modes, which according to Lemma \thmref{4.8}
have scattering coefficient of order $\lambda$ and therefore may be
grouped with the difference between dual modes 
and their constant-coefficient approximants (of smaller order by a
factor of $e^{-\theta|y|}$) in the term $\Cal{O}(\lambda)$.
\myqed

\medskip
{\bf 5.4. Extended spectral theory.} 
Finally, we cite without proof the extended spectral theory of
[ZH] that we shall need in order to establish necessity of condition
($\Cal{D}$).  In the case that ($\Cal{D}$) holds, this reduces to
simple relations \eqnref{res1}. 

Let $C^\infty_\exp$ denote the space of $C^\infty$ functions decaying exponentially 
in all derivatives at some sufficiently high rate.  

\proclaim{Definition \thmlbl{proj}}  Let $L$ be a linear ordinary differential
operator with bounded, $C^\infty$ coefficients (so that
$L:C^\infty_\exp \to C^\infty_\exp$), and let 
$G_\lambda$ denote the Green's function of $L-\lambda I$. 
Further, let $\Omega$ be an open, simply connected domain intersecting
the resolvent set of $L$, on which $G_\lambda$ has a (necessarily unique) 
meromorphic extension.  Then, for $\lambda_0\in \Omega$, we define the
{\bf effective eigenprojection} 
${\cal P}_{\lambda_0} :C^\infty_\exp \to C^\infty$
by
$$
{\cal P}_{\lambda_0}f(x)= 
\int_{-\infty}^{+\infty} P_{\lambda_0}(x,y) f(y)\,dy,
$$
where
$$
P_{\lambda_0}(x,y)= \Residue_{\lambda_0}G_\lambda(x,y)
\eqnlbl{P}
$$
and $\Residue_{\lambda_0}$ denotes residue at $\lambda_0$.
We will refer to $P_{\lambda_0}(x,y)$ as the {\bf projection kernel}.
Likewise, we define the {\bf effective eigenspace} $\Sigma'_{\lambda_0}(L)$ by
$$
\Sigma'_{\lambda_0}(L)= \Range ( {\cal P}_{\lambda_0}),
$$
and the {\bf effective point spectrum} $\sigma_p'(L)$ of $L$ in $\Omega$
to be the set of $\lambda \in \Omega$ such that 
$\text{\rm dim } \Sigma'_{\lambda_0}(L)\ne 0$.
\endproclaim

The above definition is the natural one from the point of view of the
spectral resolution of the identity,
$I= \int_\Gamma (L-\lambda I)^{-1} d\lambda$, 
hence also from the point of view of asymptotic behavior of solutions
of PDE (see Section 6, below).
Away from the essential spectrum of $L$, it corresponds with the usual definition 
[Kat,Y].

\proclaim{Definition \thmlbl{ascent}}  Let $L$, $\Omega$, $\lambda_0$ be as above,
and $K$  be the order of the pole of $(L-\lambda I)^{-1}$ at $\lambda_0$.
{} For $\lambda_0\in \Omega$, and $k$ any integer,  we define
${\cal Q}_{\lambda_0,k}: C^\infty_\exp \to C^\infty$ by
$$
{\cal Q}_{\lambda_0,k}f(x) = 
\int_{-\infty}^{+\infty} Q_{\lambda_0,k}(x,y)f(y) \, dy, 
$$
where
$$ Q_{\lambda_0,k}(x,y)=
\Residue_{\lambda_0}(\lambda-\lambda_0)^{k} G_\lambda(x,y).
$$
{} For $0\le k \le K$, we define the {\bf effective eigenspace of ascent $k$} by
$$
\Sigma'_{\lambda_0,k}(L)= \Range ( {\cal Q}_{\lambda_0,K-k})
$$
\endproclaim

With the above definitions, we obtain the following, modified Fredholm Theory.

\proclaim{Proposition \thmlbl{3.11} [ZH]} Let $L$, $\lambda_0$,
$\Omega$ be as in Definition \thmref{proj},
and $K$  be the order of the pole of $G_\lambda$
at $\lambda_0$.
Then,
\medskip
\item{(i)} 
The operators ${\cal P}_{\lambda_0}$, ${\cal Q}_{\lambda_0,k}:
C^\infty_\exp \to C^\infty$ 
are $L$-invariant, with
$$
{\cal Q}_{\lambda_0,k+1}= (L-\lambda_0 I) {\cal Q}_{\lambda_0,k}
=  {\cal Q}_{\lambda_0,k}(L-\lambda_0 I)
\eqnlbl{504.9}
$$
for all $k\ne -1$, and
$$
{\cal Q}_{\lambda_0,k}= 
(L-\lambda_0 I)^{k} {\cal P}_{\lambda_0}
\eqnlbl{504.91}
$$
for $k\ge 0$.
\medskip
\item{(ii)} 
The effective eigenspace of ascent $k$ satisfies
$$\Sigma'_{\lambda_0, k} (L) = (L -\lambda_0 I) \Sigma'_{\lambda_0, k+1} (L). 
\eqnlbl{506}
$$
for all $0\le k\le K$, with
$$\{0\} = \Sigma'_{\lambda_0, 0} (L) 
\subset \Sigma'_{\lambda_0, 1} (L) \subset \cdots \subset 
\Sigma'_{\lambda_0, K} (L) = \Sigma'_{\lambda_0} (L). \eqnlbl{505}$$
Moreover, each containment in \eqnref{505} is strict.
\medskip
\item{(iii)} On ${\cal P}_{\lambda_0}^{-1} (C^\infty_\exp)$,  
${\cal P}_{\lambda_0}$, ${\cal Q}_{\lambda_0,k}$ all commute
($k\ge 0$),
and ${\cal P}_{\lambda_0}$ is a projection.
More generally, ${\cal P}_{\lambda_0}f=f$ for any $f\in \Sigma_{\lambda_0}(L: C^\infty_{exp})$, 
hence $\Sigma_{\lambda_0,k}(L: C^\infty_{exp})\subset \Sigma'_{\lambda_0,k}(L)$ 
for all $0\le k\le K$.

\medskip
\item{(iv)} The {\it multiplicity\/} of the eigenvalue $\lambda_0$, 
defined as $\dim \Sigma'_{\lambda_0} (L)$, 
is finite and bounded by $Kn$.  Moreover,
for all $0\le k \le K$,
$$\dim \Sigma'_{\lambda_0, k} (L) = \dim \Sigma'_{\lambda^*_0, k} (L^*).
 \eqnlbl{507.5}$$
Further, the projection kernel can be expanded as
$$
P_{\lambda_0} = \sum_j \varphi_j(x)\pi_j(y),
\eqnlbl{expand}
$$
where $\{ \varphi_j \}$, 
$\{\pi_j \}$ are bases for 
$\Sigma'_{\lambda_0}(L)$, $\Sigma'_{\lambda_0^*}(L^*)$,
respectively.

\medskip
\item{(v)} (Restricted Fredholm alternative) {For $g\in C^\infty_\exp$}, 
$$
(L-\lambda_0 I)f=g
$$ 
is soluble in $C^\infty$ {\bf if} and soluble in $C^\infty_\exp$ {\bf only if}
$Q_{\lambda_0,K-1}\, g=0$, 
or equivalently 
$$
g\in \Sigma'_{\lambda^*_0, 1} (L^*)^\perp.
\eqnlbl{507.6}$$
\endproclaim
\bigskip

{\bf Proof:} See [ZH].  \myqed

{\bf Remarks:}
{} For $\lambda_0$ in the resolvent set, 
${\cal P}_{\lambda_0}$ agrees with the standard definition,
hence $\Sigma'_{\lambda_0, k} (L)$
agrees with the usual $L^p$ eigenspace of generalized eigenfunctions of 
ascent $\le k$, for all $p < \infty$, since
$C^\infty_{exp}$ is dense in $L^p$, $p>1$,  
and  $\Sigma'_{\lambda_0, k} (L)$ is closed.
In the context of stability of traveling waves (more generally,
whenever the coefficients of $L$ exponentially approach 
constant values at $\pm\infty$),
$\Sigma'_{\lambda_0, k} (L)$ lies 
{\it between} the $L^p$ subspace 
$\Sigma_{\lambda_0, k}(L)$ and the corresponding
$L^p_{\Loc}$ subspace. 

The modification from standard Fredholm Theory
is primarily that here conclusions (iii) and (v) apply on
restricted domains.  This restriction is clearly necessary,
since the various expressions occurring in their statements
would otherwise not be defined.
In this regard, note that ${\cal P}_{\lambda_0}$ is, strictly speaking,
not a projection in the case that $\lambda\in \sigma_{ess}(L)$, 
since its domain does not then match its range.
However,  it maintains the projective structure 
\eqnref{expand}, somewhat justifying our abuse of notation.
\medskip

Now, suppose further (as holds in the case under consideration) that,
on $\Omega$:
(i) $L$ is a nondegenerate operator of some order $s$,
i.e. the coefficient matrix of the principal
part of the symbol is nonsingular for all $x$\footnote{
This assumption was made in [ZH] for convenience in calculation;
the result holds in considerably more generality (for example
for systems with real viscosity).
};
and, (ii) there exists an analytic choice of bases for stable/unstable
manifolds at $+\infty$/$-\infty$ of the associated eigenvalue equation.

In this case, we can define an analytic Evans
function $D(\lambda)$ as in \eqnref{localevans},
which we assume does not vanish identically.
It is easily verified that the order to which $D$
vanishes at any $\lambda_0$ is independent of the
choice of analytic bases $\Phi^\pm$.
The following result, established in [ZH], generalizes to
the extended spectral framework the
standard result of Gardner and Jones [GJ.1--2] in the classical
setting.

\proclaim{Proposition \thmlbl{3.12} [ZH]} Let $L$, $\lambda_0$ be
as above.  Then, 
\medskip

(i) $\dim \Sigma_{\lambda_0}'(L)$ is equal to
the order $d$ to which the Evans function $D_L$
vanishes at $\lambda_0$.

(ii)
$P_{\lambda_0}(L) = \sum_j \varphi_j(x)\pi_j(y)$,
where $\varphi_j$ and $\pi_j$ are in
$\Sigma'_{\lambda_0}(L)$ and
$\Sigma'_{\lambda_0}(L^*)$, respectively,
with ascents summing to $\le K+1$, where 
$K$ is the order of the pole of $G_\lambda$
at $\lambda_0$.
\endproclaim


{\bf Remark \thmlbl{regularity2}.}
Evidently, the theory of this section goes through also, with
suitable modifications regarding regularity, 
in the case (as for relaxation shocks) 
that the coefficients of $L$ have finite regularity $C^s_{\exp}$, 
i.e. decay exponentially in up to $s$ derivatives, in which case
$L:C^s_{\exp}\to C^0_{\exp}$.

\bigskip

\newsection {\bf Section 6. Pointwise Green's function bounds.}
\sectionnumber=6 
\theoremnumber=0
\equationnumber=0
\smallskip
\TagsOnLeft

In this section, we prove the main results of the paper,
establishing pointwise bounds on the Green's function $G$
and consequently sufficiency and necessity of condition
($\Cal{D}$) for linearized orbital stability.
\medskip

{\bf 6.1. Proof of Theorem \thmref{greenbounds}.}
We will establish Theorem \thmref{greenbounds} as a corollary
of the following, more detailed result:

\proclaim{Proposition \thmlbl{aux2}}
Assuming (H0)--(H4) and the spectral stability criterion ($\CalD$)
(equivalently, (D1)--(D2)),
the Green's function $G(x,t;y)$ associated with \eqnref{linearized}
may be decomposed as in Theorem \thmref{greenbounds}, where, for
$y\le 0$:
$$
\aligned
R(x,t;y)&= 
\Cal{O}(e^{-\eta(|x-y|+t)})\\
&+\sum_{k=1}^n 
\Cal{O} \left( (t+1)^{-1/2} e^{-\eta x^+} 
+e^{-\eta|x|} \right) 
(t+1)^{-1/2}e^{-(x-y-a_k^{*-} t)^2/Mt} \\
&+
\sum_{a_k^{*-} > 0, \, a_j^{*-} < 0} 
\chi_{\{ |a_k^{*-} t|\ge |y| \}}
\Cal{O} ((t+1)^{-1}
e^{-(x-a_j^{*-}(t-|y/a_k^{*-}|))^2/Mt}
e^{-\eta x^+}, \\
&+
\sum_{a_k^{*-} > 0, \, a_j^{*+}> 0} 
\chi_{\{ |a_k^{*-} t|\ge |y| \}}
\Cal{O} ((t+1)^{-1}
e^{-(x-a_j^{*+} (t-|y/a_k^{*-}|))^2/Mt}
e^{-\eta x^-}, \\
\endaligned
\eqnlbl{Rbounds}
$$
$$
\aligned
R_x(x,t;y)&= 
\sum_{j=1}^J \Cal{O}(e^{-\eta t})\delta_{x-\bar a_j t}(-y) 
+
\sum_{j=1}^J \Cal{O}(e^{-\eta(|x-y|+t)})\\
&+\sum_{k=1}^n 
\bold {O} \left( (t+1)^{-1} e^{-\eta x^+} 
+e^{-\eta|x|} \right) 
(t+1)^{-1/2}e^{-(x-y-a_k^{*-} t)^2/Mt} \\
&+
\sum_{a_k^{*-} > 0, \, a_j^{*-} < 0} 
\chi_{\{ |a_k^{*-} t|\ge |y| \}}
\bold {O} ((t+1)^{-3/2} 
e^{-(x-a_j^{*-}(t-|y/a_k^-|))^2/Mt}
e^{-\eta x^+} \\
&+
\sum_{a_k^{*-} > 0, \, a_j^{*+} > 0} 
\chi_{\{ |a_k^{*-} t|\ge |y| \}}
\bold {O} ((t+1)^{-3/2} 
e^{-(x-a_j^{*+}(t-|y/a_k^{*-}|))^2/Mt}
e^{-\eta x^-}, \\
\endaligned
\eqnlbl{Rxbounds}
$$
and
$$
\aligned
R_y(x,t;y)&= 
\sum_{j=1}^J \Cal{O}(e^{-\eta t})\delta_{x-\bar a_j t}(-y) 
+
\sum_{j=1}^J \Cal{O}(e^{-\eta(|x-y|+t)})\\
&+\sum_{k=1}^n 
\Cal{O} \left( (t+1)^{-1/2} e^{-\eta x^+} 
+e^{-\eta|x|} \right) 
(t+1)^{-1}e^{-(x-y-a_k^{*-} t)^2/Mt} \\
&+
\sum_{a_k^{*-} > 0, \, a_j^{*-} < 0} 
\chi_{\{ |a_k^{*-} t|\ge |y| \}}
\Cal{O} ((t+1)^{-e/2} 
e^{-(x-a_j^{*-}(t-|y/a_k^{*-}|))^2/Mt}
e^{-\eta x^+} \\
&+
\sum_{a_k^{*-} > 0, \, a_j^{*+} > 0} 
\chi_{\{ |a_k^{*-} t|\ge |y| \}}
\Cal{O} ((t+1)^{-3/2} 
e^{-(x-a_j^{*+}(t-|y/a_k^{*-}|))^2/Mt}
e^{-\eta x^-}, \\
\endaligned
\eqnlbl{Rybounds}
$$
for some $\eta$, $M>0$, where $\bar a_j$ are as in Theorem
\thmref{greenbounds},  $x^\pm$ denotes the positive/negative
part of $x$,  and
indicator function $\chi_{\{ |a_k^{*-}t|\ge |y| \}}$ is 
one for $|a_k^{*-}t|\ge |y|$ and zero otherwise.
Symmetric bounds hold for $y\ge0$.
\endproclaim

{\bf Proof.}
Our starting point is the representation
$$
G(x,t;y) = {1\over 2\pi i} \text{\rm P.V.}\int_{\eta-i\infty}^{\eta + i\infty}
e^{\lambda t}G_\lambda (x,y)d \lambda,
\eqnlbl{inverseLT2}
$$
valid by Proposition \thmref{inverseLT} for any $\eta$ sufficiently large.

\medskip
{\bf Case I. $|x-y|/t$ large.}
We first treat the simpler case that $|x-y|/t\ge S$,
$S$ sufficiently large, for which $G$ is known a priori
to vanish, by the general property of finite propagation speed
for systems of hyperbolic type.
More precisely, $G\equiv 0$ for $x\not \in [z_1(y,t),z_J(y,t)]$,
where $z_j(\cdot, \cdot)$ are the standard hyperbolic characteristic paths
defined in \eqnref{char}, $a_1\le \cdots\le a_J$ denoting the eigenvalues
of $A:=df(\bar u,\bar v)$.
This may be seen, for example, by standard energy estimates as in, e.g.,
[Ev];
however, it is instructive to recover this result by direct calculation.

Suppose, then, that $x\not \in [z_1,z_J]$.  
Fixing $x$, $y$, $t$, set $\lambda= \eta + i\xi$.
Applying Proposition \thmref{highbounds}, we obtain, for
$\eta>0$ sufficiently large, the decomposition
$$
G(x,t;y) = 
{1\over 2\pi i} \text{\rm P.V.}\int_{\eta-i\infty}^{\eta + i\infty}
e^{\lambda t}H_\lambda (x,y)d \lambda
+
{1\over 2\pi i} \text{\rm P.V.}\int_{\eta-i\infty}^{\eta + i\infty}
e^{\lambda t}\Theta_\lambda (x,y)d \lambda,
\eqnlbl{split1}
$$ 
where
$$
H_{\lambda}(x,y)=
\cases
   -\sum_{j=K+1}^{J} 
a_j(y)^{-1}
e^{-\int_y^x \lambda/a_{j} (z)\, dz}r_j(x) \tilde \zeta_j(x,y) l_j^{t}(y)
\qquad &x>y,\\
   \sum_{j=1}^K
a_j(y)^{-1}
e^{-\int_y^x \lambda/a_j (z)\, dz}r_j(x) \tilde \zeta_j(x,y) l_j^{t}(y)
    \qquad &x<y,\\
                 \endcases 
\eqnlbl{nshighfreq2}
 $$
and
$$
\Theta_\lambda(x,y)= 
\lambda^{-1}B(x,y;\lambda)+ \lambda^{-1}(x-y)C(x,y;\lambda)+ 
\lambda^{-2}D(x,y;\lambda)
\eqnlbl{errorterms2}
$$
with $\zeta_j$, $B$, $C$, $D$ as defined in Proposition \thmref{highbounds}.

For definiteness, take $x>y$. Then, the first term of \eqnref{split1}
contributes to integral \eqnref{inverseLT2} the explicitly 
evaluable quantity
$$
\aligned
{1\over 2\pi i} \text{\rm P.V.}&\int_{\eta-i\infty}^{\eta + i\infty}
e^{\lambda t}H_\lambda (x,y)d \lambda =
   \sum_{j=K+1}^J 
\left( a_j(y)^{-1} \text{\rm P.V.}\int_{-\infty}^{+\infty} 
e^{ i\xi \big( t- \int_y^x(1/a_j(z) )\, dz \big)}\, d\xi \right)\\
&\qquad \times
e^{ \eta \big( t- \int_y^x(1/a_j(z) )\, dz \big)}
r_j(x) \tilde \zeta_j(x,y) l_j^{t}(y)\\
&= \sum_{j=K+1}^J 
\left( a_j(y)^{-1}\delta \big( t- \int_y^x(1/a_j(z) )\, dz \big)  \right)
r_j(x) \tilde \zeta_j(x,y) l_j^{t}(y)\\
&=
   \sum_{j=K+1}^J \delta_{x-\bar a_j t}(-y) 
r_j(x) \zeta_j(y,t) l_j^{t}(y)\\
&=H(x,t;y),\\
\endaligned
\eqnlbl{deltacalc}
$$
where $\bar a_j$, $\zeta_j$, and $H$ 
are as defined in Theorems \thmref{greenbounds} and \thmref{auxgreenbounds}.
Note: in the final equality of \eqnref{deltacalc} we have used
the general fact that
$$
f_y(x,y,t) \delta\big(f(x,y,t)\big)=\delta_{h(x,t)}(y),
\eqnlbl{generalfact}
$$
provided $f_y\ne 0$, where $h(x,t)$ is defined by $f(x,h(x,t),t)\equiv 0$.
This term clearly vanishes for $x$ outside $[z_1(y,t),z_J(y,t)]$, as claimed.

Similar calculations show that the second, error term
in \eqnref{split1} also vanishes.
For example, the term $e^{\lambda t}\lambda^{-1}B(x,y;\lambda)$
contributes
$$
{1\over 2\pi i} 
   \sum_{j=K+1}^K 
\left(\text{\rm P.V.}\int_{-\infty}^{+\infty} 
(\eta+ i\xi)^{-1}
e^{ i\xi \big( t- \int_y^x(1/a_j(z) )\, dz \big)}\, d\xi \right)
e^{ \eta \big( t- \int_y^x(1/a_j(z) )\, dz \big)}
b_j^+(x,y).
\eqnlbl{bcalc}
$$
The factor
$$\text{\rm P.V.}
\int_{-\infty}^{+\infty} 
(\eta+ i\xi)^{-1}
e^{ i\xi \big( t- \int_y^x(1/a_j(z) )\, dz \big)}\, d\xi ,
$$
though not absolutely convergent,
is integrable and uniformly bounded as a principal value integral,
for all $\eta \in \BbbR$, by explicit computation.  
On the other hand,
$$
e^{ \eta \big( t- \int_y^x(1/a_j(z) )\, dz \big)} \le
Ce^{-\eta \, d(x,[z_1(y,t),z_J(y,t)])/ \min_{j,x}|a_j(x)|}
\to 0
\eqnlbl{alphabound2}
$$
as $\eta \to +\infty$, for each $K+1\le j \le J$, since $a_j>0$
on this range.
Thus, taking $\eta \to \infty$,
we find that the product \eqnref{bcalc} goes to zero, giving the result.

Similarly, the contributions of terms $e^{\lambda t}C(x,y;\lambda)$
and $e^{\lambda t}D(x,y;\lambda)$ split into the product of a convergent,
uniformly bounded integral in $\xi$, a 
(constant) factor depending only on $(x,y)$, and a factor $\alpha(x,y,t,\eta)$
going to zero as $\eta\to 0$ at rate \eqnref{alphabound2}.  
Thus, taking $\eta \to +\infty$, we find
that each of these terms vanishes, and thus $G$ vanishes also, as claimed.

\medskip
{\bf Case II. $|x-y|/t$ bounded.}
We now turn to the critical case that $|x-y|/t\le S$ for some
fixed $S$.  In this regime, note that any contribution of order
$e^{\theta t}$, $\theta>0$, may be absorbed in the residual (error)
term $R$ (resp. $R_x$, $R_y$); we shall use this observation repeatedly.

\medskip
{\bf Decomposition of the contour.}
We begin by converting contour integral \eqnref{inverseLT2} into
a more convenient form decomposing high, intermediate, and low
frequency contributions.

\proclaim{Observation \thmlbl{decomp2}}
Assuming (H0)--(H4), 
there holds the representation
$$
\aligned
G(x,t;y)&= I_a + I_b + II_a + II_b\\
&:=
{1\over 2\pi i}\text{\rm P.V.} \int_{\eta-i\infty}^{\eta + i\infty}
e^{\lambda t}H_\lambda (x,y)d \lambda\\
&+{1\over 2\pi i} \text{\rm P.V.} \left(\int_{-\eta_1-i\infty}^{-\eta_1 - iR}
+ \int_{-\eta_1+iR}^{-\eta_1 + i\infty}\right)
e^{\lambda t}(G_\lambda-H_\lambda) (x,y)d \lambda\\
&+
{1\over 2\pi i} \oint_{\Gamma}
e^{\lambda t}G_\lambda (x,y)d \lambda
-{1\over 2\pi i} \int_{-\eta_1-iR}^{-\eta_1 + iR}
e^{\lambda t}H_\lambda (x,y)d \lambda,
\endaligned
\eqnlbl{contours}
$$
$$
\Gamma:=
[-\eta_1-iR, \eta -iR] \cup
[\eta-iR, \eta +iR] \cup
[\eta+iR, -\eta_1 +iR],
\eqnlbl{Gammacontour}
$$
for any $\eta>0$ such that \eqnref{inverseLT2} holds, $R$ sufficiently
large, and $-\eta_1<0$ as in \eqnref{Lambda}.
(Note that we have not here assumed ($\Cal{D}$)).
\endproclaim

{\bf Proof.}
We first observe that, by Proposition \thmref{highbounds}, $L$
has no spectrum on the portion of 
$\Gamma:=\{\lambda: \R \lambda \ge -\eta_1<0\}$ 
lying outside of the rectangle 
$$
\Cal{R}:=
\{\lambda: -\eta_1\le \R \lambda \le \eta, \, -R\le \Im \lambda \le R \}
\eqnlbl{rectangle}
$$
for $\eta>0$, $R>0$ sufficiently large, hence $G_\lambda$ is analytic
on this region.  Since, also, $H_\lambda$ is analytic on the whole
complex plane, contours involving either one of these contributions
may be arbitrarily deformed within $\Gamma \setminus \Cal{R}$
without affecting the result, by Cauchy's Theorem.

Recalling, further, that 
$$
|G_\lambda-H_\lambda| =\Cal{O}(\lambda^{-1})
$$
by Proposition \thmref{highbounds}, we find that
$$
{1\over 2\pi i} \text{\rm P.V.}\int_{\eta-i\infty}^{\eta + i\infty}
e^{\lambda t}(G_\lambda-H_\lambda) (x,y)d \lambda
$$
may be deformed to
$$
{1\over 2\pi i} \text{\rm P.V.}\oint_{\partial(\Lambda \setminus \Cal{R})}
e^{\lambda t}(G_\lambda-H_\lambda) (x,y)d \lambda
$$
Noting, finally, that 
$$
 +{1\over 2\pi i} \oint_{\partial \Cal{R}}
e^{\lambda t}H_\lambda (x,y)d \lambda=0,
$$
by Cauchy's Theorem, we obtain, finally, 
$$
\aligned
{1\over 2\pi i} &\text{\rm P.V.}\int_{\eta-i\infty}^{\eta + i\infty}
e^{\lambda t}G_\lambda (x,y)d \lambda =\\
&{1\over 2\pi i} \int_{\eta-i\infty}^{\eta + i\infty}
e^{\lambda t}H_\lambda (x,y)d \lambda
+
{1\over 2\pi i} \text{\rm P.V.}\int_{\eta-i\infty}^{\eta + i\infty}
e^{\lambda t}(G_\lambda-H_\lambda) (x,y)d \lambda\\
&={1\over 2\pi i} \text{\rm P.V.}\int_{\eta-i\infty}^{\eta + i\infty}
e^{\lambda t}H_\lambda (x,y)d \lambda
+
{1\over 2\pi i} \text{\rm P.V.}\oint_{\partial(\Lambda \setminus \Cal{R})}
e^{\lambda t}(G_\lambda-H_\lambda) (x,y)d \lambda\\
& \qquad +{1\over 2\pi i} \oint_{\partial \Cal{R}}
e^{\lambda t}H_\lambda (x,y)d \lambda,
\endaligned
\eqnlbl{splitcalc}
$$
from which the result follows by combining the second
and third contour integrals along their common edges $\Gamma$.
\myqed

{\bf Remark \thmlbl{move}.}
We have already shown in \eqnref{deltacalc} that 
${1\over 2\pi i} \text{\rm P.V.}\int_{\eta-i\infty}^{\eta + i\infty}
e^{\lambda t}H_\lambda (x,y)d \lambda =H(x,t;y)$ independent of $\eta$.
Thus, we may deform the the entire contour integral
$$
{1\over 2\pi i} \text{\rm P.V.}\int_{\eta-i\infty}^{\eta + i\infty}
e^{\lambda t}G_\lambda (x,y)d \lambda 
$$
to
$$
{1\over 2\pi i} \text{\rm P.V.}\int_{\partial(\Lambda \setminus \Cal{R})}
e^{\lambda t}G_\lambda (x,y)d \lambda 
$$
to obtain \eqnref{contours} by a different route.

\proclaim{Observation \thmlbl{decompf}}
Assuming ($\Cal{D}$) together with (H0)--(H4), 
we may replace \eqnref{contours} by
$$
G(x,t;y)= I_a + I_b + II_{\tilde a} + II_b+ III,
$$
where $I_a$, $I_b$ and $II_b$ are as in \eqnref{contours}, and
$$
II_{\tilde a}:=
{1\over 2\pi } \left(\int_{-\eta_1-iR}^{-\eta_1 - ir/2}
+ \int_{-\eta_1+ir/2}^{-\eta_1 + iR}\right)
e^{\lambda t}G\lambda(x,y)d xi ,
$$
$$
III:=
{1\over 2\pi i} \oint_{\tilde\Gamma}
e^{\lambda t}G_\lambda (x,y)d \lambda,
\eqnlbl{contoursf}
$$
and
$$
\tilde\Gamma:=
[-\eta_1-ir/2, \eta -ir/2] \cup
[\eta-ir/2, \eta +ir/2] \cup
[\eta+ir/2, -\eta_1 +ir/2],
\eqnlbl{tildeGammacontour}
$$
for any $\eta$, $r>0$, and $\eta_1$ sufficiently small with respect to $r$.
\endproclaim

{\bf Proof.}  By assumption ($\Cal{D}$), $L$ has no spectrum on the
region between
contour $\Gamma$ and the union of contour $\tilde \Gamma$ and 
the contour of term $II_{\tilde a}$, hence $G_\lambda$ is analytic
on that region, and
$$
II_a=II_{\tilde a}+III
$$
by Cauchy's Theorem, giving the result.
\myqed

Using the final decomposition \eqnref{contoursf},
we shall estimate in turn the high frequency contributions $I_a$ and $I_b$,
the intermediate frequency contributions $II_{\tilde a}$ and $II_b$,
and the low frequency contributions $III$.

\medskip
{\bf High frequency contribution.}
We first carry out the straightforward estimatation of
the high-frequency terms $I_a$ and $I_b$.
The principal term $I_a$ has already been computed in \eqnref{deltacalc}
to be $H(x,t;y)$.
Likewise, calculations similar to those of 
\eqnref{bcalc}--\eqnref{alphabound2}
show that the error term 
$$
\aligned
I_b&=
{1\over 2\pi i} \text{\rm P.V.}
\left(\int_{-\eta_1-i\infty}^{-\eta_1 -R}
\int_{-\eta_1+ iR}^{-\eta_1 +i\infty} \right)
e^{\lambda t}\Theta_\lambda (x,y)d \lambda\\
&=
{1\over 2\pi i} \text{\rm P.V.}
\left(\int_{-\eta_1-i\infty}^{-\eta_1 -R}
\int_{-\eta_1+ iR}^{-\eta_1 +i\infty} \right)
e^{\lambda t}
\big(\lambda^{-1}B(x,y;\lambda)\\
&\qquad + \lambda^{-1}(x-y)C(x,y;\lambda)+ 
\lambda^{-2}D(x,y;\lambda)\big)
d \lambda\\
\endaligned
$$ 
is time-exponentially small.

For example, for $x>y$ the term $e^{\lambda t}\lambda^{-1}B(x,y;\lambda)$
contributes
$$
\aligned
   \sum_{j=K+1}^K 
{1\over 2\pi i} \text{\rm P.V.}
\left( \int_{-\infty}^{-R} \int_{R}^{+\infty} \right)
&(-\eta_1+ i\xi)^{-1}
e^{ i\xi \big( t- \int_y^x(1/a_j(z) )\, dz \big)}\, d\xi \\
&\qquad \qquad \times
e^{ -\eta_1 \big( t- \int_y^x(1/a_j(z) )\, dz \big)}
b_j^+(x,y),
\endaligned
\eqnlbl{bcalc}
$$
where
$$
{1\over 2\pi i} \text{\rm P.V.}
\left( \int_{-\infty}^{-R} \int_{R}^{+\infty} \right)
(\eta+ i\xi)^{-1}
e^{ i\xi \big( t- \int_y^x(1/a_j(z) )\, dz \big)}\, d\xi <\infty,
\eqnlbl{conv}
$$
and
$$
e^{ \eta_1 \int_y^x(1/a_j(z) )\, dz }b(x,y) 
\le
C_1 e^{ \eta_1 \int_y^x(1/a_j(z) )\, dz -\theta |x-y|}
\le C_2,
\eqnlbl{alphabound3}
$$
for $\eta_1$ sufficiently small.
This may be absorbed in the first term of $R$, \eqnref{R}.

Likewise, the contributions of terms $e^{\lambda t}C(x,y;\lambda)$
and $e^{\lambda t}D(x,y;\lambda)$ split into the product of a convergent,
uniformly bounded integral in $\xi$, a 
bounded factor analogous to \eqnref{alphabound3}, and the
factor $e^{-\eta_1 t}$, giving the result.

{\it Derivative bounds.}  Derivatives
$(\partial/\partial x)I_b$ (resp. $(\partial/\partial y)I_b$)
may be treated in identical fashion using \eqnref{xder}
(resp. \eqnref{yder}, to show that they are absorbable
in the estimates given for $R_x$ (resp. $R_y$).
We point out that the integral arising from term 
$B^0_x$ (resp. $B^0_y$) of \eqnref{xder} (resp. \eqnref{yder}) 
corresponds to the first, delta-function term in \eqnref{Rxbounds}
(resp. \eqnref{Rybounds}), while the integral arising from term
$(x-y)C^0_x$ (resp. $(x-y)C^0_y$) vanishes, as the product
of $(x-y)$ and a delta function $\delta_z(y)$ with $z\ne x$.
\medskip

\medskip
{\bf Intermediate frequency contribution.}
Error term $II_b$ is time-exponentially small for $\eta_1$ sufficiently
small, by the same calculation as in \eqnref{bcalc}--\eqnref{alphabound3},
hence negligible.
Likewise, term $II_{\tilde a}$ by the basic estimate \eqnref{uKdecay}
is seen to be time-exponentially small of order $e^{-\eta_1 t}$ 
for any $\eta_1>0$ sufficiently
small that the associated contour lies in the resolvent set of $L$.

\medskip
{\bf Low frequency contribution.}
It remains to estimate the low-frequency term $III$,
which is of essentially the same form as
the low-frequency contribution analyzed in [ZH,Z.2,Z.4] in the strictly
parabolic case, in that the contour is the same and the resolvent kernel
$G_\lambda$ satisfies identical bounds in this regime.
Thus, we may conclude from these previous analyses that $III$ gives
contribution $E+S+R$, as claimed, exactly as in the strictly parabolic
case.  
For completeness, we indicate the main features of the argument here.

{\it Case $t\le 1$}. 
First observe that estimates in the short-time regime $t\le 1$ are trivial, 
since then $|e^{\lambda t}G_\lambda(x,y)|$ is uniformly bounded on
the compact set $\tilde\Gamma$, and we have $|G(x,t;y)|\le C\le e^{-\theta t}$
for $\theta>0$ sufficiently small.  But, likewise, $|E|$ and $S$
are uniformly bounded in this regime, hence time-exponentially decaying.
As observed previously, all such terms are negligible, being absorbable
in the error term $R$.  Thus, we may add $E+S$ and subtract $G$ to
obtain the result
\medskip

{\it Case $t\ge 1$}. Next, consider the critical (long-time) regime
$t\ge 1$.
For definiteness, take $y\le x\le 0$; the other two cases are similar.
Decomposing
$$
 G(x,t;y)=
{1\over 2\pi i} \oint_{\tilde\Gamma} e^{\lambda t}E_\lambda (x,y)d \lambda
+{1\over 2\pi i} \oint_{\tilde\Gamma} e^{\lambda t}S_\lambda (x,y)d \lambda
+{1\over 2\pi i} \oint_{\tilde\Gamma} e^{\lambda t}R_\lambda (x,y)d \lambda,
\eqnlbl{dec}
$$
with $E_\lambda$, $S_\lambda$, and $R_\lambda$ as defined in Proposition
\thmref{lowbounds},
we consider in turn each of the three terms on the righthand side.

\medskip

{\it $E_\lambda$ term.}
Let us first consider the dominant term
$$
 {1\over 2\pi i} \oint_{\tilde\Gamma} e^{\lambda t}E_\lambda (x,y)d \lambda,
\eqnlbl{dominant}
$$
which, by \eqnref{Elambda}, is given by
$$
\sum_{a_k^{*-} > 0}
[c^{j,0}_{k,-}]
\frac{\partial}{\partial \delta_j}
\pmatrix
\bar u^\delta(x) \\
\bar v^\delta(x)
\endpmatrix
L_k^{*-t}
\alpha_k(x,t;y),
\eqnlbl{separation}
$$
where
$$
\alpha_k(x,t;y):=
 {1\over 2\pi i} \oint_{\tilde\Gamma}
\lambda^{-1} 
e^{\lambda t}
e^{(\lambda/a^{*-}_k - \lambda^2 \beta^{*-}_k/{a^{*-}_k}^3 )y}
\, d\lambda.
\eqnlbl{dalpha}
$$

Using Cauchy's Theorem, we may move the contour $\tilde\Gamma$ to obtain
$$
\aligned
\alpha_k(x,t;y)&=
{1\over 2\pi }
\text{\rm P.V. }\int_{-r/2}^{+r/2}
(i\xi)^{-1} 
e^{i\xi t}
e^{(i\xi/a^{*-}_k + \xi^2 \beta^{*-}_k/{a^{*-}_k}^3 )y}
\, d\xi\\
&\quad +{1\over 2\pi i}
\left( \int_{-\eta_1-ir/2}^{-ir/2}
+ \int_{ir/2}^{-\eta_1+ir/2} \right)
\lambda^{-1} 
e^{\lambda t}
e^{(\lambda/a^{*-}_k - \lambda^2 \beta^{*-}_k/{a^{*-}_k}^3 )y}
\, d\lambda \\
&\quad +
\frac{1}{2}\text{\rm Residue }_{\lambda=0}\, e^{\lambda t}
e^{(\lambda/a^{*-}_k - \lambda^2 \beta^{*-}_k/{a^{*-}_k}^3 )y},\\
\endaligned
\eqnlbl{moved}
$$
or, rearranging and evaluating the final, residue term:
$$
\aligned
\alpha_k(x,t;y)&=
\left({1\over 2\pi }
\text{\rm P.V. }\int_{-\infty}^{+\infty}
(i\xi)^{-1} 
e^{i\xi(t+ y/a^{*-}_k)}
e^{\xi^2 (\beta^{*-}_k/{a^{*-}_k}^3 )y}
\, d\xi
+\frac{1}{2} \right)\\
&\quad -{1\over 2\pi }
\left(\int_{-\infty}^{-r/2} +
\int_{r/2}^{+\infty} \right)
(i\xi)^{-1} 
e^{i\xi(t+ y/a^{*-}_k)}
e^{\xi^2 (\beta^{*-}_k/{a^{*-}_k}^3 )y}
\, d\xi\\
&\quad +{1\over 2\pi i}
\left( \int_{-\eta_1-ir/2}^{-ir/2}
+ \int_{ir/2}^{-\eta_1+ir/2} \right)
\lambda^{-1} 
e^{\lambda t}
e^{(\lambda/a^{*-}_k - \lambda^2 \beta^{*-}_k/{a^{*-}_k}^3 )y}
\, d\lambda. \\
\endaligned
\eqnlbl{moved}
$$

The first term in \eqnref{moved} may be explicitly evaluated to give 
$$
\errfn\left(\frac{y+a_k^{*-}t}{\sqrt{4\beta_k^{*-}|y/a_k^{*-}|}}\right),
\eqnlbl{yversion}
$$
where
$$
\errfn(z):=
{1\over 2\pi } \int_{-\infty}^{z} e^{-y^2}dy,
\eqnlbl{errdef}
$$
whereas the second and third terms are clearly time-exponentially
small for $t\le C|y|$ and $\eta_1$ sufficiently small relative to $r$
(see detailed discussion of similar calculations below, under
{\it $R_\lambda$ term}).  In the trivial case 
$t\ge C|y|$, $C>0$ sufficiently large, we can simply move the
contour to $[-\eta_1-ir/2,-\eta_1+ir/2]$ to obtain  (complete)
residue $1$ plus a time-exponentially small error corresponding to 
the shifted contour integral, which result again agrees with 
\eqnref{yversion} up to a time-exponentially small error.

Expression \eqnref{yversion} may be rewritten as
$$
\errfn\left(\frac{y+a_k^{*-}t}{\sqrt{4\beta_k^{*-}|y/a_k^{*-}|}}\right)
\eqnlbl{tversion}
$$
plus error
$$
\aligned
\errfn&\left(\frac{y+a_k^{*-}t}{\sqrt{4\beta_k^{*-}|y/a_k^{*-}|}}\right)-
\errfn\left(\frac{y+a_k^{*-}t}{\sqrt{4\beta_k^{*-}t}}\right)\\
&\sim
\errfn'\left(\frac{y+a_k^{*-}t}{\sqrt{4\beta_k^{*-}t}}\right)
\left(
-\frac{4\beta_k^{*-}}{2}(y+a_k^{*-}t)^2(4\beta_k^{*-}t)^{-3/2}
\right)\\
&=
\Cal{O}(t^{-1}e^{ (y+a_k^{*-}t)^2/Mt}),
\endaligned
\eqnlbl{terror}
$$
for $M>0$ sufficiently large, and similarly for $x$- and $y$-derivatives.
Multiplying by 
$$
[c^{j,0}_{k,-}]
\frac{\partial}{\partial \delta_j}
\pmatrix
\bar u^\delta(x) \\
\bar v^\delta(x)
\endpmatrix
L_k^{*-t}
=\Cal{O}(e^{-\theta|x|}),
$$
we find that term \eqnref{tversion} gives contribution
$$
[c^{j,0}_{k,-}]
\frac{\partial}{\partial \delta_j}
\pmatrix
\bar u^\delta(x) \\
\bar v^\delta(x)
\endpmatrix
L_k^{*-t}
\errfn\left(\frac{y+a_k^{*-}t}{\sqrt{4\beta_k^{*-}|y/a_k^{*-}|}}\right),
\eqnlbl{tfinal}
$$
whereas term \eqnref{terror} gives a contribution absorbable in $R$ (resp.
$R_x$, $R_y$); see Remark 6.5 below for detailed discussion
of a similar calculation.

Finally, observing that
$$
[c^{j,0}_{k,-}]
\frac{\partial}{\partial \delta_j}
\pmatrix
\bar u^\delta(x) \\
\bar v^\delta(x)
\endpmatrix
L_k^{*-t}
\errfn\left(\frac{y-a_k^{*-}t}{\sqrt{4\beta_k^{*-}t}}\right)
\eqnlbl{smallterm}
$$
is time-exponentially small for $t\ge 1$, since $a_k^{*-}>0$,
$y<0$, and 
$(\partial/\partial \delta_j)(\bar u^\delta,\bar v^\delta)\le Ce^{-\theta|x|}$,
$\theta>0$, we may subtract and add this term to \eqnref{tfinal} to
obtain a total of $E(x,t;y)$ plus terms absorbable in $R$ (resp. $R_x$, $R_y$).

{\bf Remark \thmlbl{discrepancy}.}
There is a slight discrepancy between the form of $E$ used here, and the form
$$
\aligned
E(x,t;y)&:= 
\sum_{a_k^- > 0}
[c^{j,0}_{k,-}]
\frac{\partial}{\partial \delta_j}
\pmatrix
\bar u^\delta(x) 
\\
\bar v^\delta(x)
\endpmatrix
L_k^{*-t}
\left(\errfn \left(\frac{x-y+a_k^{*-}t}{\sqrt{4\beta_k^{*-}t}}\right)\right.\\
&\qquad -\errfn\left.\left(\frac{x-y-a_k^{*-}t}{\sqrt{4\beta_k^{*-}t}}\right)\right),\\
&=\sum_{a_k^- > 0}
[c^{j,0}_{k,-}]
\frac{\partial}{\partial \delta_j}
\pmatrix
\bar u^\delta(x) 
\\
\bar v^\delta(x)
\endpmatrix
L_k^{*-t}
\left(\errfn \left(\frac{-x+y+a_k^{*-}t}{\sqrt{4\beta_k^{*-}t}}\right)\right.\\
&\qquad -\errfn\left.\left(\frac{-x+y-a_k^{*-}t}{\sqrt{4\beta_k^{*-}t}}\right)\right),\\
\endaligned
\eqnlbl{Eprevious}
$$
given in [Z.2];
however, this is only cosmetic, since the two expressions are equal
modulo terms absorbable in the error terms $R$ (resp. $R_x$, $R_y$).
To see this, first note that the difference is time-exponentially small
outside the critical region 
$$
x\in [-C\log (t+1), +C\log (t+1)],
$$
for $C>0$ sufficiently large.
Then, estimating 
$$
\aligned
\left| \errfn \left(\frac{-x+y\pm a_k^{*-}t}{\sqrt{4\beta_k^{*-}t}}\right)
-
\errfn \left(\frac{y\pm a_k^{*-}t}{\sqrt{4\beta_k^{*-}t}}\right) \right|
&\le
|x| \,
\left| \errfn' \left(\frac{-x+y\pm a_k^{*-}t}{\sqrt{4\beta_k^{*-}t}}\right) \right|\\
&\le
|x| \,
t^{-1/2}e^{{(x-y\mp a_k^{*-}t)^2}/{4\beta_k^{*-}t}},\\
\endaligned
$$
and recalling that 
$(\partial/\partial \delta_j)(\bar u^\delta , \bar v^\delta )=
\Cal{O}(e^{-\theta|x|})$,
we find that the error is bounded by
$$
C|x| e^{-\theta|x|} t^{-1/2}e^{{(x-y-a_k^{*-}t)^2}/{4\beta_k^{*-}t}}
\le
C_2 e^{-\theta|x|/2} t^{-1/2}e^{{(x-y-a_k^{*-}t)^2}/{4\beta_k^{*-}t}},
\eqnlbl{errorbound}
$$
and thus absorbable in $R$.  The calculation for the $y$-derivative is
similar.
We have chosen to present the result in the form used here, since
it is more straightforward to prove, and is more convenient
for the final stability analysis (in fact, this was the form
that was actually used in the stability analysis of [Z.2], as well).
It was presented as \eqnref{Eprevious} in [Z.2] in order to capture
the parabolic property of almost-finite propagation speed; here,
we have dealt with propagation speed separately, 
which seems the preferable approach.

\medskip
{\bf Remark \thmlbl{balanced}.}
As described in [Z.2], the time-exponentially small terms 
\eqnref{smallterm} that we have (artificially)
included in $E$ are useful for linear and nonlinear
stability analyses, serving to regularize behavior at $t=0$.
Note, also, that the resulting ``balanced'' form of $E$ is the one 
that appears naturally in the explicit Green's function formula 
for viscous shock solutions of Burgers equation;
see equation (2.13), [Z.2], and surrounding example.
\medskip
{\it $S_\lambda$ term.}
Next, consider the second-order term 
$$
 {1\over 2\pi i} \oint_{\tilde\Gamma} e^{\lambda t}S_\lambda (x,y)d \lambda,
\eqnlbl{second}
$$
which, by \eqnref{Slambda2}, is given by
$$
\sum_{a_k^{*-}>0}R_k^{*-}  {L_k^{*-}}^t
\alpha_k(x,t;y)
+ 
\sum_{a_k^{*-} > 0, \,  a_j^{*-} < 0} 
[c^{j,-}_{k,-}]R_j^{*-}  {L_k^{*-}}^t
\alpha_{jk}(x,t;y)
\eqnlbl{separationS}
$$
where
$$
\alpha_k(x,t;y):=
 {1\over 2\pi i} \oint_{\tilde\Gamma}
e^{\lambda t}
e^{ (-\lambda/a^{*-}_k + \lambda^2 \beta^{*-}_k/{a^{*-}_k}^3 )(x-y)}
\, d\lambda.
\eqnlbl{dalphak}
$$
and
$$
\alpha_{jk}(x,t;y):=
 {1\over 2\pi i} \oint_{\tilde\Gamma}
e^{\lambda t}
e^{(-\lambda/a^{*-}_j + \lambda^2 \beta^{*-}_j/{a^{*-}_j}^3 )x
+(\lambda/a^{*-}_k - \lambda^2 \beta^{*-}_k/{a^{*-}_k}^3 )y}
\, d\lambda.
\eqnlbl{dalphak}
$$
Similarly as in the treatment of the $E_\lambda$ term, just above,
by deforming the contour $\tilde\Gamma$ to
$$
\Gamma':=
[-\eta_1-ir/2,  -ir/2] \cup
[-ir/2, +ir/2] \cup
[+ir/2, -\eta_1 +ir/2],
\eqnlbl{Gammaprime}
$$
these may be transformed modulo time-exponentially decaying terms
to the elementary Fourier integrals
$$
\aligned
{1\over 2\pi }
\text{\rm P.V. }\int_{-\infty}^{+\infty}
e^{i\xi(t- (x-y)/a^{*-}_k)}
&e^{\xi^2 (-\beta^{*-}_k/{a^{*-}_k}^3 )(x-y)}
\, d\xi
\\
&\qquad = 
(4\pi \beta_k^-t)^{-1/2} e^{-(x-y-a_k^{*-}t)^2 / 4\beta_k^{*-}t} 
\endaligned
\eqnlbl{alphak2}
$$
and
$$
\aligned
{1\over 2\pi }
\text{\rm P.V. }\int_{-\infty}^{+\infty}
e^{i\xi\left( (t- |y/a^{*-}_k|)-x/a^{*-}_j\right)}
&e^{-\xi^2 \left((\beta^{*-}_k/{a^{*-}_k}^3 )|y|
+ (\beta^{*-}_j/{a^{*-}_j}^3 )|x|\right) }
\, d\xi
\\
&\qquad = 
(4\pi \bar\beta_{jk}^{*-} t)^{-1/2} e^{-(x-z_{jk}^{*-})^2 / 
4\bar\beta_{jk}^{*-} t}, 
\endaligned
\eqnlbl{alphak2}
$$
respectively, where $\bar\beta_{jk}-^{+-}$ and $z_{jk}^*$ 
are as defined in \eqnref{barbeta} and \eqnref{zjk}.
These correspond to the first and third terms in expansion
\eqnref{S}, the latter of which has an additional factor
$e^{-x}/(e^x+ e^{-x})$.
Noting that the second and fourth terms of \eqnref{S} are
time-exponentially small for $t\ge 1$, $y\le x\le 0$, and
that
$$
\aligned
|(4\pi \bar\beta_{jk}^{*-} t)^{-1/2} &e^{-(x-z_{jk}^{*-})^2 / 
4\bar\beta_{jk}^{*-} t}
(1- e^{-x}/(e^x+ e^{-x}))|\\
&\quad \le
|(4\pi \bar\beta_{jk}^{*-} t)^{-1/2} e^{-(x-z_{jk}^{*-})^2 / 
4\bar\beta_{jk}^{*-} t}
e^{-\theta |x|}
\endaligned
$$
for some $\theta>0$, so is absorbable in error term $R$, we find that the
total contribution of this term, modulo terms absorbable in $R$,
is $S$.

\medskip
{\it $R_\lambda$ term.}
Finally, we briefly discuss the estimation of error term
$$
 {1\over 2\pi i} \oint_{\Gamma} e^{\lambda t}R_\lambda (x,y)d \lambda,
\eqnlbl{error}
$$
which decomposes into the sum of integrals involving the various terms
of $R^E_\lambda$ and $R^S_\lambda$ 
given in \eqnref{RElambda} and \eqnref{RSlambda2}.
Since each of these are separately analytic, they may be
split up and estimated sharply via the {\it Riemann saddlepoint method}
(method of steepest descent), as described at great length in [ZH,HoZ.2].
That is, for each summand 
$\alpha_\lambda(x,y)\sim e^{\beta(\lambda)x+\gamma(\lambda)y}$ in 
$R_\lambda$ we deform the contour $\tilde \Gamma$ to a new, 
contour in $\Lambda$ that is a mini-max contour for
the modulus 
$$
m_\alpha(x,y,\lambda):=|e^{\lambda t}\alpha_\lambda(x,y)|=
e^{\R\lambda t + \R \beta x + \R \gamma y},
$$
passing through an appropriate saddlepoint/critical 
point of $m_\alpha(x,y,\cdot)$: necessarily lying on the
real axis, by the underlying complex symmetry resulting
from reality of operator $L$.

Since terms of each type appearing in $R$ have been sharply
estimated in [ZH], we shall omit the details, only describing
two sample calculations to illustrate the method:

{\it Example 1: $e^{-\theta|x-y|}$}. 
Contour integrals of form
$$
 {1\over 2\pi i} \oint_{\tilde\Gamma} e^{\lambda t}
e^{-\theta|x-y|}
d \lambda,
\eqnlbl{e1}
$$
arising through the pairing of fast forward and fast dual modes,
may be deformed to
$$
 {1\over 2\pi i} \int_{-\eta_1-ir/2}^{-\eta_1+ir/2} e^{\lambda t}
e^{-\theta|x-y|}
d \lambda 
\eqnlbl{d1}
$$
and estimated as
$\Cal{O}(e^{-\eta_1 t -\theta |x-y|})$, a negligible, time-exponentially
decaying contribution.
\myqed
\medskip

{\it Example 2:} 
$\lambda^r e^{ (-\lambda/a^{*-}_k + \lambda^2 \beta^{*-}_k/{a^{*-}_k}^3 )(x-y)}$.
Contour integrals of form
$$
 {1\over 2\pi i} \oint_{\tilde \Gamma} e^{\lambda t}
\lambda^r e^{ (-\lambda/a^{*-}_k + \lambda^2 \beta^{*-}_k/{a^{*-}_k}^3 )(x-y)}
d \lambda,
\eqnlbl{e2}
$$
$a_k^{*-}>0$,
arising through the pairing of slow forward and slow dual modes,
may be deformed to contour
$$
\Gamma':=
[-\eta_1-ir/2, \eta_* -ir/2] \cup
[\eta_*-ir/2,\eta_* +ir/2] \cup 
[\eta_*+ir/2, -\eta_1 +ir/2],
\eqnlbl{Gammaprime2}
$$
where saddlepoint $\eta_*$ is defined as
$$
\eta_*(x,y,t)
:=\cases
\frac{\ratio}{p} &\text{if $|\frac{\ratio}{p}|\le \varepsilon$} \\
\pm \varepsilon &\text{if $\frac{\ratio}{p} \gtrless \varepsilon$, }\endcases
\eqnlbl{123}
$$
with
$$
\ratio := \frac{x-y-a^-_k t}{2 t}, \quad p:=\frac{\beta_k^-(x-y)}{(a^-_k)^2 t}>0,
\eqnlbl{111.5}
$$
and $\eta_1$, $\varepsilon >0$ are chosen sufficiently small with
respect to $r$, to yield, modulo time-exponentially decaying terms, 
the estimate
$$
e^{-(x-y-a_k^-t)^2/4\beta^{*-}t}
 \int_{-\infty}^{+\infty}\Cal{O}(|\eta_*|^r+|\xi|^r) e^{-\theta \xi^2 t}\,d\xi
=
\Cal{O}(t^{-(r+1)/2} e^{-(x-y-a_k^-t)^2/Mt})
\eqnlbl{d2}
$$
if $|\frac{\ratio}{p}|\le \varepsilon$, and
$$
e^{-\varepsilon t/M}
 \int_{-\infty}^{+\infty}\Cal{O}(|\eta_*|^r+|\xi|^r) e^{-\theta \xi^2 t}\,d\xi
=
\Cal{O}(t^{-(r+1)/2} e^{-\eta t}).
\eqnlbl{d3}
$$
if $|\frac{\ratio}{p}|\ge \varepsilon$, $\theta>0$.
In either case, the main contribution lies along the central portion
$[\eta_*-ir/2,\eta_* +ir/2]$ of contour $\Gamma'$.

To see this, note that $|x-y|$ and $t$ are comparable when
$|\ratio/p|$ is bounded, whence the evident
spatial decay $e^{-\theta\xi^2 |x-y|}$ of the integrand along contour
$\lambda=\eta_*+i\xi$ may be converted 
to the $e^{-\xi^2 t}$ decay displayed in \eqnref{d2}--\eqnref{d3};
likewise, temporal growth in $e^{\lambda t}$ 
on the horizontal portions of the contour,
of order $\le e^{(|\eta_1|+|\varepsilon|)t}$,
is dominated by the factor $e^{-\theta \xi^2}= e^{-\theta r^2/4}$,
provided $|\eta_1|+|\varepsilon|$ is sufficiently small with respect
to $r^2$.
For $|\ratio/p|$ large, on the other hand, $\eta_*$ is uniformly negative,
and also $|x-y|$ is negligible with respect to $t$, whence the estimate
\eqnref{d3} holds trivially.

We point out that $\eta_*$ is easily determined, as the minimal point
on the real axis of the quadratic function
$$
f_{x,y,t}(\lambda):=
\lambda t (-\lambda/a^{*-}_k + \lambda^2 \beta^{*-}_k/{a^{*-}_k}^3 )(x-y),
\eqnlbl{quad}
$$
the argument of the integrand of \eqnref{e2}.
\myqed
\medskip

Other terms may be treated similarly:
All ``constant-coefficient'' terms $\phi_k \tilde \phi_{k*}$
or $\psi_k\tilde \psi_{k*}$ are of either the form treated
in Example 1 (fast modes) or in example 2 (slow modes).
Scattering pairs involving slow forward and slow dual modes
from different families (i.e.,
terms with coefficients $M_{jk}$, $d^\pm_{jk}$)
may be treated similarly as in Example 2.
Scattering pairs involving two fast modes are of the form already
treated in Example 1, since both modes of a scattering pair are
decaying.
Scattering pairs involving one fast and one slow mode
may be factored as the product of a term of the form treated in example 2
and a term that is uniformly exponentially decaying in either $x$ or $y$;
factoring out the exponentially decay, we may treat such terms as in
Example 2.

{\bf Remark \thmlbl{contourchoice}.}
In [ZH,HoZ.2], the saddlepoint contours followed were hyperbolae
through $\eta_*$, rather than the simpler, vertical contours used here.
Any reasonable contour through the saddlepoint will suffice,
yielding the same estimates.

\medskip
{\bf Scattering coefficients.}
The relation \eqnref{scattering} may now be
deduced, a posteriori, from conservation of mass in
the $u$ variable of the linearized flow \eqnref{linearized}.
For, by inspection, all terms save $E$ and $S$ in the
decomposition \eqnref{ourdecomp} of $G$ decay in $L^1$,
whence these terms must, time-asymptotically, 
carry exactly the mass in $u$ of the initial perturbation 
$\delta_y(x)I_{n}$, i.e., 
$$
\lim_{t\to \infty}
(I_n,0)\int_{-\infty}^{+\infty}\big(E(x,t;y) + S(x,t;y)\big) \, dx
=
I_{n+r}.
\eqnlbl{massbalance}
$$
Taking $y\le 0$ for definiteness,
and right-multiplying both sides of \eqnref{massbalance} by 
$(r_k^{*-t} , 0)^t $, we thus obtain the result
from definitions \eqnref{E}--\eqnref{S}.

Moreover, rewriting \eqnref{scattering} as
$$
(r_1^{*-},\dots,r_{n-i_-}^{*-}, 
r_{i_+}^{*+},\dots,r_n^{*+}, 
m_1,\dots,m_\ell)
\pmatrix 
[c_{k,-}^{1,-}]\\
\vdots\\
 [c_{k,-}^{n-i_-,-}]\\
 [c_{k,-}^{i_+,+}]\\
\vdots\\
 [c_{k,-}^{n,+}]\\
 [c_{k,-}^{1,0}]\\
\vdots\\
 [c_{k,-}^{\ell,0}]\\
\endpmatrix
=
r_k^{*-},
\eqnlbl{Deltaversion}
$$
where
$$
m_j:=\int_{-\infty}^{+\infty}(\partial/\partial \delta_j)
\bar u^\delta(x) \, dx,
\eqnlbl{mdef}
$$
we find that, under assumption ($\Cal{D}$), 
it uniquely determines the scattering
coefficients $[c_{k,-}^{j,\pm}]$.  For, the determinant
$$
\det
(r_1^{*-},\dots,r_{n-i_-}^{*-}, 
r_{i_+}^{*+},\dots,r_n^{*+}, 
m_1,\dots,m_\ell)
\eqnlbl{determinant}
$$
of the matrix on the righthand side of \eqnref{Deltaversion}
is exactly the inviscid stability coefficient $\Delta$
defined in [ZS,Z.4], hence is nonvanishing by the equivalence
of ($\Cal{D}$2) and (D2) 
(recall discussion of Section 1.2, just below Theorem \thmref{D}).
Note that \eqnref{determinant} reduces to the Liu--Majda determinant
\eqnref{liumajda} in the Lax case, $\ell=1$, with $\bar u^\delta(x)$ 
parametrized as $\bar u^\delta(x):=\bar u(x-\delta_1 x)$, for which
$m_1$ reduces to $[u]$.

Relation \eqnref{pi} follows from the observation that
$$
\aligned
P_\lambda(x,y)&:=
\text{\rm Residue }_{\lambda=0} G_\lambda(x,y)\\
&=
\sum_{j=1}^\ell (\partial/\partial \delta_j)
\pmatrix \bar u^\delta \\ \bar v^\delta \endpmatrix(x)
\, \pi_j,\\
\endaligned
$$
with $\pi_j$ defined by the first, or second expressions appearing
in \eqnref{pi}, according as $y$ is less than or equal to/greater
than or equal to zero.
For, the extended spectral theory of Section 5.4 then implies that
$\pi_j$ are the left effective eigenfunctions associated with
right eigenfunctions 
$(\partial/\partial \delta_j)(\bar u^{\delta } , \bar v^{\delta } )$.
Alternatively, \eqnref{pi} may be deduced 
by linear-algebraic manipulation directly from \eqnref{Deltaversion}
and its counterpart for $y\ge 0$ (the same identity with
$k,-$ everywhere replaced by $k,+$).

{\bf Remark \thmlbl{pidef2}.}
As pointed out in [LZ.2], the left eigenfunctions (constant vectors) 
$\pi_j=(\pi_{j,u},0)$ may themselves be determined by the relation
$$
\pi_{j,u}^t
(r_1^{*-},\dots,r_{n-i_-}^{*-}, 
r_{i_+}^{*+},\dots,r_n^{*+}, 
m_1,\dots,m_\ell)
=
(0, e_j^t),
\eqnlbl{Deltapi}
$$
from which \eqnref{pi} is easily deduced,
where $e_j$ here denotes the $j$th standard basis element in $\BbbR^\ell$.
This may be deduced a priori through the extended spectral theory of
Section 5.4;
it would interesting also
to derive \eqnref{Deltaversion} a priori, directly
through the relations \eqnref{Mplus}--\eqnref{dplusminus}.
\medskip

\myqed

{\bf Remark \thmlbl{uc} (The undercompressive case).}
In the undercompressive case, the result of Lemma \thmref{4.7}
is false, and consequently the estimates of Lemma \thmref{4.8} 
do not hold.
This fact has the implication that shock dynamics are not
governed solely by conservation of mass,
as in the Lax or overcompressive case,
but by more complicated dynamics of front interaction; 
for related discussion, see [LZ.2,Z.2].
At the level of Proposition \thmref{lowbounds}, it means
that the simple representations of $E_\lambda$ and $S_\lambda$
in terms of slow dual modes alone (corresponding to equilibrium
characteristics that are incoming to the shock) are no longer
valid in the undercompressive case, and there appear new terms 
involving rapidly decaying dual modes $\sim e^{-\theta|y|}$ related
to inner layer dynamics.
Though precise estimates can nonetheless be carried out,
we have not found a similarly compact representation of the
resulting bounds as that of the Lax/overcompressive case,
and so we shall not state them here.
We mention only that this rapid variation in the $y$-coordinate precludes
the $L^p$ stability arguments used here and in [Z.2], requiring 
instead detailed pointwise bounds as in [HZ.2,Z.5].
See [LZ.1--2,ZH,Z.2,Z.5] for further discussion of this interesting case.
\medskip

{\bf 6.2. Inner layer dynamics.}  Similarly as in [ZH], we
now investigate dynamics of the inner shock layer,
in the case that ($\Cal{D}$) does not necessarily hold,
in the process establishing the necessity of ($\Cal{D}$)
for linearized orbital stability.

\proclaim{Proposition \thmlbl{point}}
Let (H0)--(H4) hold.  
Then, there exists $\eta>0$ such that,
for $x$, $y$ restricted to any bounded
set, and $t$ sufficiently large,
$$
G(x,y;t)=
\sum_{\lambda \in \sigma_p'(L) \cap \{Re(\lambda)\ge 0\}}
e^{\lambda t}
\sum_{k\ge 0} t^{k} (L-\lambda I)^k  P_{\lambda}(x,y)
+ \Cal{O}(e^{-\eta t}),
\eqnlbl{behavior}
$$
where $ P_{\lambda}(x,y)$ is the effective projection kernel
described in Definition \thmref{proj}, and $\sigma'_-(L)$
the effective point spectrum.
\endproclaim

{\bf Proof.}  Similarly as in the proof of Proposition \thmref{aux2},
decompose $G$ into terms $I_a$, $I_b$, $II_a$, and $II_b$
of \eqnref{contours}.  Then, the same argument (Case II. $|x-y|/t$ bounded) 
yields that $I_a=H(x,t;y)$,
while $I_b$ and $II_b$ are time-exponentially small.
Observing that $H=0$ for $x$, $y$ bounded, and $t$ sufficiently large,
we have reduced the problem to the study of 
$$
II_a:=
{1\over 2\pi i} \oint_{\Gamma}
e^{\lambda t}G_\lambda (x,y)d \lambda,
\eqnlbl{reduced}
$$
where, recall,
$$
\Gamma:=
[-\eta_1-iR, \eta -iR] \cup
[\eta-iR, \eta +iR] \cup
[\eta+iR, -\eta_1 +iR].
\eqnlbl{Gammacontour}
$$
Note that the high-frequency bounds of Proposition \thmref{highbounds}
imply that $L$ has no point spectrum in $\Omega$ outside of the rectangle
$\Cal{R}$ enclosed by $\Gamma \cup [-\eta_1-iR, -\eta +iR] $.
Choosing $\eta_1$ sufficiently small, therefore, 
we may ensure that no effective point spectrum lies within the strip 
$-\eta_1\le \R \lambda < 0$, by compactness of $\Cal{R}$ together with
the fact that effective eigenvalues are isolated from one another,
as zeroes of the analytic Evans function.

Recall that $G_\lambda$ is meromorphic on $\Omega$.  
Thus, we may express \eqnref{reduced}, using Cauchy's Theorem, as
$$
{1\over 2\pi i} \int_{-\eta_1-iR}^{-\eta_1 + iR}
e^{\lambda t}G_\lambda (x,y)d \lambda
+
\text{\rm Residue }_{\lambda\in \Cal{R}} \,
e^{\lambda t}G_\lambda(x,t;y).
\eqnlbl{move}
$$ 
By Definition \thmref{proj} and Proposition \thmref{3.11}, 
$$
\eqalign{
\Residue_{\lambda \in 
\{\R \lambda \ge 0\}} 
&e^{\lambda t}G_\lambda (x,y)\cr
&=
\sum_{\lambda_0 \in \sigma_p'(L) \cap 
\{\R \lambda \ge 0\}} 
e^{\lambda_0 t}
\Residue_{\lambda_0}  e^{(\lambda -\lambda_0)t} G_\lambda(x,y) \cr
&=
\sum_{\lambda_0 \in \sigma_p'(L) \cap 
\{\R \lambda \ge 0\}} 
e^{\lambda_0 t}
\sum_{k\ge 0} (t^k/k \factorial) 
\Residue_{\lambda_0} (\lambda-\lambda_0)^k G_{\lambda}(x,y) \cr
&=
\sum_{\lambda_0 \in \sigma_p'(L) \cap 
\{\R \lambda \ge 0\}} 
e^{\lambda_0 t}
\sum_{k\ge 0} (t^k/k \factorial)  Q_{\lambda_0,k}(x,y) \cr 
&=
\sum_{\lambda_0 \in \sigma_p'(L) \cap 
\{\R \lambda \ge 0\}} 
e^{\lambda_0 t}
\sum_{k\ge 0} (t^k/k\factorial) (L-\lambda_0 I)^k P_{\lambda_0}(x,y)). \cr
}
\eqnlbl{q}
$$

On the other hand, for $x$, $y$ bounded and $t$ sufficiently large,
$t$ dominates $|x|$ and $|y|$ and we obtain from Propositions
\thmref{intbounds} and \thmref{lowbounds}  that
$$ |G_\lambda (x,t;y)|\le C 
$$
for all $\lambda \in \Omega$, for $\eta_1$ sufficiently small, hence
$$
{1\over 2\pi i} \int_{-\eta_1-iR}^{-\eta_1 + iR}
e^{\lambda t}G_\lambda (x,y)d \lambda
\le
2CRe^{-\eta_1 t}.
\eqnlbl{small}
$$
Combining \eqnref{move},  \eqnref{q}, and \eqnref{small}, we obtain the
result.
\myqed

\proclaim{Corollary \thmlbl{nec}}
Let (H0)--(H4) hold.  Then, ($\Cal{D}$) is necessary for
linearized orbital stability with respect to compactly supported
initial data, as measured in any $L^p$ norm, $1\le p\le \infty$.
\endproclaim

{\bf Proof:}
{}From \eqnref{behavior}, we find that $(\bar u,\bar v)(\cdot)$
is linearly orbitally stable only if $P_{\lambda}=0$
for all $\lambda \in \{\R \lambda \ge 0\} \setminus \{0\}$
and $\Range \Cal{P}_{0}=\text{\rm Span } 
\{ (\partial/\partial \delta_j)(\bar u^\delta,\bar v^\delta \}$.
By Propositions \thmref{3.11} and \thmref{3.12}, this is equivalent to \stab.
\myqed

\bigskip
\newsection {\bf Section 7. Stability.}
\sectionnumber=7
\theoremnumber=0
\equationnumber=0
\smallskip
\TagsOnLeft

We now establish our remaining results on linearized and
nonlinear stability.
\medskip

{\bf 7.1. Linearized stability.}  Linearized orbital stability 
follows immediately from the pointwise bounds we have established.
Similarly as in [Z.2], define the {\it linear instantaneous projection}:
$$
\aligned
\varphi(x,t)&:=
\int_{-\infty}^{+\infty} E(x,t;y)U_0(y)\, dy
\\
&=:\delta(t) (\partial/\partial \delta)
\pmatrix \bar u^\delta(x)\\ \bar v^\delta(x) \endpmatrix,\\
\endaligned
\eqnlbl{linproj}
$$
where $U_0$ is the initial data for the linearized evolution equation
\eqnref{linearized}.
The coordinates $\delta\in \Bbb{R}^{1\times \ell}$ 
may be expressed, alternatively, as
$$
\delta(t)= 
\int_{-\infty}^{+\infty} e(y,t)v_0(y)\, dy,
$$
where
$$
E(x,t;y)=:
(\partial/\partial \delta)
\pmatrix \bar u^\delta(x)\\ \bar v^\delta(x) \endpmatrix
e(y,t),
\eqnlbl{eE}
$$
i.e., 
$$
e(y,t):=
\sum_{a_k^{*-}}
\left(\errfn\left(\frac{y+a_k^{*-}t}{\sqrt{4\beta_k^{*-}t}}\right)
-\errfn \left(\frac{y-a_k^{*-}t}{\sqrt{4\beta_k^{*-}t}}\right)\right)
L_k^{*-}
\eqnlbl{e}
$$
for $y\le 0$, and symmetrically for $y\ge 0$.

Then, the solution $U$ of \eqnref{linearized} satisfies
$$
U(x,t)-\phi(x,t)=
\int_{-\infty}^{+\infty}(H + \tG)(x,t;0)U_0(y)\, dy,
\eqnlbl{convolution}
$$
where 
$$
\tG:=S+R
\eqnlbl{tG}
$$
is the regular part and $H$ the singular part of
the time-decaying portion of the Green's function $G$.

\proclaim{Lemma \thmlbl{2.05}}  Under assumptions  (H0)--(H4) and ($\Cal{D}$),
there hold: 
$$
|\int_{-\infty}^{+\infty} \tilde G(\cdot,t;y)f(y)dy|_{L^p}
\le C (1+t)^{-\frac{1}{2}(1-1/r)} |f|_{L^q},
\eqnlbl{tGbounds}
$$
$$
|\int_{-\infty}^{+\infty} \tilde G(\cdot,t;y)(0,I_r)^tf(y)dy|_{L^p}
\le C (1+t)^{-\frac{1}{2}(1-1/r)-1/2} |f|_{L^q},
\eqnlbl{tGprojbounds}
$$
$$
|\int_{-\infty}^{+\infty} \tilde G_y(\cdot,t;y)f(y)dy|_{L^p}
\le C (1+t)^{-\frac{1}{2}(1-1/r)-1/2} |f|_{L^q}
+ Ce^{-\eta t}|f|_{L^p},
\eqnlbl{tGybounds}
$$
and
$$
|\int_{-\infty}^{+\infty} H(\cdot,t;y)f(y)dy|_{L^p}
\le Ce^{-\eta t} |f|_{L^p},
\eqnlbl{Hbounds}
$$
for all $t\ge 0$, some $C$, $\eta>0$, for any
$1\le r\le p$ and $f\in L^q$ (resp. $L^p$), where $1/r+1/q=1+1/p$.
\endproclaim

{\bf Proof.}  
Bounds \eqnref{tGbounds}--\eqnref{tGybounds} follow by
the Hausdorff-Young inequality together
with bounds \eqnref{S} and \eqnref{Rbounds}--\eqnref{Rybounds},
precisely as in [Z.2].
Bound \eqnref{Hbounds} follows by direct computation and
the fact that particle paths $ z_j(y,t)$
satisfy uniform bounds 
$$
1/C \le |(\partial/\partial y)z_j| < C,
$$
for all $y$, $t$, by the fact that characteristic speeds $a_j(x)$ converge
exponentially as $x\to \pm \infty$ to constant states.
\myqed

\proclaim{Corollary \thmlbl{suff}}
Let (H0)--(H4) hold.  Then, ($\Cal{D}$) is sufficient for
$L^1\cap L^p\to L^p$ linearized orbital stability, for any $p>1$.
More precisely, 
for initial data $U_0\in L^1\cap L^p$, the solution $U(x,t)$
of the linearized evolution equations \eqnref{linearized} 
satisfies the linear decay bounds:
$$
|U(\cdot, t)-\varphi(\cdot, t)|_{L^p}\le 
C(1+t)^{-\frac{1}{2}(1-1/p)}(|U_0|_{L^1}+ |U_0|_{L^p}).
\eqnlbl{lindecaybound}
$$
\endproclaim

{\bf Proof.}
Immediate, from \eqnref{convolution} 
and bounds \eqnref{tGbounds} and \eqnref{Hbounds}, with $q=p$.
\myqed

{\bf Proof of Theorem \thmref{D}.}
Theorem \thmref{D} now follows immediately
from Proposition \thmref{nec} together with Corollary \thmref{suff}.
\myqed

\medskip

{\bf 7.2. Nonlinear stability (Proof of Theorem \thmref{nonlin}).}  
We conclude by establishing nonlinear stability of relaxation profiles 
for simultaneously symmetrizable and discrete kinetic models,
as described in Theorem \thmref{nonlin}.

\medskip
{\bf 7.2.1. Simultaneously symmetrizable models.}
We first establish the result of Theorem \thmref{nonlin} for 
relaxation profiles of simultaneously symmetrizable relaxation models,
under the weak shock assumption
$$
|(\bar u,\bar v)'|_{L^\infty}\le \varepsilon,
\eqnlbl{weak}
$$
$\varepsilon>0$ sufficiently small with respect to the parameters of
system \eqnref{general}.
Before beginning, we recall some needed preliminaries regarding the class
of simultaneously symmetrizable systems:

\medskip
{\bf Definition \thmlbl{symm}.} A relaxation model \eqnref{general}
is called {\it simultaneously symmetrizable}, if there exists
a smooth $A^0(u,v)$ symmetric, positive definite, such that
$A^0(u,v)A(u,v)$ and $A^0(u,v)Q(u,v)$ are symmetric, and 
$A^0(u,v)Q(u,v)$ is negative semidefinite, where, as in the introduction,
$A:=(df,dg)^t$ and $Q:=(0,dq)^t$.
\medskip

\proclaim{Lemma \thmlbl{skew} ([SK])}
Assuming simultaneous symmetrizability, condition (H3)
is equivalent to either of:
\medskip
(K1) \quad 
There exists a smooth skew-symmetric matrix $K(u)$ such that
$$
\text{\rm Re }(KA - A^0Q)(u) \ge \theta>0.
\eqnlbl{skew}
$$
$A^0$ as in Definition \thmref{symm}.

\medskip  
(K2)  \quad
There is no eigenvector of $A(u)$ lying in the kernel of $Q(u)$.
\endproclaim

{\bf Proof.}
These and other useful equivalent formulations are established in [SK].
\myqed

Next, we recall the basic iteration scheme of [Z.2];
for precursors of this scheme,  see [Go.2,K.1--2,LZ.1--2,ZH,HZ.1--2]. 
Restricting now to Lax-type shocks,
define the nonlinear perturbation
$$
U(x,t):= \pmatrix u \\v \endpmatrix (x+\delta(t),t)-
\pmatrix \bar u \\ \bar v \endpmatrix(x),
\eqnlbl{28}
$$
where $\delta(t)$ (estimating shock location) is to be determined
later; for definiteness, fix $\delta(0)=0$.  Then,
$$
U_t-LU=N_1(U,U)_x+ (0,I_r)^t N_2(U,U)
+\dot \delta (t)(\BU_x + U_x),
\eqnlbl{29}
$$
$$
\bar U:= \pmatrix \bar u\\ \bar v \endpmatrix,
\eqnlbl{BU}
$$
where
$$
N_j(U,U)=\Cal{O}(|U|^2)
\eqnlbl{29.1}
$$
so long as $|U|$ remains bounded.
By Duhamel's principle, and the fact that
$$
\int^\infty_{-\infty}G(x,t;y)\BU_x(y)dy=e^{Lt}\BU_x(x)=\BU_x(x),
\eqnlbl{stationary}
$$
we have 
$$
\aligned
& U(x,t)=\int^\infty_{-\infty}G(x,t;y)U_0(y)dy\\
&-\int^t_0 \int^\infty_{-\infty} G_y(x,t-s;y)(N_1(U,U)+
\dot \delta U) (y,s)dyds\\
&+\int^t_0 \int^\infty_{-\infty} G(x,t-s;y)
(0,I_r)^tN_2(U,U)(y,s)dyds\\ 
&+ \delta (t)\BU_x.\\
\endaligned
\eqnlbl{30}
$$

Defining, by analogy with the linear case,
the {\it nonlinear instantaneous projection}:
$$
\aligned
\varphi(x,t)
&:= -\delta(t)\BU_x\\
&:= \int^\infty_{-\infty}{E}(x,t;y)U_0(y) dy\\
&-\int^t_0 \int^\infty_{-\infty}E_y(x,t-s;y)(N_1(U,U)+
\dot \delta U)(y,s)dy,\\
\endaligned
\eqnlbl{proj}
$$
or equivalently, the {\it instantaneous shock location}:
$$
\aligned
\delta (t)
&=-\int^\infty_{-\infty}e(y,t) U_0(y)dy\\
&+\int^t_0\int^{+\infty}_{-\infty} e_{y}(y,t-s)(N_1(U,U)+
\dot \delta U)(y,s) dy ds,\\
\endaligned
\eqnlbl{31}
$$
where $E$, $e$ are defined as in \eqnref{E}, \eqnref{e},
and recalling \eqnref{tG}, we thus
obtain the {\it reduced equations}:
$$
\aligned
U(x,t)
&=\int^\infty_{-\infty} (H+\tG)(x,t;y)U_0(y)dy\\
&+\int^t_0
\int^\infty_{-\infty}H(x,t-s;y)(N_1(U,U)_x+
(0,I_r)^tN_2(U,U) 
+\dot \delta U_x)(y,s)dy \, ds\\
&-\int^t_0
\int^\infty_{-\infty}\tG_y(x,t-s;y)(N_1(U,U)+
\dot \delta U)(y,s) dy \, ds\\
&+\int^t_0
\int^\infty_{-\infty}\tG(x,t-s;y)(0,I_r)^tN_2(U,U)
dy \, ds,\\
\endaligned
\eqnlbl{32}
$$
and, differentiating \eqnref{31} with respect to $t$,
$$
\aligned
\dot \delta (t)
&=-\int^\infty_{-\infty}e_t(y,t) U_0(y)dy\\
&+\int^t_0\int^{+\infty}_{-\infty} e_{yt}(y,t-s)(N_1(U,U)+
\dot \delta U)(y,s) dy ds.\\
\endaligned
\eqnlbl{33}
$$

{\bf Notes:}  
In deriving \eqnref{33}, we have used the fact 
that $e_y (y,s)\rightharpoondown 0$ as $s \to 0$, as the 
difference of approaching heat kernels, in evaluating the 
boundary term
$$
\int^{+\infty}_{-\infty} e_y (y,0)(N_1(U,U)+
\dot \delta U)(y,t)dy=0.
\eqnlbl{35}
$$
(Indeed, $|e_y(\cdot, s)|_{L^1}\to 0$, see Remark 2.6, below).
We have also 
used the fact that $E(0,I_r)\equiv 0$,
to eliminate the term corresponding to nonlinear source $(0,I_r)N_2$
in \eqnref{proj}--\eqnref{31}. 
Moreover, we have in the first place used the fact that the
dimension $\ell$ of the stationary manifold is one, a property 
of Lax and undercompressive shocks, in order to factor out
$\BU(x)$ in the first line of \eqnref{proj};
for this reason, our analysis is inherently limited to 
the case of Lax shocks (recall that the favorable Green's function 
bounds of Theorem \thmref{greenbounds} are available for Lax and 
overcompressive shocks).
\medskip

The defining relation $ \delta (t)\bar u_x:= -\varphi$ in \eqnref{proj}
can be motivated heuristically by
$$
\eqalign{
\tilde U(x,t)-\varphi(x,t) \sim U&= \pmatrix u\\v\endpmatrix (x+\delta(t),t) 
- \pmatrix \bar u\\ \bar v\endpmatrix(x) \cr
&\sim \tilde U(x,t) + \delta (t)\bar U_x(x),}
$$
where $\tilde U$ denotes the solution of the linearized
perturbation equations, and $\bar U$ the background profile.  
Alternatively, it can be thought
of as the requirement that the instantaneous projection 
of the shifted (nonlinear) perturbation variable $U$ be zero, [HZ.1--2].

\proclaim{Lemma \thmlbl{3} [Z.2]}  The kernel ${e}$ satisfies
$$
|{e}_y (\cdot, t)|_{L^p},  |{e}_t(\cdot, t)|_{L^p} 
\le C t^{-\frac{1}{2}(1-1/p)},
\eqnlbl{36}
$$
$$
|{e}_{ty}(\cdot, t)|_{L^p} 
\le C t^{-\frac{1}{2}(1-1/p)-1/2},
\eqnlbl{37}
$$
for all $t>0$.  Moreover, for $y\le 0$ we have the pointwise bounds
$$
|{e}_y (y,t)|, |{e}_t (y,t)| \le Ct^{-\frac{1}{2}} e^{-\frac{(y+a_-t)^2}{Mt}},
\eqnlbl{38}
$$
$$
|{e}_{ty} (y,t)| \le C t^{-1} e^{-\frac{(y+a_-t)^2}{Mt}}, 
\eqnlbl{39}
$$
for $M>0$ sufficiently large (i.e. $>4b_\pm$),
and symmetrically for $y\ge 0$.
\endproclaim

{\bf Proof.}  
For definiteness, take $y \le 0$. 
Then, \eqnref{e} gives
$$
{e}_y(y,t)=
\left(\frac{1}{u_+-u_-}\right)
\left(K(y+a_-t,t)-K(y-a_-t,t)\right),
\eqnlbl{41}
$$
$$
{e}_t(y,t)=
\left(\frac{1}{u_+-u_-}\right)
\left( (K+K_y)(y+a_-t,t)-(K+K_y)(y-a_-t,t)\right),
\eqnlbl{41}
$$
$$
{e}_{ty}(y,t)=
\left(\frac{1}{u_+-u_-}\right)
\left( (K_y+K_{yy})(y+a_-t,t)-(K_y+K_{yy})(y-a_-t,t)\right),
\eqnlbl{42}
$$
where
$$
K(y,t):= \frac{e^{-y^2/4b_-t}}{\sqrt{4\pi b_-t}}
\eqnlbl{43}
$$
denotes an appropriate heat kernel.  The pointwise bounds
\eqnref{38}--\eqnref{39} follow immediately for $t\ge 1$ 
by properties of the heat kernel, in turn yielding \eqnref{36}--\eqnref{37} 
in this case.  The bounds for small time $t\le 1$ follow from estimates
$$
\aligned
|K_y (y+a_-t,t)-K_y (y-a_-t,t)|
&=|\int^{y-a_-t}_{y+a_-t}K_{yy}(z,t)dz| \\
&\le Ct^{-3/2}\int^{y-a_-t}_{y+a_-t} e^{\frac{-z^2}{Mt}} dz\\
&\le Ct^{-1/2}e^{-\frac{(y+a_-t)^2}{Mt}},\\
\endaligned
\eqnlbl{44}
$$
and, similarly,
$$
\aligned
|K_{yy}(y+a_-t,t)-K_{yy}(y-a_-,t)|
&= |\int^{a_-t}_{-a_-t}K_{yyy}(z,t)dz|\\
&\le Ct^{-2}\int^{y-a_-t}_{y+a_-t} e^{\frac{-z^2}{Mt}} dz,\\
&\le Ct^{-1}e^{-\frac{(y+a_-t)^2}{Mt}}.\\
\endaligned
\eqnlbl{45}
$$
The bounds for $|{e}_y|$ are again immediate.
Note that we have taken crucial account of cancellation in
the small time estimates of $e_t$, $e_{ty}$.
\qed

{\bf Remark \thmlbl{2.6}:}
For $t\le 1$, a calculation analogous to that of
\eqnref{44} yields
$
|e_y(y,t)|\le 
Ce^{-\frac{(y+a_-t)^2}{Mt}},
$
and thus $|e(\cdot,s)|_{L^1}\to 0$ as $s\to 0$.

\medskip

With these preparations, we are ready to carry out our analysis:
\medskip

{\it Preliminary result.} To illustrate the method,
we first carry out the analysis completely for the 
simplest case of discrete kinetic models with
the minimal regularity assumptions 
$$
|U_0|_{L^1\cap H^1}\le \zeta_0,
\eqnlbl{minreg}
$$
(still keeping the weak shock assumption \eqnref{weak}).
That is, we establish the second assertion of Theorem \thmref{nonlin}.
Define
$$
\aligned
\zeta(t)
&:= \sup_{0\le s \le t}
\Big( \, |U(\cdot, s)|_{L^2}(1+s)^{1/4}
+ |\dot \delta (s)|(1+s)^{\frac{1}{2}}
+ |\delta (s)|\\
&\qquad
+ |U_x(\cdot, s)|_{L^2}(1+s)^{-1/4} 
+ |U(\cdot, s)|_{L^\infty}  
\, \Big).\\
\endaligned
\eqnlbl{48}
$$
We shall establish:

{\it Claim.} For all $t\ge 0$ for which a solution exists 
with $\zeta$ uniformly bounded by some fixed, sufficiently small constant,
there holds
$$
\zeta(t) \leq C_2(\zeta_0 + \zeta(t)^2).
\eqnlbl{claim}
$$
\medskip
{}From this result, it follows by continuous induction that,
provided $\zeta_0 < 1/4C_2$,  
there holds 
$$
\zeta(t) \leq 2C_2\zeta_0
\eqnlbl{bd}
$$
for all $t\geq 0$ such that $\zeta$ remains small.
By standard short-time theory/continuation, we find that
the solution (unique, in this regularity class) in fact remains in 
$L^2\cap H^1$ for all $t\ge 0$,
with bound \eqnref{bd}, yielding existence and the claimed $L^p$ bounds.
Thus, it remains only to establish the claim above.
\medskip

{\it Proof of Claim.}
We must establish bounds on $|U|_{L^2}$,
$|\dot \delta|$, $|\delta|$,
$|U_x|_{L^2}$, and $|U|_{L^\infty}$.
Provided that we can establish the others, the
$L^\infty$ bound follows by Sobolev estimate
$$
|U|_{L^\infty}\le
|U|_{L^2}^{1/2} |U_x|_{L^2}^{1/2}.
$$

Likewise, $|U_x|_{L^2}$ may be estimated in by now routine fashion
using a nonstandard energy estimate of the type 
formalized by Kawashima [Kaw].  The origin of this approach goes back
to [Ka,MN] in the context of gas dynamics; see, e.g., [HoZ.2] for further
discussion/references.

We first perform a standard, ``Friedrichs-type'' 
estimate for symmetrizable hyperbolic systems.
Expanding and differentiating \eqnref{29}, and using the key fact
that $N_1(U,U)\equiv 0$ in the case of discrete kinetic models,
we obtain
$$
{U_x}_t-(AU)_{xx}-(QU)_x=(0,I_r)^t N_2(U,U)_x
+\dot \delta (t)(\BU_{xx} + U_{xx}),
\eqnlbl{kderiveq}
$$
where
$$
N_2(U,U)_x= \Cal{O}(|U|^2+ |U||U_x|)
\eqnlbl{Nderiv}
$$
provided $|U|_{L^\infty}$ remains uniformly bounded as we have assumed. 
Let $A^0=A^0(x)$ denote the simultaneous symmetrizer of $A$, $Q$.
Taking the $L^2$ inner product of $\tilde A^0(x)U_x$ against \eqnref{kderiveq}, 
we obtain after rearrangement/integration by parts, and applications of
Young's inequality, the energy estimate
$$
\frac{1}{2}\langle A^0 U_x,U_x\rangle_t 
\le \langle U_x, A^0 QU_x \rangle
 + \varepsilon |U_x|_{L^2}^2 + 
\varepsilon^{-1}\Cal{O}(|U|_{L^2}^2 + |\dot\delta(t)|^2),
\eqnlbl{est1}
$$
where $\varepsilon>0$ is as in \eqnref{weak}.  Here, we have freely
used assumption \eqnref{weak} and the consequent higher order
bound $|\BU''|\le C\varepsilon$, plus the bounds $|\BU_x|_{L^p}$,
$|\BU_{xx}|_{L^p}\le C$.

Next, we perform a nonstandard, ``Kawashima-type'' derivative estimate.
Let $K=K(x)$ denote the skew-symmetric matrix described in
Lemma \thmref{skew}, associated with $A$, $Q$.
Taking the $L^2$ inner product of $U_x$ against $K$ times the undifferentiated
equation
$$
{U}_t-(AU)_{x}-QU=(0,I_r)^t N_2(U,U)
+\dot \delta (t)(\BU_{x} + U_{x}),
\eqnlbl{kderiveq}
$$
and noting that (integrating by parts, and using skew-symmetry of $K$)
$$
\aligned
\frac{1}{2}\langle U_x,KU \rangle_t
&=
 \frac{1}{2}\langle U_x,KU_t \rangle+
\frac{1}{2}\langle U_{xt},KU \rangle\\
&=
\frac{1}{2}\langle U_x,KU_t \rangle+
-\frac{1}{2}\langle U_{t},KU_x \rangle
-\frac{1}{2}\langle  U_{t},K_x U \rangle\\
&=
\langle U_x,KU_t \rangle
+\frac{1}{2}\langle U,K_x U_t \rangle,\\
\endaligned
\eqnlbl{parts}
$$
we obtain after rearrangement the auxiliary energy estimate:
$$
\frac{1}{2}\langle U_x,KU \rangle_t \le
-\langle  U_x , KA U_x\rangle
 + \varepsilon |U_x|_{L^2}^2 + 
\varepsilon^{-1}\Cal{O}(|U|_{L^2}^2 + |\dot\delta(t)|^2).
\eqnlbl{est2}
$$
Adding \eqnref{est1} and \eqnref{est2}, and recalling \eqnref{skew},
we obtain, finally:
$$
\aligned
\frac{1}{2}(
\langle A^0 U_x,U_x\rangle
+ \langle U_x,KU \rangle)_t &\le
\frac{1}{2}\langle U_x, (A^0 Q - KA) U_x\rangle
+\varepsilon^{-1}\Cal{O}(|U|_{L^2}^2 + |\dot\delta(t)|^2)\\
&\le
C(|U|_{L^2}^2 + |\dot\delta(t)|^2),\\
\endaligned
\eqnlbl{estfinal}
$$
provided $\varepsilon $ is sufficiently small
that $\varepsilon |U_x|_{L^2}^2 $ may be absorbed in
the negative definite term
$\frac{1}{2}\langle U_x, (A^0 Q - KA) U_x\rangle$.

Integrating \eqnref{estfinal} from $0$ to $t$, and recalling
definition \eqnref{48}, thus yields
$$
(\langle A^0 U_x,U_x\rangle + \langle U_x,KU \rangle)|^t_0 \le
C\int_0^t (|U|_{L^2}^2 + |\dot\delta|^2)(s) \,ds.
\eqnlbl{intcalc0}
$$
Assuming, then, that claim \eqnref{claim} holds when restricted to $|U|_{L^2}$
and $|\dot \delta|$, i.e.,
$$
|U|_{L^2}(s)\le
C(\zeta_0 + \zeta(t)^2)(1+s)^{-1/4}
\eqnlbl{restrictU}
$$
and
$$
|\dot \delta|(s)\le
C(\zeta_0 + \zeta(t)^2)(1+s)^{-1/2},
\eqnlbl{restrictdeltadot}
$$
for all $0\le s\le t$, we obtain
$$
\aligned
(\langle A^0 U_x,U_x\rangle + \langle U_x,KU \rangle)|^t_0 &\le
C(\zeta(0)+ \zeta(t)^2)^2\int_0^t 
(1+s)^{-1/2}\, ds\\
&\le 
C(\zeta(0)+ \zeta(t)^2)^2 (1+t)^{1/2}.
\endaligned
\eqnlbl{intcalc}
$$
Using Young's inequality to bound
$$
\aligned
\langle U_x,KU \rangle(t) &\le
\epsilon |U_x|_{L^2}^2(t) + C\epsilon^{-1} |U|_{L^2}^2(t)\\
&\le \epsilon |U_x|_{L^2}^2(t) + C\epsilon^{-1} 
(\zeta(0)+ \zeta(t)^2)^2 
(1+t)^{-1/2},\\
\endaligned
\eqnlbl{youngbd}
$$
and recalling, by assumption, that 
$(\langle A^0 U_x,U_x\rangle + \langle U_x,KU \rangle)(0)\le C\zeta_0^2 $,
we may rearrange \eqnref{intcalc} to obtain
$$
\langle A^0 U_x,U_x\rangle \le C(\zeta(0)+ \zeta(t)^2)
(1+t)^{1/2},
$$
which, by uniform positive definiteness of $A^0$,
verifies \eqnref{claim} for $|U_x|_{L^2}$ as well.

Thus, in order to establish the result, 
we have only to verify claim \eqnref{claim} when restricted 
to $|U|_{L^2}$, $|\dot\delta|$, and $|\delta|$:
that is, to establish \eqnref{restrictU}, \eqnref{restrictdeltadot}, and
$$
|\delta|(s)\le
C(\zeta_0 + \zeta(t)^2),
\eqnlbl{restrictdelta}
$$
for all $0\le s\le t$.
By \eqnref{32}--\eqnref{33}, and $N_1\equiv 0$, we have
$$
\aligned
\big|U\big|_{L^2}(t)
&\le \big|\int^\infty_{-\infty} (H+\tG)(x,t;y)U_0(y)\big|_{L^2}dy\\
&+\big|\int^t_0
\int^\infty_{-\infty}H(x,t-s;y)
(0,I_r)^tN_2(U,U) dy\,ds\big| \\
&+\big|\int^t_0
\int^\infty_{-\infty}H(x,t-s;y)
\dot \delta U_x(y,s)dy \, ds\big|_{L^2}\\
&+\big|\int^t_0
\int^\infty_{-\infty}\tG_y(x,t-s;y)
\dot \delta U(y,s) dy \, ds\big|_{L^2}\\
&+\big|\int^t_0
\int^\infty_{-\infty}\tG(x,t-s;y)(0,I_r)^tN_2(U,U)
dy \, ds\big|_{L^2},\\
&=: I_a + I_b + I_c + I_d+ I_e,\\
\endaligned
\eqnlbl{reduced}
$$
$$
\aligned
\big|\dot \delta\big|(t)
&=\big|\int^\infty_{-\infty}e_t(y,t) U_0(y)dy\big|\\
&+\big|\int^t_0\int^{+\infty}_{-\infty} e_{yt}(y,t-s)
\dot \delta U(y,s) dy ds\big|\\
&=: II_a + II_b.\\
\endaligned
\eqnlbl{deltadot}
$$
and
$$
\aligned
\big|\delta\big|(t)
&=\big|\int^\infty_{-\infty}e(y,t) U_0(y)dy\big|\\
&+\big|\int^t_0\int^{+\infty}_{-\infty} e_{y}(y,t-s)
\dot \delta U(y,s) dy ds\big|\\
&=: III_a + III_b.\\
\endaligned
\eqnlbl{delta}
$$
We estimate each term in turn, following the approach of [Z.2].

The linear term $I_a$  satisfies bound
$$
I_a \le  C\zeta_0 (1+t)^{-1/4},
\eqnlbl{Ia}
$$ 
as already shown in
the proof of Corollary \thmref{suff}.
Likewise, applying the bounds of Lemma \thmref{2.05} together
with definition \eqnref{48}, we have:
$$
\aligned
I_b&= \big|\int^t_0
\int^\infty_{-\infty}H(x,t-s;y)
(0,I_r)^tN_2(U,U) dy\,ds\big|_{L^2} \\
&\le
C\int_0^t e^{-\eta (t-s)}|U|_{L^\infty}|U|_{L^2}(s) ds\\
&\le
C\zeta(t)^2
\int_0^t e^{-\eta (t-s)}(1+s)^{-1/4}ds\\
&\le
C\zeta(t)^2 (1+t)^{-1/4},\\
\endaligned
\eqnlbl{Ib}
$$
$$
\aligned
I_c&= \big|\int^t_0
\int^\infty_{-\infty}H(x,t-s;y)
\dot \delta U_x(y,s)dy \, ds\big|_{L^2}\\
&\le
C\int_0^t e^{-\eta (t-s)}|\dot\delta||U_x|_{L^2}(s) ds\\
&\le
C\zeta(t)^2
\int_0^t e^{-\eta (t-s)}(1+s)^{-1/4}ds\\
&\le
C\zeta(t)^2 (1+t)^{-1/4},\\
\endaligned
\eqnlbl{Ic}
$$
$$
\aligned
I_d&= 
\big|\int^t_0
\int^\infty_{-\infty}\tG_y(x,t-s;y)
\dot \delta U(y,s) dy \, ds\big|_{L^2}\\
&\le
C\int_0^t \big(1+(t-s)\big)^{-1/2}|\dot\delta||U|_{L^2}(s) ds\\
&\le
C\zeta(t)^2
\int_0^t \big(1+(t-s)\big)^{-1/2}
(1+s)^{-3/4}ds\\
&\le
C\zeta(t)^2 (1+t)^{-1/4},\\
\endaligned
\eqnlbl{Id}
$$
and
$$
\aligned
I_e&= 
\big|\int^t_0
\int^\infty_{-\infty}\tG(x,t-s;y)(0,I_r)^tN_2(U,U)
dy \, ds\big|_{L^2},\\
&\le
C\int_0^t \big(1+(t-s)\big)^{-3/4}|U|_{L^2}^2(s) ds\\
&\le
C\zeta(t)^2
\int_0^t \big(1+(t-s)\big)^{-3/4}
(1+s)^{-1/2}ds\\
&\le
C\zeta(t)^2 (1+t)^{-1/4},\\
\endaligned
\eqnlbl{Ie}
$$
Summing bounds \eqnref{Ia}--\eqnref{Ie}, we obtain \eqnref{restrictU},
as claimed.

Similarly, applying the bounds of Lemma \thmref{3} together
with definition \eqnref{48}, we find that
$$
\aligned
II_a&=
\big|\int^\infty_{-\infty}e_t(y,t) U_0(y)dy\big|\\
&\le
|e_t(y,t)|_{L^\infty}(t) |U_0|_{L^1}\\
&\le
C\zeta_0 (1+t)^{-1/2}
\endaligned
\eqnlbl{IIa}
$$
and 
$$
\aligned
II_b&=
\big|\int^t_0\int^{+\infty}_{-\infty} e_{yt}(y,t-s)
\dot \delta U(y,s) dy ds\big|\\
&\le
\int^t_0
|e_{yt}|_{L^2}(t-s) |\dot\delta||U|_{L^2}(s) ds\\
&\le
C\zeta(t)^2 
\int^t_0
(t-s)^{-3/4} (1+s)^{-3/4} ds\\
&\le
C\zeta(t)^2 (1+t)^{-1/2},
\endaligned
\eqnlbl{IIb}
$$
while
$$
\aligned
III_a&=
\big|\int^\infty_{-\infty}e(y,t) U_0(y)dy\big|\\
&\le
|e(y,t)|_{L^\infty}(t) |U_0|_{L^1}\\
&\le
C\zeta_0 
\endaligned
\eqnlbl{IIIa}
$$
and 
$$
\aligned
III_b&=
\big|\int^t_0\int^{+\infty}_{-\infty} e_{y}(y,t-s)
\dot \delta U(y,s) dy ds\big|\\
&\le
\int^t_0
|e_{y}|_{L^2}(t-s) |\dot\delta||U|_{L^2}(s) ds\\
&\le
C\zeta(t)^2 
\int^t_0
(t-s)^{-1/4} (1+s)^{-3/4} ds\\
&\le
C\zeta(t)^2. 
\endaligned
\eqnlbl{IIIb}
$$
Summing \eqnref{IIa}--\eqnref{IIb} and
\eqnref{IIIa}--\eqnref{IIIb},
we obtain \eqnref{restrictdeltadot} and \eqnref{restrictdelta},
as claimed.

This completes the proof of the claim, and the result.
\myqed

{\bf Remark \thmlbl{SX}.} The hypotheses $U_0\in L^1\cap H^1$ 
made in this subsection on initial perturbations are similar
to but slightly weaker than those for the stability result announced 
in [S,SX.2] for weak profiles of the discrete kinetic Broadwell model.
We obtain sharp rates of decay in $L^2$ and smallness
in $L^\infty$, with nonsharp rates of decay in $L^p$, $2<p<\infty$,
whereas they assert decay in $L^p$, all $p\ge 2$, with no rates.  
(Recall that spectral stability holds for weak Broadwell shocks,
by the results of [KM], so that we obtain a complete stability
result in this case).
\medskip

{\it The general case.} 
We next sharpen our analysis for simultaneously symmetrizable
systems under the strengthened regularity assumption 
$$
|U_0|_{L^1\cap H^2}\le \zeta_0,
\eqnlbl{strenreg}
$$
and assuming that $q\in C^3$ for discrete kinetic models as well
as general systems.
That is, we establish the first assertion of Theorem \thmref{nonlin}.

Similarly as in \eqnref{48}, define
$$
\aligned
\zeta(t)
&:= \sup_{0\le s \le t}
\Big( \, |U(\cdot, s)|_{H^1}(1+s)^{1/4}
+ |\dot \delta (s)|(1+s)^{1/2}
+ |\delta (s)|\\
&\qquad
+ |U_{xx}(\cdot, s)|_{L^2}(1+s)^{-1/4} 
+ |U(\cdot, s)|_{L^\infty} (1+s)^{1/2} 
+ |U_x(\cdot, s)|_{L^\infty}  
\, \Big).\\
\endaligned
\eqnlbl{zeta2}
$$
We shall establish, as before, the claim \eqnref{claim}, from
which the desired sharp $L^2$ and $L^\infty$ bounds immediately follow,
and thus, by interpolation, the claimed sharp bounds in 
all $L^p$ such that $2\le p\le \infty$.

\medskip
{\it Proof of claim.}
As in the previous proof, we first observe that the bound on
$|U_x|_{L^\infty}$ follows by Sobolev embedding, provided that
we can establish the claimed $H^1$ and $H^2$ bounds on $U$.

As before, we next show by energy estimates that the $H^2$ estimate follows
provided we can establish the claimed bounds on $|U|_{H^1}$,
$|\dot\delta|$, and $|U|_{L^\infty}$.
However, we must now take a bit more care in this step, 
due to the presence of the higher-order
nonlinearity term $N_1(U,U)_x$ in \eqnref{29} in the general case.
To avoid this term, we rewrite \eqnref{29} more strategically,
using the quadratic Leibnitz relation
$$
A_2 U_2- A_1 U_1= A_2(U_2-U_1) + (A_2-A_1)U_1,
\eqnlbl{nLeib}
$$
in the quasilinear form
$$
{U}_t-\tilde AU_{x}-\tilde QU=M_1(U)\BU_x+ (0,I_r)^t M_2(U)\BU
+\dot \delta (t)(\BU_{x} + U_{x}),
\eqnlbl{strat}
$$
where 
$$ 
\aligned
\tilde A&:=(df,dg)^t(u,v)(x+\delta(t),t),\\
\tilde Q&:=(0,dq)^t(u,v)(x+\delta(t),t),\\
\endaligned
\eqnlbl{tA}
$$
and
$$ 
\aligned
M_1(U)&=\Cal{O}(|U|):= \tilde A(x,t)- A(x),\\
M_2(U)&=\Cal{O}(|U|):= \tilde Q(x,t)- Q(x),\\
\endaligned
\eqnlbl{Mj}
$$
$A(x)$ and $Q(x)$ defined as usual as $(df,dg)^t(\bar u, \bar v)$
and $(0,dq)^t(\bar u,\bar v)$,
or, differentiating twice:
$$
{U_{xx}}_t-(\tilde AU_{x})_{xx}-(\tilde QU)_{xx}=
(M(U)\BU_x)_{xx}+ (0,I_r)^t (M_2(U)\BU)_{xx}
+\dot \delta (t)(\BU_{xxx} + U_{xxx}).
\eqnlbl{gderiveq}
$$

Let $\tilde A^0$ denote the simultaneous symmetrizer of $\tilde A$, $\tilde Q$.
Taking the $L^2$ inner product of $\tilde A^0(x)U_{xx}$ 
against \eqnref{gderiveq}, 
we obtain after rearrangement/integration by parts, and several
applications of Young's inequality, the energy estimate
$$
\langle \tilde A^0 U_{xx},U_{xx}\rangle_t 
\le \langle U_{xx},\tilde A^0 \tilde QU_{xx}\rangle
 + \varepsilon |U_{xx}|_{L^2}^2 + 
\varepsilon^{-1}\Cal{O}(|U|_{H^1}^2 + |\dot\delta(t)|^2),
$$
where $\varepsilon>0$ is as in \eqnref{weak}.  Here, we have freely
used assumption \eqnref{weak} and the consequent higher order
bounds $|\BU''|\le C\varepsilon$, $|\BU'''|\le C\varepsilon$, 
plus the bounds $|\BU_x|_{L^p}$, $|\BU_{xx}|_{L^p}$.
$|\BU_{xxx}|_{L^p}\le C$, and also used the fact that
$(|U|_{L^\infty}+ |U_x|_{L^\infty})\le C\zeta(t)$
is assumed to remain small, to bound the new term
$$
\langle \tilde A^0_t U_{xx},U_{xx}\rangle
=\Cal{O}(|U_{xx}|_{L^2}|U_t|_{L^\infty})
=\Cal{O}(|U_{xx}|_{L^2}(|U|_{L^\infty}+ |U_x|_{L^\infty})
\eqnlbl{newterm}
$$
arising from the fact that $\tilde A^0$ is now time-dependent.

Similarly, taking the $L^2$ inner product of $U_{xx}$
against $\tilde K $ times the singly differentiated equation
$$
{U_{x}}_t-(\tilde AU_{x})_{x}-(\tilde QU)_{x}=
(M(U)\BU_x)_{x}+ (0,I_r)^t (M_2(U)\BU)_{x}
+\dot \delta (t)(\BU_{xx} + U_{xx}),
\eqnlbl{g1deriveq}
$$
where $\tilde K$ denotes 
the skew-symmetric matrix described in
Lemma \thmref{skew}, associated with $\tilde A$, $\tilde Q$,
and proceeding as before, we obtain the auxiliary energy estimate:
$$
\frac{1}{2}\langle \tilde U_{xx},KU_x \rangle_t \le
-\langle  U_{xx}, \tilde KA U_{xx}\rangle
 + \varepsilon |U_{xx}|_{L^2}^2 + 
\varepsilon^{-1}\Cal{O}(|U|_{H^1}^2 + |\dot\delta(t)|^2).
\eqnlbl{est2}
$$
From here, we may proceed as in the previous case to obtain the
desired bound on $|U|_{H^2}$, 
assuming the claimed bounds on $|U|_{H^1}$ and $|\dot \delta|$.

But, the bounds on $|U|_{H^1}$, $|\dot \delta|$, and $|\delta|$
follow exactly as in the previous case, once we observe that
$\tilde G_x$ breaks up into the sum of a time-exponentially decaying
delta-function term $H_1$ analogous to $H$ and a kernel
$\tilde G_1$ obeying the same bounds as $\tilde G$.
This term therefore gives the same bounds as in the $L^2$
estimates of the previous case.  The term involving $H$,
on the other hand, may be estimated 
(see comment just beneath (7.75), below) as
$$
\aligned
\big|\int^t_0
&\int^\infty_{-\infty}H(x,t-s;y)
\Cal{O}(|U_{xx}||U|+ |U_x|^2)(y,s)dy \, ds\big|_{L^2}\\
&\le
C\int_0^t e^{-\eta (t-s)}||(|U_{xx}|_{L^2}|U|_{L^\infty}
+ |U_x|_{L^2}|U_x|_{L^\infty})
(s) ds\\
&\le
C\zeta(t)^2
\int_0^t e^{-\eta (t-s)} (1+s)^{-1/4} ds\\
&\le
C\zeta(t)^2 (1+t)^{-1/4},\\
\endaligned
$$
again giving the desired result.

Finally the claimed bound on $|U|_{L^\infty}$ follows, for terms
involving the regular kernel $\tilde G$, by the exactly the computations 
used to estimate $|\dot \delta|$, and for terms involving the singular
kernel $H$, using the new bound $|U_x|_{L^\infty}\le C\zeta(t)$ coming 
from definition \eqnref{zeta2}, and proceeding similarly as in the
estimation of $U_{H^1}$.
We omit the details, which are tedious but straightforward.
\myqed

\medskip
{\bf 7.2.2. Discrete kinetic models.}
Finally, we briefly indicate the extension of our arguments to
general $L^p$, $1\le p\le \infty$, and strong relaxation profiles, 
for discrete kinetic models with sufficiently regular coefficients
$q\in C^3$, not necessarily simultaneously symmetrizable, and
perturbations that are sufficiently small in $W^{3,1}\cap W^{3,\infty}$.
That is, we establish the third, and last assertion of Theorem \thmref{nonlin},
completing the proof of the theorem.

The main idea in this subsection is to take advantage of the special
property of discrete kinetic models that nonlinear terms, in the original,
unshifted equations \eqnref{general}, appear without derivatives.
Likewise, the basic, unshifted perturbation variable 
$$
\tilde U:=(u,v)(x,t)- (\bar u,\bar v)(x)
\eqnlbl{tU}
$$
satisfies an equation with this same property:
$$
\tilde U_t -(A\tilde U)_x -Q\tilde U= (0,I_r)^t\tilde N_2(\tilde U,\tilde U),
\eqnlbl{tUeq}
$$
where $\tilde N(\tilde U,\tilde U)=\Cal{O}(|\tilde U|^2)$, and, more generally,
$$
|(d/dx)^k \tilde N(\tilde U,\tilde U)|\le C
\sum_{0\le \alpha + \beta \le k: \, \alpha,\beta \ge 0}
|(d/dx)^\alpha \tilde U (d/dx)^\beta \tilde U|,
\eqnlbl{xtUbd}
$$
$k=0,\dots,3$,
provided $|\tilde U|_{W^{1,\infty}}$ remains bounded, 
as we shall assume it does throughout our argument.
(Note: a straightforward Leibnitz calculation shows that,
provided coefficients $q$ are sufficiently regular,
bound \eqnref{xtUbd} extends to $k=3K+2$, provided
$|\tilde U|_{W^{K,\infty}}$ remains bounded).

Using this observation, and the previously observed fact that $x$-derivative
bounds on the Green's function, modulo time-exponentially decaying
singular terms (i.e., delta-functions and their derivatives), are essentially
independent of the order of differentiation, we may use 
Green's function bounds instead of energy estimates to bound derivatives
of the solution.  That is, the integral equations for $\tilde U$ close,
with no loss of derivatives, in the special case of discrete kinetic models.
On the other hand, the variable $\tilde U$ is not expected to decay, since
it includes a possible stationary component corresponding to 
translation of the original profile.  Thus, the bounds we obtain by this
technique are quite crude, namely 
$$
|U_{xxx}|_{L^p} \le 
C|\tilde U_{xxx}|_{L^p}\sim t,
\eqnlbl{order}
$$
where $U$ as usual denotes the shifted perturbation variable
$U(x,t):=(u,v)^t(x+\delta(t),t)-(\bar u,\bar v)^t(x)$.
However, this is sufficient for our argument, since as we have
seen in previous subsections, $|(d/dx)^k U|_{L^p}$ bounds improve as
$k$ decreases, by order $t^{-1/2}$ for each successive derivative.
At level $|U_x|_{L^p}$, we therefore obtain a uniform $\Cal{O}(1)$
bound, which, as in the proof of the previous case, is sufficient
to establish sharp bounds in all $L^p$, $1\le p\le \infty$.
Note that the resulting order-one $|U|_{W^{1,\infty}}$ bound,
which implies a corresponding bound on $|\tilde U|_{W^{1,\infty}}$,
is consistent with the assumption under which we obtained bound
\eqnref{xtUbd}, an important point.

{\bf Remark \thmlbl{reason}.} The equations for the more favorable
(i.e., decaying) shifted perturbation $U$ feature the nonlinear forcing 
term $\dot \delta(t)U_x$, which does involve a derivative.  Thus, the
integral equations for $U$ do not close, involving a loss of one derivative,
and so we 
cannot directly estimate $|U_{xxx}|_{L^p}$ 
by this approach.
\medskip

We now describe the argument in somewhat further detail.
Define
$$
\aligned
\zeta(t)
&:= \sup_{0\le s \le t}
\Big( \, 
|U(\cdot, s)|_{L^\infty}(1+s)^{1/2}
+|U(\cdot, s)|_{L^1}
+ |\dot \delta (s)|(1+s)^{\frac{1}{2}}
+ |\delta (s)|\\
&
+ |U_{x}(\cdot, s)|_{L^1\cap L^\infty}
+ |U_{xx}(\cdot, s)|_{L^1\cap L^\infty}(1+s)^{-1/2} 
+ |U_{xxx}(\cdot, s)|_{L^1\cap L^\infty}(1+s)^{-1} 
\, \Big).\\
\endaligned
\eqnlbl{zeta3}
$$
As usual, our goal is to establish 
$$
\zeta(t) \leq C_2(\zeta_0 + \zeta(t)^2)
\eqnlbl{repeatclaim}
$$
for all $t\ge 0$ for which a solution exists with $\zeta(t)$ sufficiently
small, from which we may conclude as before both global existence and
the bound
$$
\zeta(t)\le C\zeta_0
\eqnlbl{repeatbd}
$$
for $\zeta_0$ sufficiently small.
From \eqnref{repeatbd} and the definition \eqnref{zeta3}, we then
immediately obtain the claimed sharp bounds in $L^1$ and $L^p$,
and, by interpolation, sharp bounds for all $1\le p\le \infty$,
completing the proof.

{\bf Observation \thmlbl{highderiv}.}
Assuming $q\in C^3$, (H0)--(H4), and ($\Cal{D}$), there holds
for any $0\le k\le 3$ the decomposition:
$$
(d/dx)^k(G-E)= \tilde G_k +  \sum_{j=0}^k H_j,
\eqnlbl{higherdecomp}
$$
where $E$ is as defined in \eqnref{E};
kernel $\tilde G_k$ is a measurable function in
$C^1(y)$, satisfying the same bounds asserted for $\tilde G$
in Lemma \thmref{2.05}; and 
$$
H_k(x,t;y)= 
\mu_k(x,t;y) (d/dy)^k
\eqnlbl{measure}
$$
where $\mu_k(x,t;\cdot)$ are measures satisfying
$$
|\mu_k(x,t;\cdot|= \Cal{O}(e^{-\eta t}).
\eqnlbl{measurebd}
$$
\medskip

{\bf Proof.} The regularity $q\in C^3$ yields $Q\in C^2$ for
the coefficients of the  linearized operator $L$, whence we
obtain differentiability of order $C^3(x)$ (resp. $C^3(y)$)
for solutions of the eigenvalue (resp. adjoint eigenvalue) ODE
from which the resolvent kernel $G_\lambda$ is constructed:
more precisely, the functions $V(x;\lambda)$ (resp. $\tilde V(y;\lambda)$)
from which modes $e^{\mu{\lambda}x}V(x;\lambda)$ 
(resp. $e^{-\mu{\lambda}x}\tilde V(x;\lambda)$) are constructed are
uniformly bounded in $C^3(x)$ (resp. $C^3(y)$).

With this observation, the claimed bounds follow in straightforward
fashion, by the same procedure used to estimate $G$, $G_x$, $G_y$
in previous sections.  Indeed, the proof is somewhat simpler, since
we only require modulus bounds here, so need not carry out the detailed
bookkeeping of previous sections.
That bounds for the regular part $\tilde G_k$ do not degrade with
increasing $k$ is clear from the elementary calculation
$$
(d/dx)e^{\mu{\lambda}x}V(x;\lambda)=
\mu e^{\mu{\lambda}x}V(x;\lambda) +
e^{\mu{\lambda}x}(d/dx)V(x;\lambda),
$$
together with the fact that, modulo negligible error terms, 
bounds on the regular part come entirely from a bounded set
of the complex plane, on which $\mu$ is bounded.
Likewise, the singular terms $H_k$ derive from large-$\lambda$
bounds by explicit computation, similarly as in previous sections.
\myqed

We first apply Observation \thmref{highderiv} in crude fashion
to obtain the crucial bound $|U|_{W^{3,1}\cap W^{3,\infty}}$
that will be used to close our iteration;
this step replaces the energy estimates of the previous subsections.
Clearly, it is sufficient to bound instead 
$|\tilde U|_{W^{3,1}\cap W^{3,\infty}}$, $\tilde U$ defined as in
\eqnref{tU} above, using
$$
U(x-\delta(t),t)-\tilde U(x,t)=
\BU(x) - \BU(x-\delta(t)) \sim \delta(t)\BU(x)
\eqnlbl{equiv}
$$
and assuming the desired order one bound on $\delta(t)$.

By \eqnref{tUeq}, Observation \thmref{highderiv}, and Duhamel's
formula, we have the representation:
$$
\aligned
\tilde U_{xxx}(x,t)
&=\int^\infty_{-\infty} \Big(
(E_{xxx}+\tilde G_3)(x,t;y)\tilde U_0(y)
+\sum_{k=0}^3 H_k(x,t;y) (d/dy)^k \tilde U_0(y)
\Big)dy\\
&+\int^t_0
\int^\infty_{-\infty}
(E_{xxx}+\tilde G_3)(x,t-s;y)
(0,I_r)^t\tilde N_2(\tilde U,\tilde U)(y,s) 
dy \, ds\\
&
+\sum_{k=0}^3 
\int^t_0
\int^\infty_{-\infty}
H_k(x,t;y) (d/dy)^k 
(0,I_r)^t\tilde N_2(\tilde U,\tilde U)(y,s) 
dy \, ds,\\
\endaligned
\eqnlbl{trep}
$$
from which we readily obtain the desired bound using the estimates
of Observation \thmref{highderiv}, together with lower derivative bounds
$$
|(d/dx)^k\tilde U |_{L^p\cap L^1}(s)\le 
\cases
C\zeta(t)(1+s)^{(k-1)/2} & 1\le k\le 2,\\
C\zeta(t) & k=0,\\
\endcases
$$
$0\le s\le t$,
which follow from the $U$ bounds of definition \eqnref{zeta3} 
by an argument similar to that given in \eqnref{equiv} above.

For example, we may bound the term
$$
\big|\int^t_0
\int^\infty_{-\infty}
H_3(x,t;y) 
\Cal{O}(|\tilde U||\tilde U_{xxx}|)(y,s) 
dy \, ds,
\big|_{L^1\cap L^\infty}
$$
arising in the expansion via \eqnref{xtUbd} of
the worst-case singular term
$$
\int^t_0
\int^\infty_{-\infty}
H_3(x,t;y) (d/dy)^3
(0,I_r)^t\tilde N_2(\tilde U,\tilde U)(y,s) 
dy \, ds,
\eqnlbl{worst}
$$
by
$$
\aligned
C\int^t_0 e^{-\eta (t-s)}
|\tilde U|_{L^\infty}|\tilde U_{xxx}|_{L^1\cap L^\infty}(y,s) 
dy \, ds
&\le
C\zeta(t)^2 
\int^t_0 e^{-\eta (t-s)} (1+s)ds\\
&\le
C\zeta(t)^2 (1+t),
\endaligned
$$
and similarly for all other terms arising in \eqnref{worst},
since by homogeneity of our derivative estimates the bound on the
nonlinear term depends only on the total number of derivatives.

Likewise, the worst-case regular term
$$
\big|
\int^t_0
\int^\infty_{-\infty}
E_{xxx}(x,t-s;y)
(0,I_r)^t\tilde N_2(\tilde U,\tilde U)(y,s) 
dy \, ds
\big|_{L^1\cap L^\infty}
\eqnlbl{worstreg}
$$
may be bounded using the triangle inequality by
$$
\aligned
C\int^t_0 |E_{xxx}|_{L^1\cap L^\infty}(t-s)
|\tilde N|_{L1}(s) ds
&\le
C\zeta(t)^2 
\int^t_0 ds\\
&\le
C\zeta(t)^2 (1+t),\\
\endaligned
$$
completing the verification of the claimed bound 
for $|U_{xxx}|_{L^p\cap L^\infty}$.

From this highest-derivative bound, the lower-derivative bounds on 
the shifted perturbation $U$ may be obtained in routine 
fashion from the corresponding integral representations for $U$,
using the bounds of Observation \thmref{highderiv} and estimates
precisely analogous to those of the previous subsections.  We
omit these lengthy but straightforward calculations.
This completes the proof of the claim, and the Theorem
\myqed

{\bf Remark \thmlbl{noopt}.} In this subsection, we have made
little attempt to minimize regularity requirements.  It seems
likely that the assumption $W^{3,1}\cap W^{3,\infty}$
on initial perturbation could be substantially reduced.
\bigskip

\specialhead
A.\quad Appendices:  Auxiliary Results/Calculations.
\endspecialhead

\specialhead
Appendix A1.\quad Structure and existence of Relaxation Profiles.
\endspecialhead
\sectionnumber=8
\theoremnumber=0
\equationnumber=0
\smallskip
\TagsOnLeft

In this appendix, we prove the results cited in the introduction
concerning structure and (small-amplitude) existence of relaxation profiles.

{\bf Proof of Lemma \thmref{structure}.}
First, note that matrix \eqnref{redcoeff} can have no nonzero purely
imaginary eigenvalues $i\xi$, by (H3).
Observing that
$$
\aligned
\det
(q_u,q_v)
\pmatrix 
f_u & f_v\\
g_u & g_v\\
\endpmatrix ^{-1}
\pmatrix 
0\\
I_r\\
\endpmatrix
&=
\det
\pmatrix 
f_u & f_v\\
q_u & q_v\\
\endpmatrix
\pmatrix 
f_u & f_v\\
g_u & g_v\\
\endpmatrix ^{-1}\\
&=
\det df^* \det q_v 
\det
\pmatrix 
f_u & f_v\\
g_u & g_v\\
\endpmatrix ^{-1}\\
\endaligned
$$ 
is nonzero, by (H1), (H2), and \eqnref{qv}, we find that $g+$ is a hyperbolic
rest point of \eqnref{gode}, from which \eqnref{expdecay} follows.
\myqed

{\bf Proof of Lemma \thmref{connection}.}
In contrast to the analogous result of [MP] in the viscous case,
which can be phrased as $i_\pm=d_\pm$, $n=r$, {\it independently}
for each state $u_\pm$, the present result must be stated
jointly, and in fact depends in an essential way on the fact that the
two states are joined by an admissible connecting profile.
Moreover, the proof emerges in a natural way from the investigation of
stability properties through the study of the linearized eigenvalue ODE.
A similar argument was sketched in Remark 2.3, [ZH], for the parabolic case.

Precisely, recall from Lemma \thmref{frozen}, Section 3,
that the $(n+r)$-dimensional eigenvalue ODE satisfies the 
{\it consistent splitting} hypothesis at spatial plus/minus infinity,
as a consequence of (H3): that is, the dimensions of the stable/unstable
subspaces of the ``frozen,'' limiting coefficient matrices at
plus and minus infinity have common values, for all frequencies
$\R \lambda>0$.
On the other hand, by Lemma \thmref{zeroexpansion}, Section 5, 
the eigenvalues of the frozen coefficient matrices
at $\lambda=0$ consist of $r$ ``fast'' modes $\R \mu_j\ne 0$,
corresponding to the eigenvalues of the reduced traveling wave
ODE linearized at $(u_\pm,v_\pm)$,
and $n$ ``slow'' modes $\mu_j \equiv 0$ corresponding to the
$n$ constants of motion of the full traveling wave ODE.
At plus infinity, the former contribute $d_+$ decaying values,
while the latter bifurcate as $\R \lambda$ increases into 
$i_+$ growing modes and $(n-i_+)$ decaying modes, accounting for
a total of $d_++(n-i_+)$ negative real part eigenvalues of the coefficient
matrix for small $\lambda$ with $\R \lambda > 0$.
Likewise, at minus infinity, we find that there are $d_-+(n-i_-)$ positive real
part eigenvalues for small $\lambda$ with $\R \lambda >0$.
Invoking consistent splitting, we have that
$d_++(n-i_+)+ d_-+(n-i_-)= n+r$, yielding the result.
\myqed

{\bf Remark.}
In the strictly parabolic case, the argument above yields instead
that $d_\pm + (n-i_\pm)=n$, or $d_\pm=i_\pm$, the result of Majda 
and Pego [MP].
\medskip

{\bf Existence of weak profiles.}
For completeness, we present also a general result on existence
of small-amplitude profiles, encompassing and slightly extending
that of [YoZ]; in particular, we do not require
the ``genuine coupling condition'' (2.5) of [YoZ],
which was induced by the particular choice of coordinates
made there.
As in [YoZ], our analysis is based on a center manifold
construction analogous to the one of [MP] in the viscous, strictly
parabolic case; however, our argument is somewhat simpler than
that of [YoZ],
and reveals more clearly the analogy to [MP].

\proclaim{Proposition \thmlbl{existence}}
Let (H0)--(H1),
\smallskip

\noindent($\tilde{ H2}$) \quad $\displaystyle{
 \sigma\left( df^*(u_\pm)\right)}$ real and distinct,
\smallskip
\noindent and (H3) hold for $(u,v)$ 
in a neighborhood of an equilibrium state
$(u_\pm,v_\pm)_0:=(u_0,v^*(u_0))$,
let $s_0:=a_p(u_0)$, where $a_p^*(u)$
denotes the $p$th eigenvalue $a_p^*(u)$
(characteristic speed) of the equilibrium flux Jacobian $A^*:=df^*(u)$,
and let $\Cal{U}$ be a sufficiently small neighborhood of $(u_0,v^*(u_0))$.
Then, for each equilibrium states 
$(u_\pm,v^*(u_\pm),s)$ in a sufficiently small neighborhood of
$(u_0,v_0,s_0)$, such that the Rankine--Hugoniot condition 
$ f^*(u_+)-f^*(u_-)= s(u_+-u_-)$
and the noncharacteristic condition $a_p^*(u_\pm)\ne s$
for the equilibrium system \eqnref{hyp} are satisfied, 
i.e., $(u_\pm,s)$ forms a noncharacteristic
shock triple for \eqnref{hyp}, there exists a 
traveling wave connection \eqnref{traveling wave}
lying in $\Cal{V}$
if and only if the triple $(u_\pm,s)$ satisfies the
Liu--Oleinik entropy condition [L.3] for the equilibrium
system \eqnref{hyp}.

Such a local connecting orbit, if it exists, is unique.
Moreover, for triples that are uniformly noncharacteristic
in the sense that $|a_p^*(u_\pm)- s|/|u_+-u_-|\ge \theta$,
$\theta>0$ a fixed constant (note: this holds for all triples
when characteristic $a_p^*$ is genuinely nonlinear at $u_0$),
the local profiles satisfy uniform bounds
$$
(\bar u- u_\pm, \bar v-v_\pm)(z) \le C|u_+-u_-|
e^{-|a_p^*(u_\pm)-s||z|/C}
\eqnlbl{expbd}
$$
for $z\ge 0$, $z\le 0$, respectively.
\endproclaim

{\bf Proof.}
For simplicity of exposition,
we shall establish the result with $(u_-, v^*(u_))$ held fixed,
without loss of generality taking $(u_-, v^*(u_),s_0)=(0,0,0)$,
and taking $z\le 0$ only in \eqnref{expbd};
The full result then follows, with $z$ still held nonpositive
in \eqnref{expbd}, by compactness, and the observation
that all of our estimates are uniform in the various parameters
of the problem; finally, we may extend \eqnref{expbd} to $z\ge 0$
by symmetry with respect to $u_-$ and $u_+$.
Alternatively,  one may as in [MP] work with the ``translated variable''
$(u,v)-(u_-,v_-)$, with $(u_-,v_-)$ carried as an additional
$(n+r)$-dimensional parameter, to obtain the full result through
a single, larger center manifold construction.

Including now the new parameter $s$ as a supplementary variable,
we may write the traveling wave ODE as
$$
\aligned
(f(u,v)-su)'&=0,\\
(g(u,v)-sv)'&= q(u,v),\\
s'&= 0.\\
\endaligned
\eqnlbl{s ode}
$$  
The map 
from $(u,v)$ to $(\tilde f,\tilde g):= (f-su,g-sv)$
(with $s$ held fixed)
has Jacobian
$$
\pmatrix 
f_u-s & f_v  \\
g_u & g_v-sv  \\
\endpmatrix 
\eqnlbl{jacob}
$$
at $(u_0,v_0)=(0,0)$ $s_0=0$, hence is locally invertible
by the noncharacteristic assumption (H1). 

Thus, we may rewrite \eqnref{s ode} in the tilde variables as
$$
\aligned
\tilde f'&=0,\\
\tilde g'&= q(u,v)(\tilde f,\tilde g),\\
s'&= 0,\\
\endaligned
\eqnlbl{tilde ode}
$$  
or, eliminating $\tilde f\equiv 0$, in the convenient form
$$
\pmatrix 
\tilde g\\
s \\
\endpmatrix'
=
\pmatrix 
q(u,v)(0,\tilde g)\\
0 \\
\endpmatrix.
\eqnlbl{g ode}
$$

Linearizing \eqnref{g ode} about base point 
$(\tilde g_0,s_0):=(g(0,0), 0)$, and using \eqnref{jacob},
we obtain after a brief calculation the linearized equations
$$
\pmatrix 
\tilde g\\
s \\
\endpmatrix'
=
\pmatrix 
M
& 0\\
0 & 0\\
\endpmatrix
\pmatrix 
\tilde g\\
s \\
\endpmatrix,
\eqnlbl{lingode}
$$
where
$$
M:=
\pmatrix 
q_u & q_v\\
\endpmatrix 
\pmatrix 
f_u & f_v\\
g_u & g_v\\
\endpmatrix ^{-1}
\pmatrix 
0_n \\ I_r\\
\endpmatrix
_{|_{u,v=0,0}}.
\eqnlbl{M}
$$
Evidently, the center subspace of the coefficient matrix of
\eqnref{lingode} is the direct sum of $(0_r,1)$ and $(C,0)$,
where $C$ is the center subspace of $M$.

{\it Claim.} The center subspace $C$ of $M$ is one-dimensional, 
consisting of its kernel.
More precisely,
it is spanned by ${g^*_u r^*_p}(0)$, where $r^*_p$ is the right 
eigenvector of the equilibrium flux $A^*$ associated
with the principal characteristic speed $a_p^*$,
and $g_u^*$ analogously to $f^*_u$ is defined
as $g_u - q_v^{-1}q_u g_v$, the variation of $g$ along the equilibrium
manifold: equivalently,
$$
\pmatrix 
f_u & f_v \\
g_u & g_v \\
\endpmatrix ^{-1}
\pmatrix
0\\
g_u^*r^*_p\\
\endpmatrix
=
\pmatrix r^*_p \\
-q_v^{-1}q_u r^*_p. \\
\endpmatrix
\eqnlbl{uvcoord}
$$
Likewise, the left kernel of $M$ is spanned by $-l_p^* f_v q_v^{-1}(0)$,
where $l_p^*$ is the left zero eigenvector of $f^*_u$ dual to $r_p^*$.

\medskip
{\it Proof of Claim.}
Clearly, the eigenvalues of $M$ are contained among the
eigenvalues of
$$
M_1:=
\pmatrix 
0& 0\\
q_u & q_v\\
\endpmatrix 
\pmatrix 
f_u & f_v\\
g_u & g_v\\
\endpmatrix ^{-1},
\eqnlbl{M1}
$$
which has no nonzero pure imaginary eigenvalues, by (H3).
Thus, $C$ consists entirely of the zero eigenspace of $M$.

Next, supposing that $x$ lies in the kernel of $M$, we
find that 
$$
\pmatrix 
f_u & f_v \\
g_u & g_v \\
\endpmatrix ^{-1}
\pmatrix
0\\
x\\
\endpmatrix
$$
lies in the kernel of $(q_u,q_v)$, hence must be
of form $(r, -q_v^{-1}q_u r)^t$.
From the first coordinate of the equation
$$
\pmatrix 
f_u & f_v \\
g_u & g_v \\
\endpmatrix 
\pmatrix
r\\
-q_v^{-1}q_u r\\
\endpmatrix
=\pmatrix 0\\x\endpmatrix
$$
we find that $r$ must lie in the one-dimensional kernel of $f^*_u(0)$, whence
we obtain \eqnref{uvcoord} and the claimed description of the kernel 
of $M$.

A direct computation then shows that
$$
\aligned
-l_p^*f_v q_v^{-1}M&=
-l_p^*
\pmatrix f_v q_v^{-1}q_u & f_v\\
\endpmatrix
\pmatrix
f_u  & f_v\\
g_u & g_v\\
\endpmatrix ^{-1}
\pmatrix
0_n \\
I_r\\
\endpmatrix
\\
&=
l_p^*
\pmatrix f^*_u & 0 \\
\endpmatrix
\pmatrix
f_u  & f_v\\
g_u & g_v\\
\endpmatrix ^{-1}
\pmatrix
0_n \\
I_r\\
\endpmatrix
 =0,
\\
\endaligned
$$
i.e., $\tilde x:=-l_p^*f_v q_v^{-1}$ satisfies $\tilde xM=0$,
and more generally,
$$
\aligned
-l_p^*f_v q_v^{-1}
\pmatrix q_u & q_v\\
\endpmatrix
\pmatrix
f_u  & f_v\\
g_u & g_v\\
\endpmatrix ^{-1}
&=
-l_p^*
\pmatrix f_v q_v^{-1}q_u & f_v\\
\endpmatrix
\pmatrix
f_u  & f_v\\
g_u & g_v\\
\endpmatrix ^{-1}\\
&=
-l_p^*
\Big(
\pmatrix -f^*_u  & 0\\
\endpmatrix
+ \pmatrix f_u & f_v\\
\endpmatrix
\Big)
\pmatrix
f_u  & f_v\\
g_u & g_v\\
\endpmatrix ^{-1}
\\
&=
-l_p^* \pmatrix I_n & 0_r \endpmatrix,
\endaligned
\eqnlbl{useful}
$$
from which we may conclude that $\tilde x\ne 0$.
This verifies the asserted description of the left kernel of $M$.

Observing that
$$
\tilde x x =
-l_p^*f_v q_v^{-1}g_u^*r_p^*
=
l_p^*B^*r^*_p,
$$
where $B^*$ is the Chapman--Enskog viscosity defined in \eqnref{B}
(Note: in the final equality, we have used the fact that $f^*_u r_p^*=0$),
we thus find that a necessary and sufficient condition for existence of a 
generalized eigenvector $y$, $My=x$, 
 is that 
$$
\tilde x x=
l_p^*B^*r^*_p = 0.
\eqnlbl{visc cond}
$$
But, the low-frequency
Fourier expansion in the following appendix (A2) shows that the stability
condition (H3) implies $l_p^* B^* r_p^*>0$, and so no such generalized 
eigenvector occurs.  
\myqed


{\bf Remark.}  The conclusions of the Claim may be obtained in a
lengthier but more transparent fashion by the observation that
zero eigenvalues/generalized eigenvalues of $M$ correspond
with those of
$$
\pmatrix 
f_u & f_v \\
q_u & q_v \\
\endpmatrix 
\pmatrix 
f_u & f_v \\
g_u & g_v \\
\endpmatrix ^{-1}.
$$
\medskip
\myqed

With these preparations, the result now follows similarly as in [MP].
By standard center manifold construction, there exists a two-dimensional
manifold tangent to the center subspace of $M$, which contains all
locally bounded orbits of \eqnref{g ode}: in particular, all rest
points $(\tilde g_\pm,s)$.
For each fixed $s$, this reduces to a one-dimensional fiber
$$
\tilde g=\tilde g(\theta),
\eqnlbl{param}
$$
$\theta \in \BbbR$, containing
both of the rest points $\tilde g_\pm$.
Thus, there is a connection between them if and only if there
is no rest point of \eqnref{g ode} lying between them.

By the Lax structure theorem [La], locally,
all rest points lie on the $p$th Hugoniot curve through $u_-=0$
of the equilibrium system, mapped through $g(u,v^*(u))- s v^*(u)$
to the $\tilde g$ coordinates, and both Hugoniot curve and fiber
are smooth curves lying approximately in the $g^*_u r^*_p(0)$
direction. Thus, the ordering of the rest points must be the 
same along the center manifold fiber as along the Hugoniot curve.
But, the analysis of [L.3] shows that two noncharacteristic
rest points are consecutive along the Hugoniot curve if and only
if the Liu--Oleinik admissibility condition or its reverse is
satisfied.
Thus, it remains only to check that the sense of a connection
must agree with that of the (forward) Liu--Oleinik condition,
or, equivalently, that there hold the one-sided 
{\it Lax characteristic condition} $a^*_p(u_-)> s$; for details, see [L.3].

To see that the one-sided Lax condition holds, we may use invariance
of \eqnref{param} under the flow of \eqnref{g ode} to obtain
$$
h(\theta) \tilde g_\theta= q(u,v)(0,\tilde g),
\eqnlbl{pre}
$$
where
$$
\theta'=h(\theta)
$$
describes the reduced flow along the center manifold.
Differentiating \eqnref{pre} with respect to $\theta$ at $\theta=0$,
we obtain
$$
h_\theta (0)\tilde g_\theta(0)= M(s) \tilde g_\theta(0),
$$
where
$$
\aligned
M(s)&:=
\pmatrix
q_u & q_v
\endpmatrix
\pmatrix
f_u-su & f_v\\
g_u & g_v -sv\\
\endpmatrix^{-1}_{|_{u,v=0,0}}
\pmatrix
0_n \\
I_r \\
\endpmatrix
\\
&= 
M +
sM_2
+ \Cal{O}(s^2)
\\
\endaligned
\eqnlbl{Mpert}
$$
is a one-parameter perturbation of $M$ above, with
$$
M_2:=
 \pmatrix
q_u & q_v
\endpmatrix
\pmatrix
f_u & f_v\\
g_u & g_v \\
\endpmatrix^{-2}_{|_{u,v=0,0}}
\pmatrix
0_n \\
I_r \\
\endpmatrix.
\eqnlbl{M2}
$$
That is, $h_\theta(0)$, $\tilde g_\theta(0)$
are a right eigenvalue, eigenvector pair for $M(s)$,
bifurcating from the values $0$, $g^*_u r_p^*(0)$ for
$M$ at $s=0$.
Since we have already ascertained that $0$ is an isolated
eigenvalue of $M$, with left and right eigenvectors 
$l_p^* f_v q_v^{-1}(0)$ and
$g_u^*r_p^*(0)$, respectively, we may conclude by standard
matrix perturbation theory [Kat] that $h_\theta(0)$ varies
with respect to $s$ as
$$
\aligned
h_\theta(s)&= s 
\frac{l_p^* f_v q_v^{-1}(0) M_2 g_u^*r_p^*(0)}
{l_p^* f_v q_v^{-1}(0)  g_u^*r_p^*(0)}
 + \Cal{O}(s^2)\\
&= s \frac{l_p^* r_p^*(0)}
{l_p^* B^* r_p^*(0)}
 + \Cal{O}(s^2)\\
&=  \frac{s}
{l_p^* B^* r_p^*(0)}
 + \Cal{O}(s^2),\\
\endaligned
\eqnlbl{hpert}
$$
where in the second equality
we have used \eqnref{useful} and \eqnref{uvcoord} in simplifying
the numerator of the first-order coefficient.
With \eqnref{visc cond}, this gives 
$$
\sgn h_\theta(0) = \sgn s= \sgn (s-a_p^*(u_-))
$$ 
for $s$ sufficiently
small, verifying that $\theta=0$ is a repelling point of the reduced
flow if and only if the one-sided Lax condition holds.
At the same time, it gives the bound \eqnref{expbd} for $z\le 0$
by straightforward estimates on the scalar system \eqnref{reduced}.
This completes the proof of the Proposition.
\myqed

\bigskip

\specialhead
Appendix A2. Expansion of the Fourier Symbol.
\endspecialhead

In this appendix, we carry out the Taylor expansions
about $\xi=0$ and $\xi=\infty$ of the Fourier
symbol $-i\xi A(x) + Q(x)$ of the frozen, constant-coefficient
operator at a given $x$.

\medskip
{\bf High frequency expansion.} By straightforward matrix perturbation
theory [Kat], the first-order expansion at infinity 
of $P(i\xi)$ with respect to $(i\xi)^{-1}$ is just
$$
P=-i\xi A + R \diag \{\eta_1,\dots, \eta_{J}\}L,
$$
$$
\diag \{\eta_1,\dots, \eta_{J}\}:= LQR,
$$
where $L$, $R$ are composed of left, right eigenvectors of $A$,
as described in the introduction.
\medskip

{\bf Low frequency expansion.}  We next carry out the expansion of $P$
in $\xi$ about zero.  Recall that
$$
Q:=
\pmatrix
0&0\\q_u &q_v
\endpmatrix,\quad 
A:= \pmatrix
f_u &f_v\\ g_u &g_v
\endpmatrix .
\eqnlbl{6}
$$

{\it Claim.} To second order, ''slow'' dispersion relations
$$
\lambda(i\xi) =\sigma(Q-i\xi A), \quad \lambda(0)=0
\eqnlbl{7}
$$
are given by
$$
\lambda_j(\xi )= -i\xi a^*_j - \beta^*_j\xi^2+ \dots,
\quad j=1,\cdots n,
\eqnlbl{8}
$$
with corresponding eigenvectors
$$
R_j^*(i\xi)=
\pmatrix
r^*_j & \\
-q_v^{-1} q_u r^*_j
\endpmatrix
+ \dots,
$$
and 
$$  
\beta^*_j:= l^*_j {B^*} r^*_j,
\eqnlbl{9}
$$
$$
{B^*}:= f_v\Big[g_u-g_vq^{-1}_v q_u+
f^*_u q^{-1}_vq_u\Big],
\eqnlbl{10}
$$
where $a^*_j$, $r_j^*$, $l_j^*$ are eigenvalues and associated
right and left eigenvectors of $A^*:=f^*_u$.
The remaining $r$ ``fast'', or strictly stable roots are $\sim \sigma(q_v)$.

{\bf Remark.}
This expansion may be viewed as a rigorous justification of
the Chapman--Enskog expansion in the constant-coefficient case.
\medskip

{\bf Proof.}
Setting 
$$
(Q-i\xi A-(-i\xi \alpha_j-\beta_j \xi^2)\cdots)
(V_0 +V_1 \xi+V_2 \xi^2+\cdots)=0,
$$
and collecting terms of successive orders in $\xi$, we obtain:

(Order $1$): 
$$
QV_0=0 \qquad 
\Rightarrow V_0=
\pmatrix
r_j & \\
-q_v^{-1} q_u r_j
\endpmatrix
\eqnlbl{11}
$$

\medskip

(Order $\xi$): 
$$
QV_1-i(A-\alpha_j)V_0=0,
\eqnlbl{12}
$$
$\Rightarrow$ (looking at first coordinate)
$$
(A_*-\alpha_j)r_j=0,\qquad A^*:= f^*_u=f_u-f_v q^{-1}_v q_u 
\eqnlbl{13}
$$
$\Rightarrow$
$\alpha_j$, $r_j$ are eigenvalue $a_j^*$, right eigenvector $r_j^*$ of $A^*$.

Further, comparing with the second coordinate in \eqnref{12}, we have
$$
{1 \over \tau} (q_u,q_v)V_1= i(g_u,g_v -a^*_j)V_0,
\eqnlbl{15}
$$
giving 
$$
V_1=
\pmatrix
s & \\
q^{-1}_v(-q_u s+\tau i\Big[g_u-(g_v-a^*_j)q^{-1}_vq_u\Big]
r_j)
\endpmatrix.
\eqnlbl{16}
$$
\medskip

(Order $\xi^2$):
$$
QV_2 -i(A-\alpha_j)V_1+\beta_jV_0=0.
$$
Examining the first coordinate, we have
$$
\aligned
\beta_jr_j
&= i(f_u-a^*_j,f_v)V_1\\
&= i(A^*-a^*_j)s-\tau f_vq^{-1}_v[g_u -
(g_v-a^*_j)q^{-1}_v q_u]r^*_j.
\endaligned
\eqnlbl{17}
$$
Taking inner product with $l^*_j$, left e--vector of $A_*$, we obtain
$$
\aligned
\beta^*_j
&=-\tau l^*_j f_v q^{-1}_v[g_u-(g_v-a^*_j)q^{-1}_v q_u]r^*_j\\
&=-\tau l^*_j f_v q^{-1}_v[g_u-g_vq^{-1}_vq_uf^*_u]r^*_j,
\endaligned
\eqnlbl{18}
$$
completing the result.
\myqed

{\bf Remark.}  The ``strict asymptotic parabolicity'' condition 
$$
\R\beta_j^*>0, \quad j=1,\dots, n,
\eqnlbl{apar}
$$
reduces in the case $n=r=1$ treated by Liu [L.2] 
to the classical subcharacteristic condition [Wh];
however, these are not equivalent.
As discussed in [Yo.4], the relation between the subcharacteristic
condition and stability is quite analogous
to that between the CFL condition and stability in numerical analysis.
In the symmetrizable case,
condition \eqnref{apar} is equivalent (see [Kaw,Ze.2]) to
$Q\le 0$ and no vector in $ker(Q)$ is an eigenvector of $A$.

\specialhead
Appendix A3. The Gap and Tracking Lemmas.
\endspecialhead

For completeness, we present here two technical lemmas
referred to throughout the text,  relating behavior of
variable- and constant-coefficient ODE in, respectively, 
the {\it asymptotically constant}
and the {\it slowly varying} coefficient case.
The resulting bounds are crucial in our study of solutions
of the generalized eigenvalue equations,
in the former case on bounded domains in frequency space 
$\lambda$, and in the latter as $|\lambda|\to \infty$.

\specialhead
A3.1.\quad The Gap Lemma.
\endspecialhead
A common problem arising in 
the asymptotic study of eigenvalue ODE
is to relate
behavior near $x=\pm \infty$ of solutions of an asymptotically 
constant-coefficient eigenvalue equation to that of solutions of the 
corresponding limiting, constant-coefficient equations, in a manner that 
is smooth with respect to spectral parameters.
More generally, consider a general ODE with parameter
$$
W'=\BbbA(\lambda, x)W,
\eqnlbl{3.1}
$$
``$\prime$'' denoting $(d/dx)$,
where the coefficient $\BbbA$ is $C^0$ in the evolution variable $x$ and
analytic (resp. $C^r$) in parameter $\lambda$, and
converges as $x\to \pm \infty$ to limiting values $\BbbA_\pm$.
It is well known (see [Co], Thm. 4, p. 94) that, provided that
$$
\int^{\pm\infty}_0 | \BbbA - \BbbA_\pm | dx < +\infty,
\eqnlbl{3.2}
$$
there is a one-to-one correspondence between the normal modes $V^{\pm}_j 
e^{\mu_j^\pm x}$ of the asymptotic systems
$$
W'=\BbbA_\pm(\lambda )W,
\eqnlbl{3.3}
$$
where $V^\pm_j$, $\mu^\pm_j$ are eigenvector and 
eigenvalue of $\BbbA_-$ (alternatively, 
$V^\pm_j x^{\ell} e^{\mu_j^\pm x}$, if 
$V^\pm_j$ is a generalized eigenvector of order $\ell$) and certain 
solutions $W^\pm_j$ of \eqnref{3.1} having the same asymptotic behavior,
i.e., 
$$
W^\pm_j (\lambda, x) = V^\pm_j e^{\mu^\pm_j x} (1 + o(1))\  
{\text as}\  x \to \pm \infty
\eqnlbl{3.4}
$$
(alternatively, $W^\pm_j (\lambda, x) = V^\pm_j x^\ell
e^{\mu^\pm_j x}(1 + o(1))$. That is, the flows near $\pm \infty$ of 
\eqnref{3.1} and \eqnref{3.3} are homeomorphic.

Such a correspondence is of course highly nonunique, since \eqnref{3.4} determines 
$W^\pm_j$ only up to faster decaying modes. However, provided that $\R(\mu^\pm_j)$ 
is {\it strictly separated} from all other $\R(\mu^\pm_k)$, i.e.,
that there is a {\it 
spectral gap}, the choice defined in [Co], Theorem 4 by fixed point iteration 
is in fact analytic (resp. $C^r$)
in $\lambda$, as the uniform limit of an analytic (resp. $C^r$) sequence of 
iterates. The argument breaks down at points $\lambda_0$ 
where $\R(\mu_j) = \R(\mu_k)$ for some $k \not= j$, since in this 
case $(\R(\mu_j) - \R(\mu_k))$ does not have a 
definite sign, and the definition of the fixed point iteration is 
determined by the signs of all $(\R(\mu_j) - \R(\mu_k))$.

The purpose of the present section is to point out that analyticity (resp.
smoothness) in $\lambda$ can be recovered 
in the absence of a spectral gap, by virtually the same argument as in [Co] if 
we substitute for \eqnref{3.2} the stronger hypothesis:
$$
|\BbbA - \BbbA_\pm| = \bfO (e^{-\alpha|x|})\ {\text as\ } x \to
\pm\infty. 
\eqnlbl{3.5}
$$

This observation is a special case of the ``Gap Lemma of [GZ], also
proved independently in [KS].  The original version was phrased in
terms of the projectivized flow associated with \eqnref{3.1}.  Here,
we give an alternative statement and derivation directly in terms of
\eqnref{3.1}, a form more convenient for our needs.  

\proclaim{Proposition \thmlbl{3.1}} In \eqnref{3.1}, let $\BbbA$ be 
$C^{0}$ in $x$ and analytic (resp. $C^r$) in $\lambda$, with 
$|\BbbA-\BbbA_- (\lambda)| = 
\bfO(e^{-\alpha|x|})$ as $x \to - \infty$ for $\alpha>0$, and
$\balpha < \alpha $.  
If $V^-(\lambda)$ is an 
eigenvector of $\BbbA_-$ with eigenvalue $\mu(\lambda)$, both 
analytic (resp. $C^r$) in $\lambda$, 
then there exists a solution $W(\lambda, x)$ of $\eqnref{3.1}$ of
form 
$$
W(\lambda, x) = V (x;\lambda ) e^{\mu(\lambda) x},
$$
where $V$ (hence $W$) is $C^{1}$ in $x$ and 
{\it locally analytic (resp. $C^r$)}  in $\lambda$ and 
for each $j=0,1,\ldots$ satisfies 
$$
(\partial/\partial \lambda)^j V(x;\lambda )=(\partial/\partial
\lambda )^j V^-(\lambda ) +
\bfO (e^{-\balpha|x|}|V^- (\lambda)|),\quad x < 0,
\eqnlbl{3.6}
$$
Moreover, if $\R\mu(\lambda)> \R\Tmu (\lambda) - \alpha $ 
for all (other) eigenvalues $\Tmu$ of $\BbbA_-$, then $W$ is uniquely 
determined by \eqnref{3.6}, and \eqnref{3.6} holds for $\balpha 
= \alpha $.
\endproclaim

{\bf Proof.} Setting $W(x)=e^{\mu x}V(x)$, we can rewrite $W'=\BbbA
W$ as 
$$
\aligned
V'
&= (\BbbA_- - \mu I)V+\theta V,\\
\theta
&:= (\BbbA - \BbbA_-)=\bfO(e^{-\alpha|x|}),
\endaligned \eqnlbl{3.7}
$$
and seek a solution $V(x;\lambda )\to V^-(x)$ as $x \to \infty$.

Set 
$$
\balpha =\balpha_1< \alpha _1 < \alpha _2  <\balpha_2< \alpha .  
\eqnlbl{alphas}
$$
Fixing a base point $\lambda_0$, we can define on some neighborhood
of $\lambda_0$ to the complementary $\BbbA_-$-invariant projections
$P(\lambda)$ and $Q(\lambda)$ where $P$ projects onto the direct sum
of all eigenspaces of $\BbbA_-$ with eigenvalues $\Tmu$ 
satisfying 
$$
\R(\Tmu) < \R(\mu) + \alpha_2,
\eqnlbl{3.9}
$$
and $Q$ projects onto the direct sum of the remaining eigenspaces,
with eigenvalues satisfying 
$$
\R(\Tmu) \ge \R(\mu) + \alpha _2 > \R(\mu) 
+\alpha _1.
\eqnlbl{3.9}
$$ 
By basic matrix perturbation theory (eg. [Kat]) it follows that 
$P$ and $Q$ are analytic (resp. $C^r$) in a neighborhood of $\lambda_0$,  
with
$$
\aligned
\left|e^{(\BbbA_- - \mu I)x} P \right| &= \bfO (e^{\balpha_2 x}),
\quad x>0, \\
\left|e^{(\BbbA_- - \mu I)x} Q \right| &= \bfO (e^{\balpha_1 x}), 
\quad x<0.
\endaligned \eqnlbl{3.10}
$$

Thus, for $M>0$ sufficiently large, the map $\CalT$ defined by 
$$
\CalT V(x) 
= V^- + \int^x_{-\infty} e^{(\BbbA_- - \mu I)(x-y)} P
\theta (y) V(y) dy 
- \int^{-M}_x e^{(\BbbA_- - \mu I)(x-y)} Q \theta (y) V(y) dy
\eqnlbl{3.11}
$$
is a contraction on $L^\infty(-\infty, -M]$;
for, applying \eqnref{3.10}, we have
$$
\aligned
\left|\CalT V_1 - \CalT V_2 \right|_{(x)} 
&\le \bfO(1) |V_1 - V_2|_\infty 
\bigg(\int^x_{-\infty} e^{\balpha_2(x-y)} e^{\alpha y} dy 
+ \int^{-M}_{x} e^{\balpha _1(x-y)}e^{\alpha y}dy\bigg)\\
&= \bfO(1) |V_1 - V_2|_\infty 
\bigg(e^{\balpha_2 x}e^{(\alpha-\balpha_2)y}|^x_{-\infty}+
e^{\balpha_1 x}e^{(\alpha-\balpha_1)y}|^{-M}_{x}\bigg) \\
&= \bfO(1) |V_1 - V_2|_\infty (e^{\alpha x} + e^{\balpha x})
\\
&= \bfO(1) |V_1 - V_2|_\infty e^{-\balpha M} <
\frac{1}{2}.
\endaligned \eqnlbl{3.12}
$$

By iteration, we thus obtain a solution $V \in L^\infty (-\infty, -M]$ of $V = 
\CalT V$ with $V=\bfO(|V^-|)$; since $\CalT$ clearly preserves analyticity
(resp. smoothness), $V(\lambda, x)$  is 
analytic (resp. $C^r$) in $\lambda$ as the uniform limit of analytic (resp.
smooth) iterates (starting with $V_0=0$). 
Differentiation shows that $V$ is a bounded solution of
$V=\CalT V$ iff it is a bounded solution of \eqnref{3.7}. 
Further, taking $V_1=V$, $V_2=0$ in \eqnref{3.12}, we obtain from
the second to last equality that 
$$
|V-V^-| = |\CalT(V) - \CalT(0)| = \bfO(1) e^{\balpha x} 
|V| = \bfO(e^{\balpha x})|V^-|,
\eqnlbl{3.13}
$$
giving \eqnref{3.6} for $j=0$.  Derivative bounds, $j>0$, follow by
standard interior estimates, or, alternatively, by differentiating
\eqnref{3.11} with respect to $\lambda $ and repeating the same
argument.  Analyticity (resp. smoothness), and the bounds \eqnref{3.6},  
extend to $x<0$ by standard analytic dependence for the initial value
problem at $x=-M$. 

Finally, if $\R(\mu(\lambda))> \R(\Tmu(\lambda))- 
\frac{\alpha}{2}$ for all 
other eigenvalues, then $P=I$, $Q=0$, and $V=\CalT V$ 
must hold for any $V$ satisfying \eqnref{3.6}, 
by Duhamel's principle. Further, the only term appearing in
\eqnref{3.12} is the first integral, giving bound \eqnref{3.13} for
$\balpha = \alpha $. \myqed

Proposition \thmref{3.1} extends also to subspaces of solutions.  This can be
seen most easily by associating to a $k$-plane of solutions, Span
$\{W_1(x),\ldots , W_k(x)\}$, the corresponding $k$-form $\eta = W_1
\wedge \cdots \wedge W_k$.  The equations \eqnref{3.1} induce a
linear flow
$$
\eta ' =\BbbA^{(k)}(x;\lambda )\eta,
\eqnlbl{3.14}
$$
on the space of $k$-forms via the Leibnitz rule,
$$
\BbbA^{(k)}(W_1\wedge \cdots \wedge W_k)=
(\BbbA W_1 \wedge \cdots \wedge W_k)
+\cdots + (W_1 \wedge \cdots \wedge \BbbA W_k).
\eqnlbl{3.15}
$$
The evolution of the $k$-plane of solutions of \eqnref{3.1} is
clearly determined by that of $\eta (x;\lambda )$.  It is easily seen
that for a given (constant) matrix $\BbbA$, the eigenvectors of
$\BbbA^{(k)}$ are of form $V_1 \wedge \cdots \wedge V_k$, where Span
$\{V_1,\cdots , V_k\}$ is an invariant subspace of $\BbbA$, and that
the corresponding eigenvalue is the trace of $\BbbA$ on that
subspace. 

\proclaim{Definition \thmlbl{3.2}}  Let $C=\Span\ \{V_{k+1},\ldots , V_N\}$
and $E=\Span\ \{V_1,\ldots , V_k\}$ be complementary
$\BbbA$-invariant subspaces.  We define their {\bf spectral gap} to
be the difference $\beta $ between the real part of the eigenvalue of
minimal real part of $\BbbA$ restricted to $C$ and the real part of
the eigenvalue of maximal real part of $\BbbA$ restricted to $E$.
\endproclaim

If $\eta$ is a $k$-form associated with an $A$-invariant subspace
$E$ as in the definition above, then the spectral gap $\beta $ is the
minimum difference between the real part of the eigenvalue $\mu$ of
$\BbbA^{(k)}$ associated with $\eta$ and the real part of the
eigenvalue associated with any other eigenvector of $\BbbA^{(k)}$.
Combining these observations with the result of Proposition \thmref{3.1}, we
obtain a complete version of the Gap Lemma of [GZ]:

\proclaim{Corollary \thmlbl{3.3} (The Gap Lemma)}  
Let $\BbbA(x;\lambda )$ be $C^{0}$ 
in $x$, analytic (resp. $C^r$) in $\lambda $, with $\BbbA(x;\lambda )\to
\BbbA_\pm(\lambda )$ as $x \to \pm \infty$ at exponential rate
$e^{-\alpha |x|}$, $\alpha >0$, and let $\eta^-(\lambda )$ and $\zeta^-$
be analytic (resp. $C^r$) $k$ and $n-k$-forms associated to complementary
$\BbbA_-(\lambda )$-invariant subspaces $C^-$ and $E^-$ as in
Definition 3.2, with arbitrary spectral gap $\beta $, and let
$\tau_{C^-}(\lambda)$ be analytic, where 
$\tau_{C^-}(\lambda)$ is
 the trace of $\BbbA_-(\lambda)$ restricted to $C^-$.
Then, there exists a solution $\CalW(\lambda ,x)$ of \eqnref{3.14} of
form
$$
\CalW (\lambda ,x)=\eta(\lambda ,x)e^{\tau_{C^-}x}
$$
where $\eta$ (hence $\CalW$) is $C^{1}$ in $x$ locally
analytic (resp. $C^r$) in $\lambda $, and for each $j=0,1,\ldots$ satisfies
$$
(\partial/\partial \lambda)^j \eta(x;\lambda )=(\partial/\partial
\lambda )^j \eta^-(\lambda ) +
\bfO (e^{-\balpha|x|}|\eta^- (\lambda)|),\quad x < 0,
\eqnlbl{3.16}
$$
for all $\balpha < \alpha $.  Moreover, if $\beta >-\alpha $, then
$\eta$ is uniquely determined by \eqnref{3.16} and \eqnref{3.16}
holds for $\balpha = \alpha $.
\endproclaim

{\bf Remark \thmlbl{achoice}.}
Note that the construction of Corollary \thmref{3.3} determines
at the same time analytic subspaces $C(x;\lambda)$, $E(x;\lambda)$
asymptotic to $C^-$, $E^-$,
and thereby corresponding complementary analytic projections
$P_C(x;\lambda)$, $P_E(x;\lambda)$.  Fixing $x_0$ and $\lambda_0$,
and applying projections $P_C(x_0,\lambda)$, $P_E(x_0,\lambda)$
to fixed bases of $C(x_0,\lambda_0)$, $E(x_0,\lambda_0)$, we
may thus obtain a {\it locally analytic choice} of individual
solutions $W_j^-(x;\lambda)$ spanning 
$C(x;\lambda)$, $E(x;\lambda)$, in a neighborhood of any $\lambda_0$.  
Moreover, if the eigenvalues associated with $C^-$ have 
real parts lying strictly between $\eta_1$ and $\eta_2$, 
then the standard homeomorphism
between solutions of limiting, constant coefficient equations and
their variable-coefficient counterparts (described, e.g. in Section 3)
yields that the flow on $C(\lambda,\cdot)$ decays/grows exponentially 
with $|x|$ on $x\le 0$, with rate uniformly bounded between $\eta_1$ and $\eta_2$.
Alternatively, such an analytic choice may
be obtained in more straightforward fashion by a standard
''stable-manifold type'' construction, as described in [Z.1],
from which exponential decay/growth may be deduced at the same time.
\medskip

{\bf Remark \thmlbl{noinfo}.}
Bounds \eqnref{3.6}, \eqnref{3.16} assert no information about the 
direction of $V(\cdot;\lambda)$, $\eta(\cdot;\lambda)$ at $x=0$.  
A review of our fixed-point
construction reveals that in general no such information is available, 
i.e., $\Cal{T}$ is a contraction only on $L^\infty(-\infty, -M]$, for $M>0$
sufficiently large, and not on $L^\infty(-\infty,0]$. Specifically,
the second integral in definition \eqnref{3.11} becomes for $x>-M$
an expansive term with expansion coefficient growing exponentially
in $x+M$ with rate at least $\balpha_1$ and possibly arbitrarily
larger; at $x=0$, this becomes at least order $e^{\balpha_1 M}$,
and typically greater than one.

In particular, we obtain from the bounds of the Gap Lemma no lower
bound on the Evans function $|D(\lambda)|$, i.e., no information
on existence or nonexistence of point spectrum at $\lambda$.
Indeed, this could not be so, since we know a priori in many
cases that $D(0)=0$ even though the Gap Lemma estimates are uniformly
valid up to $\lambda=0$.  To obtain lower bounds on the Evans function,
for example in the high-frequency limit,
we will use instead an alternative construction described in the 
following subsection, based on the {\it projectivized flow} of
the eigenvalue ODE.  By factoring out variations in modulus, this
allows more accurate estimates of direction; in particular,
the mapping analogous to $\Cal{T}$ in projectivized coordinates
((8.54) below) will be seen to be a {\it global} $L^\infty$-contraction 
on the whole real line.
\footnote{
In this regard, we point out a key error in the proof of Lemma 5
of [R], specifically of lower bound (20) asserted on the Evans function
(bounds (18)--(19) clearly extend from $x_0$ to $x=0$ by Abel's formula).
Namely, the claimed bounds on operator $N$ of (34), following directly
from contractivity of the mapping $T_\lambda$ (our $\Cal{T}$, above)
are valid only on the interval $[x_0,\infty)$ where $T_\lambda$ is
a contraction, and not all of $[0,+\infty)$ as claimed.
The error results from an incorrect citation of the Gap Lemma construction
of [GZ,ZH],
stating that mapping $T_\lambda$ itself is a contraction on $[0,\infty)$).
As discussed in Remark \thmref{noinfo} above, such information cannot
be obtained by the Gap Lemma construction, but rather requires the 
projectivized version described below
(indeed, the required lower bound in [R] is a trivial consequence of
the machinery developed here, in Section 4 and Appendix A3.2 below).
}
\medskip

\specialhead
A3.2.\quad The Tracking Lemma.
\endspecialhead

Another general situation that arises in the asymptotic study of 
eigenvalue equations is an ODE
$$
w'=(\BbbA(x,\delta )+\delta \Theta(x,\delta ))w,\quad w \in \BbbC^N,
\eqnlbl{32}
$$
with a small parameter $\delta \to 0$, satisfying
$$
|\BbbA |+ |\Theta|\le C
\eqnlbl{33}
$$
and
$$
|\BbbA'| \le C \delta,
\eqnlbl{A'}
$$
where ``$^\prime$'' denotes differentiation with respect to
$x$: that is, an ODE with slowly varying coefficients.
This situation arises in the limit as the frequency
rather than the spatial variable goes to infinity,
after rescaling to a length scale on which the resulting
rapid oscillations in the solution have period of order one
with $\delta \to 0$ as $|\hbox{frequency}|\to \infty$.
Thus, it is dual to the previous case.

\medskip
{\bf A3.2.1. The basic estimate.}
To introduce ideas, we first present an analysis of the ``standard''
case treated in [GZ,ZH,Z.1,Z.4]. 
Suppose in addition to \eqnref{33}--\eqnref{A'} that 
the spectrum of the matrix $\BbbA(x, \delta)$
divides into two spectral groups, 
$$
\alpha _1(x,\delta), \dots \alpha _\ell(x,\delta) \le \underline{\alpha}(x) < 
\balpha(x) \le \alpha _{\ell+1}(x,\delta), \dots \alpha _N(x,\delta),
\eqnlbl{34}
$$
where $<$ denotes ordering with respect to {real parts},
with a {\it uniform spectral gap}
$$
0< 2\eta \le
\balpha(x) - \underline{\alpha}(x);
\eqnlbl{etagap}
$$
in the setting of asymptotic eigenvalue equations,
this is achieved by restricting to an appropriately small subset of
the region of consistent splitting,
or normal set, of the operator $L$ under investigation.
By standard matrix perturbation theory [Kat], therefore, we have 
$\BbbA(\cdot)$-invariant projections 
$P(x)$ and $Q(x)$ onto the eigenspaces associated 
with $(\alpha_1,\dots,\alpha_\ell)$ and $(\alpha_{\ell+1},\dots \alpha_v)$, 
respectively, satisfying
$$
|P'|+|Q'|\le C_2 \delta,\ |P|, |Q| \le C;
\eqnlbl{35}
$$
we assume that these bounds hold uniformly over all $-\infty <x<+\infty$
(this can be guaranteed, for example, by the assumption that $\BbbA$ varies
within a compact set, as is the case in the applications we have in mind).

Under these general assumptions, we will show that the 
``stable''/``unstable'' manifolds of solutions of 
\eqnref{32}, decaying at $+\infty$/$-\infty$ with rates 
$\sim e^{\underline {\alpha}x}$/$e^{\balpha x}$,
respectively, approximately {\it track} the corresponding
subspaces of the principal coefficient matrix $\BbbA(x,\delta)$ 
as they vary with $x$, lying
always within angle $\Cal{O}(\delta/\eta)$ and decreasing/increasing with
uniform rate 
$\sim e^{\tilde{\underline{\alpha}}x}$/$e^{\tilde {\balpha}x}$
for any
$\tilde{\underline{\alpha}}>\max_x \underline{\alpha}$,
$\tilde{\balpha}<\min_x\balpha$.

\medskip

{\bf Remark \thmlbl{generalgap}.}  The most common applications
of this result 
concern the situation that 
$$
\underline{\alpha}(x)\equiv
\underline{\alpha}<0<\balpha \equiv \balpha(x),
$$
in which case the manifolds described are truly stable/unstable
and $\balpha=\min_x \balpha$, $\underline{\alpha}=\max_x \underline{\alpha}$.
However, neither our results nor their useful application are 
limited to this case.
\medskip

We will follow an approach based on energy estimates/invariant cones,
though other approaches are certainly possible;
see [AGJ,GZ], for example, for a revealing alternative formulation 
of Proposition 8.22 (just below) in terms of projectivized differential forms.
For further discussion, we refer the reader to [ZH], Section 7.

\proclaim{Proposition \thmlbl{6}}  For $C$ sufficiently large, $\delta/\eta $ 
sufficiently small, the cone
$$
\BbbK_-(x):= \big\{w: \frac{|P(x)w|}{|Q(x)w|}\le C \delta/\eta \big\}
$$
is {\bf positively invariant} under the flow of \eqnref{32} (i.e.,
invariant in forward time $x \nearrow$), and {\bf exponentially
attracting} on 
$\BbbJ_-(x):= \{w:\frac{|P(x)w|}{|Q(x)w|}\le \eta/C\delta \}$. 
Here, $C$ depends only on
the bounds \eqnref{33}--\eqnref{A'} and \eqnref{35}
and $\eta$ as defined in \eqnref{etagap} measures spectral gap.  
\endproclaim

The motivation behind Proposition \thmref{6} is clear:  the $P$-component of
$w$ is in some sense growing exponentially slower than the
$Q$-component, by the amount of the spectral gap.  To quantify this
observation, we use the following standard linear algebraic result:

\proclaim{Lemma \thmlbl{7}}
For matrices $M$ s.t. $(M) \le C$ and $\R \delta (M)\ge \eta > 0$,
and for any $0 < \etaT < \eta$,
there exists a $C^\infty$ choice of coordinate transformation $S(M)$
such that $\tildeM:=SMS^{-1}$ is {\bf real positive definite},
satisfying 
$$
\R(\tildeM):= \frac{1}{2}(\tildeM + \tildeM^*)\ge \etaT>0.
$$
\endproclaim

{\bf Proof.}  See, e.g., Proposition 
A.9, p. 361 of [St].\myqed

It is worthwhile to consider the import of the
Lemma in the context of a constant coefficient ODE
$$
w'=Mw,
\eqnlbl{36}
$$
namely, the existence of a coordinate change $\tildeW=Sw$ for which
the flow of \eqnref{36}, now of form
$$
\tildeW '=\tildeM \tildeW,
\eqnlbl{37}
$$
strictly increases $|\tildeW|$.  For,
$$
\aligned
(\frac{1}{2} |\tildeW|^2)' = \frac{1}{2} 
(\langle\tildeW, \tildeM \tildeW \rangle +
\langle M \tildeW, \tildeW\rangle)
&= \langle \tildeW, \R \tildeM \tildeW\rangle\\
&\ge \frac{\eta}{2} |\tildeW|^2.
\endaligned \eqnlbl{38}
$$
Related topics are Lyapunov theory and the Kreiss matrix theorem.

\medskip
{\bf Proof of Proposition \thmref{6}.}  By rescaling if necessary by
$y=e^{-1/2(\underline{\alpha}+ \balpha)x}W$, 
we can reduce without loss of generality to the symmetric case
$$
\alpha _1, \dots , \alpha _\ell \le -\eta < 0 < \eta \le \alpha
_{\ell + 1}, \dots , \alpha _N.
\eqnlbl{38.1}
$$
By \eqnref{A'}, there exists a change of coordinates $w\to S(x)w$
reducing \eqnref{32} further to the case that $\BbbA(x,\delta )$ 
has block diagonal form,
$$
\BbbA(x,\delta ) =
\pmatrix 
M_1 &0\\
0 &M_2\\
\endpmatrix .
\eqnlbl{39}
$$
For, by standard matrix theory [Kat], \eqnref{35} implies the
existence of bases spanning the left and right eigenspaces 
associated with projections $P$, $Q$, and depending smoothly on $P$, $Q$; 
it follows that these bases, and thus the corresponding diagonalizing
transformation $S$ have derivatives of order $\delta$, whence
error term $S^{-1}S'w$ may be absorbed in $\Theta$.
Finally, using Lemma \thmref{7}, we can reduce by a further
coordinate change $w\to Tw$ to the case that
$$
\R M_1 < -\etaT, \quad \R M_2 > \etaT.
\eqnlbl{40}
$$

Expressing \eqnref{32} coordinate-wise, we have
$$
w_j'=M_j w_j + \delta \Theta_j w,
\eqnlbl{43}
$$
where $|\Theta_j|\le C $,
from which we obtain the growth/decay estimates
$$
\aligned
|w_j |'
&= \R \langle w_j/|w_j|,w_j'\rangle
\\
&=\langle w_j/|w_j|,\R(M_j)w_j\rangle +
\R \langle w_j/|w_j|, \delta \Theta_jw \rangle \\
&\zigzag \mp \etaT|w_j| + C \delta |w|.
\endaligned \eqnlbl{44}
$$
Defining $r:=\frac{|w_1|}{|w_2|}$, we thus find,
after some simplification, that
$$
\aligned
r'
&=\frac{|w_2|' |w_1| - |w_2||w_1|'}
{|w_2|^2}\\
&< \frac{-2\etaT|w_2| |w_1| + C \delta
|w|^2}{|w_2|^2}. \\
&= -2\etaT r + C \delta (1+r^2).\\
\endaligned \eqnlbl{45}
$$

It is thus clear that $\BbbK_-:=\{w:r\le \frac{\tildeC \delta
}{4\etaT}\}$ is invariant, provided $\delta $ is small enough that
$\frac{C \delta }{4\etaT}\le 1$.  More generally, we have,
for $\frac{C \delta }{4\etaT}\le r \le
\frac{\etaT}{4C \delta }$, that $r' \le -\etaT r$.
Thus, $\BbbK_-$ is exponentially attracting on $\BbbJ_-:=\{w:r\le
\frac{\etaT}{4C \delta}\}$.\myqed

\proclaim{Corollary \thmlbl{8}}  For $C$ sufficiently large, $\delta/\eta $ 
sufficiently small,
solutions in 
$$
\BbbJ_-(x):=\big\{w:\frac{|P(x)w|}{|Q(x)w|}\le \eta/C\delta \big\}
$$
increase as $x \to +\infty$ at exponential rate 
$e^{\tilde \alpha x}$, for any 
$\tilde \alpha< \liminf_{x\to +\infty}
 \balpha$ ($\balpha$ as defined in \eqnref{34}).
\endproclaim

{\bf Proof.}  By Proposition \thmref{6}, solutions in $\BbbJ_-(x)$ eventually
enter cone 
$$
\BbbK_- := \{w: \frac{|P(x)w|}{|Q(x)w|}\le C \delta/\eta \}.
$$
Thus, \eqnref{44}, translated back to original coordinates, gives
$$
\aligned
|w_2|'
&>\buildrel \approx \over \alpha |w_2| - (C^2 \delta^2/\eta)
|w_2|\\ 
&\ge  \tilde \alpha |w_2|,
\endaligned \eqnlbl{46}
$$
for any $\buildrel \approx \over \alpha< \alpha$, 
$\tilde \alpha< \buildrel \approx \over \alpha- (C^2 \delta^2/\eta) $, 
from which the result follows.\myqed 

\proclaim{Corollary \thmlbl{9} (The Tracking Lemma)}  
For $\delta/\eta $ sufficiently small,
solutions $w^+$/$w^-$ of \eqnref{32} decaying at $+\infty$/$-\infty$ 
at rate $e^{ \tilde{\tilde{\underline{\alpha}}} x}$/$e^{
\tilde{\tilde{\balpha}} x}$
lie always within cones
$$
\BbbK_-(x):= \{w:\frac{|P(x)w|}{|Q(x)w|}\le C\delta/\eta \},
\eqnlbl{47}
$$
$$
\BbbK_+(x):= \{w:\frac{|Q(x)w|}{|P(x)w|}\le C\delta/\eta \},
\eqnlbl{48}
$$
respectively, for any $\tilde{\tilde{\underline{\alpha}}}
<\liminf_{x\to +\infty} {\balpha}$,
$\tilde{\tilde{\balpha}}> \limsup_{x\to -\infty} {\underline{\alpha}}$.
Moreover, there hold the uniform decay/growth rates
$$
\aligned
\frac{|w^+(x)|}{|\bar w^+(y)|} &\le C e^{\tilde{ \underline{\alpha}}|x-y|}, \\
\frac{|w^-(x)|}{|\bar w^-(y)|} &\ge C^{-1} e^{\tilde{ \balpha} |x-y|}, \\
\endaligned
\eqnlbl{uniformrates}
$$
for all $x>y$, and symmetrically for $x<y$, for any
$\tilde{ \underline{ \alpha}}> \max_x \underline{\alpha}$, 
$\tilde {\balpha}< \min_x \balpha$. (Note: $C$ depends in part upon the
choice of $\tilde {\underline{\alpha}}$ and $\tilde{\balpha}$).
\endproclaim

{\bf Proof.}  Without loss of generality, consider the $+\infty$ case.  
Bounds \eqnref{47}--\eqnref{48} follow immediately by contradiction 
from the previous Corollary, whence \eqnref{uniformrates} follows
from \eqnref{46}.\myqed

{\bf Remark \thmlbl{constcase}.}
In the ``standard'' case that $\underline{\alpha}<0<\balpha$,
$\underline{\alpha}$ and $\balpha$ constant, 
the statement of Corollary \thmref{9} considerably simplifies;
in particular, the conclusion applies to solutions merely
{\it decaying} at $+\infty$/$-\infty$, with no specified rate.
In the case of varying $\underline{\alpha}$ and $\balpha$, 
the rates given in \eqnref{uniformrates} can clearly be sharpened to
$e^{\int_y^x \tilde{\underline{\alpha}}(z)dz}$ and
$e^{\int_y^x \tilde{\balpha}(z)dz}$, without any change in the argument. 

\medskip

{\bf A3.2.2. Extensions.}
(i) {\it The case $\eta\to 0$}. 
The treatment above is essentially that given in [ZH,Z.1,Z.4], 
and suffices for the sectorial operators studied there.
It suffices also for the nonsectorial, dispersive--diffusive 
operators studied in [HZ.2], or even purely dispersive operators,
since in either case a uniform spectral gap is maintained up to the 
essential spectrum boundary.
However, it does not directly apply to the nonsectorial
operators arising in the study of relaxation or real viscosity models, 
for which $\underline{\alpha}$, $\balpha$, and $\eta$
typically vary with $\delta$, with $\eta$, $\delta/\eta \to 0$ 
as $\delta \to 0$, and derivative $A'$ is typically of order
$\eta >> \delta$.

Here, we only point out that all conclusions of the
above subsection remain valid also in this more general setting, provided
that we can find a coordinate change reducing \eqnref{32} to form
\eqnref{39}--\eqnref{40}.
Such a reduction may in fact be carried out under rather general
circumstances, as discussed in Section 4, above: for example,
if the distance between {\it all} eigenvalues of $A$ remains uniformly 
bounded from below even as the difference between their real parts 
goes to zero (as occurs for strictly hyperbolic relaxation models),
or if the distance between {\it groups}
$\{\alpha _1, \dots \alpha _\ell\}(x,\delta)$ and
$\{\alpha _{\ell+1}, \dots \alpha _N\}(x,\delta)$ 
remains uniformly bounded from below, while at the same time
$|A'|\le \eta/C$ with $C>0$ sufficiently large, 
(as occurs for general relaxation models, away from the shock layer).
\medskip

(ii) {\it Contraction mapping formulation.} 
The analysis of this subsection can equivalently 
be phrased as an integral equation/contraction mapping argument 
giving existence of stable/unstable manifolds at the same
time that estimates are obtained.  This has the advantage
of slightly further generality: for example, one can obtain
complete {short-time} theory/pointwise bounds (which depend
only on high-frequency estimates) by this technique with no 
assumptions other than smoothness on the coefficients of $L$.
Moreover, it disassociates regularity in $\lambda$ from 
smoothness in $x$, allowing high accuracy in $\lambda$
with minimal regularity assumptions on coefficients. 
For example, in the application to relaxation systems in Section 4,
it allows us to terminate our diagonalization process
at the first order, with $\delta \sim O(|\lambda|^{-2})$,
$\eta \sim |\lambda|^{-1}$, while still achieving the necessary
regularity in $\lambda$ to order $\lambda^{-1}$, with error
estimates to order $|\lambda|^{-2}$, whereas with the invariant
cone approach we should have had to diagonalize to one further order,
$\delta \sim \Cal{O}(|\lambda|^{-3}$, to obtain the needed estimates.
This would have increased our smoothness requirements 
in hypothesis (H0) to $f$, $g$, $q\in C^4$ for general, strictly hyperbolic
relaxation systems, $C^3$ for strictly hyperbolic 
discrete kinetic systems.

Specifically,  after first reducing to form \eqnref{39}--\eqnref{40},
where $M_j$ are assumed further to be {\it diagonal},
we may consider the ``lifted equations'' to the space of 
$\ell$- (resp. $(N-\ell)$-) forms, as in \eqnref{3.14}, Section A3.1,
in this case appearing as
$$
\eta ' =\BbbA^{(\ell)}(x;\lambda )\eta+ \delta \Bbb{\Theta} \eta,
\eqnlbl{lifted}
$$
where $\Cal{E}_1:=e_1\wedge \cdots\wedge e_\ell$ is the unique eigenvector
of $\BbbA^{(k)}$ corresponding to the eigenvalue of minimum real part, 
which moreover is separated from all other eigenvalues with spectral
gap $2\eta$, $e_j$ as usual denoting the standard $j$th basis element
in $\BbbC^n$.

Enumerating the standard basis elements of the space of $\ell$-forms
as $\Cal{E}_j$, and defining $\eta_j$ to be the coordinate
representation of $\eta$ with respect to this basis, we obtain
from structure \eqnref{39}--\eqnref{40} equations for $\eta_j'$
with the same, block-diagonal structure, where the first block,
$\Cal{M}_1$, corresponding to $\eta_1$ is scalar, with $\R \Cal{M}_1<-\eta$,
and (by a straightforward calculation using definition \eqnref{3.15})
the second block again satisfies $\R \Cal{M}_2 >\eta$.
Thus, defining $\hat {\eta}_j:=\eta_j/\eta_1$ for $j\ge 2$, 
$\hat{\eta}:=(\hat{\eta}_2, \dots)$,
we obtain projectivized equations of form 
$$
{\hat \eta}= (\Cal{M}_2-\Cal{M}_1 I)(x) \hat{\eta} + \delta \Cal{Q}(\hat{\eta},x)
\eqnlbl{etahat}
$$
where 
$$
\R (\Cal{M}_2-\Cal{M}_1I) < -2\eta
\eqnlbl{posM}
$$
for all $j\ge 2$ and $\Cal{Q}(r,x)$ is a quadratic polynomial in $r$,
with uniformly bounded coefficients,
similarly as in \eqnref{45}.

Define, similarly as in Section A3.1, the map
$$
\Cal{T}\hat{\eta}(x) :=
\delta \int^\infty_x \Cal{F}^{z\to x}
\Cal{Q}(\hat{\eta}(z),z) \, dz,
\eqnlbl{T}
$$
where $\Cal{F}^{z\to x}$ denotes the flow from $z$ to $x$ of the
approximate equation
$$
{\hat \eta}= (\Cal{M}_2-\Cal{M}_1 I)(x) \hat{\eta}.
\eqnlbl{approxetahat}
$$
From \eqnref{posM}, we obtain the bound 
$$
|\Cal{F}^{z\to x}| \le C e^{-2\eta (z-x)}
\eqnlbl{flowbound}
$$
for all $z>x$, similarly as in \eqnref{44}.
From this, together with the quadratic form of $\Cal{Q}$, 
we readily find that $\Cal{T}$ is a contraction on any ball in $L^\infty$,
with constant $C\delta/\eta$,
provided that $\delta/\eta$ is sufficiently small.
From this observation, we obtain both existence of the desired
fixed point $\hat{\eta}=\Cal{T}(\hat{\eta})$ and the bound
$$
|\hat{\eta}-\Cal{T}(0)|=
|\Cal{T}(\hat{\eta})-\Cal{T}(0)|\le (C\delta/\eta) |\hat{\eta}|,
$$
yielding $|\hat{\eta}|\le C_2 |\Cal{T}(0)|\le C_3 \delta/\eta$: the same
rate of tracking as obtained before using the invariant cone approach.
On the other hand, we have at the same time a means to generate 
the solution $\hat{\eta}$ to any desired order in $\delta/\eta$, by fixed point
iteration: in particular, we have the first-order expansion 
$$
\hat{\eta}=\Cal{T}(0)+ \Cal{O}(|\delta/\eta|^2),
\eqnlbl{series}
$$
where the first order term is explicitly given as an (finite) exponential
integral against the constant term of $\Cal{Q}$.

{\bf Remark \thmlbl{layercalc}.} 
Note that block-diagonal structure of $\BbbA$ is necessary in order
that the approximate flow \eqnref{approxetahat} be linear; in general, 
it could be an arbitrary quadratic function, in which case 
Duhamel's principle would not apply.
\medskip

{\bf Remark \thmlbl{graph}.} 
The estimate \eqnref{series} in terms of 
$\ell$- (resp. $(N-\ell)$)
forms yields a corresponding estimate in terms of subspaces: that is,
a description of the stable manifold at $+\infty$ (resp.
unstable manifold at $-\infty$) as a graph over $w_2\equiv 0$,
smooth in $x$ and analytic in $\delta/\eta$. 	
\medskip

{\bf Remark \thmlbl{layercalc}.} 
Clearly, the global decay bound \eqnref{flowbound} suffices for the above 
construction.  The pointwise decay afforded by negativity
of $\Cal{M}_2-\Cal{M}_1I$ is not strictly necessary, and may
be allowed to fail by order $\alpha$ on intervals with total length of order
$1/\alpha$.  This observation allows us to treat strong relaxation shocks,
for which negativity may fail by order $\eta$ on a shock layer of 
order $1/\eta$.  
Likewise, in the tracking construction of A3.2.1 and A3.2.2(i),
negativity/positivity of $M_1$/$M_2$ may be allowed to fail by
order $\alpha$ on intervals with total length of order $1/\alpha$,
during which solutions may move to a larger cone (by a bounded
multiple),  which may be
chosen to be invariant on ``good'' intervals.
%

\medskip

{\bf A3.2.3. The reduced flow.}
Finally, suppose that we have successfully carried out
the reduction \eqnref{39} of \eqnref{32} to block-diagonal
form in such a way that \eqnref{33} remains valid; 
that we have verified positivity 
condition \eqnref{40} (except perhaps on 
a finite collection of intervals of total length order $1/\alpha$,
where positivity may be allowed to fail to order $\alpha$, see
Remark \thmref{layercalc}), or by some other means established
uniform exponential decay/growth 
$$
|\Cal{F}_j^{y\to x}|\lessgtr Ce^{\mp \theta \eta |x-y|},
\eqnlbl{unexp}
$$
$\theta>0$,
of flows $\Cal{F}_j$ associated with
the approximate (decoupled) equations
$$
w_j'=M_jw_j;
\eqnlbl{approxdec}
$$
and that $\delta/\eta \to 0$ as $\delta\to 0$,
so that the entire stable/unstable manifold construction 
described in Sections A3.2.2 (i)--(ii) can be carried out
(alternatively, the estimates of Sections A3.2.1-A3.2.2(i)).
Then, the ``reduced flow'' restricted to stable/unstable
manifolds takes the simple form
$$
w_j'= M_jw_j + \delta \Theta_{jj} w_j + 
(\delta^2/\eta) \Cal{R}_j,
\eqnlbl{redflow}
$$
where $\Theta_{jj}$, $j=1,2$ are the diagonal blocks of error $\Theta$,
and terms $\Cal{R}_j=\Cal{O}(1)$ corresponding to off-diagonal blocks,
arise through the basic estimate $|w_{\tilde j}|/|w_j|\le \delta/\eta$,
$\tilde j\ne j$.

\proclaim{Proposition \thmlbl{redest}}
The flows (i.e., solution operator) $\bar{\Cal{F}}_j^{y\to x}$ of the
reduced equations \eqnref{redflow} satisfy
$$
\bar{\Cal{F}}_j^{y\to x}=\Cal{F}_j^{y\to x} + (\delta/\eta)\Cal{E}_j(x,y,\delta)
+ \Cal{O}(\delta^2/\eta^2)e^{-\tilde \theta \eta |x-y|},
\eqnlbl{redexp}
$$
for $x\gtrless y$, for any $0<\tilde \theta<\theta$,
where 
$$
\Cal{E}_1:= 
\eta \int^{x}_y \Cal{F}_1^{z\to x}
\Theta_{11}(z,\delta){\Cal{F}}_1^{y\to z}
\, dz  
\eqnlbl{Edef}
$$
satisfies the uniform exponential decay estimate
$$
|\Cal{E}_j(x,y,\delta)| \le 
C_1\eta|x-y| e^{-\theta \eta |x-y|}
\le C_2
e^{-\tilde{ \theta} \eta |x-y|},
\eqnlbl{undecay}
$$
is smooth in $x$ and $y$ and as smooth in $\delta$ as is $\Theta_{11}$,
and $\Cal{F}_j^{y\to x}$ are the approximate flows associated with 
\eqnref{approxdec}
$C_j$ and $\Cal{O}(\cdot)$, depending only on $\tilde \theta$
and the bounds \eqnref{unexp}, are of order $C(\theta-\tilde \theta)^{-1}$,
where $C$ is the constant of \eqnref{unexp}.
\endproclaim

{\bf Proof.}
Restrict to the case of the stable manifold, $j=1$; the case
$j=2$ may be treated in symmetric fashion.
Fixing $y$, define similarly as in previous subsections the map
$$
\Cal{T}\bar{\Cal{F}}_1^{y\to x} :=
{\Cal{F}}_1^{y\to x} +
\int^{x}_y \Cal{F}_1^{z\to x}
\big(\delta \Theta_{11} + (\delta^2/\eta)\Cal{R}_1\big)(z,\delta)
\bar{\Cal{F}}_1^{y\to z}
\, dz,
\eqnlbl{T}
$$
where $\Cal{F}_1^{z\to x}$ denotes the flow from $z$ to $x$ of the
approximate equation \eqnref{approxdec}.
Then, from the decay estimate \eqnref{unexp}, it is readily established
as in the argument of Section A3.2.2(ii) that
$\Cal{T}$ is a contraction on $L^\infty[y,+\infty)$, with constant
of order $\delta/\eta$, establishing the existence of a unique fixed-point
solution.  Carrying out two iterations of the fixed
point iteration, we obtain the estimate
$$
\aligned
\bar{\Cal{F}}_1^{y\to x} &=
\Cal{T}^2(0) + \Cal{O}(\delta^2/\eta^2)\\
&=
{\Cal{F}}_1^{y\to x} +
\int^{x}_y \Cal{F}_1^{z\to x}
\big(\delta \Theta_{11} + (\delta^2/\eta)\Cal{R}_1\big)(z,\delta)
{\Cal{F}}_1^{y\to z}
\, dz  + \Cal{O}(\delta^2/\eta^2)\\
&=
{\Cal{F}}_1^{y\to x} +
\delta \int^{x}_y \Cal{F}_1^{z\to x}
\Theta_{11}(z,\delta){\Cal{F}}_1^{y\to z}
\, dz  + \Cal{O}(\delta^2/\eta^2),\\
\endaligned
\eqnlbl{lastcalc}
$$
yielding \eqnref{redexp} with $\tilde \theta=0$ and \eqnref{Edef}.

The term $\Cal{E}_1$ may be estimated as 
$$
\aligned
|\Cal{E}_1|&\le 
\eta \int^{x}_y |\Cal{F}_1^{z\to x}|
|\Theta_{11}(z,\delta)|\,|{\Cal{F}}_1^{y\to z}| 
\, dz
\le
C\eta \int^{x}_y e^{-\theta \eta |x-z|}
e^{-\theta \eta |y-z|}
\, dz\\
&=
C\eta|x-y| e^{-\theta \eta |x-y|}
\le C \eta (|x-y| e^{-(\theta-\tilde \theta)\eta |x-y|})
e^{-\tilde{ \theta} \eta |x-y|},\\
&\le C(\theta-\tilde \theta)^{-1}
e^{-\tilde{ \theta} \eta |x-y|},\\
\endaligned
$$
for any $0<\tilde \theta<\theta$, and likewise the term
involving $\Cal{R}_1$ in the second line of \eqnref{lastcalc}
may be estimated as
$\Cal{O}(\delta^2/\eta^2) e^{-\tilde{ \theta} \eta |x-y|})$.
By a similar calculation, we can show that
$\Cal{T}$ is actually a contraction in the exponentially
weighted norm 
$$
\|f\|_{\tilde \theta}:= 
\|f(\cdot)e^{\tilde \theta \eta (\cdot-y)} \|_{L^\infty[y,+\infty)},
$$
yielding the improved bound 
$\Cal{O}(\delta^2/\eta^2)e^{-\tilde \theta \eta |x-y|}$,
on the third, $\Cal{O}(\cdot)$ term in the second line of \eqnref{lastcalc}.
Specifically, we have
$$
\aligned
\|\Cal{T}\Cal{F}_1- &\Cal{T}\Cal{F}_2\|_{\tilde \theta}\\
&\le
e^{\tilde \theta \eta |x-y|}
\|\Cal{F}_1- \Cal{F}_2\|_{\tilde \theta}
\delta \int^{x}_y |\Cal{F}_1^{z\to x}|
|(\Theta_{11}+ (\delta/\eta)\Cal{R}_1)(z,\delta)|\,e^{-\tilde \eta |z-y|}
\, dz\\
&\le
C\delta \|\Cal{F}_1- \Cal{F}_2\|_{\tilde \theta}
e^{\tilde \theta \eta |x-y|}
\int^{x}_y e^{-\theta \eta |x-z|}
e^{-\tilde \theta \eta |y-z|}
\, dz\\
& =
C\delta \|\Cal{F}_1- \Cal{F}_2\|_{\tilde \theta}
\int^{x}_y e^{-(\theta - \tilde \theta) \eta |x-z|} \, dz\\
&\le
C(\theta-\tilde \theta)^{-1}
(\delta/\eta) \|\Cal{F}_1- \Cal{F}_2\|_{\tilde \theta}.
\endaligned
$$
Combining these observations, we obtain \eqnref{redexp} as claimed.
\myqed

{\bf Remark \thmlbl{diagrmk}.}  Represention \eqnref{redexp} becomes
particularly simple in the case that $M_j$ are {\it diagonal},
and $\Cal{F}_j^{y\to x}=e^{\int_y^x M_j(z)} dz$.
This occurs, for example,
in the case of {\it strictly hyperbolic} relaxation systems,
or in the case of {\it scalar} (third and) higher order equations 
considered in [HZ.2].
However, for more general applications, it is important that we
do not limit ourselves to this case.  For {\it nonstrictly hyperbolic}
relaxation systems, and for higher order {\it systems},
$M_j$ are block diagonal, with blocks of form
$$
m_{j,k}= \beta_{j,k} + \Cal{O}(\delta),
$$
$\beta_{j,k}$ scalar.
Thus, it is possible again to express the dominant (typically
highly oscillatory) effects in each block as a simple exponential
integral factor multiplying a higher-order, coupled flow. 
We describe such a calculation in detail in Section 4, above.

Note that we do not require in Proposition \thmref{redest}
any spectral gap between different modes (blocks) of  $M_j$
in order to obtain an approximately block-diagonal flow.
This is important, for example, in the case of dispersive--diffusive
equations considered in [HZ.2].
Indeed, Proposition \thmref{redest} and the contraction mapping
construction of Section A3.2.2(ii) together repair an omission in the
analogous constructions carried out in [HZ.2]:
specifically in the proof of Lemma 3.2, [HZ.2],
wherein a nonexistent gap was implicitly assumed.
\footnote{
More precisely, by mistake,  only estimates for a solo mode were carried out, 
rather than for conjugate modes as stated.  
Note that error bounds of type \eqnref{redexp}--\eqnref{undecay} 
indeed fit into the claimed bounds of Lemma 3.2, [HZ.2], 
validating the later conclusions of that paper.
}

\bigskip

\Refs
\medskip\noindent
[AGJ] J. Alexander, R. Gardner and C.K.R.T. Jones,
{\it A topological invariant arising in the analysis of
traveling waves}, J. Reine Angew. Math. 410 (1990) 167--212.
\medskip\noindent
[AR] A. Aw and M. Rascle, 
{\it Resurrection of "second order" models of traffic flow,}
SIAM J. Appl. Math. 60 (2000), no. 3, 916--938 (electronic).
\medskip\noindent
[BE] A.A. Barmin and S.A. Egorushkin,
{\it Stability of shock waves,}
Adv. Mech. 15 (1992) No. 1--2, 3--37.
\medskip\noindent
[BGDV] G.I. Barenblatt, J. Garcia-Azorero, A. De Pablo, J.L. Vazquez, 
{\it Mathematical model of the non-equilibrium
water-oil displacement in porous strata,} 
Appl. Anal. 65 (1997), no. 1-2, 19--45. 
\medskip\noindent
[BV] G.I. Barenblatt and A.P. Vinnichenko, 
{\it Nonequilibrium filtration of nonmixing fluids,}
(Russian) Adv. in Mech. 3 (1980), no. 3, 35--50.
\medskip\noindent
[BSZ] S. Benzoni--Gavage, D. Serre, and K. Zumbrun,
{\it Alternate Evans functions and viscous shock waves,}  
SIAM J. Math. Anal. 32 (2001), no. 5, 929--962 (electronic).
\medskip\noindent
[Br.1] L. Q. Brin, {\it Numerical testing of the stability of viscous
shock waves}, Ph.D. dissertation, Indiana University, May 1998.
\medskip\noindent
[Br.2] L. Q. Brin,
{\it Numerical testing of the stability of viscous shock waves,}
to appear, Math. Comp.
\medskip\noindent
[BR]
G. Boillat and T. Ruggeri, 
{\it On the shock structure problem for hyperbolic system 
of balance laws and convex entropy,} 
Contin.  Mech. Thermodyn. 10 (1998), no. 5, 285--292. 
\medskip\noindent
[B] F. Bouchut
{\it Construction of BGK models with a family of kinetic entropies 
for a given system of conservation laws,} 
J. Statist. Phys. 95 (1999), no. 1-2, 113--170. 
\medskip\noindent
[CaL] R.E. Caflish and T.-P. Liu,
{\it Stability of shock waves for the Broadwell equations,}
Comm. Math. Phys. 114 (1988) 103--130.
\medskip\noindent
[Ce] C. Cercignani,
{\it The Boltzmann equation and its applications,}
Applied Mathematical Sciences, 67. Springer-Verlag, New York (1988)
xii+455 pp. ISBN: 0-387-96637-4.
\medskip\noindent
[CLL] G.-Q. Chen, D.C. Levermore, and T.-P. Liu, 
{\it Hyperbolic conservation laws with 
stiff relaxation terms and entropy,} Comm. Pure
Appl. Math. 47 (1994), no. 6, 787--830. 
\medskip\noindent
[CL] E.A. Coddington and N. Levinson, {\it Theory of ordinary
differential equations,} McGraw-Hill Book Company, Inc., 
New York-Toronto-London (1955) xii+429 pp.
\medskip\noindent
[Co] W. A. Coppel,
{\it Stability and asymptotic behavior of differential equations},
D.C. Heath and Co., Boston, MA (1965) viii+166 pp. 
\medskip\noindent
[CH.1] R. Courant and D. Hilbert, 
{\it Methods of mathematical physics. Vol. I,}
Interscience Publishers, Inc., New York, N.Y., 1953. xv+561 pp.
\medskip\noindent
[CH.2] R. Courant and D. Hilbert, 
{\it Methods of mathematical physics. Vol. II:
Partial differential equations,}
(Vol. II by R. Courant.) Interscience Publishers (a division
of John Wiley \& Sons), New York-Lon don 1962 xxii+830 pp. 
\medskip\noindent
[Do] J. Dodd,
{\it Convection stability of shock profile solutions of a modified
KdV--Burgers equation},
Thesis, University of Maryland (1996).
\medskip\noindent
[Dre] W. Dreyer, 
{\it Maximisation of the entropy in nonequilibrium,}
J. Phys. A 20 (1987), no. 18, 6505--6517. 
\medskip\noindent
[Er] J. J. Erpenbeck,
{\it Stability of step shocks,} Phys. Fluids 5 (1962) no. 10, 1181--1187.
\medskip\noindent
[E.1] J.W. Evans,
{\it Nerve axon equations: I. Linear approximations,}
Ind. Univ. Math. J. 21 (1972) 877--885.
\medskip\noindent
[E.2] J.W. Evans,
{\it Nerve axon equations: II. Stability at rest,}
Ind. Univ. Math. J. 22 (1972) 75--90.
\medskip\noindent
[E.3] J.W. Evans,
{\it Nerve axon equations: III. Stability of the nerve impulse,}
Ind. Univ. Math. J. 22 (1972) 577--593.
\medskip\noindent
[E.4] J.W. Evans,
{\it Nerve axon equations: IV. The stable and the unstable impulse,}
Ind. Univ. Math. J. 24 (1975) 1169--1190.
\medskip\noindent
[Ev] L.C. Evans, 
{\it Partial differential equations,} Graduate Studies in Mathematics, 
19. American Mathematical Society, Providence, RI,
1998. xviii+662 pp. ISBN: 0-8218-0772-2. 
\medskip\noindent
[FreZe] H. Freist\"uhler and Y. Zeng,
{\it Shock profiles for systems of balance
laws with relaxation,} preprint (1998).
\medskip\noindent
[Fr] A. Friedman,  
{\it Partial differential equations of parabolic type},
Prentice-Hall, Englewood Cliffs, NY (1964), Reprint Ed. (1983).
\medskip\noindent
[GJ.1] R. Gardner and C.K.R.T. Jones,
{\it A stability index for steady state solutions of
boundary value problems for parabolic systems},
J. Diff. Eqs. 91 (1991), no. 2, 181--203. 
\medskip\noindent
[GJ.2] R. Gardner and C.K.R.T. Jones,
{\it Traveling waves of a perturbed diffusion equation
arising in a phase field model}, 
Ind. Univ. Math. J. 38 (1989), no. 4, 1197--1222.
\medskip\noindent
[GZ] R. Gardner and K. Zumbrun,
{\it The Gap Lemma and geometric criteria for instability
of viscous shock profiles}, 
Comm. Pure Appl.  Math. 51 (1998), no. 7, 797--855. 
\medskip\noindent
[G] P. Godillon,
{\it Linear stability of shock profiles for systems of conservation laws 
with semi-linear relaxation,} (English. English summary) 
Phys. D 148 (2001), no. 3-4, 289--316. 
\medskip\noindent
[Go.1] J. Goodman, {\it Nonlinear asymptotic stability of viscous
shock profiles for conservation laws,}
Arch. Rational Mech. Anal. 95 (1986), no. 4, 325--344.
\medskip\noindent
[Go.2] J. Goodman, 
{\it Remarks on the stability of viscous shock waves}, in:
Viscous profiles and numerical methods for shock waves 
(Raleigh, NC, 1990), 66--72, SIAM, Philadelphia, PA, (1991).
\medskip\noindent
[He] D. Henry,
{\it Geometric theory of semilinear parabolic equations},
Lecture Notes in Mathematics, Springer--Verlag, Berlin (1981),
iv + 348 pp.
\medskip\noindent
[HoZ.1] D. Hoff and K. Zumbrun, 
{\it Multi-dimensional diffusion waves for the Navier-Stokes equations of
compressible flow,} Indiana Univ. Math. J. 44 (1995), no. 2, 603--676.
\medskip\noindent
[HoZ.2] D. Hoff and K. Zumbrun, 
{\it Pointwise Green's function bounds for multi-dimensional 
scalar viscous shock fronts,}  preprint (2000).
\medskip\noindent
[H.1] P. Howard,
{\it Pointwise estimates on the Green's function 
for a scalar linear convection-diffusion equation,}
J. Differential Equations 155 (1999), no. 2, 327--367.
\medskip\noindent
[H.2] P. Howard,
{\it Pointwise methods for stability of a scalar conservation law,}
Doctoral thesis (1998).
\medskip\noindent
[H.3] P. Howard,
{\it Pointwise Green's function approach to 
stability for scalar conservation laws,}
Comm. Pure Appl. Math. 52 (1999), no. 10, 1295--1313.
\medskip\noindent
[HZ.1] P. Howard and K. Zumbrun, 
{\it A tracking mechanism for one-dimensional 
viscous shock waves,} preprint (1999).
\medskip\noindent
[HZ.2] P. Howard and K. Zumbrun,
{\it Pointwise estimates for dispersive-diffusive shock waves,}
to appear, Arch. Rational Mech. Anal. 
\medskip\noindent
[Hu] J. Humpherys, 
{\it Stability of Jin--Xin relaxation shocks,} preprint (2000).
\medskip\noindent
[JX] S. Jin and Z. Xin, 
{\it The relaxation schemes for systems of conservation laws 
in arbitrary space dimensions,}
Comm. Pure Appl.  Math. 48 (1995), no. 3, 235--276.
\medskip\noindent
[J] C.K.R.T. Jones,
{\it Stability of the travelling wave solution of the FitzHugh--Nagumo system},
 Trans. Amer. Math. Soc.  286 (1984), no. 2, 431--469.
\medskip\noindent
[JGK] C. K. R. T. Jones, R. A. Gardner, and T. Kapitula,
{\it Stability of travelling waves
for non-convex scalar viscous conservation laws},
Comm. Pure Appl. Math. 46 (1993) 505--526.
\medskip\noindent
[Ka] Ya. Kanel,
{\it On a model system of equations of of one-dimensional
gas motion,} Diff. Eqns. 4 (1968) 374--380.
\medskip\noindent
[K.1] T. Kapitula,
{\it Stability of weak shocks in $\lambda$--$\omega$ systems},
Indiana Univ. Math. J. 40 (1991), no. 4, 1193--12.
\medskip\noindent
[K.2] T. Kapitula,
{\it On the stability of travelling waves in weighted $L^\infty$ spaces},
J. Diff. Eqs. 112 (1994), no. 1, 179--215.
\medskip\noindent
[KS] T. Kapitula and B. Sandstede,
{\it Stability of bright solitary-wave solutions 
to perturbed nonlinear Schrödinger equations,} Phys. D 124
(1998), no. 1-3, 58--103.
\medskip\noindent
[Kaw] S. Kawashima,
{\it Systems of a hyperbolic--parabolic composite type,
 with applications to the equations of magnetohydrodynamics},
thesis, Kyoto University (1983).
\medskip\noindent
[KM] S. Kawashima and A. Matsumura,
{\it Asymptotic stability of traveling wave solutions of systems 
for one-dimensional gas motion,} Comm.
Math. Phys. 101 (1985), no. 1, 97--127.
\medskip\noindent
[Kat] T. Kato,
{\it Perturbation theory for linear operators},
Springer--Verlag, Berlin Heidelberg (1985).
\medskip\noindent
[LU] O.A. Ladyzenskaja and N.N. Ural'tseva,
{\it On linear and quasij-linear equations
and systems of elliptic and parabolic types,} 1963 Outlines Joint Sympos.
Partial Differential Equations (Novosibirsk, 1963) 146--150 Acad. Sci. USSR
Siberian Branch, Moscow.
\medskip\noindent
[LSU] O.A. Ladyzenskaja, V.A. Solonnikov, and N.N. Ural'tseva,
{\it Linear and quasi-linear equations of parabolic type,}
Translations of Math. Monographs 23, AMS, Providence, RI (1968).
\medskip\noindent
[La] P.D. Lax,
{\it Hyperbolic systems of conservation laws and the mathematical theory of shock waves},
Conference Board of the Mathematical Sciences Regional Conference
Series in Applied Mathematics, No. 11. Society for Industrial and Applied Mathematics,
Philadelphia, Pa., 1973. v+48 pp.
\medskip\noindent
[LP] Lax, Peter D.; Phillips, Ralph S. Scattering theory. Second
edition. With appendices by Cathleen S. Morawetz and Georg Schmidt. 
Pure and Applied Mathematics, 26. Academic Press, Inc., 
Boston, MA, 1989. xii+309 pp. ISBN: 0-12-440051-5. 
\medskip\noindent
[Lev]
D.C. Levermore, 
{\it Moment closure hierarchies for kinetic theories,}
(English. English summary) 
J. Statist. Phys. 83 (1996), no. 5-6, 1021--1065. 
\medskip\noindent
[Le] E.E. Levi, {\it Sulle equazioni lineari totalmente ellittiche
alle derivate parziali,}
Rend. Circ. Mat. Palermo 24 (1907) 275--317.
\medskip\noindent
[Li] T. Li, 
{\it Global solutions and zero relaxation limit for a traffic flow model,}
SIAM J. Appl. Math. 61 (2000), no. 3, 1042--1061 (electronic).
\medskip\noindent
[Liu] H. Liu,
{\it Asymptotic stability of relaxation shock
profiles for hyperbolic conservation laws,}
preprint (2000).  
\medskip\noindent
[L.1] T.-P. Liu,
{\it Pointwise convergence to shock waves for viscous conservation laws},
Comm. Pure Appl. Math. 50 (1997), no. 11, 1113--1182.
\medskip\noindent
[L.2] T.-P. Liu, 
{\it Hyperbolic conservation laws with relaxation,}
Comm. Math. Phys. 108 (1987), no. 1, 153--175. 
\medskip\noindent
[L.3] T.-P. Liu, 
{\it The Riemann problem for general systems of conservation laws,}
J. Differential Equations 18 (1975), 218--234.
\medskip\noindent
[LZ.1] T.P. Liu and K. Zumbrun,
{\it Nonlinear stability of an undercompressive shock for complex
Burgers equation,} Comm. Math. Phys. 168 (1995), no. 1, 163--186.
\medskip\noindent
[LZ.2] T.P. Liu and K. Zumbrun,
{\it On nonlinear stability of general undercompressive viscous shock waves,}
Comm.  Math. Phys. 174 (1995), no. 2, 319--345.
\medskip\noindent
[LZe] T.-P. Liu and Y. Zeng, 
{\it Large time behavior of solutions for general 
quasilinear hyperbolic--parabolic systems of conservation laws},
AMS memoirs 599 (1997).
\medskip\noindent
[LW] M.J. Lighthill and G.B. Whitham,
{\it On kinematic waves. II. A theory of traffic flow on long crowded roads,}
Proc. Roy. Soc. London. Ser. A. 229 (1955), 317--345. 
\medskip\noindent
[M.1], A. Majda,
{\it The stability of multi-dimensional shock fronts -- a
new problem for linear hyperbolic equations,}
Mem. Amer. Math. Soc. 275 (1983).
\medskip\noindent
[M.2], A. Majda,
{\it The existence of multi-dimensional shock fronts,}
Mem. Amer. Math. Soc. 281 (1983).
\medskip\noindent
[M.3] A. Majda,
{\it Compressible fluid flow and systems of conservation laws in several
space variables,} Springer-Verlag, New York (1984), viii+ 159 pp.
\medskip\noindent
[MP] A. Majda and R. Pego,
{\it Stable viscosity matrices for systems of conservation laws},
J. Diff. Eqs. 56 (1985) 229--262.
\medskip\noindent
[MPl] Marchesin-Plohr, multiphase, TODO.
\medskip\noindent
[Ma] C. Mascia,
{\it Renormalization Method for the Chapman-Enskog Expansion for
Linear Relaxation Systems,} preprint (2001).
\medskip\noindent
[MaN] C. Mascia and R. Natalini, 
{\it $L\sp 1$ nonlinear stability of traveling waves for a 
hyperbolic system with relaxation,} J. Differential
Equations 132 (1996), no. 2, 275--292.
\medskip\noindent
[MN] A. Matsumura and K. Nishihara, 
{\it On the stability of travelling wave solutions of a one-dimensional 
model system for compressible viscous gas,}
Japan J. Appl. Math. 2 (1985), no. 1, 17--25.
\medskip\noindent
[MR] I. M\"uller and T. Ruggeri,
{\it Rational extended thermodynamics,}
Second edition, with supplementary chapters by H. Struchtrup 
and Wolf Weiss, Springer Tracts in Natural Philosophy, 37, 
Springer-Verlag, New York (1998) xvi+396 pp. ISBN: 0-387-98373-2.
\medskip\noindent
[N] R. Natalini, 
{\it Recent mathematical results on hyperbolic relaxation problems,}
TMR Lecture Notes (1998).
\medskip\noindent
[NT] R. Natalini and A. Tesei,
{\it On the Barenblatt model for non-equilibrium two phase flow 
in porous media,} Arch. Ration. Mech. Anal.
150 (1999), no. 4, 349--367.
\medskip\noindent
[OZ] M. Oh and K. Zumbrun,
{\it Stability of periodic solutions of viscous conservation laws:
Pointwise Green's function bounds,} Preprint (2001).
\medskip\noindent
[Pa] A. Pazy, {\it Semigroups of linear operators and applications 
to partial differential equations,} Applied Mathematical Sciences, 44, 
Springer-Verlag, New York-Berlin, (1983) viii+279 pp. ISBN: 0-387-90845-5.
\medskip\noindent
[PI] T.  Platkowski and R. Illner, 
{\it Discrete velocity models of the Boltzmann equation: 
a survey on the mathematical aspects of the
theory,} SIAM Rev. 30 (1988), no. 2, 213--255. 
\medskip\noindent
[Ro] P. Rosenau, {\it Extending hydrodynamics via the 
regularization of the Chapman-Enskog expansion,}
Phys. Rev. A (3) 40 (1989), no. 12, 7193--7196. 
\medskip\noindent
[R] F. Rousset, 
{\it Zero mass stability of strong relaxation shock profiles,}
Preprint (2001).
\medskip\noindent
[RS]  M. Reed and B. Simon, 
{\it Methods of modern mathematical physics. I-IV},
Academic Press, Inc. [Harcourt Brace Jovanovich,
Publishers], New York--London, 1980.
\medskip\noindent
[Sat] D. Sattinger,
{\it On the stability of waves of nonlinear parabolic systems},
Adv. Math. 22 (1976) 312--355.
\medskip\noindent
[SK] Y Shizuta and S. Kawashima,
{\it Systems of equations of hyperbolic--parabolic
type with applications to the discrete Boltzmann equation,}
Hokkaido Math. J. 14 (1984) 435--457.
\medskip\noindent
[Sl.1] M. Slemrod, 
{\it A renormalization method for the Chapman-Enskog expansion,}
Phys. D 109 (1997), no. 3-4, 257--273. 
\medskip\noindent
[Sl.2] M. Slemrod, 
{\it Renormalization of the Chapman-Enskog expansion: isothermal 
fluid flow and Rosenau saturation,} J. Statist. Phys. 91 (1998),
no. 1-2, 285--305.
\medskip\noindent
[St] J.C. Strikwerda, {\it Finite difference schemes
and partial differential equations,}
(Chapman and Hall, New York (1989) xii+ 386 pp.
\medskip\noindent
[S] A.  Szepessy, 
{\it On the stability of Broadwell shocks,}
Nonlinear evolutionary partial differential equations 
(Beijing, 1993), 403--412, 
AMS/IP Stud. Adv. Math., 3, Amer. Math. Soc., Providence, RI, 1997. 
\medskip\noindent
[SX.1] A. Szepessy and Z. Xin,
{\it Nonlinear stability of viscous shock waves,}
Arch. Rat. Mech. Anal. 122 (1993) 53--103.
\medskip\noindent
[SX.2] A. Szepessy and Z. Xin,
{\it Stability of Broadwell shocks,}
unpublished manuscript.
\medskip\noindent
[Ti] E.C. Titschmarsh,
{\it Eigenfunction expansion associated with second-order
differential equations I-II,} Oxford (1946), Revised Ed. (1958).
\medskip\noindent
[Wh] G.B. Whitham, 
{\it Linear and nonlinear waves,}
Reprint of the 1974 original, Pure and Applied Mathematics,
A Wiley-Interscience Publication, John Wiley \& Sons, Inc., 
New York, 1999. xviii+636 pp. ISBN: 0-471-35942-4. 
\medskip\noindent
[Yo.1] W.-A. Yong, 
{\it Singular perturbations of first-order 
hyperbolic systems,} PhD Thesis,
Universit\"at Heidelberg (1992).
\medskip\noindent
[Yo.2] W.-A. Yong, 
{\it Singular perturbations of first-order 
hyperbolic systems,} (English. English summary) 
Nonlinear hyperbolic problems: theoretical, applied, 
and computational aspects (Taormina, 1992), 597--604, 
Notes Numer. Fluid Mech., 43, Vieweg, Braunschweig, 1993. 
\medskip\noindent
[Yo.3] W.-A. Yong, 
{\it Singular perturbations of first-order hyperbolic systems 
with stiff source terms,} J. Differential Equations 155 (1999),
no. 1, 89--132. 
\medskip\noindent
[Yo.4] W.-A. Yong, 
{\it Basic properties of hyperbolic relaxation systems,}
Lecture notes, TMR Summer School Lectures, Kochel am See (1999).
\medskip\noindent
[YoZ] W.-A. Yong and K. Zumbrun,
{\it Existence of relaxation shock profiles 
for hyperbolic conservation laws,}
SIAM J. Appl. Math. 60 (2000), no. 5, 1565--1575 (electronic).
\medskip\noindent
[Y] K. Yosida, {\it Functional analysis,} 
Reprint of the sixth (1980) edition, Classics in Mathematics, 
Springer-Verlag, Berlin, 1995, xii+501 pp. ISBN: 3-540-58654-7.
\medskip\noindent
[Ze.1] Y. Zeng,
{\it $L^1$ asymptotic behavior of compressible, isentropic, 
viscous $1$-d flow,}
Comm. Pure Appl. Math. 47 (1994) 1053--1092.
\medskip\noindent
[Ze.2] Y. Zeng,
{\it Gas dynamics in thermal nonequilibrium 
and general hyperbolic systems with relaxation,}
Arch. Ration. Mech. Anal.  150 (1999), no. 3, 225--279.
\medskip\noindent
[ZH] K. Zumbrun and P. Howard,
{\it Pointwise semigroup methods and stability of viscous shock waves,}
Indiana Mathematics Journal V47 (1998) no. 4, 741--871.
\medskip\noindent
[Z.1] K. Zumbrun,  {\it Stability of viscous shock waves},
Lecture Notes, Indiana University (1998).
\medskip\noindent
[Z.2] K. Zumbrun, {\it Refined Wave--tracking and Nonlinear 
Stability of Viscous Lax Shocks}, to appear, Math. Anal. Appl. (2001);
preprint (1999).
\medskip\noindent
[Z.3] K. Zumbrun, {\it Refined Green's function bounds
and scattering behavior of viscous shock waves,} 
in preparation.
\medskip\noindent
[Z.4] K. Zumbrun, {\it Multidimensional stability of
planar viscous shock waves}, TMR Summer School Lectures:
Kochel am See, May, 1999, to appear,
Birkhauser's Series: Progress in Nonlinear Differential
Equations and their Applications (2001), 207 pp.
\medskip\noindent
[Z.5] K. Zumbrun, {\it Stability of general undercompressive shocks
of viscous conservation laws,}
in preparation.
\medskip\noindent
[ZS] K. Zumbrun and D. Serre,
{\it Viscous and inviscid stability of multidimensional 
planar shock fronts,} Indiana Univ. Math. J. 48 (1999), no. 3,
937--992.
\endRefs
\enddocument